\documentclass[english,onefignum,onetabnum]{siamart190516}

\usepackage[T1]{fontenc}
\usepackage{xcolor}
\usepackage{array}
\usepackage{float}
\usepackage{calc}
\usepackage{mathrsfs}
\usepackage{multirow}
\usepackage{amsmath}
\usepackage{amssymb}
\usepackage{graphicx}
\usepackage{wasysym}
\usepackage{multicol}

\makeatletter

\providecommand{\tabularnewline}{\\}
\floatstyle{ruled}
\newfloat{algorithm}{tbp}{loa}
\providecommand{\algorithmname}{Algorithm}
\floatname{algorithm}{\protect\algorithmname}

\numberwithin{equation}{section}
\numberwithin{figure}{section}

\usepackage{algorithm,algpseudocode}
\usepackage{array}

\makeatother

\usepackage{babel}
\begin{document}

\title{\textsc{A simple geometric method for navigating  the energy landscape of centroidal voronoi tessellations}}

\author{Ivan Gonzalez\thanks{Department of Mathematics and Statistics,
McGill University, Montreal, Quebec, H3A0B9
  (\email{ivan.gonzalez@mail.mcgill.ca}).}
\and Rustum Choksi\thanks{Department of Mathematics and Statistics,
McGill University, Montreal, Quebec, H3A0B9
 \newline  (\email{rustum.choksi@mcgill.ca}, \url{http://www.math.mcgill.ca/rchoksi/}).}
\and Jean-Christophe Nave\thanks{Department of Mathematics and Statistics,
McGill University, Montreal, Quebec, H3A0B9 
 \newline (\email{jean-christophe.nave@mcgill.ca}, \url{https://www.math.mcgill.ca/jcnave/}).}
}

\maketitle

\begin{abstract}
Finding optimal (or low energy) centroidal Voronoi tessellations (CVTs) on a 2D domain is a challenging problem. One must navigate an energy landscape whose desirable critical points  have sufficiently small basins of attractions
that they are inaccessible with Monte-Carlo initialized gradient descent methods. 
We present a simple deterministic method for efficiently navigating the energy landscape in order to {\it access} these  low energy CVTs.  The  method has two parameters and is based upon
each generator moving away from the closest neighbor 
by a certain distance. We give a statistical analysis of the performance of this hybrid method comparing with 
the results of a large number of runs for both Lloyd's method and state of the art quasi-Newton methods. Stochastic alternatives are also considered.
\end{abstract}



\section{Introduction}

\begin{figure}[]
\centering

\begin{tabular}{cc}
\includegraphics[width=0.25\columnwidth]{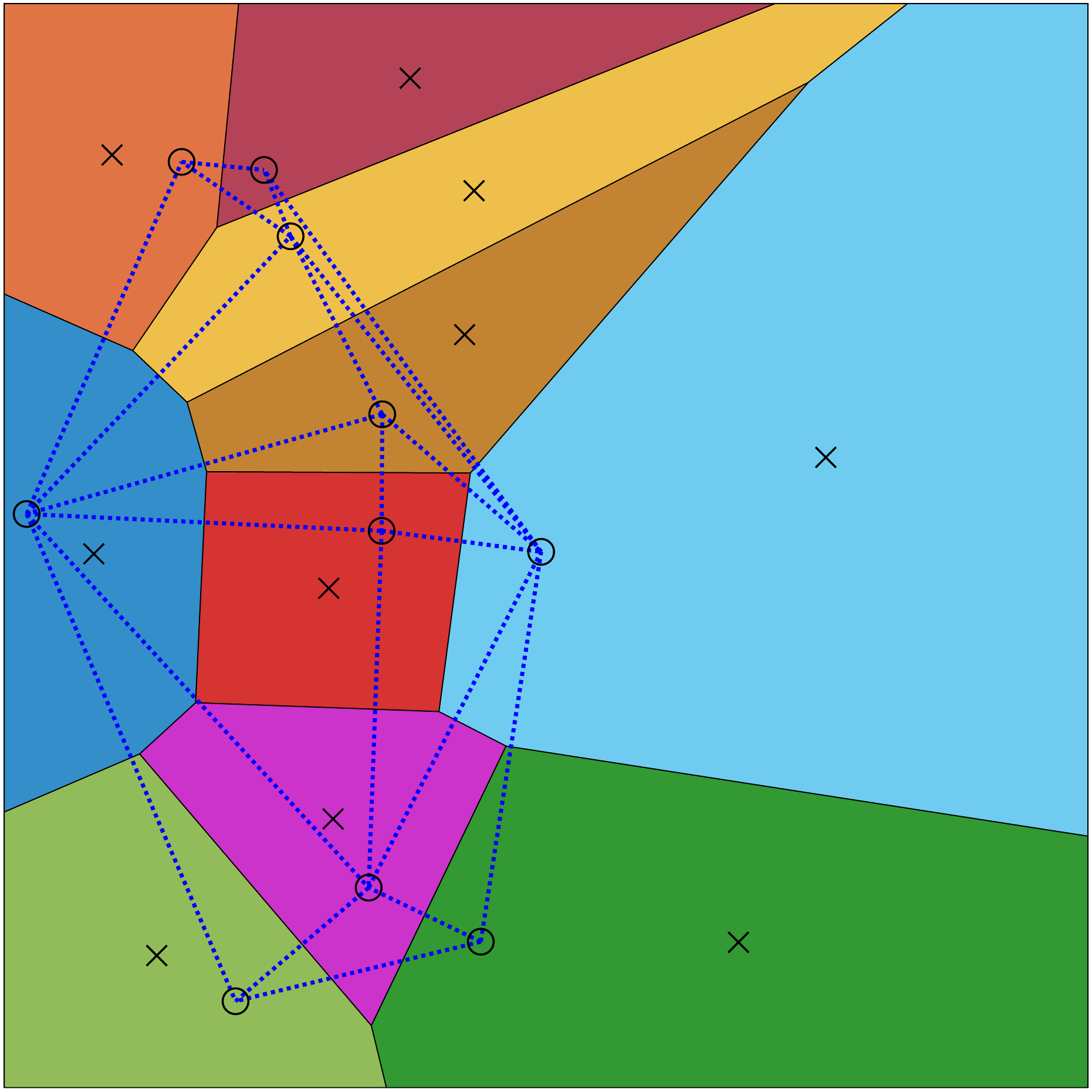} & \includegraphics[width=0.25\columnwidth]{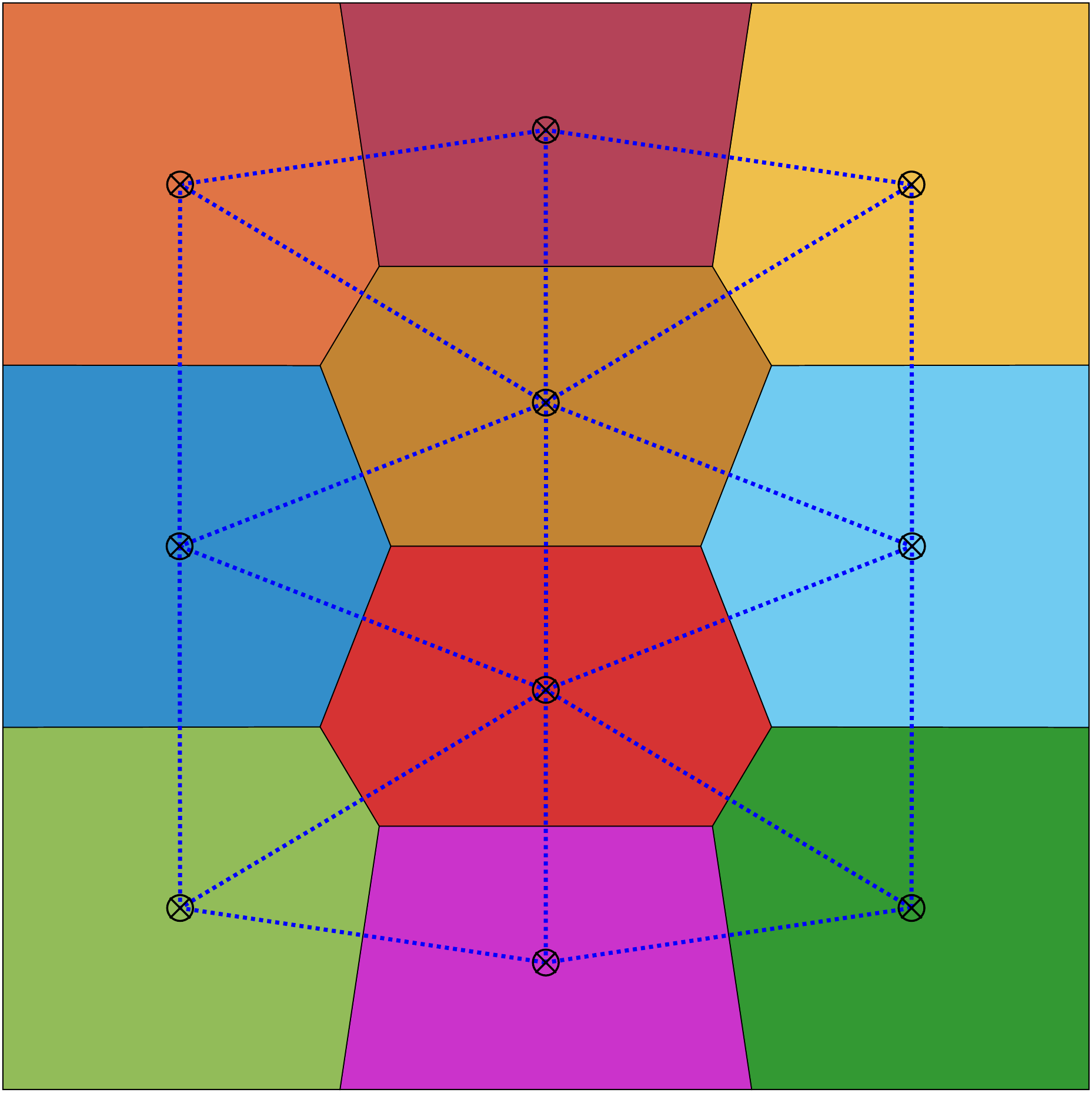}\tabularnewline
\footnotesize{(a) $\mathcal{V}(\mathbf{X})$}& \multicolumn{1}{c}{\footnotesize{(b) CVT} }\tabularnewline
 \includegraphics[width=0.25\columnwidth]{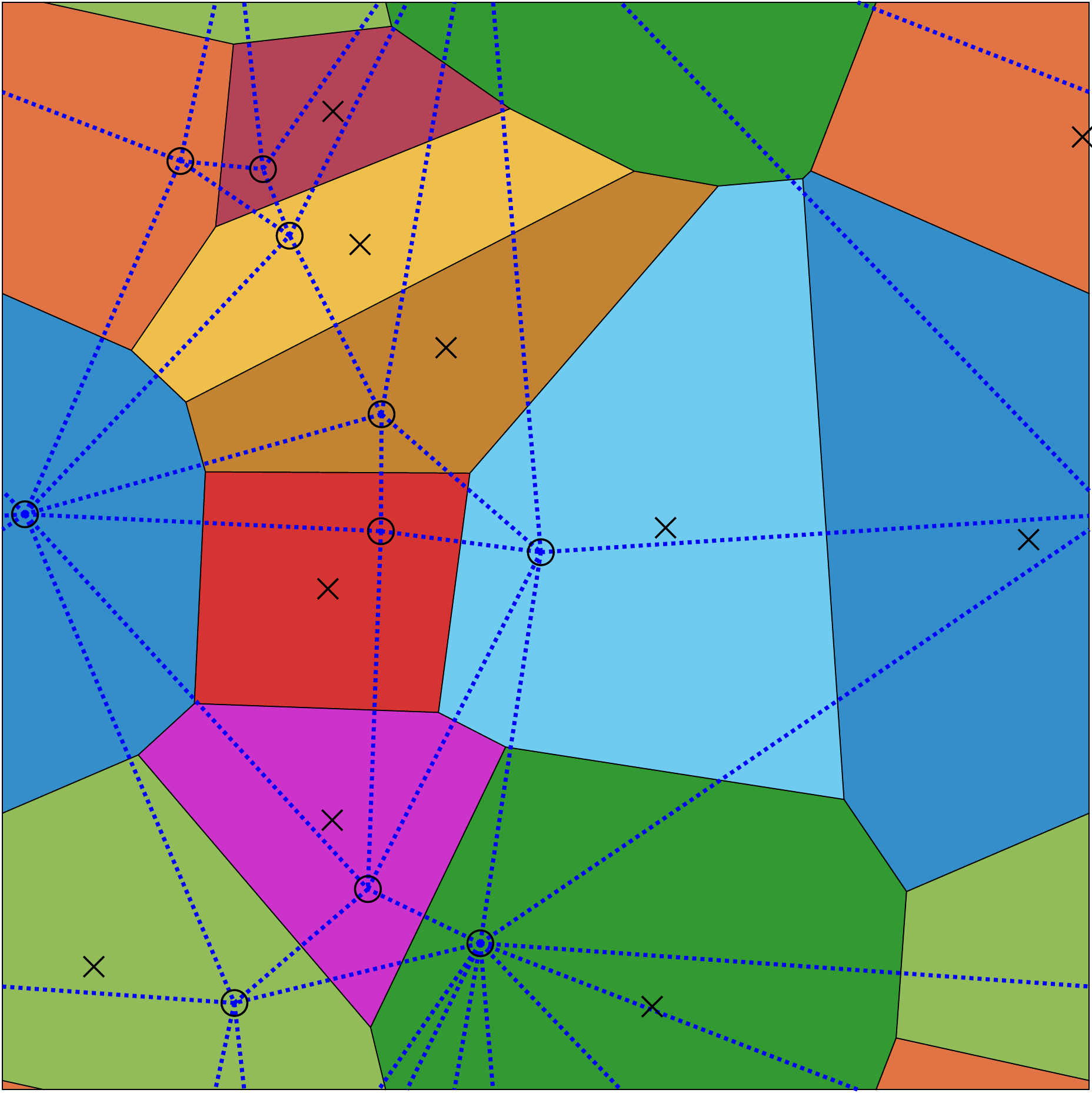} & \includegraphics[width=0.25\columnwidth]{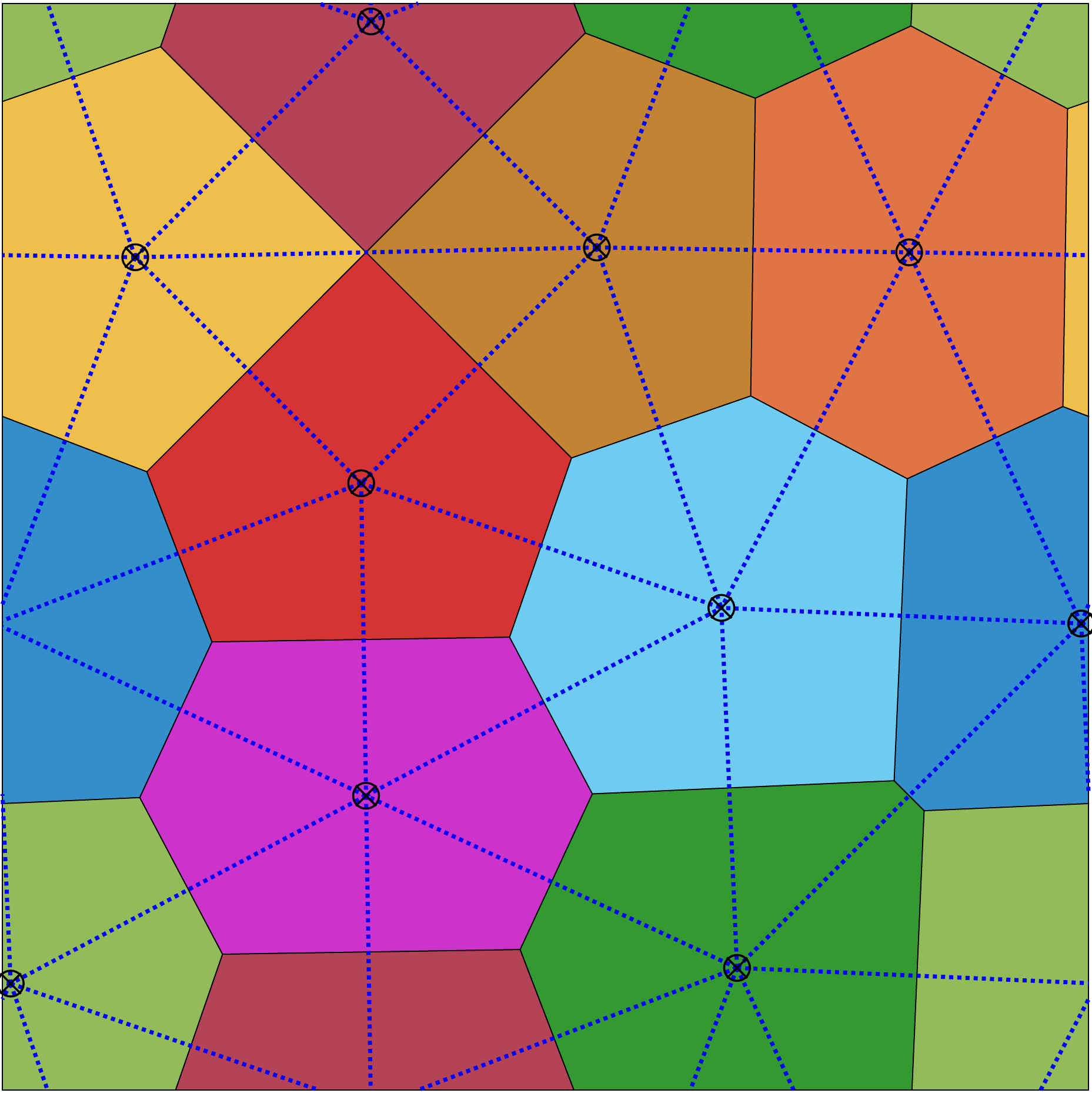}\tabularnewline
\footnotesize{(c) periodic-$\mathcal{V}(\mathbf{X})$} & \footnotesize{(d) periodic-CVT} \tabularnewline
\end{tabular}
\center

\caption{Two examples of the 2D framework for Voronoi tessellations with $N=10$,
the top row considers a bounded square $\Omega$
while the bottom row shows $\Omega$ as a square flat torus. \protect \\
Sets of generators $\mathbf{X}:=\{x_{i}\}_{i=1}^{N}$ are marked as
``$\circ$'' and centroids $c_{i}$ of the corresponding cells $V_{i}$
are marked as ``$\times$''. \protect \\
(a) Generic collection $\mathbf{X}$ and associated
Voronoi tessellation. (b) A centroidal Voronoi tessellation, i.e. generators $x_{i}$ and
respective centroids $c_{i}$ coincide $\forall\,i$. \protect \\
(c) The same sampling of generators found in (a) to emphasize
the changes in centroids and in the connectivity of the tessellation
when on the torus. (d) A periodic centroidal Voronoi
tessellations (PCVT). \protect \\
The dual graphs of the boundary sets $\cup_{i\leq N}\partial V_{i}$
(i.e. the Delaunay triangulations) are shown in dotted lines. Notice the superior regularity of the tessellation and of its corresponding
triangulation in both centroidal cases.\label{fig:Example Framework in 2D}}
\end{figure}

A fundamental problem in information theory and 
discrete geometry is known, respectively, as optimal quantization and optimal centroidal Voronoi tessellations (CVT); see \cite{CVT_AppAndAlg_Du,AdvancesOnCVT_AppAndAlg_Du,Quantization_Monograph} and \cite{Okabe} for an overview of the many concrete applications. Let us present the problem in its simplest form where in the underlying density is assumed to be uniform. Consider a bounded  $\Omega \subset  \mathbb{R}^d$ with a collection of $N$ distinct points  $\mathbf{X}:=\{x_{i}\}_{i=1}^{N} \subset \Omega$ referred to as \textit{generators}. The collection of points $\mathbf{X}$ gives rise to a 
Voronoi tessellation $\mathcal{V}(\mathbf{X})=\{V_i\}_{i=1}^N$ of $\Omega$ where
\begin{equation}
V_{i}(\mathbf{X})  \, = \, \{ y\in\Omega\,\big\vert\,||x_{i}-y||<||x_{j}-y||,\,\, \forall j\neq i\}\quad  i=1,...,N
\label{eq: def Voronoi cell}
\end{equation}
In other words,  Voronoi cell $V_i$  contains the points of $\Omega$ closer to $x_i$ than to any other generator.
We define for any $\mathbf{X}$ the non-local energy  
\begin{equation}\label{Vor-energy}
F(\mathbf{X})\,=\, \int_\Omega {\rm dist}^2(y,\mathbf{X})\, d y\, =\, \sum_{i=1}^N \int_{V_i} ||y-x_i||^2\, dy\,=:\sum_{i=1}^{N}F_{i}, 
\end{equation}
and minimize $F$ over the possible positions $\mathbf{X}$ of the generators. 
As seen later in \textsection \ref{sec:2}, criticality of this energy gives rise to a CVT; that is,  a placement of the generators $\{x_{i}\}_{i=1}^{N}$ such that 
they are exactly the centroids of their associated Voronoi cell $V_{i}$, i.e. 
\[
 x_{i}=c_i:=\frac{1}{|V_{i}|}\int_{V_{i}}y\,dy\qquad  i=1,...,N.
 \]
In the context of information theory and particularly in vector quantization \cite{Quantization_Monograph}, the set $\mathbf{X}$ is viewed as a {\it quantizer} that discretely models data; here, uniformly  distributed over $\Omega$. The {\it quantization error} is given by $F(\mathbf{X})$ with the {\it optimal quantizer} being the minimizer $\mathbf{X}^*$ with the least error (alternatively {\it the CVT with lowest energy}). 
Figure \ref{fig:Example Framework in 2D} contains some  visual examples to illustrate these notions. 

The energy (\ref{Vor-energy}) has a wealth of critical points (CVTs) \cite{Global_MCM_Lu,PCVTS_Du_Zhang} and  {\it low energy} CVTs have tiny basins of attractions making them difficult or impossible to find via gradient based descent with random initializations.  Finding optimal or at least low energy CVTs with desirable geometric properties is of fundamental importance in many concrete applications, c.f. \cite{CVT_AppAndAlg_Du,AdvancesOnCVT_AppAndAlg_Du,Quantization_Monograph}, for example: 1) spatial optimization, 2) mesh generation and numerical analysis, 3) vector quantization.
To this end,  one does have a benchmark for the optimal geometry in the limit $N \to \infty$ wherein we dispense with shape and boundary effects. Indeed, Gersho's conjecture \cite{OptiQuant_Gersho} addresses the periodic structure of the optimal quantizer as $N \to \infty$. 
The conjecture is completely solved in 2D wherein the optimal Voronoi cell is the regular hexagon, corresponding  to  generators on a triangular lattice. However, to date 
it remains open in 3D wherein the belief is that 
the optimal Voronoi {\it period-cell} is the truncated octahedron, corresponding to generators on a body centred cubic (BCC) lattice; see \cite{Gersho3D_Lattice_Barnes, BoundsGersho_Choksi,On3DGersho_Du}. 

The purpose of this article is to present and assess in 2D a simple deterministic method for efficiently navigating the energy landscape in order to {\it access} low energy CVTs which are otherwise inaccessible with Monte Carlo initializations coupled to gradient based descent methods (i.e. taking optimal results amongst hundreds of thousands to millions of gradient based descents on randomly sampled initial configurations). 
The proposed method has two parameters, namely the \textit{preconditioning number} $K$ and the \textit{probing number} $Q$, and is based upon a very simple scheme: 
each generator \textbf{moving away from the closest neighbor} (\textit{MACN})
by a certain distance (displacement); i.e., the basic \textit{MACN} scheme moves each generator in the opposite direction from its closest neighbor in $\mathbf{X}$. For  this scheme we consider  two choices of displacement:

\smallskip
\begin{itemize}
\item The individual displacement $||x_i-c_i||$ of $x_i$ to the centroid of its Voronoi cell. Our scheme using this displacement is labelled as the \textit{MACN-c} step.

\smallskip

\item A fixed displacement $\delta$ for all generators where
\begin{equation}
\delta \,: = \, \frac{1}{4}\sqrt{\frac{|\Omega|}{N}}, 
\label{eq: original delta}
\end{equation}
is set in terms of the average distance between generators
in the approximate regular hexagonal lattice. Our scheme with this distance is labelled as the  \textit{MACN-$\delta$} step.\\

\end{itemize}

\noindent Our method is then a \textit{coupling} procedure repeating  three steps. Starting from a random initialization (placement of generators in $\Omega$): 
\medskip
\begin{itemize} 
\item[Step \, 1] Iterate the \textit{MACN-c} step $K$ times as a \textbf{preconditioning};
\smallskip
 \item[Step \, 2] Implement Lloyd's, or any other deterministic method surveyed in \textsection \ref{sec:3},  minimizing (\ref{Vor-energy}) to a CVT;
 \smallskip
 \item[Step \, 3] A single \textit{MACN}-$\delta$ step working as a geometric-based \textbf{"annealing"} to disrupt the CVT from Step 2;
 \smallskip
 \item[Step \, 4] Repetition of Steps 1-3 a total number of $Q$ times, ending at Step 2.\\
 
\end{itemize}

Step 1 does result in a configuration which is \textit{close} to a CVT. However, this is not the point: it results in a configuration which lies in the basin of an {\it energetically desirable} CVT. 
By {\it energetically desirable} we mean two things: (i) it has low energy in the sense that its energy is comparable with the optimal result of any standard gradient-based descent algorithm assessed over a ``large" number of runs (for us with $N \sim 1000$, a large number of runs is on the order of $100,000$); (ii) the same holds true for the measures of regularity described in \textsection \ref{sec:2}. 
Step 2 achieves this {\it energetically desirable} CVT. 
Step 3  breaks away from this basin to another basin which can contain a more optimal CVT. In the end, one run of our method will have probed a number $Q$ of 
potentially different basins of attraction and their corresponding CVTs.

The choice of $\delta$  in (\ref{eq: original delta}) is subtle: it is large enough in order to change basins but sufficiently small in order to not lose the desired regularity achieved thus far. In 2D,  
the optimal configuration is partial to $N$ regular hexagonal Voronoi cells. In this case, $\delta$ can be thought of as half the distance from the generator to its Voronoi cell boundary plane. Indeed, in the case of regular hexagonal Voronoi cells, the 6 neighboring generators are equally distant. In this scenario, Step 3 chooses one of the 6 closest generators according to some pre-established tie breaking rule.

The probing number $Q$ encapsulates a degree of freedom common to all coupling methods (see \textsection \ref{sec:3}). The preconditioning number $K$ is chosen to remain constant over the $Q$ stages. This simple choice of constant $K$ is shown to produce excellent results on 2D domains. Moreover, preliminary implementation with a constant $K$ has also produced similar results on the sphere. The adaptation of the \textit{MACN} method to the sphere and possibly other manifolds is, however, outside the scope of the present paper and will be address in subsequent work.
On the other hand, one could devise several strategies --with higher degrees of freedom-- for a generalized sequence $\{K_q\}_{q=0}^{Q-1}$; for example taking an initial value $K_0$ combined with a suitable decaying profile adapted to both present and past history of the energy along the dynamics. 
Such variants could potentially improve convergence and energy efficiency of the method. It is, however, difficult to find consistent strategies that will systematically render desirable results across different problem sizes.  
Our crude tuning using a constant preconditioning number does render desirable results across different domains and problem sizes. (c.f. \textsection \ref{sec:5}).

\begin{figure}[H]
\begin{centering}
\begin{tabular}{cc}
\includegraphics[width=0.3\columnwidth]{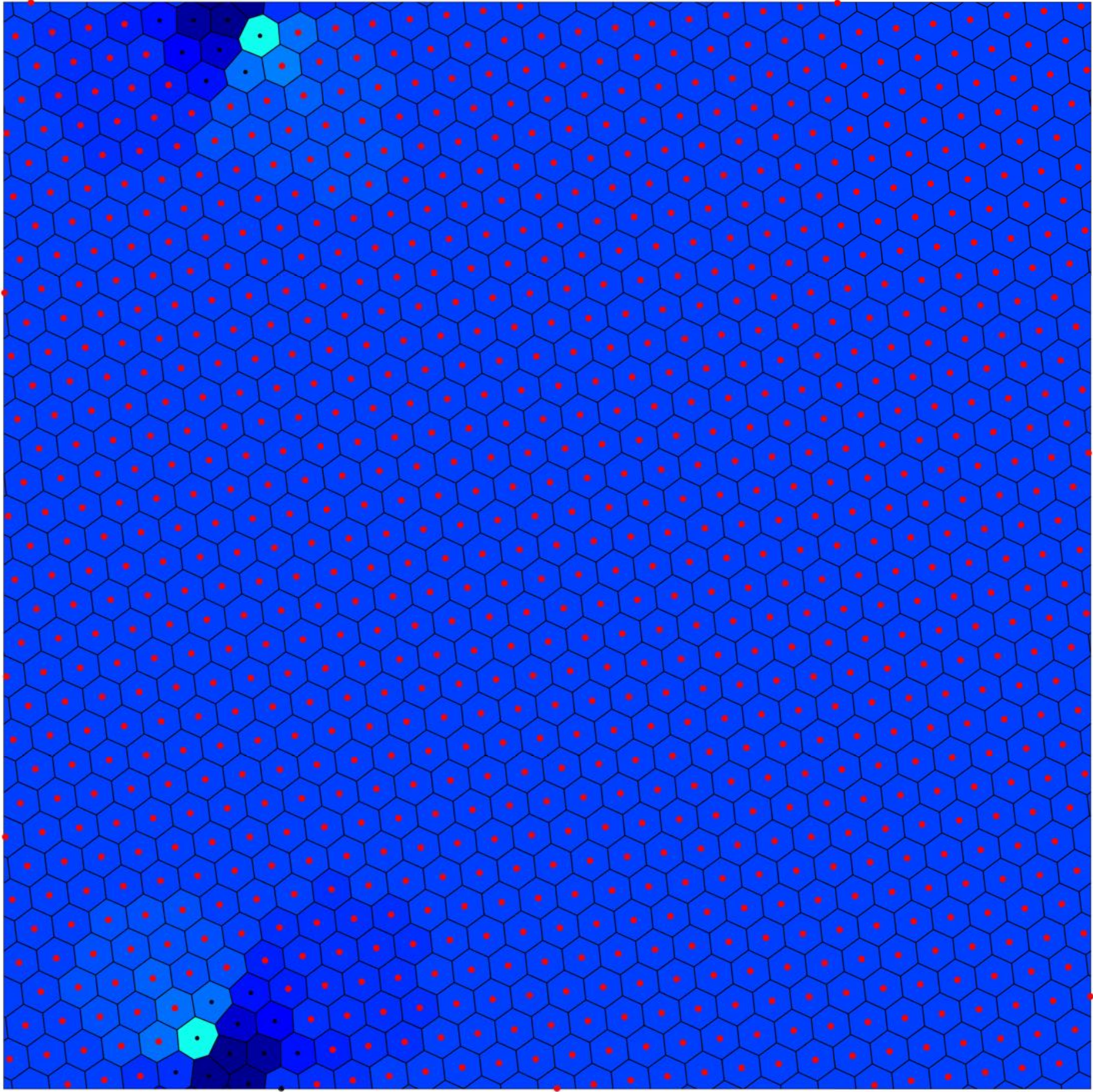} & \quad\includegraphics[width=0.3\columnwidth]{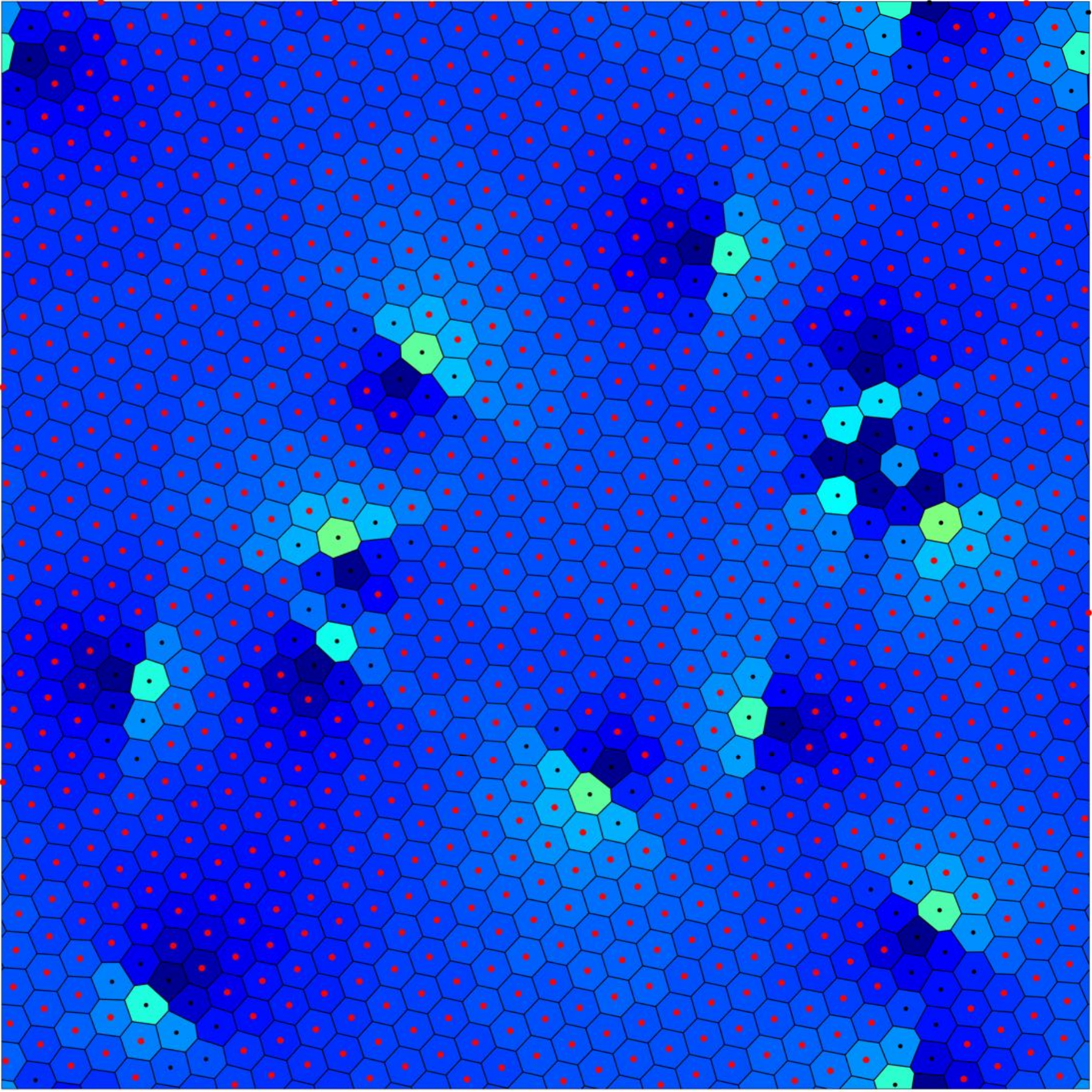}\tabularnewline
\footnotesize{lowest energy PCVT sampled} & \footnotesize{lowest energy PCVT sampled}\tabularnewline
 \footnotesize{with our hybrid method} & \footnotesize{with gradient based algorithms}
\end{tabular}
\par\end{centering}
\centering{}\caption{Case $N=1000$ on the square torus; comparative performance of our hybrid method with $Q=10$ (left)
vs. the lowest energy configuration obtained during our statistical
sampling of the landscape with
selected gradient based methods (right). 
The PCVT on the left is energetically 548\% closer to the non-achievable
regular hexagonal lattice than the one on the right. Full statistical
detail on this comparison is found in Figure \ref{fig: data N1000} and Table \ref{tab: N=1000}.
The color of each $V_i$ scales with $F_i$ to indicate local energy contribution, the reader is referred to \textsection \ref{sec:2} and subsequent Figures for further detail on the color map.}
\label{fig: First Hand comparison Hybrid vs GL-PLBFGS}
\end{figure}

We implement and assess our method on two choices of 2D domains, both having periodic boundary conditions in order to dispense with boundary effects. To start we work on a primary domain for which the ground state is known:  the periodic regular hexagon which can be tessellated into $N$ regular hexagons provided $N$ is suitably chosen \cite{PerimMinim_PentaTiling_Chung}. We then work on a primary domain that does not permit a perfect regular tilling: the flat square torus. Here there is always frustration due to the size effects, and it is surprisingly unclear as to the true nature of the lowest energy state. 
With the number of generators $N$ taken from the range 1000 to 4000,  we show that 
our hybrid method with $Q \leq 10$, implemented with less than two dozens initializations, readily finds states with far lower energy 
(and other metrics of optimal regularity) than the ones accessible with Lloyds' method, or any state of the art gradient 
based descent method,  assessed over hundreds of thousands random initializations, see Figure \ref{fig: First Hand comparison Hybrid vs GL-PLBFGS}. In \textsection \ref{sec:5} and \textsection \ref{sec:6} we present the full details of this comparison, emphasizing the role of statistics for assessing ours and other numerical methods. \\

The motivation and scope of our hybrid method is twofold. First, there is the direct goal of generating low energy CVTs which has an impact in many applications such as the ones already presented above. In contrast to methods which are based upon initial sampling or building on regular hexagons, our method is based only on basic structures of the energy (centroids and distance functions). Moreover, on domains like the square torus, optimality is subtle with regard to the presence of non-regular hexagons and defects (non-hexagonal cells): indeed, there is probably not even one perfectly regular hexagon in the Voronoi diagram of the optimal energy configuration.

There is a second scope to our work.  Probing non-convex and non-local energy landscapes is a fundamental problem in physics and applied mathematics.  
Deterministic strategies for navigating such landscapes are far and few. The CVT energy (\ref{Vor-energy}) is perhaps the simplest non-trivial example of such a landscape; while it is finite dimensional with a simple geometric characterization of criticality (namely a CVT), it is challenging (even in 2D) to navigate the landscape of CVTs.  As such, it perfects an ideal problem to address deterministic navigation. 
Our results demonstrate that while the CVT energy landscape on the square torus with $N\sim 1000$ is  indeed complex, our hybrid algorithm is able to efficiently navigate it with only 9-10 deterministic  "annealing" steps. 

We emphasize that our \textit{MACN} algorithm it designed for navigating the CVT energy landscape. We make no claim that our algorithm has any direct applications for navigating general non-convex energy landscapes. In fact, even for the closest energetic models of crystallographic particle interactions \cite{Crystallization_Conjecture}, it is unclear whether or not our work has any application.

In light of our two scopes we would like to point out that there is another method able to successfully probe the landscape, namely the global Monte Carlo method from \cite{Global_MCM_Lu} (cf. \textsection \ref{sec:3}). 
However, the novelty and change of paradigm here is that  our method, motivated entirely by geometry, 
 probes the energy in a  completely deterministic way.
Moreover,  our method shows a probing improvement compared to the global Monte Carlo. Precisely, the number of computed  PCVTs needed to reach low energy configurations is lower with the \textit{MACN} algorithm than with this alternative method. We elaborate on this comparison in \textsection \ref{sec:6}. 
\medskip

Finally let us remark on the empirical nature of our work. While we give a heuristic rationale for our \textit{MACN} steps, the precise nature of the distance chosen is based in part on empirical tries. For the \textit{MACN}-$\delta$ step, we experimented with other choices of distance, for example the intrinsic length-scale of the Voronoi cell (c.f. \textsection \ref{sec:6}). Our choice is one with sensibly the most effective performance. 
Overall, while we only have heuristics to explain certain aspects of our method, we feel the remarkable results justify its presentation and discussion here.

The paper is organized as follows: in \textsection \ref{sec:2} we give a brief description of periodic centroidal Voronoi tessellations (PCVTs). We also discuss, besides the normalized energy, certain natural {\it local measures of regularity}, one of which is novel. Then in \textsection \ref{sec:3}, we survey methods to generate and improve CVTs. We present our hybrid algorithm in \textsection \ref{sec:4} and then an analysis for its predictions in \textsection \ref{sec:5}. Later in \textsection \ref{sec:6} we discuss the advantages and disadvantages of our method compared to stochastic alternatives and then finish with closing remarks and future directions in \textsection \ref{sec:7}.\\

\section{Periodic CVTs and Regularity Measures\label{sec:2}}

For the remaining of the paper we will consider 2D torii spaces; namely polygons $\Omega$, with opposite and pairwise identifiable sides, that can be periodically extended to the plane. With the enforcement of these boundary conditions follow the definitions of a periodic-VT (PVT) and periodic-CVT (PCVT) by using the metric $||.||_{\mathbb{T}}$ that is inherited from the Euclidian metric in the above definitions (\ref{eq: def Voronoi cell}) and (\ref{Vor-energy}) respectively. See Figure \ref{fig:Example Framework in 2D} (c) and (d) for a concrete example.
In addition, the dual graph of the set  $\cup_{i\leq N}\partial V_i$ known as the periodic-Delaunay triangulation (PDT) will be of importance. The PDT has an edge $e_{ij}$ linking $x_i$ and $x_j$ iff they are neighbors in the PVT, in other words the index sets $\mathscr{N}_{i}:=\{j\ne i\,|\,\partial V_{j}\cap\partial V_{i}\neq\emptyset\};\;i=1,...N$ are intrinsic to the triangulation. 
We refer the reader to \cite{Okabe} for an introductory treatment on Delaunay
Triangulations and other variants of the Voronoi Tessellation and to \cite{PCVTS_Du_Zhang,2D_PCVTS_Levy,3D_DT_Caroli} for a detailed definition of the periodic framework.
For further detail on a popularized use of the Voronoi/Delaunay duality see for example \cite{Linear_SuperConv_SimpliMeshes,On3DGersho_Du}. There the notion of \textit{optimal triangulation} used in mesh generation and numerical analysis is defined as the dual of a CVT/PCVT as opposed to directly optimizing the connectivity of the triangulation according to some criteria \cite{VariationalTetrahedraMeshing,OptimalDelaunay}; thus tackling a continuous problem rather than a combinatorial one.

\subsection{Centroidal characterizations\label{subsec:2.2}}

In addition to the geometrical property $x_i=c_i\;\forall i$, a PCVT admits a variational characterization via $F:\mathcal{D}\to\mathbb{R}_{+}$  given by (\ref{Vor-energy}), here $\mathcal{D}:=\{\mathbf{X}\in\mathbb{R}^{2N}|\,x_{i}\neq x_{j}\;\forall\,j\neq i\,;\;x_{i}\in\Omega\;\forall\,i\}$ represents the set of non-degenerate configurations of generators in our periodic space.
Depending on the application the functional $F=\sum_{i=1}^N{F_{i}}$ is often referred to as the distortion, cost function or potential. Alternatively,
one can interpret each $F_{i}$ from a physical point of view as the trace of the $2\times2$ inertia tensor of $V_{i}$ with respect
to $x_{i}$, i.e. the resistance of $V_{i}$ to its rotation
around an axis passing through $x_{i}$ that is orthogonal to $\Omega$.

It is a well established fact in rigid mechanics (parallel axis theorem)
that the trace of this
tensor will be minimized whenever the orthogonal axis of rotation
passes through the centroid. This locally translates
to a Voronoi region being generated by its centroid
and therefore globally translates to a PCVT. Indeed, one can formally prove this: by the definition of $\mathcal{V}(\mathbf{X})$ with the
Euclidian metric we have according to \cite{GeoOptimi_Murota,NewtonLloyd_Du}
that ${\displaystyle \frac{\partial F}{\partial\mathcal{V}}(\mathbf{X},\mathcal{V}(\mathbf{X}))\equiv0}$
and thus Reynold's Transport Theorem gives
the gradient components
\[
D_{i}\,F(\mathbf{X}):=\frac{\partial F}{\partial x_{i}}(\mathbf{X})=2|V_{i}|(x_{i}-c_{i})\quad i=1,...,N
\]
which show that the gradient vanishing configurations of $F$ are
exactly PCVTs. Moreover a non-trivial second application of the theorem
yields the entries of the Hessian $D^{2}F(\mathbf{X})$ 
\begin{equation}
\begin{cases}
\frac{\partial^{2}F}{\partial x_{i}^{(m)}\partial x_{i}^{(m)}}=2|V_{i}|-\sum_{j\in\mathscr{N}_{i}}\frac{2}{||x_{i}-x_{j}\text{||}}\int_{\partial V_{i}\cap\partial V_{j}}\scriptstyle{(x_{i}^{(m)}-y^{(m)})^{2}\,dy}\\

\frac{\partial^{2}F}{\partial x_{i}^{(m)}\partial x_{i}^{(l)}}=-\sum_{j\in\mathscr{N}_{i}}\frac{2}{||x_{i}-x_{j}\text{||}}\int_{\partial V_{i}\cap\partial V_{j}}\scriptstyle{(x_{i}^{(m)}-y^{(m)})(x_{i}^{(l)}-y^{(l)})\,dy} & \scriptstyle{m\neq l}\\

\frac{\partial^{2}F}{\partial x_{i}^{(m)}\partial x_{j}^{(l)}}=\frac{2}{||x_{i}-x_{j}\text{||}}\int_{\partial V_{i}\cap\partial V_{j}}\scriptstyle{(x_{i}^{(m)}-y^{(m)})(x_{j}^{(l)}-y^{(l)})\,dy }& \scriptstyle{j\in\mathscr{N}_{i}}\\

\frac{\partial^{2}F}{\partial x_{i}^{(m)}\partial x_{j}^{(l)}}\equiv0 & \scriptstyle{j\neq i,\,j\notin\mathscr{N}_{i}}
\end{cases}\label{eq: Hessian D2E}
\end{equation}
Furthermore, $F$ has been shown to be $C^{2}(\mathcal{D})$ \cite{CVT_EnerSmooth_Levy,PCVTS_Du_Zhang} and although the Hessian is sparse it is in general not definite,
as noted in \cite{Global_MCM_Lu} this is a consequence of the energy
functional being highly non-convex with the presence of saddle points, see also \cite{PCVTS_Du_Zhang}.

\subsection{Optimal PCVT, ground state energy and regularity measures\label{subsec:2.3}}

Based upon the \textit{hexagon theorem} \cite{TothProof_by_Gruber,HexTheor_Newman,2DGersh_Toth}, and in a similar spirit to \cite{SecondMoments_Conway,HierachicalCVT_Wang}, we scale the energy functional with respect to $N\,F_{hex}$.
Here $F_{hex}$ denotes the second moment of a regular hexagon with
area $|\Omega|/N$ and a simple calculation shows that 
$$F_{hex}={\displaystyle \frac{5}{18\sqrt{3}}\frac{|\Omega|^{2}}{N^{2}}}\;.$$
The scaled energy that is thus independent of the size of the problem is
\[
E(\mathbf{X}):=\frac{F(\mathbf{X})}{N F_{hex}}=\frac{1}{N}\sum_{i=1}^{N}\frac{F_{i}}{F_{hex}}=:\frac{1}{N}\sum_{i=1}^{N}E_{i}
\]
and with the same scaling carrying out trivially to $DE(\mathbf{X})$ and $D^{2}E(\mathbf{X})$.

As an additional way to help us quantify the quality of a PCVT, and
as it has been customary in the literature, we use the fraction
of hexagonal cells:
\[
H(\mathbf{X}):=\frac{\#\{V_{i}\,,\,i=1,...,N\;|\;\#\mathscr{N}_{i}=6\}}{N}
\]
However, regularity measures such as $H$ or the Voronoi Entropy used
in \cite{Entropy_Bormashenko} (and rediscussed in \cite{Entropy_ReDiscussed_Bormashenko})
only give us information about the connectivity of the DT and make
no direct connection with the hexagon theorem requiring hexagons to
be regular to achieve the energy $E=1$. For this reason we create
a refinement on $H$ that takes into account the regularity of the hexagonal cells through
their isoperimetric ratio. We recall that the isoperimetric ratio $r$
of a polygon is the dimensionless ratio of its perimeter squared and
area. Using $r(V_{i})$ of each Voronoi cell together
with the isoperimetric ratio of the regular hexagon $r_{hex}:=8\sqrt{3}$, and for given $\epsilon>0$ small
we introduce the following regularity measure of a tessellation
$\mathcal{V}(\mathbf{X})$ :

\[
R_{\epsilon}(\mathbf{X}):=\frac{\#\{V_{i}\,,\,i=1,...,N\;|\;\#\mathscr{N}_{i}=6\;;|1-\frac{r(V_{i})}{r_{hex}}|\leq\epsilon\}}{N}
\]
That is $R_{\epsilon}$ is the fraction of cells that are
regular hexagons within an isoperimetric tolerance of $\epsilon$. Generators whose cell is
taken into account by $R_{\epsilon}$ will be depicted in red in the
PVTs of Figure \ref{fig:Example of MACN-c dynamics} and onwards.\\

Later in \textsection \ref{sec:4}-\ref{sec:6} we will rely on $E-1$
, $H$ and $R_{\epsilon}$ to measure the performance of our method and others in reaching high quality PCVTs. The data will also make the compelling case that $R_{\epsilon}$ is a more sensitive measure of
low values of $E$ than $H$ (for appropriate $\epsilon$).
Finally the statistics on these three quantities will be at the core
of the conjecture that our method acts at a semi-global scale.

\section{Generating and improving CVTs\label{sec:3}}

We divide this brief literature survey into methods available to compute CVTs and ways to enhance their quality.

\subsection{Computing CVTs\label{subsec:3.1}}

These methods rely either on the characterization $x_i=c_i\;\forall i$ and seek to solve these nonlinear equations or are of variational character on $E(\bf{X})$:\\

\noindent 1. Lloyd's method, introduced in the seminal work \cite{LloydAlg_Lloyd}, is unquestionably the most widespread method due to its simplicity. It iteratively applies the map $\mathbf{T}:\;\mathcal{D}\rightarrow\mathbf{\mathcal{D}}$
defined by $\mathbf{T}(\mathbf{X})=\mathbf{C}$ where $\mathbf{C}:=\{c_{i}\}_{i=1}^{N}$ is the collection of centroids of $\mathcal{V}(\mathbf{X})$.
We refer the reader to
\cite{CVT_AppAndAlg_Du,LloydAlgConvergence_Du,NonDegen_LloydConvergence_Emelianenko}
for properties and analysis of the map.\\
Due to the importance of the method in the remaining of this paper we provide explicit pseudo-code as a sub-routine in Algorithm \ref{alg:hybrid-algorithm Q stages}.\\

\noindent 2. McQueen's and probabilistic Lloyd's \cite{ProbabilisticMcQueen_Du},  generate small sets of random sampling points  in $\Omega$ and use these to approximate the centroid of
Voronoi regions. \\

\noindent 3. Lloyd-Newton \cite{NewtonLloyd_Du} uses Newton's root finding method on $\mathbf{S}(\mathbf{X}):=\mathbf{X}-\mathbf{T}(\mathbf{X})$ after some Lloyd steps to reach the vicinity of a root in $\mathcal{D}$. The authors also propose multilevel based extensions: one using an algebraic multigrid preconditioner combined with Gauss-Sidel iterations to solve the linear system involved in each Newton step and another using a nonlinear multigrid approach in solving $D E=0$. See also \cite{FastmultilevelCVT} for further development. However, as discussed in \cite{CVT_EnerSmooth_Levy}, unstable CVTs may be produced (i.e. saddle points of $E$).\\

\noindent 4. Newton's classical technique that minimizes Hessian-based quadratic models
of $E\in C^{2}(\mathcal{D})$. With (\ref{eq: Hessian D2E}) available, the method converges at least quadratically when coupled with a line-search ensuring the
strong Wolfe conditions \cite{NumOptimization_Nocedal}. Nonetheless
this method suffers from two downsides: i) $D^{2}E$ is often indefinite and needs
to be altered, for example by adding  a ``small''
matrix to render it SPD prior to executing
an incomplete Cholesky factorization, see \cite{ModCholeskyDecomp_McSweeney,IncomplCholesky_LimitMemo_More} for theory and low memory algorithms. ii) the Hessian is expensive to populate due to the boundary integrals. \\

\noindent 5. Quasi-Newton BFGS collection. These methods only  use $E$ and $DE$ to give an
iterative approximation of the inverse Hessian. They remain to this
day the favored methods in the literature for fast CVT/PCVT computation due to their expected super-linear convergence whenever they are coupled with a line-search
method that ensures the strong Wolfe conditions. The two families suited for medium/large scale problems that will be used in this paper are:\\

\begin{enumerate}

\item[(a)] the low memory L-BFGS($M$) in which the inverse Hessian approximation remains
sparse and is computed recursively from the $M$ previous
approximations. \\

\item[(b)] the Preconditioned-L-BFGS($M$,$T$) uses, every modulo $T$
iterations, a SPD preconditioner matrix $\widetilde{A}$ (that does not necessarily need to approximate the Hessian) with the goal of redirecting the algorithm to a more suitable
descent direction. \\
\end{enumerate}

\noindent See \cite{NumOptimization_Nocedal} for a thorough
description and analysis of the classical BFGS and L-BFGS($M$) and \cite{PLBFGS_MolecularEner_Jiang,CVT_EnerSmooth_Levy}
for the use P-L-BFGS($M$,$T$). Finally, an explicit routine of the preconditioned algorithm is provided
in Appendix \ref{App:C}. \\

\noindent 6. Non-Linear Conjugate Gradient (NLCG) methods generalize the classical CG used in quadratic programming. Several updates for the conjugate directions are available \cite{NLCGsurvey_Hager,NumOptimization_Nocedal}. One can also use a relevant preconditioner SPD matrix $\widetilde{A}$ to improve the descent-conjugate directions.\\

There are several preconditioner matrices $\widetilde{A}$, for example: the Hessian $D^2E$ itself (along with the often necessary modification to make it SPD) or a Graph Laplacian $G$ introduced in \cite{GraphLap_LongChen} whose purpose is to approximate $D\,E(\mathbf{X})=0$ to first order by the matrix equation $G\mathbf{X}=0$.\\
Originally presented for the $\Omega$-bounded case, below is our adaptation of $G$ for the periodic case; denoting by $\tau_{i,j}$ the pyramid with base $\partial V_{j}\cap\partial V_{i}$ and apex $x_{i}$
then the $N\times N$ matrix is given by 
\begin{equation}
G:=\begin{cases}
g_{ij}=-\int_{\tau_{i,j}\cup\tau_{j,i}}\rho(y)dy & \text{if }j\in\mathscr{N}_{i}\\
g_{ii}=\sum_{j\in\mathscr{N}_{i}}|g_{ij}|\\
0 & \text{otherwise}
\end{cases}\label{eq:Graph Laplacian def}
\end{equation}
However, contrary to the original
construction for $\Omega$-bounded, our adaptation
is symmetric and positive semi-definite ($G$ given by (\ref{eq:Graph Laplacian def})
is not strictly diagonal dominant and can thus be singular). As a
consequence we will need to use a modified Cholesky factorization of $G$.

\subsection{Energetic improvements\label{subsec:3.3}}

Next we survey initialisations and other methods used jointly with the above algorithms in order to improve the quality of CVTs.\\

\noindent \textit{Initialisations}
\vspace{3mm}
 
\noindent (a) A Greedy Edge-Collapsing initialization \cite{EdgeCollapseAlgo_Moriguchi}
that meshes $\Omega$ using more than $N$ vertices, then it repeatedly uses an edge-collapsing
scheme and finally the decimated vertices are employed as the initial
generators.\\

\noindent (b)  Quasi Random samplings of $\Omega$ reduce the discrepancy of the initial cloud of generators, i.e. the sampling of each site depends on the position of the others. These QR samplings use low discrepancy sequences such as Halton's, Hammersley's, Niederreiter's and Sobol's \cite{NiederreiterSampling,HammersleySampling,HaltonSampling}, see also \cite{HammerslyIni_Quinn} for CVT results.\\

\noindent \textit{Couplings}
\vspace{3mm}

\noindent (a) A Hierarchical method \cite{HierachicalCVT_Wang} that refines a CVT by cleverly inserting new generators over the DT. In this way the "regular" portions of the CVT "grow" when alternating with an energy descent method.\\

\noindent (b) An Atomic Operation method \cite{AtomicOperation_Men}, here the authors
establish three operations on ``defects'' that merge or split non-hexagonal
cells prior to minimizing the energy; the process is then repeated.\\

\noindent (c) A global Monte Carlo method \cite{Global_MCM_Lu} which applies ideas from Simulated Annealing. This method starts at a CVT and after a specific random perturbation
of generators, that is dependent on the size of each $V_{i}$, a new
CVT is obtained by a non-linear minimization method. If the new energy is lower than the previous one then the algorithm
automatically accepts that new configuration, otherwise it accepts
it only according to a transitional probability that is dependent
both on a cooling temperature and on the energy gap
between the two CVTs. \\
While the method has been proven to find the ground state in infinite time \cite{MonteCarlo_ProteinFolding}, a crucial disadvantage in practice is the number of parameters that need to be adjusted to obtain good
performances, namely: the initial temperature, the temperature decay
to zero, the perturbation amplitude and the number of iterations repeating the procedure.\\

Indeed, the above-mentioned literature shows that these initialisations and couplings yield significant energetic improvements, even more so when combined.

\section{Our hybrid algorithm\label{sec:4}}

This section describes in more detail our three stage method and subsequently provides an individual and deeper insight on the \textit{MACN-c} and \textit{MACN-$\delta$} dynamics respectively.

Given a set of generators $\mathbf{X}:=\{x_{i}\}_{i=1}^{N}$ and its associated $\mathcal{V}(\mathbf{X})$, a closest neighbor to $x_{i}$ is denoted $x_{j_{i}^{*}}$; i.e. $x_{j_{i}^{*}}$ solves $\displaystyle \min_{j\ne i}||x_{i}-x_{j}||=\min_{j\in\mathscr{N}_{i}}||x_{i}-x_{j}||$ where $\mathscr{N}_{i}$ was defined earlier to be the index set of Delaunay edge-connected neighbors to the site $x_{i}$. Notice that $x_{j_{i}^{*}}$ may not be unique however it suffices
for our implementation to simply pick one solution via a pre-established tie
breaking rule (e.g. selection of the closest candidate up to machine
precision).

Our \textbf{\textit{MACN}} scheme is a displacement of $x_i$ by a distance $d_i$ in the opposite direction to $x_{j_{i}^{*}}$, i.e. a \textit{MACN} iteration with distances $\{d_{i}\}_{i=1}^{N}$, consists on the update scheme
\begin{equation}
x_{i}\leftarrow x_{i}+d_i\,\frac{x_{i}-x_{j_{i}^{*}}}{||x_{i}-x_{j_{i}^{*}}||}\qquad i=1,...,N\label{eq: update formula MACN}
\end{equation}
 Here $d_i$ is a dummy quantity representing either the distance to the centroid $d_i:=||x_i-c_i||$ for a \textit{MACN-c} step or $d_i:=\delta=\, \frac{1}{4}\sqrt{\frac{|\Omega|}{N}}$ for a \textit{MACN-$\delta$} step --as in (\ref{eq: original delta}).\\
Explicit pseudo-code of the overall method comprising the 4 steps listed on page 3 is provided in Algorithm \ref{alg:hybrid-algorithm Q stages} for completeness.

\begin{algorithm}[H]
\begin{footnotesize}
\caption{\textit{MACN} algorithm over $Q$ stages with Lloyd subroutine\label{alg:hybrid-algorithm Q stages}}

\begin{algorithmic}
\State \textbf{Input:} 
1) initial generators $\mathbf{X}=\{x_{i}\}_{i=1}^{N}$; 2) probing number $Q$; 3) preconditioning number $K$; 4) $tol$ for convergence to a PCVT\\

\For {$q=0:Q-1$}\\

	\State I. \textit{MACN-c}:
	\For{$k=0:K-1$}
		\State 1. compute PDT($\mathbf{X}$), associated PVT($\mathbf{X}$) and centroids $\{c_i\}_{i=1}^N$
		\State 2. extract $\{\mathscr{N}_{i}\}_{i=1}^{N}$ and find $\{j_{i}^{*}\}_{i=1}^N$ solving $\min_{j\in\mathscr{N}_{i}}||x_{i}-x_{j}|| \; \text{for each}\,i$
		\State 3. set $d_i=||x_i-c_i||$ and update $\mathbf{X}$:  $x_{i}\leftarrow x_{i}\;\forall\,i$ with (\ref{eq: update formula MACN})
	\EndFor\\

	\State II. \textit{Reaching criticality}:
	\State {[}$\textbf{X}_{q}^{*}$, PDT($\textbf{X}_{q}^{*}$), PCVT($\textbf{X}_{q}^{*}$)]=Lloyd$(\textbf{X},tol)$\\

	\State III. \textit{MACN-$\delta$}:
	\If{$q<Q-1$}
		\State 1. extract $\{\mathscr{N}_{i}\}_{i=1}^{N}$ from PDT($\textbf{X}_{q}^{*}$) and find $\{j_{i}^{*}\}_{i=1}^N$ solving $\min_{j\in\mathscr{N}_{i}}||x_{i}-x_{j}|| \; \text{for each}\,i$
		\State 2. set $d_i=\delta$ and get new $\textbf{X}$ by $\delta$-perturbing $\textbf{X}_{q}^{*}$:  $x_{i}\leftarrow x_{i}\;\forall\,i$ with (\ref{eq: update formula MACN})
	\EndIf\\

\EndFor\\

\State \textbf{Output:} $\left\{\textbf{X}_{q}^{*}\right\}_{q=0}^{Q-1}$, a collection of stable local minimizers of \textit{E} and their corresponding PCVTs and PDTs.\\

\State \textbf{subroutine} [$\textbf{X}$, PDT($\textbf{X}$), $\mathcal{V}(\textbf{X})$]=Lloyd$(\textbf{X},tol)$\\
set \textit{diff}=Inf
	\While{\textit{diff > tol}}
		\State 1. compute $\mathcal{V}(\mathbf{X})$ (i.e. PVT($\mathbf{X}$)), centroids $\{c_i\}_{i=1}^N$, areas $\{|V_i|\}_{i=1}^N$ and set  \textit{diff}$=||D\,E||/N$
		\State 2. update $\mathbf{X}$:  $x_{i}\leftarrow c_{i}\;\forall\,i$
	\EndWhile
\State \textbf{end subroutine}\\

\end{algorithmic}
\end{footnotesize}
\end{algorithm}

As briefly mentioned in the Introduction, the $K$ initial \textit{MACN-c} iterations yields more evenly distributed sets of generators.  This preconditioning is common in outcome to other initialisations listed in \textsection \ref{sec:3}. Thus the success of \textit{MACN-c}+Lloyd in reaching "low" energy PCVTs is not surprising. Then by disrupting the centroidal configuration with \textit{MACN-$\delta$} and introducing the $Q-1$ supplementary repetitions we create a coupling technique involving a symbiosis of relaxation and contraction. By this we mean that the combination of our preconditioning and perturbation help Lloyd's algorithm in reaching lower energy states whilst, conversely, Lloyd's gets the system to PCVTs that our \textit{MACN} stages use as a stepping-stone to successfully probe the landscape.

We emphasize next that the sequence $\{E(\mathbf{X}_{q}^{*})\}_{q=0}^{Q-1}$ obtained with our method
may not be strictly decreasing. It is in fact quite likely that our algorithm moves to a higher energy basin of attractions from one stage to another, thus resembling Simulated Annealing in that sense (whilst remaining completely deterministic). The advantage however is, again, the low probing number needed. With $Q\leq10$ we are able to sample low energy states that are impossible or scarcely achievable by other deterministic
methods (even when using a stochastic sampling of the landscape done with a great
number of initial configurations),  c.f. \textsection \ref{sec:5}.

Concerning the complexity of our method:
most of the software used to construct VT/PVT with Euclidean distance
(e.g. the CGAL Computational Geometry Algorithms
Library --https://www.cgal.org \cite{CGAL_Library}-- for C++ or the
built-in \textsc{Matlab} function \texttt{voronoin}) rely on an
early construction of the Delaunay Triangulation. Thus one can extract $\mathscr{N}_{i}\;\,\forall\,i$
before the computation of the tessellation while keeping the same complexity.
Moreover, the remaining difference between one \textit{MACN-c} step
and one Lloyd step is the computation of the indices $j_{i}^{*}\;\forall\,i$
which simply adds a lower order term $O(N)$ to the optimal overall complexity $O(N\ln(N))$
in 2D \cite{Okabe}. Hence, our overall scheme (\ref{eq: update formula MACN}) benefits from the same low complexity
of Lloyd's algorithm as well as a comparable simple implementation. This is the reason we chose Lloyd's for our coupling rather than any other gradient based method.

\subsection{\textit{MACN-c}\label{subsec:4.1}}
\begin{figure}[H]
\centering

\begin{tabular}{ccc}
\includegraphics[width=0.23\columnwidth]{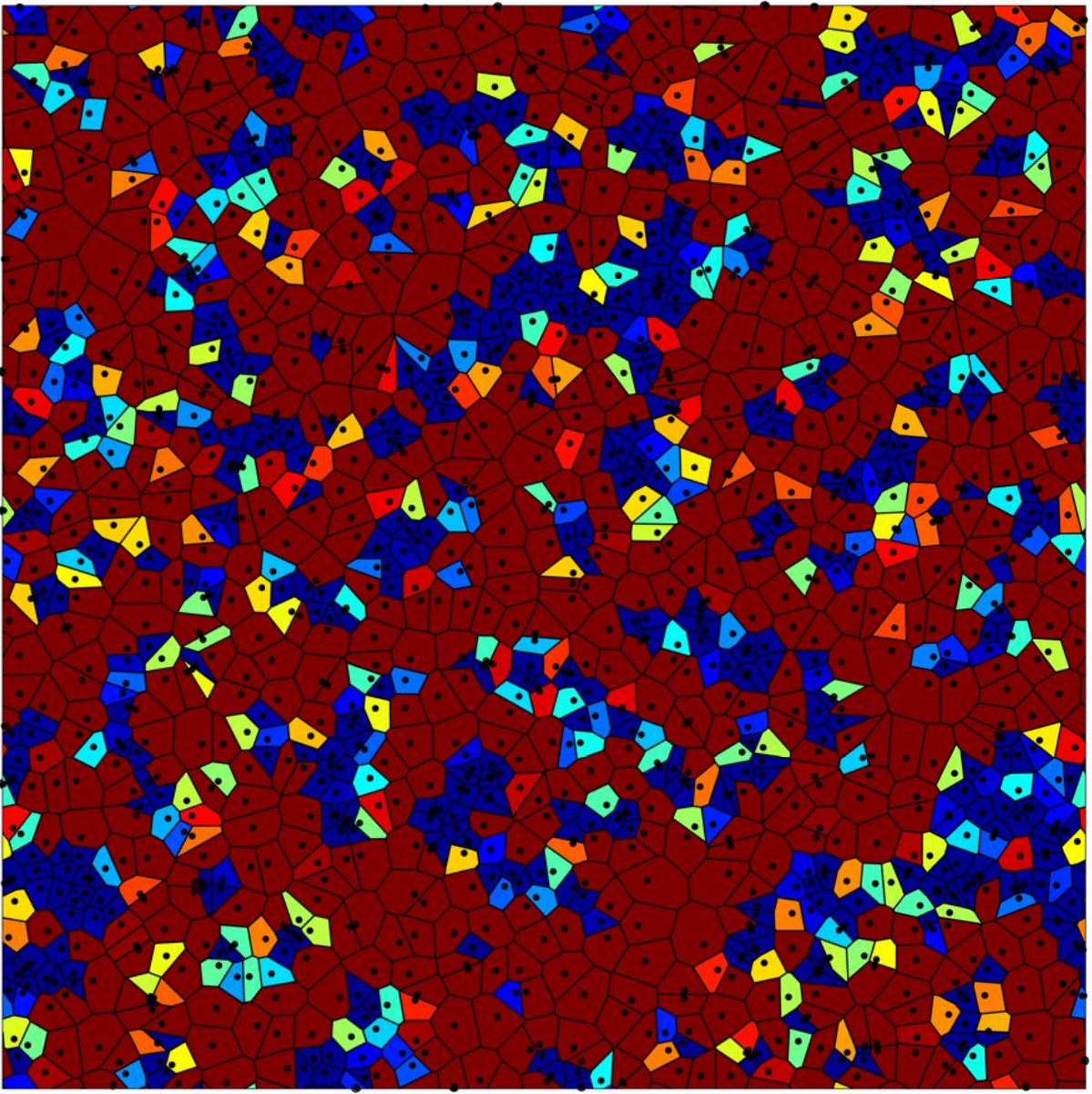} & \includegraphics[width=0.23\columnwidth]{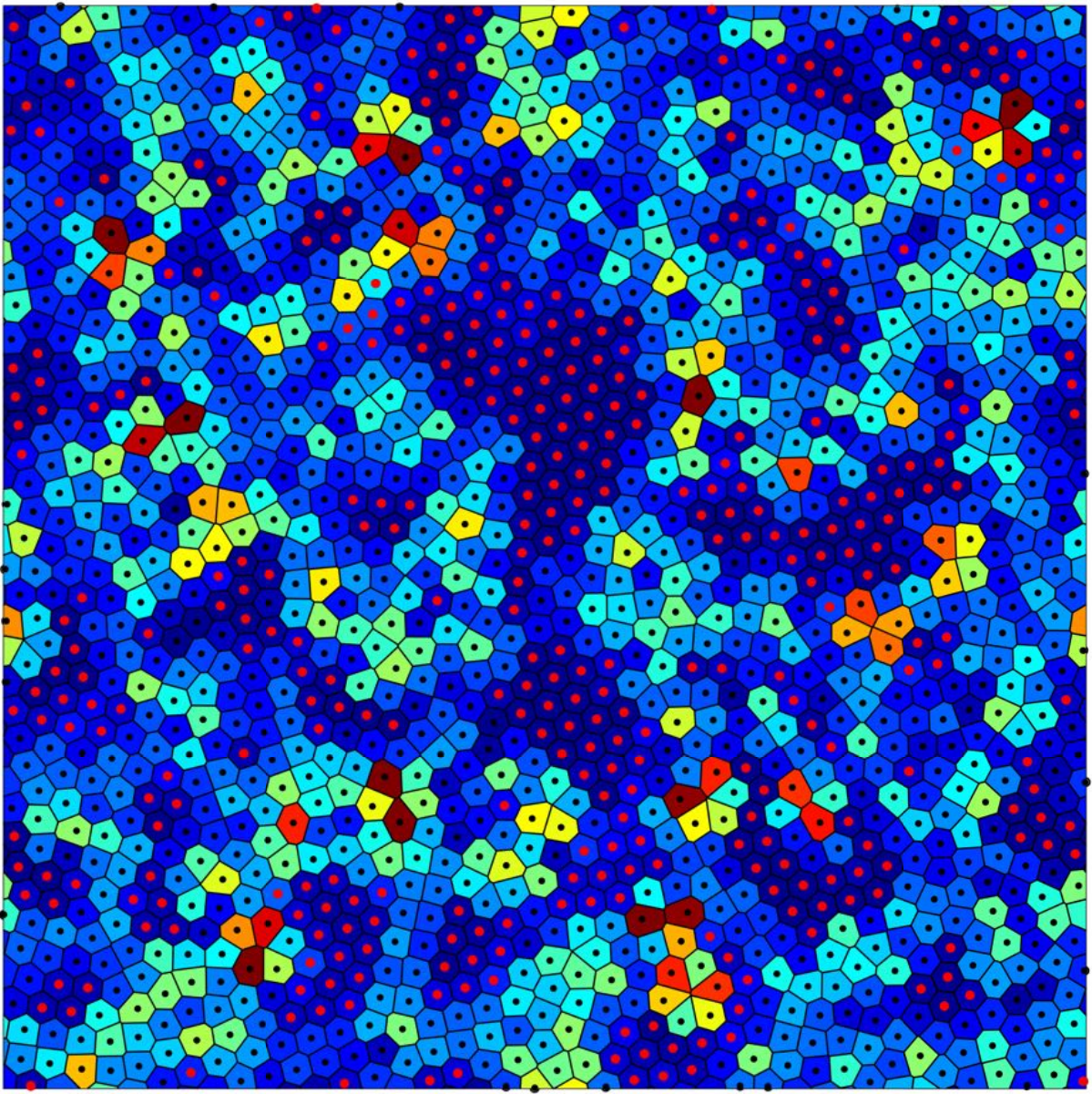} & \includegraphics[width=0.23\columnwidth]{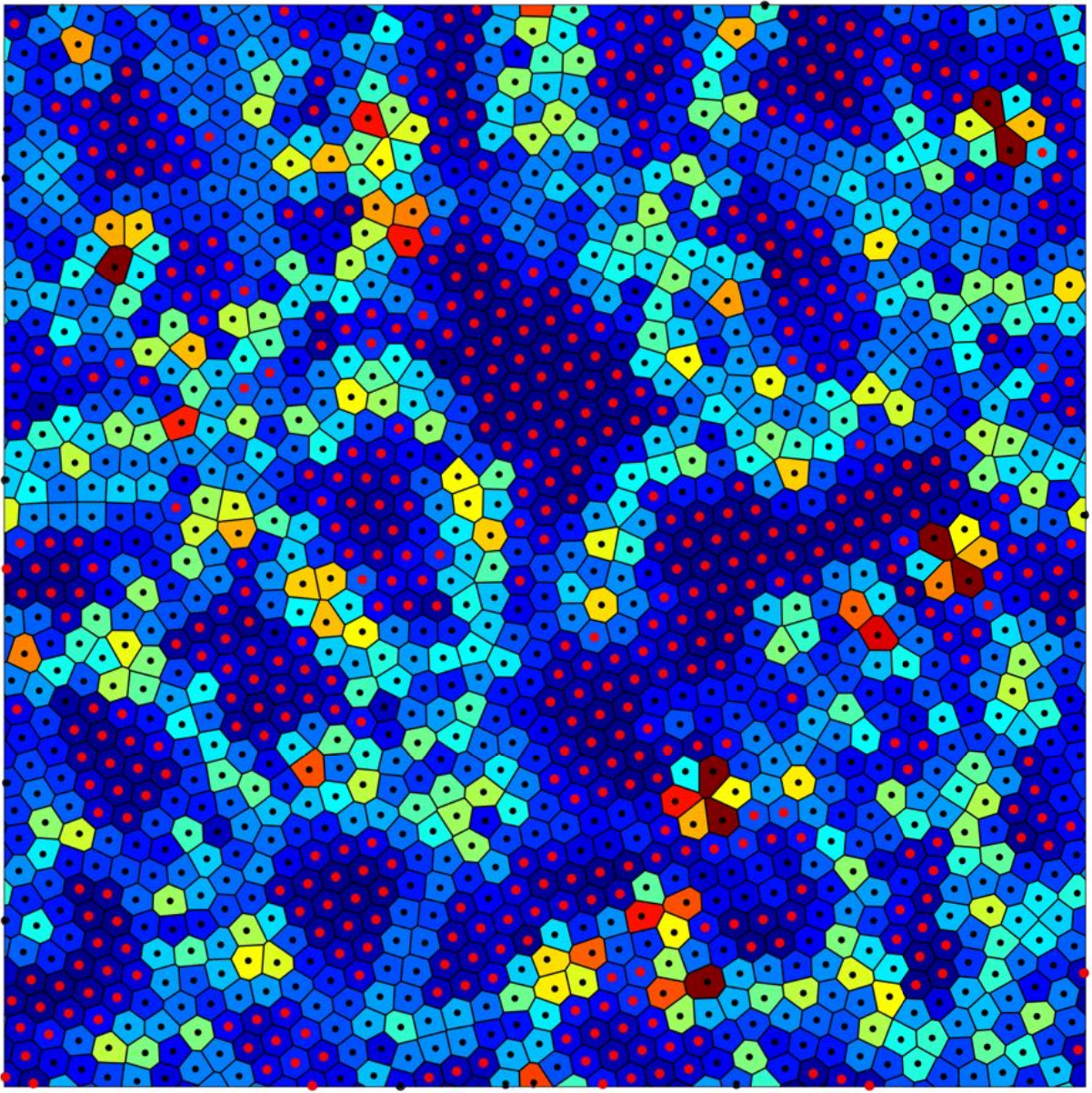}\tabularnewline
\footnotesize{(a) $k=0$} & \footnotesize{(b) $k=500$} & \footnotesize{(c) $k=1000$}\tabularnewline
\includegraphics[width=0.23\columnwidth]{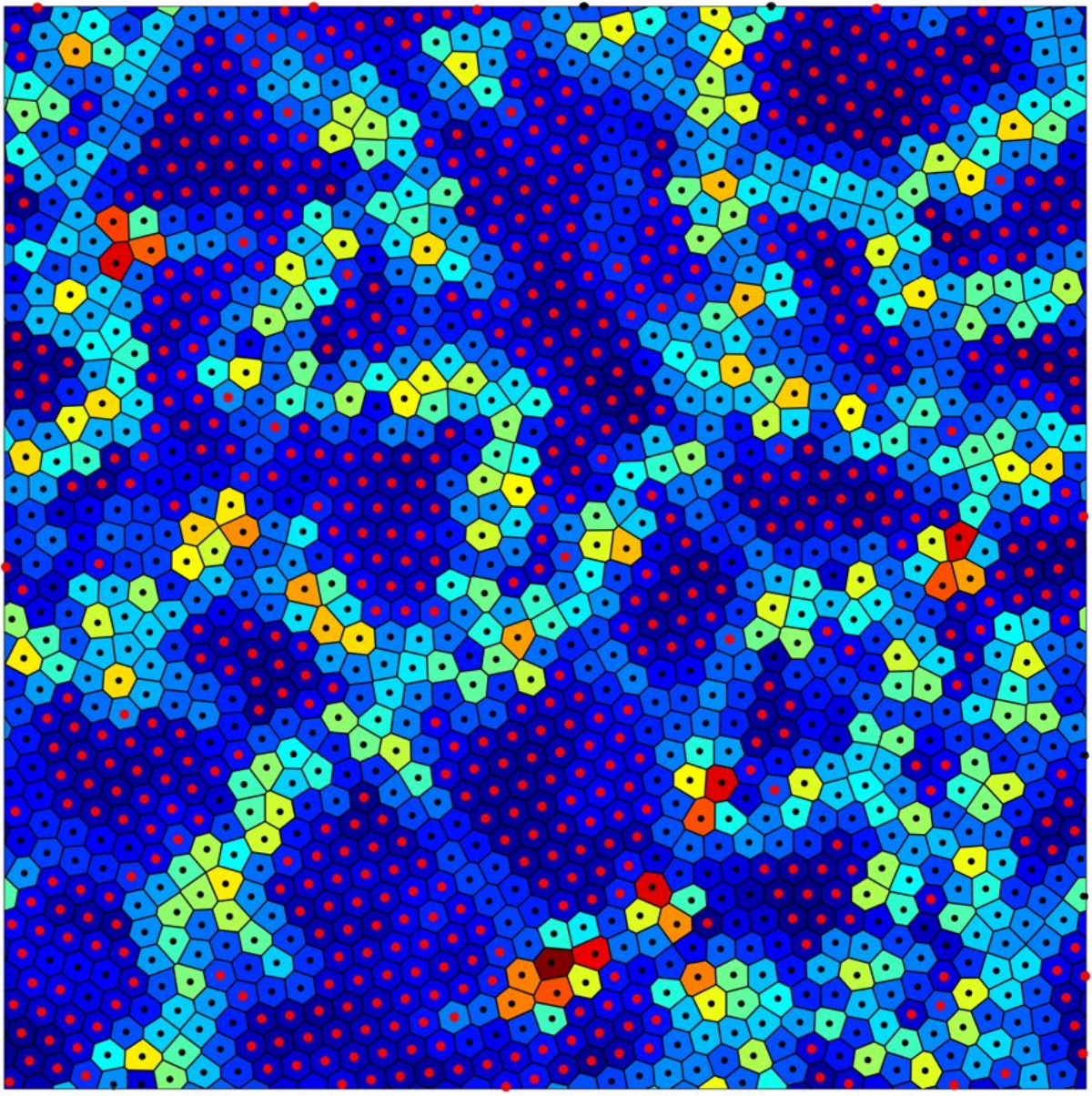} & \includegraphics[width=0.23\columnwidth]{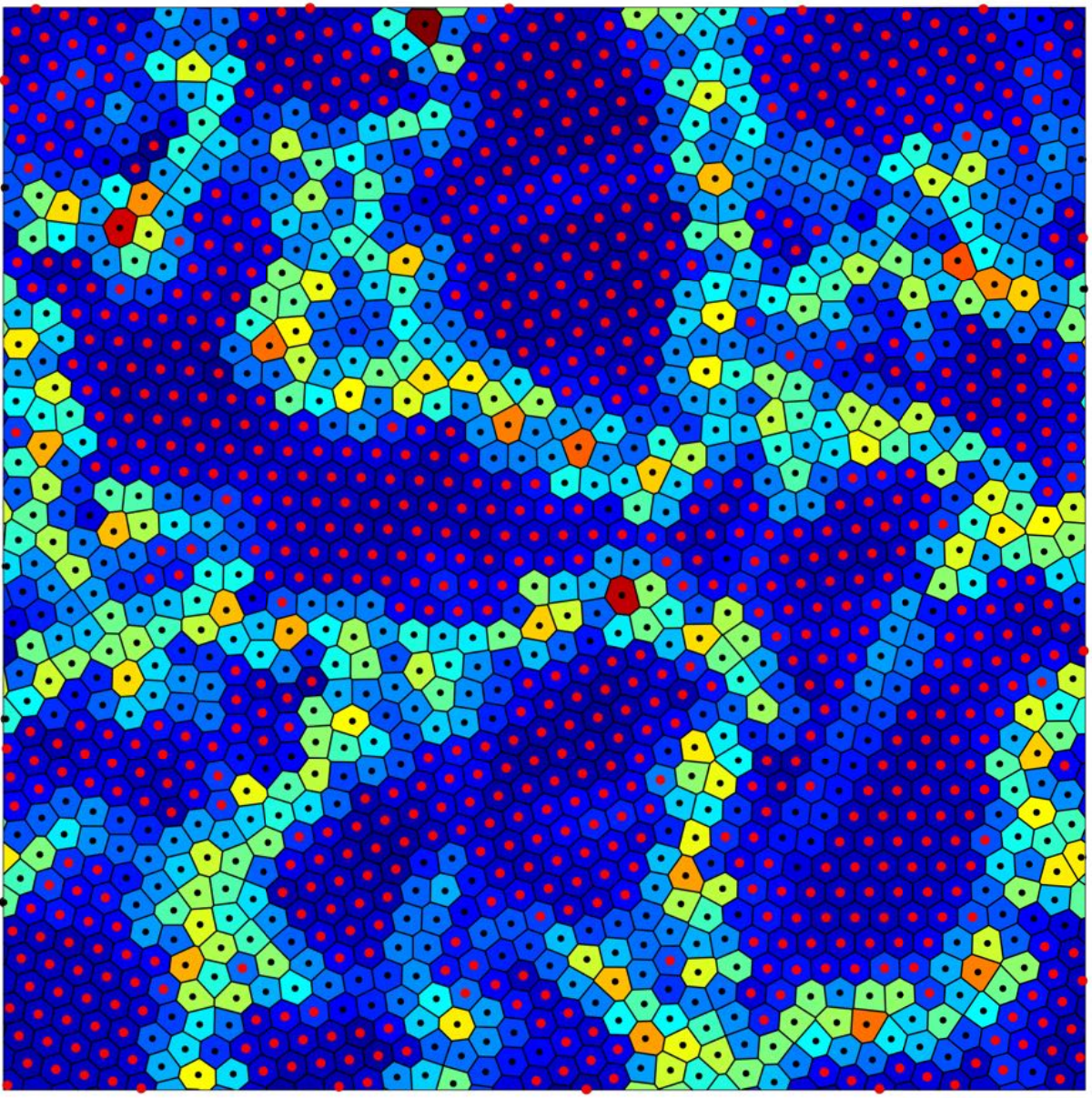} & \includegraphics[width=0.23\columnwidth]{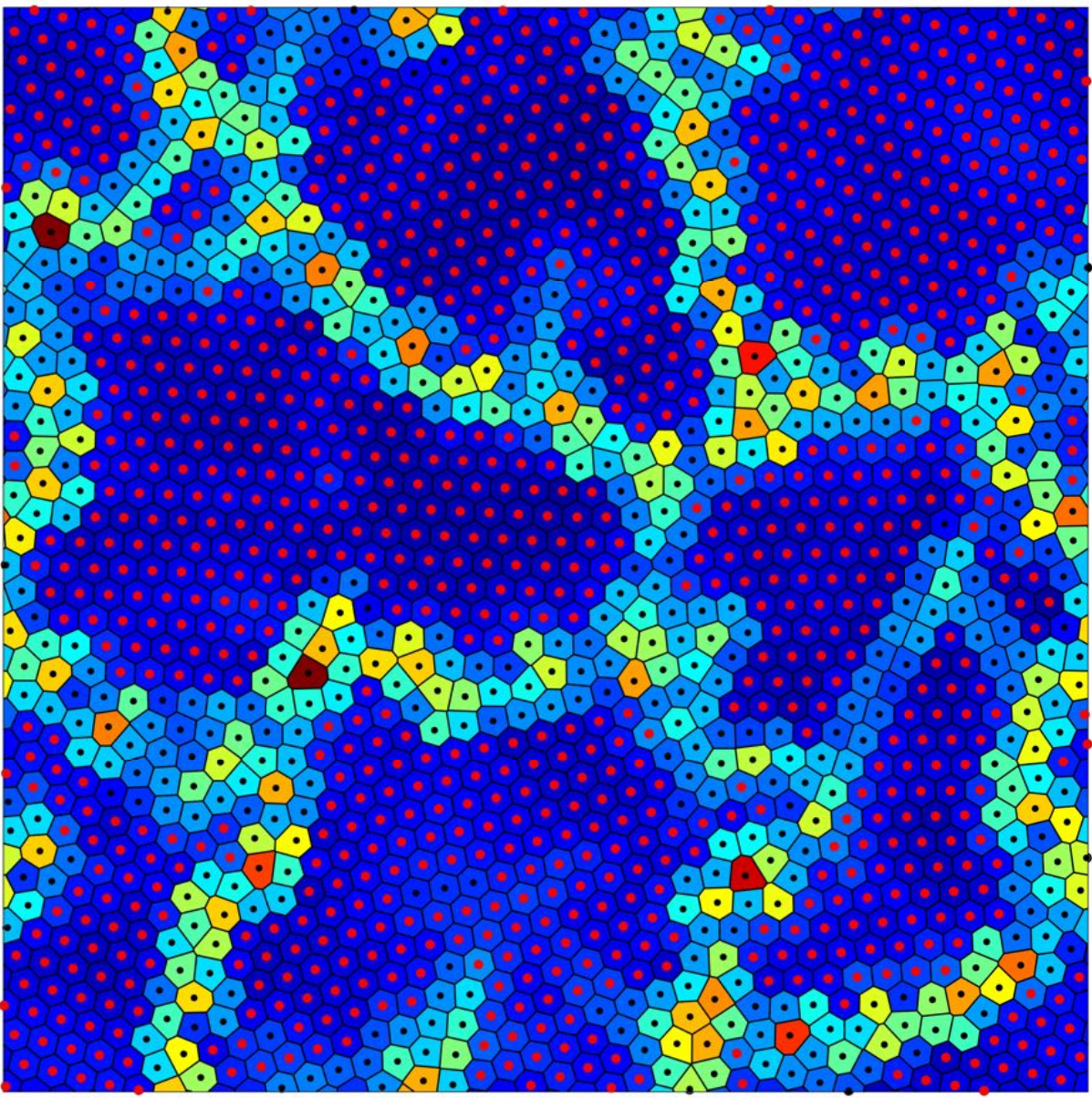}\tabularnewline
\footnotesize{(d) $k=2000$} & \footnotesize{(e) $k=4000$} & \footnotesize{(f) $k=8000$}\tabularnewline
\multicolumn{3}{c}{\includegraphics[width=0.4\columnwidth]{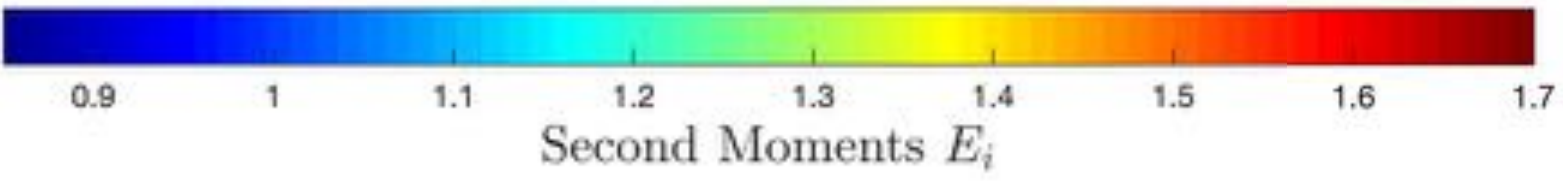}}\tabularnewline
\end{tabular}

\vspace{3mm}

\begin{tabular}{|c|c|c|c|c|c|c|}
\hline 
{\scriptsize{}$k$} & {\scriptsize{}0} & {\scriptsize{}500} & {\scriptsize{}1000} & {\scriptsize{}2000} & {\scriptsize{}4000} & {\scriptsize{}8000}\tabularnewline
\hline 
\hline 
{\scriptsize{}$E-1$} & {\scriptsize{}0.98007} & {\scriptsize{}0.02933} & {\scriptsize{}0.02598} & {\scriptsize{}0.02132} & {\scriptsize{}0.01875} & {\scriptsize{}0.01608}\tabularnewline
\hline 
{\scriptsize{}$H$ (\%)} & {\scriptsize{}29.33} & {\scriptsize{}77.86} & {\scriptsize{}81.60} & {\scriptsize{}84.80} & {\scriptsize{}87.46} & {\scriptsize{}89.60}\tabularnewline
\hline 
{\scriptsize{}$R_{\epsilon}$ (\%)} & {\scriptsize{}0.00} & {\scriptsize{}25.60} & {\scriptsize{}34.13} & {\scriptsize{}42.13} & {\scriptsize{}52.26} & {\scriptsize{}62.86}\tabularnewline

\hline

\multicolumn{7}{c}{{\footnotesize{}(g) measures specific to configurations (a)--(f)}}\tabularnewline
\end{tabular}

\vspace{4mm}

\begin{tabular}{cc}
\includegraphics[width=0.4\columnwidth, height=0.25\columnwidth]{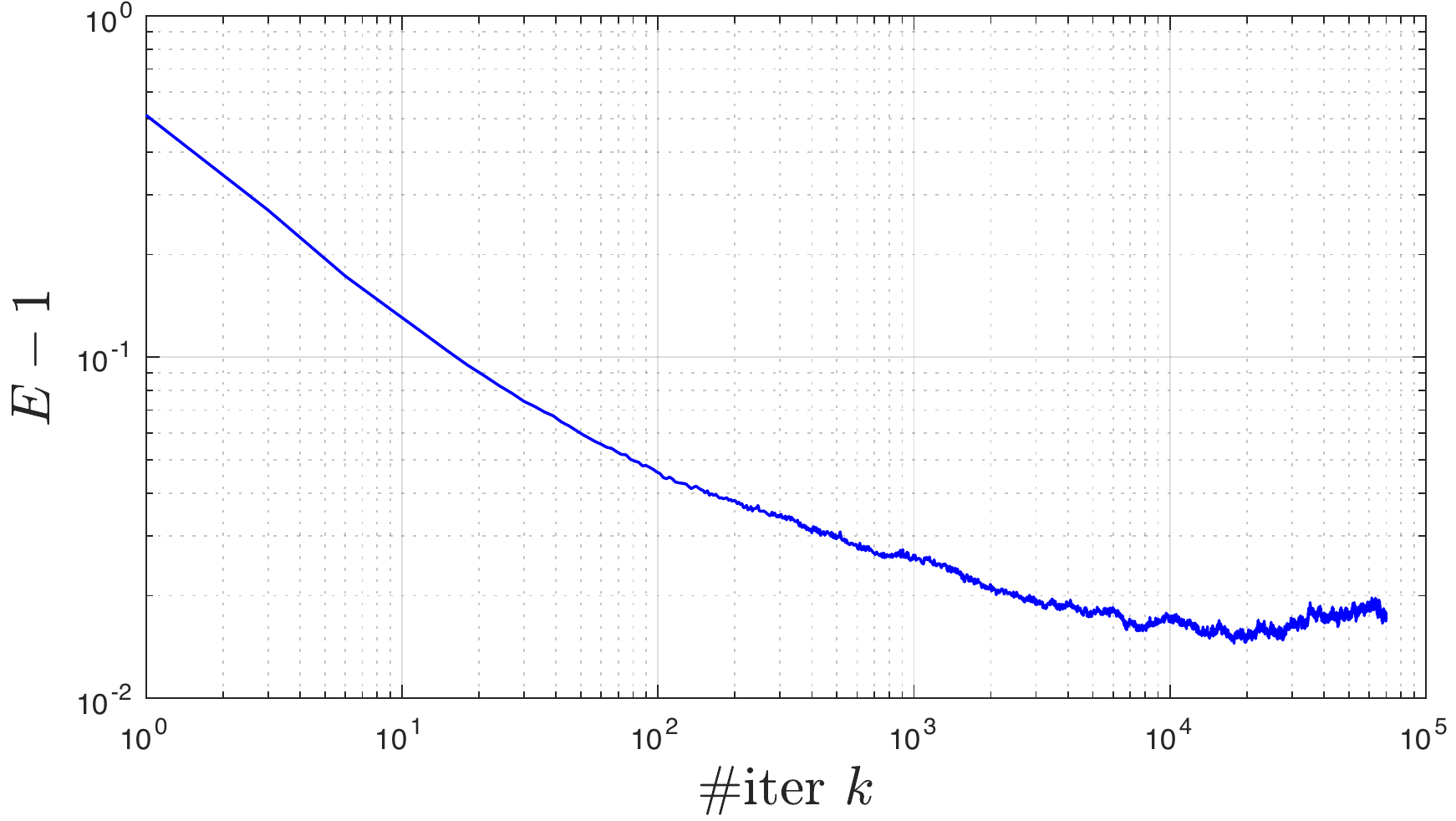} & \includegraphics[width=0.4\columnwidth, height=0.25\columnwidth]{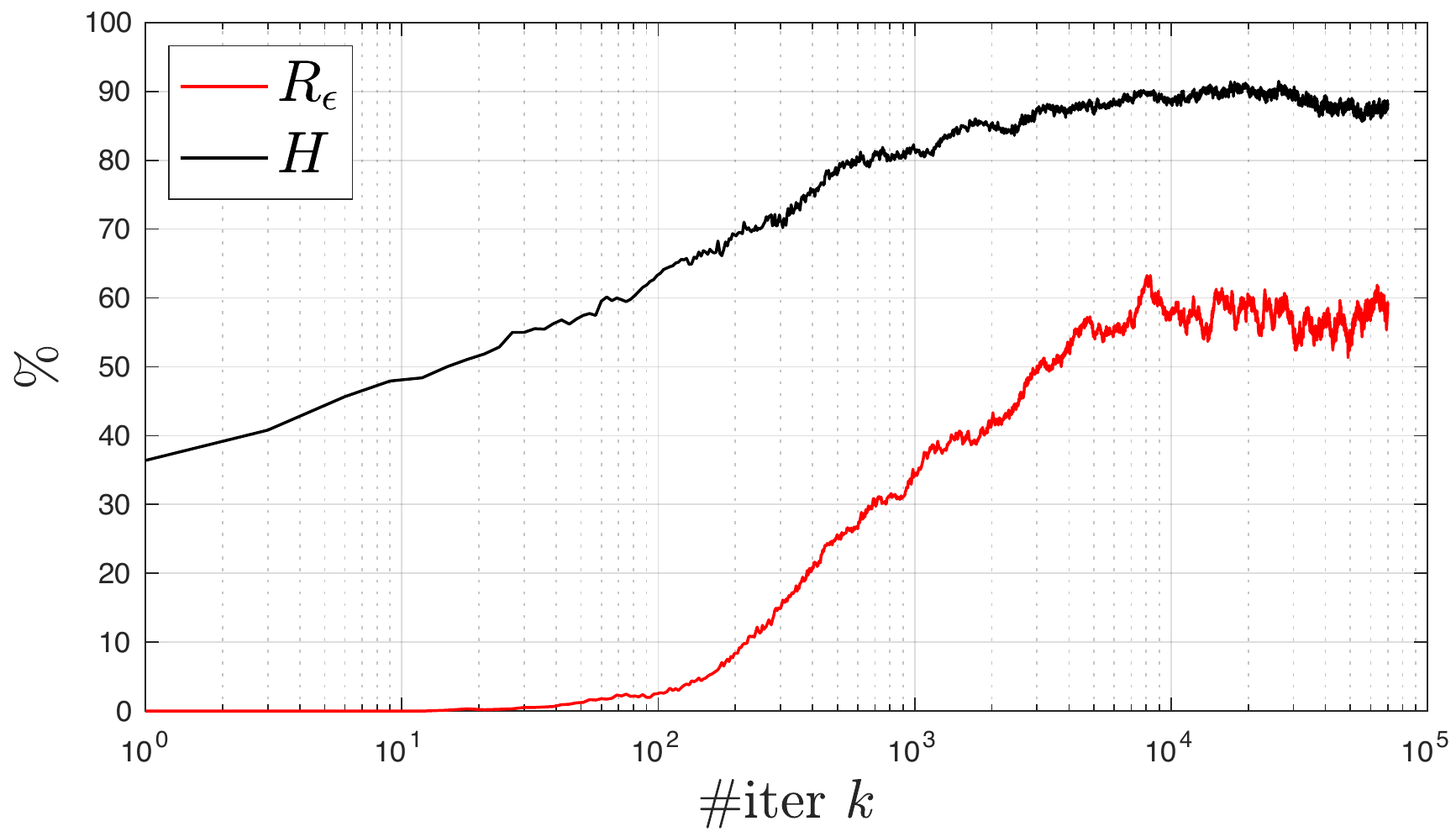} \tabularnewline
\footnotesize{(h) energy profile showing meso-stability} & \footnotesize{(i) $R_{\epsilon}$ and $H$ profiles}
\end{tabular}

\caption{Example of \textit{MACN-c} dynamics with $N=1500$ in the square torus $\Omega$.
The $k^{\text{{th}}}$ iterations displayed start from the uncorrelated
uniform random sampling of the domain shown in (a), the table contains the measures of the respective tessellations . The colouring represents $E_i$ and generators in red are those taken into account by $R_\epsilon\;(\epsilon=0.5\%)$. Finally we plot the three regularity measures against the iterate number up to $k=7\times10^{4}=K$
to show the different regimes; the system's meso-stability after $k\approx10^{4}$ is apparent.
\label{fig:Example of MACN-c dynamics}}

\end{figure}

We focus now on \textit{MACN-c} as a stand alone initialization in order to gain a first insight on the functionality of this iterative scheme.
With the concise example presented in Figure \ref{fig:Example of MACN-c dynamics} we show the general long term behavior of these dynamics over iterations $k=0,...,K-1$ with $K=7\times 10^4$.\\
The first immediate observation is the reduction of the discrepancy of the original cloud of points $\textbf{X}$.  Indeed when applying \textit{MACN-c} we obtain a significant decrease in energy, often close to two orders of magnitude when compared to a general random sampling. It is worth noting however that this scheme alone is not contractive due to the abrupt topological changes in the PDT from one iteration to another and thus will most likely fail to provide a PCVT by itself, even in the limit $K\to\infty$. \\
The second observation is that the system eventually reaches a "low" energy plateau and oscillates around it (see Figure \ref{fig:Example of MACN-c dynamics} (h)). This suggests that, in order to maximize energetic performance, the subsequent coupling with Lloyd's algorithm should be performed when the \textit{MACN-c} dynamics attain this regime. However it is difficult to estimate a priori when the system will reach such a mesostability
and even more so the required $K$ might be too large. Thus in practice $K$ (more precisely $Q\times K$) should be chosen by the user to retain tractability of the overall method we presented in Algorithm \ref{alg:hybrid-algorithm Q stages}. For this reason the presentation of our numerical results in \textsection \ref{sec:5} will start by a parameter sweep over values of $K$ that yield not necessary the best energetic results but a suitable trade-off between computation time and energy.

\subsection{\textit{MACN-$\delta$}\label{subsec:4.2}}

\begin{figure}[H]
\centering

\begin{tabular}{c}
\includegraphics[width=0.53\columnwidth]{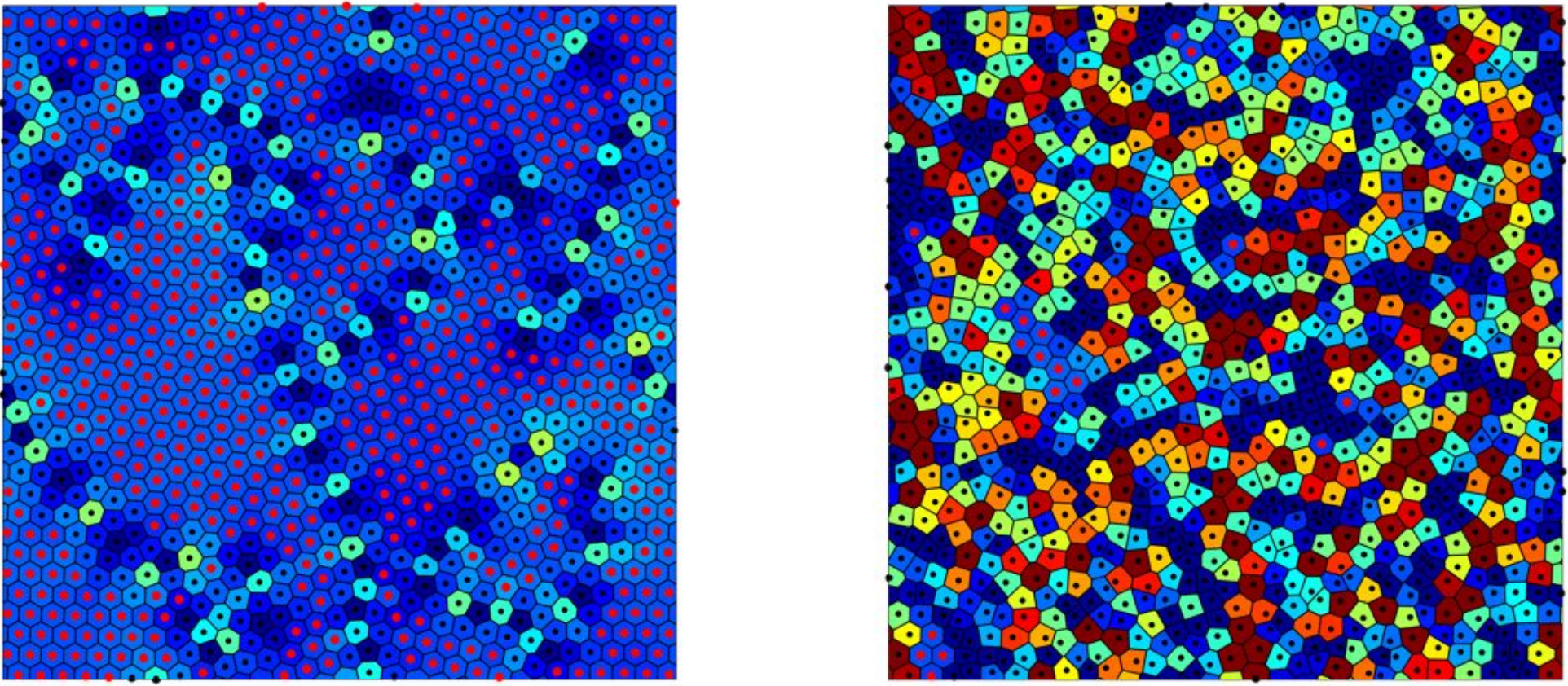}\tabularnewline
\includegraphics[width=0.4\columnwidth]{Graphics/CommonFigures/EnergyColorBar_horizontal}\tabularnewline
\end{tabular}

\caption{Example of a \textit{MACN-$\delta$} perturbation with $N=1000$ in
a square torus $\Omega$. On the left is a generic PCVT and on the right
is its resulting perturbation with $\delta=\frac{1}{4}\sqrt{|\Omega|/N}$. The colormap represents energy and hexagonal cells with red generators
are regular within the isoperimetric tolerance $\epsilon=0.5\%$.
\label{fig:example MACN-delta}}
\end{figure}

We emphasize once more the empirical nature of our choice for the distance $\delta$ (\ref{eq: original delta}) that was found among scalings of $\sqrt{|\Omega|/N}$ (the linear length-scale of the geometry in question). In particular the a priori mysterious factor $1/4$ was chosen to be in the small threshold of such scalings that simultaneously allows a change of basin of attraction while remaining small enough to preserve a "certain regularity" in the structure and topology of the tessellation.
Perhaps the strongest ascertainment is that $\delta$ equals the intrinsic length-scale of the regular hexagonal tessellation, namely
$|V|/|\partial V|$ where $V$ represents here the regular
hexagon with inscribed circle of diameter $\sqrt{|\Omega|/N}$. More
details on this latter interpretation of $\delta$ and on the global behavior of other variants of \textit{MACN} for PCVT dislocations are provided in \textsection \ref{sec:6}.\\
We finish the section with the example of Figure \ref{fig:example MACN-delta} showing a prototypical PCVT dislocation using \textit{MACN-$\delta$}.

\section{Numerical results\label{sec:5}}

A total of six Examples with different $N$'s are presented below: the first two are on the regular hexagonal torus while the remaining four are on the square torus. In the former the perfect honey comb lattice is attainable iff $N=a^{2}+ab+b^{2}\;\forall a,b\in\mathbb{N}$ \cite{PerimMinim_PentaTiling_Chung} while in the latter the ground state is unknown due to size frustration.
For conciseness of the exposition we will fix from now onwards the probing number to $Q=10$, this value will be large enough to demonstrate the power of our method vis-\`{a}-vis of the alternatives. Our study of the six cases thus begins by finding a preconditioning number $K$ offering a suitable trade-off between time and energy performance. Once $K$ is chosen we deepen the analysis of the performance of our method and then finish the section with a comparison of $R_{\epsilon}$ and $H$ as measures faithful to $E$.

Throughout the section we use uniformly uncorrelated initial configurations over $\Omega$ and our hybrid method is compared with L-BFGS(7) and P-L-BFGS(20,20), for the latter our implementation of  the pseudo code of Appendix \ref{App:C} uses the periodic adaptation of the preconditioner matrix $\widetilde{A}=G$ as described in (\ref{eq:Graph Laplacian def}). Additionally, we include results obtained with Lloyd's in our periodic set up to have a solid point of reference for past and future work since this is the only algorithm with a complete lack of tuning parameters.\\
Following collectively the results of \cite{CVT_EnerSmooth_Levy,HierachicalCVT_Wang,GraphLap_LongChen}
as well as our own implementation of some of the deterministic methods
recalled in \textsection \ref{sec:3}, we believe that the two Quasi-Newton choices of comparison
paint a good overview of the current deterministic state of the art methods: in particular, $\widetilde{A}=D^{2}E$
as well as other $(M,T)$ values were tested for the torus
but did not achieve noticeable systematic improvements in the regularity
measures. 

We further emphasize that, while the objective and main contribution of this paper is to establish a dynamical and fully deterministic way of sampling energy basins with as few parameters as possible, we tested the global Monte Carlo method from \cite{Global_MCM_Lu} on the torus; we report results later in \textsection \ref{sec:6} where the energy outcomes are comparable to ours and we discuss the comparative advantages of both methods. Furthermore, although not adapted to the periodic boundaries, quasi random samplings of $\Omega$ were also tested prior to Quasi-Newton minimization; the energies achieved out of 1,000 runs were comparable with those of the basins sampled by Lloyd's out the 10,000 runs initialized with uncorrelated distributions that are discussed below. 

Let us next introduce some notation; once $K$ is fixed, our hybrid method and its lowest sampled energy $E^{min}$ will primarily be compared with respect to the reference energy $E^{ref}$ given by the minimal energy PCVT obtained from a large batch of initial configurations using the three gradient based
methods.  More precisely, $E^{min}$ will be recorded with our method over 100 or 1,000 initial configurations (depending on the Example) while $E^{ref}$ will be the lowest energy amongst 100,000 runs for each of
L-BFGS(7) and P-L-BFGS(20,20) as well as 10,000 runs for Lloyd's.\\
We define then the following performance ratio for the sampled minimal energies. 

\begin{equation}
\tau:=\dfrac{E^{min}-1}{E^{ref}-1} \label{eq: ratio tau}
\end{equation}

We will also employ the empirical cumulative distribution functions (ECDFs) 
$f_{E-1},f_{R_{\epsilon}},f_{H}$  of our respective regularity
measures as well as the values $f_{E-1}^*,f_{R_{\epsilon}}^*,f_{H}^*$ obtained when evaluating the ECDFs of our hybrid method at $E^{ref},\,R_{\epsilon}^{ref},\,H^{ref}$ respectively; these quantities will establish the frequency of PCVTs for which our hybrid algorithm outperforms the best comparative method. Other basic statistics provided on regularity measures include averages and standard deviations taken over the designated number of runs, we denote them by $\langle \cdot \rangle$  and $\sigma_{(\cdot)} $ respectively. Finally we fix the isoperimetric tolerance to be $\epsilon\equiv0,5\%$, further detail on this choice will be discussed at the end of the section.

A final note on our energy measurements; no quadrature was involved (i.e. exact calculations were performed) and 
a tolerance was used on $||DE||/N$ guaranteeing that the values of energy listed in all the tables and figures are accurate at least up to the significant digits provided.

\subsection{Choosing $K$\label{subsec:5.1}}
\begin{figure}[H]
\centering

\begin{tabular}{cc}
\includegraphics[width=0.45\columnwidth, height=0.27\columnwidth]{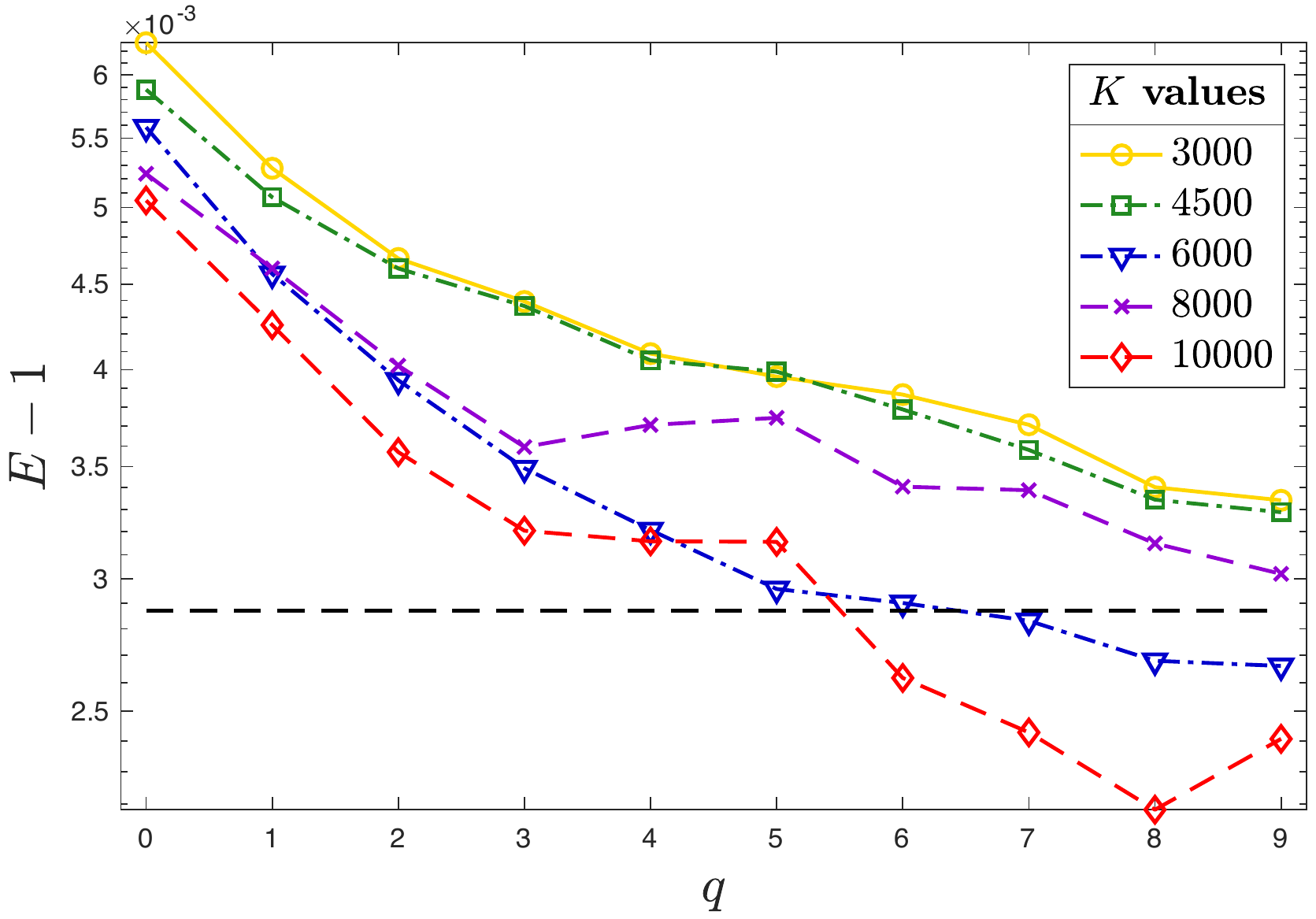}
& \includegraphics[width=0.45\columnwidth, height=0.27\columnwidth]{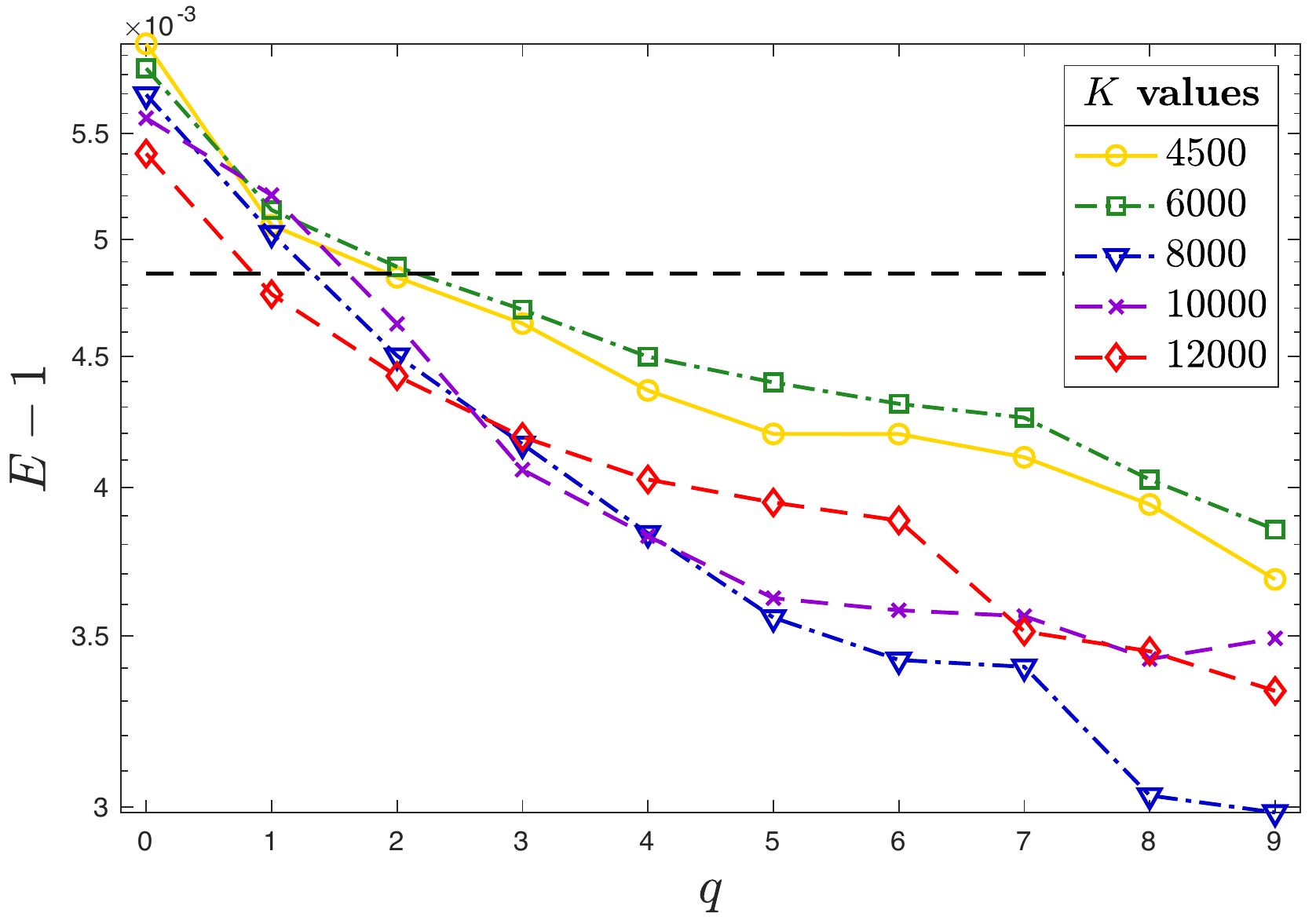} \tabularnewline
\footnotesize{(a) $N=973$ on the hexagonal torus} & \footnotesize{(b) $N=2029$ on the hexagonal torus} \tabularnewline
\includegraphics[width=0.45\columnwidth, height=0.27\columnwidth]{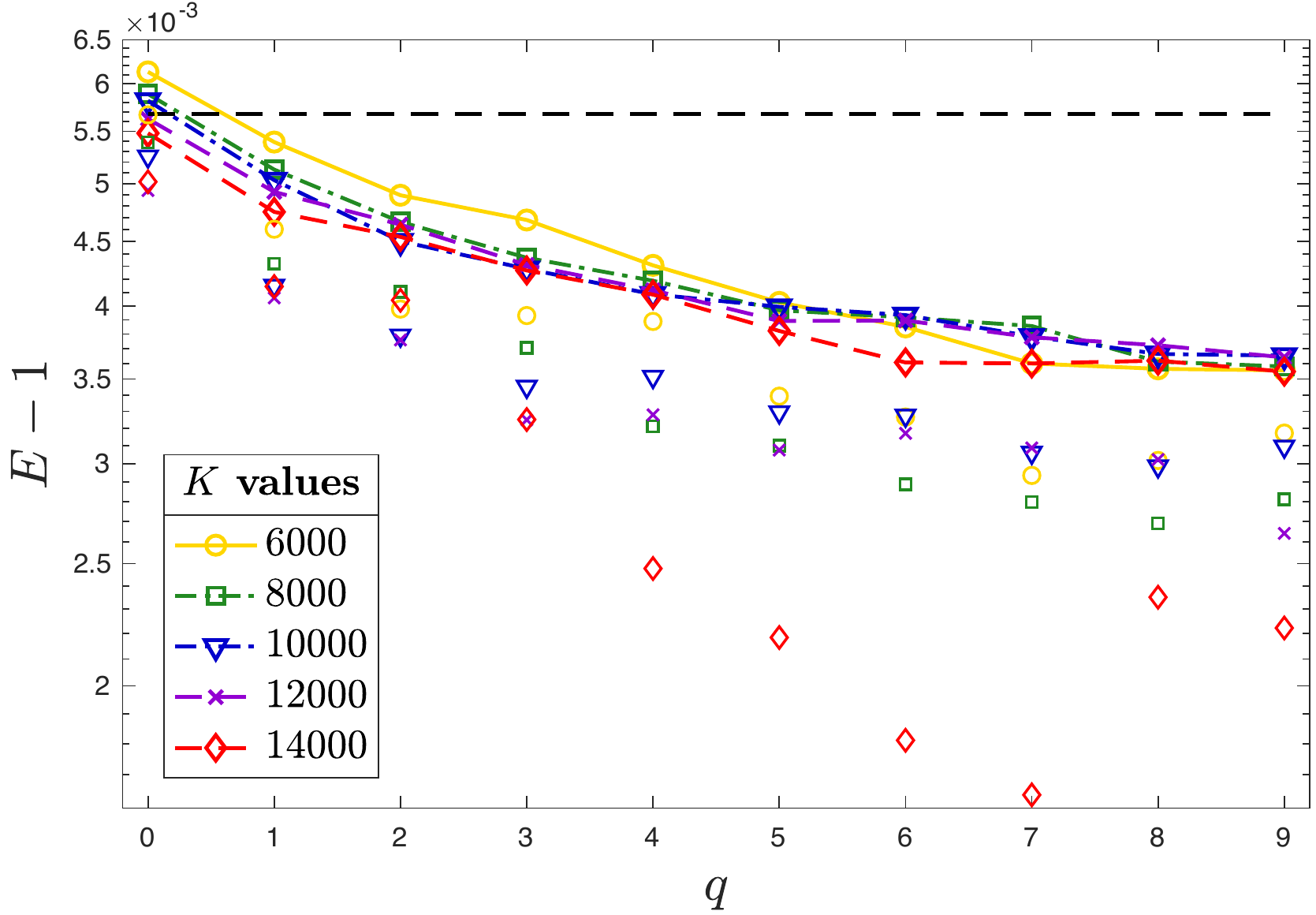}
& \includegraphics[width=0.45\columnwidth, height=0.27\columnwidth]{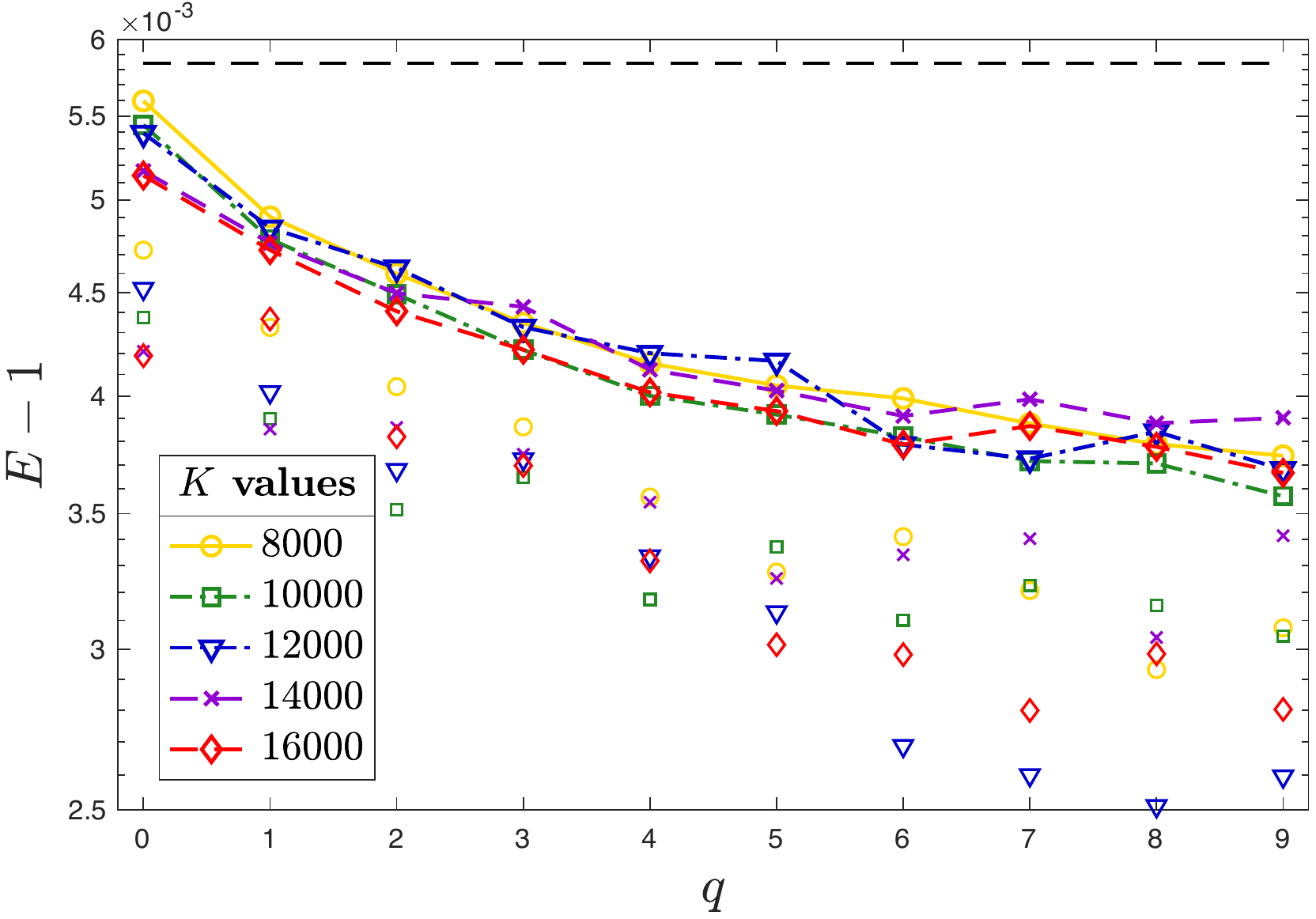}\tabularnewline
\footnotesize{(c) $N=3000$ on the square torus} & \footnotesize{(d) $N=4000$ on the square torus}

\end{tabular}
\caption{$K$-sweeps: the joint markers represent $\langle E-1 \rangle$ while the isolated ones represent minimums over  15 initial configurations. In black dotted lines are the values $E^{ref}-1$ appearing in Tables \ref{tab: N=973},\ref{tab: N=2029},\ref{tab: N=3000} and \ref{tab: N=4000} for a first comparison.\label{fig: K sweep}}
\end{figure}

The two scenarios on the hexagonal torus $\Omega$ (allowing the honey comb tiling) are with $N=973=17^{2}+17\times19+19^{2}$ and $N=2029=25^{2}+25\times27+27^{2}$. On the square torus we'll work with the values $N=\{n\times1000\}_{n=1}^{4}$. For these set ups we run our hybrid method on a reduced set of 15 initial configurations using a selected list of $K$'s, the energy results shown in Figure \ref{fig: K sweep} will allow us to choose trade-off values between time and energy. We remark on the general decrease tendency over the $Q$ hybrid stages but that
the decrease is not monotone in $K$.\\
Note as well that because the first two values of $N$ on the square torus are close to those on the hexagonal torus there is no need to run sweeps for $N=1000$ and $2000$, we'll just retain the same parameter values.

The graphs of $\langle E-1 \rangle$ in Figure \ref{fig: K sweep} (a) and (b) suggest that we pick $K=6000$ for Examples 1 \& 3 and $K=8000$ for Examples 2 \& 4. For Examples 5 \& 6 however, energy averages do not provide clear insight, we turn then to minimums from which  Figure \ref{fig: K sweep} (c) and (d) suggest we pick $K=8000$ and $K=12000$ respectively.

\subsection{Hexagonal torus $\Omega$\label{subsec:5.2}}

\paragraph*{Example 1\label{par:Example-1:}}

\begin{figure}[H]
\centering

\begin{tabular}{cc}
\includegraphics[width=0.33\columnwidth]{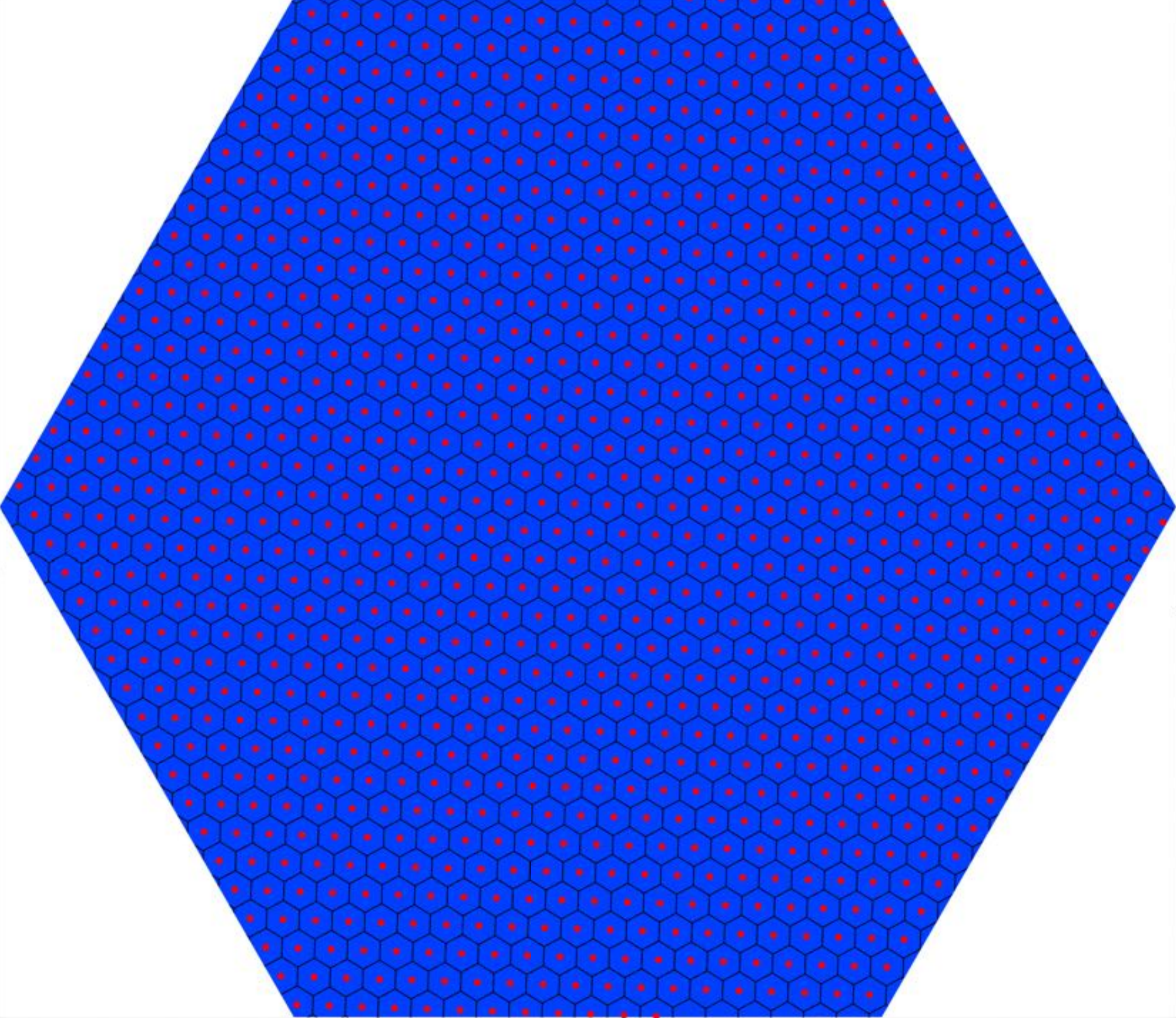} & \includegraphics[width=0.33\columnwidth]{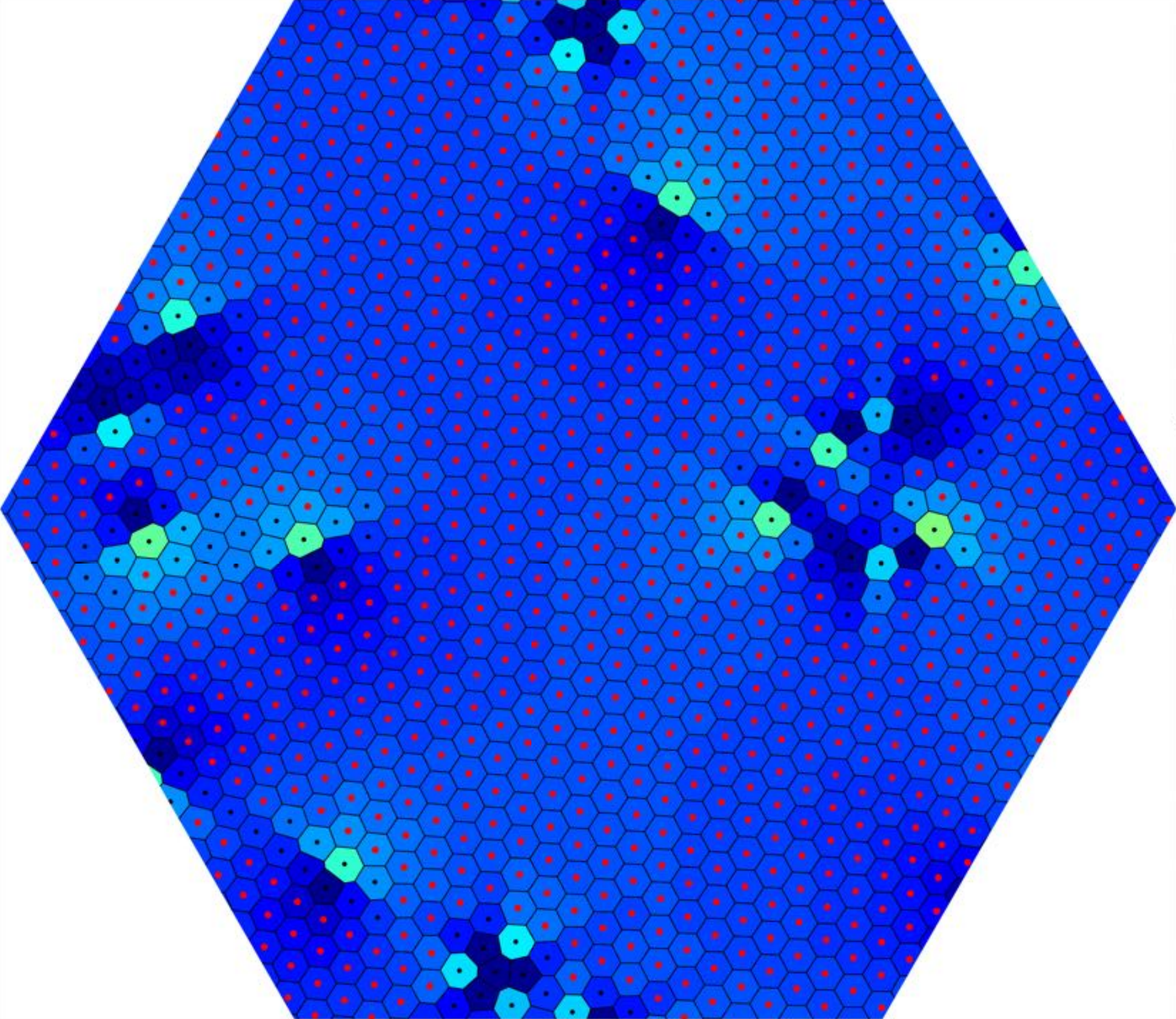}\tabularnewline
{\footnotesize{}(a) PCVT achieving $E^{min}$} & {\footnotesize{}(b) PCVT with $E^{ref}$}\tabularnewline
\end{tabular}
\begin{tabular}{cc}
\includegraphics[width=0.45\columnwidth, height=0.28\columnwidth]{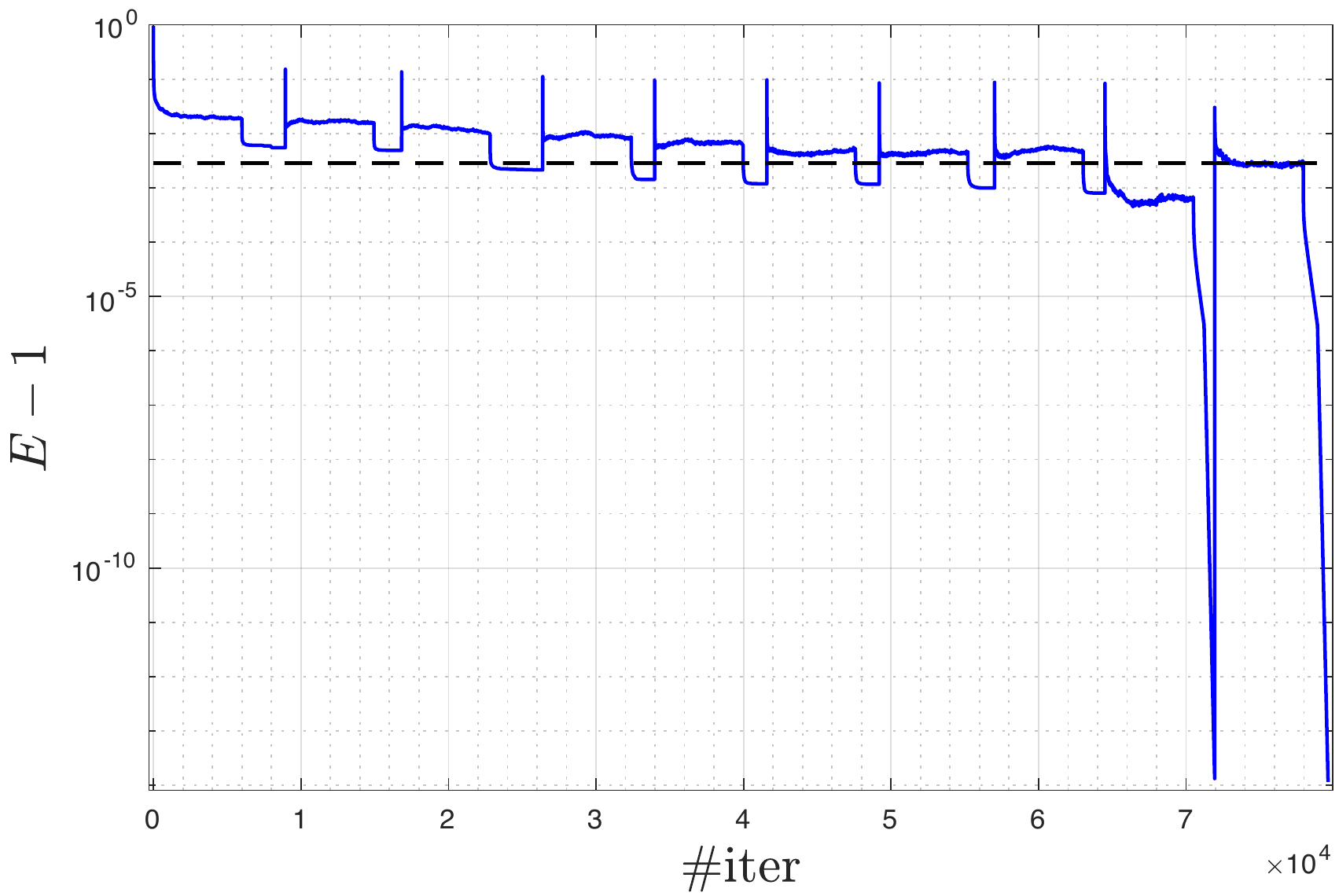} & \includegraphics[width=0.45\columnwidth, height=0.28\columnwidth]{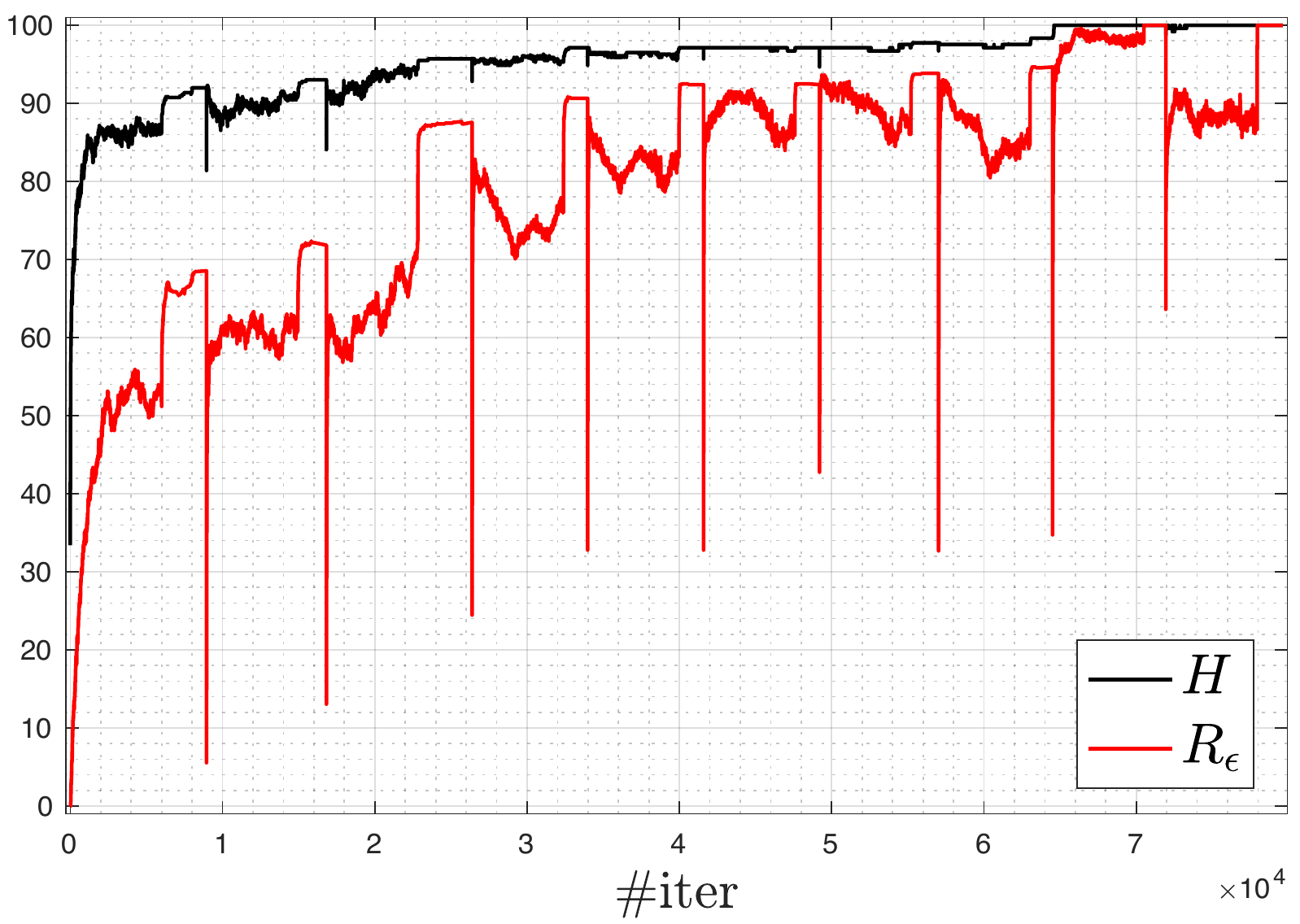}\tabularnewline
{\footnotesize{}(c) hybrid run reaching $E^{min}$ } & {\footnotesize{}(d) $R_{\epsilon},H$ profiles of (c)}\tabularnewline
\end{tabular}
\caption{Optimal PCVTs from Example 1: (a) ground state with $E^{min}-1\approx1e-14$ reached by our hybrid method. (b) configuration carrying $E^{ref}-1=0.00287$ achieved by L-BFGS(7). (c) and (d) are the measures of a hybrid run that reached $E^{min}$ amongst the larger batch of 100 runs with $K=6000$ (the dotted black line designates $E^{ref}-1$ again). \label{fig: best configurations N973}}
\end{figure}

\begin{figure}[H]
\centering

\begin{tabular}{c}
\includegraphics[width=0.55\columnwidth, height=0.25\columnwidth]{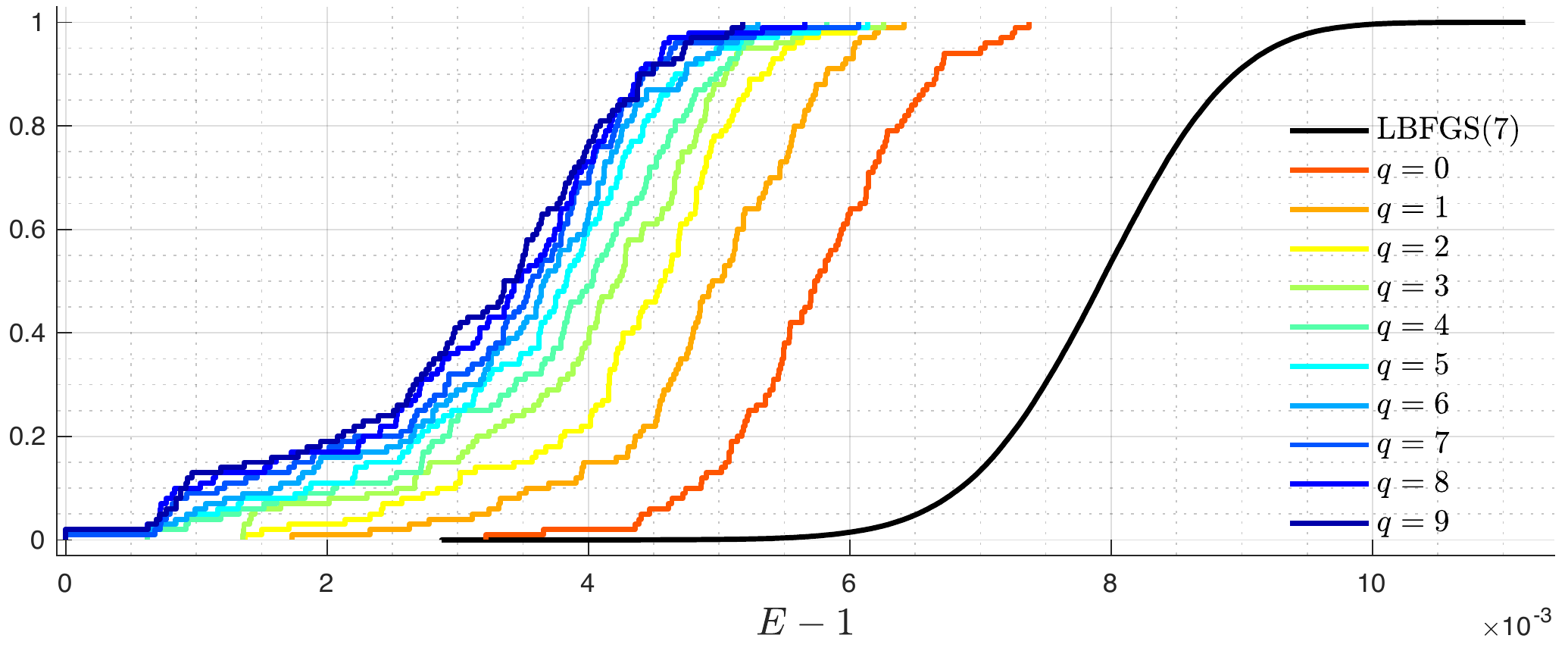}\tabularnewline
{\footnotesize{}(a) $f_{E-1}$}\tabularnewline
\end{tabular}

\begin{tabular}{cc}
\includegraphics[width=0.43\columnwidth, height=0.26\columnwidth]{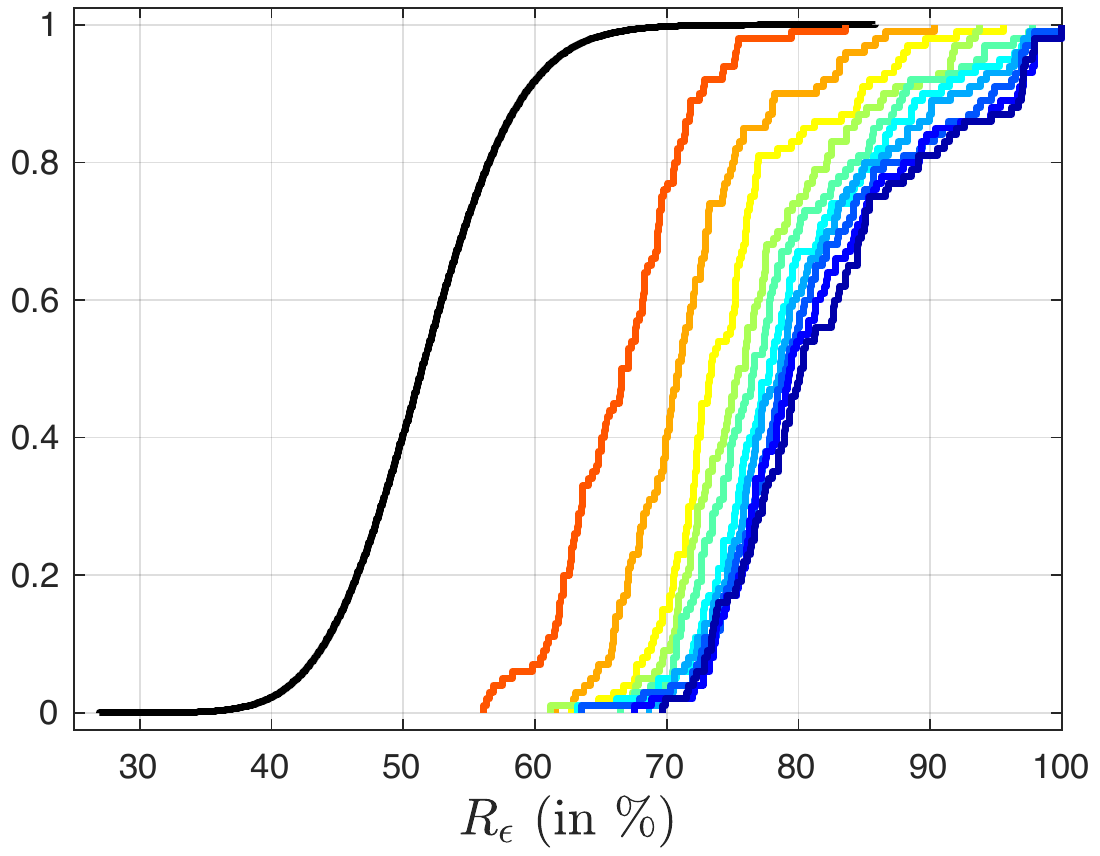} & \includegraphics[width=0.43\columnwidth, height=0.26\columnwidth]{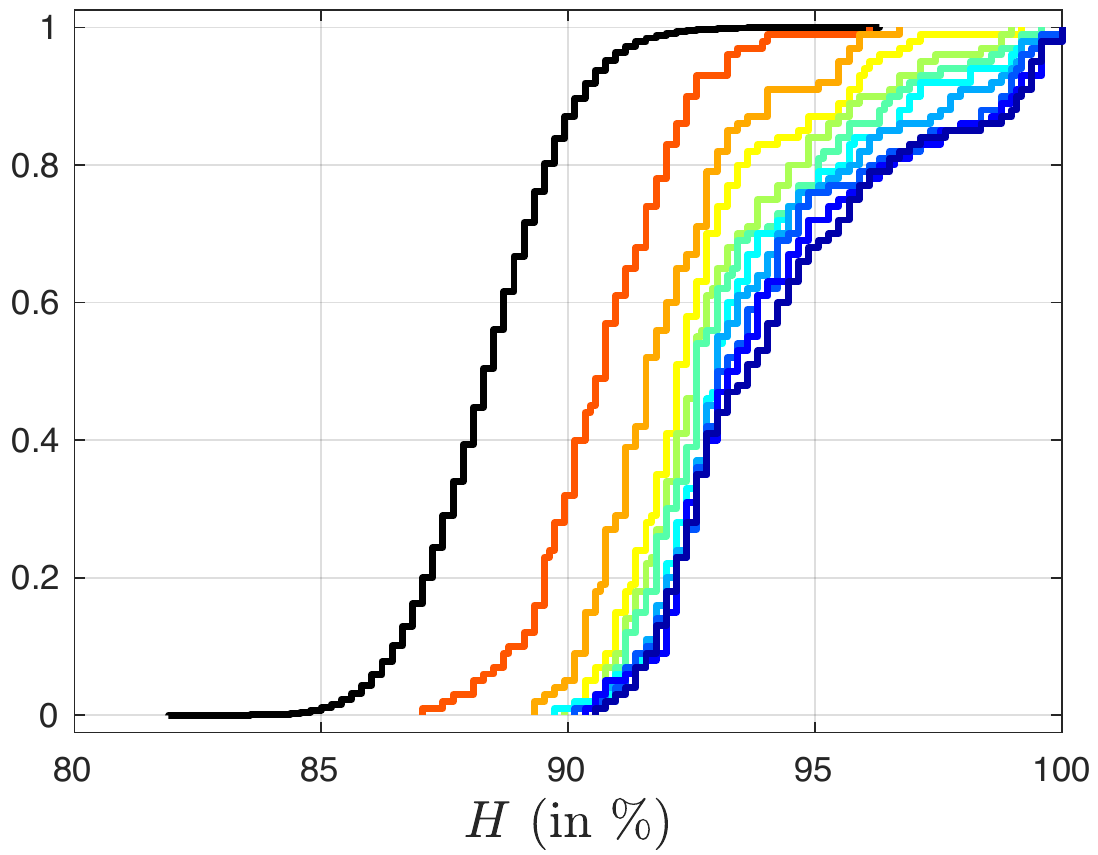}\tabularnewline
{\footnotesize{}(b) $f_{R_{\epsilon}}$} & {\footnotesize{}(c)$f_{H}$}\tabularnewline
\end{tabular}

\caption{Data from Example 1: ECDFs of the three regularity measures for 100 hybrid runs with $K=6000$ along with 100,000 runs of L-BFGS(7). \label{fig: data N973}}
\end{figure}

With $N=973$ and $K=6000$ chosen, we run our method on a larger batch of 100 initial configurations; it reached the ground state with $E^{min}-1\approx1$e-14 (highest precision allowed by our implementation) while the optimal PCVT from the comparative methods is $E^{ref}-1=0.00287$ (achieved by L-BFGS(7)), see Figures \ref{fig: best configurations N973} \& \ref{fig: data N973} as well as the statistics summary of Table \ref{tab: N=973}. In particular, Figure \ref{fig: best configurations N973} (c) shows how the symbiotic blocks of \textit{MACN}-$\delta$+\textit{MACN-c} act on probing non-PCVT configurations with energy close to (if not below) $E^{ref}$.\\
Notice at last how the values $f_{E-1}^{*},f_{R_{\epsilon}}^{*},f_{H}^{*}$ indicate that our method outperforms the others for a significant fraction of the runs. In particular, not only do we get the ratio $\tau$ (\ref{eq: ratio tau}) to be sensibly zero but; on average one needs to run
our hybrid algorithm on approximatively three uniformly sampled initial
configurations to obtain an energy lower that $E^{ref}$. In other
words, we only need $\approx30$ PCVTs so that our way of probing
the energy landscape achieves comparable results to the sampling of the basins done by L-BFGS(7) out of 100,000 runs.\\ 

\paragraph*{Example 2\label{par:Example-2:}}
For $N=2029$ we first remark from the reduced sets of runs from Figure \ref{fig: K sweep} (b) that energy averages and their deviation from $E^{ref}$ immediately compare favorably to those in Example 1, this is the first indicator that our hybrid method is less affected than the comparative methods by the increasing non-convexity of the energy with $N$. The same conclusion can be drawn from the data from a hundred runs in Tables \ref{tab: N=973} and \ref{tab: N=2029}; more precisely we have that although our hybrid method did not achieve the honey comb structure, our optimal PCVT with $E^{min}-1=0.00150$ is still far more regular and gets $\tau^{-1}\thickapprox3.2$ times closer
to the ground state than the compared methods. Furthermore, the impressive value $f_{E-1}^{*}=1.00$ achieved suggests that despite the increased non-convexity with $N$, the required number of PCVTs needed to navigate the energy
landscape, in a similar fashion as the random sampling made with
gradient based methods, has decreased to 9 (i.e. less than one full run).

\begin{table}[H]

\resizebox{1\columnwidth}{0.13\columnwidth}{
\begin{centering}
\begin{tabular}{|>{\centering}p{13mm}||cccccccccc||ccc|}
\cline{2-14} 
\multicolumn{1}{>{\centering}p{9mm}|}{} &  &  &  &  & Hybrid &  &  &  &  &  & \multicolumn{1}{c|}{\multirow{2}{*}{Lloyd}} & \multicolumn{1}{c|}{\multirow{2}{*}{LBFGS}} & \multirow{2}{*}{PLBFGS}\tabularnewline
\cline{2-11} 
\multicolumn{1}{>{\centering}p{9mm}|}{} & $\mathbf{X}_{0}^{*}$ & $\mathbf{X}_{1}^{*}$ & $\mathbf{X}_{2}^{*}$ & $\mathbf{X}_{3}^{*}$ & $\mathbf{X}_{4}^{*}$ & $\mathbf{X}_{5}^{*}$ & $\mathbf{X}_{6}^{*}$ & $\mathbf{X}_{7}^{*}$ & $\mathbf{X}_{8}^{*}$ & $\mathbf{X}_{9}^{*}$ &  &  & \tabularnewline
\hline 
{$\langle E-1\rangle$} & 0.00574 & 0.00488 & 0.00435 & 0.00401 & 0.00377 & 0.00356 & 0.00342 & 0.00328 & 0.00316 & 0.00309 & 0.00807 & 0.00790 & 0.00789\tabularnewline
{$\sigma_{E-1}$} & 0.00075 & 0.00087 & 0.00095 & 0.00105 & 0.00113 & 0.00116 & 0.00116 & 0.00121 & 0.00123 & 0.00124 & 0.00082 & 0.00082 & 0.00082\tabularnewline
{max} & 0.00737 & 0.00641 & 0.00626 & 0.00625 & 0.00582 & 0.00613 & 0.00529 & 0.00606 & 0.00565 & 0.00517 & 0.01102 & 0.01116 & 0.01114\tabularnewline
{min} & 0.00321 & 0.00173 & 0.00137 & 0.00135 & 0.00062 & \textbf{$\approx$1e-14} & \textbf{$\approx$1e-14} & \textbf{$\approx$1e-14} & \textbf{$\approx$1e-14} & \textbf{$\approx$1e-14} & 0.00383 & \textbf{0.00287} & 0.00349\tabularnewline
{$f_{E-1}^{*}$} & 0.00 & 0.04 & 0.10 & 0.15 & 0.19 & 0.23 & 0.26 & 0.28 & 0.33 & \textbf{0.36} & - - - & - - - & - - -\tabularnewline
\hline 
{$\langle R_{\epsilon}\rangle$ } & 66.71 & 71.76 & 76.94 & 76.88 & 78.22 & 79.48 & 80.22 & 81.11 & 81.81 & 82.33 & 50.42 & 51.56 & 51.63\tabularnewline
{$\sigma_{R_{\epsilon}}$} & 5.01 & 5.57 & 6.06 & 6.69 & 7.05 & 7.33 & 7.34 & 7.83 & 7.90 & 7.96 & 5.99 & 5.98 & 5.98\tabularnewline
{max} & 83.55 & 90.33 & 95.58 & 93.73 & 97.73 & \textbf{100} & \textbf{100} & \textbf{100} & \textbf{100} & \textbf{100} & 79.44 & \textbf{85.81} & 82.01\tabularnewline
{min} & 56.11 & 61.56 & 62.79 & 61.15 & 66.49 & 63.20 & 68.65 & 63.51 & 67.52 & 69.68 & 31.55 & 26.92 & 28.57\tabularnewline
{$1-f_{R\epsilon}^{*}$} & 0.00 & 0.02 & 0.08 & 0.12 & 0.15 & 0.19 & 0.20 & 0.21 & 0.24 & \textbf{0.25} & - - - & - - - & - - -\tabularnewline
\hline 
{$\langle H\rangle$ } & 90.72 & 91.98 & 92.66 & 93.11 & 93.38 & 93.69 & 93.85 & 94.08 & 94.23 & 94.35 & 88.26 & 88.38 & 88.39\tabularnewline
{$\sigma_{H}$} & 1.51 & 1.60 & 1.76 & 2.00 & 2.15 & 2.23 & 2.34 & 2.52 & 2.59 & 2.60 & 1.46 & 1.47 & 1.47\tabularnewline
{max} & 96.09 & 96.71 & 99.17 & 98.97 & 99.58 & \textbf{100} & \textbf{100} & \textbf{100} & \textbf{100} & \textbf{100} & 94.65 & \textbf{96.30} & 95.68\tabularnewline
{min} & 87.05 & 89.31 & 89.92 & 89.92 & 90.33 & 89.72 & 90.33 & 90.13 & 90.33 & 90.54 & 82.52 & 81.91 & 81.91\tabularnewline
{$1-f_{H}^{*}$} & 0.00 & 0.01 & 0.05 & 0.10 & 0.14 & 0.16 & 0.16 & 0.20 & \textbf{0.22} & 0.21 & - - - & - - - & - - -\tabularnewline
\hline 
\end{tabular}
\par\end{centering}
}
\caption{Statistics of Example 1 for $N=973$ with a 100 runs of our hybrid method using $K=6000$, the values $E^{min}-1,\,R_{\epsilon}^{min} ,\,H^{min}$ as well as $E^{ref}-1,\,R_{\epsilon}^{ref} ,\,H^{ref}$ are bold faced.\label{tab: N=973}}

\resizebox{1\columnwidth}{0.13\columnwidth}{
\begin{centering}
\begin{tabular}{|>{\centering}p{13mm}||cccccccccc||ccc|}
\cline{2-14} 
\multicolumn{1}{>{\centering}p{9mm}|}{} &  &  &  &  & Hybrid &  &  &  &  &  & \multicolumn{1}{c|}{\multirow{2}{*}{Lloyd}} & \multicolumn{1}{c|}{\multirow{2}{*}{LBFGS}} & \multirow{2}{*}{PLBFGS}\tabularnewline
\cline{2-11} 
\multicolumn{1}{>{\centering}p{9mm}|}{} & $\mathbf{X}_{0}^{*}$ & $\mathbf{X}_{1}^{*}$ & $\mathbf{X}_{2}^{*}$ & $\mathbf{X}_{3}^{*}$ & $\mathbf{X}_{4}^{*}$ & $\mathbf{X}_{5}^{*}$ & $\mathbf{X}_{6}^{*}$ & $\mathbf{X}_{7}^{*}$ & $\mathbf{X}_{8}^{*}$ & $\mathbf{X}_{9}^{*}$ &  &  & \tabularnewline
\hline 
{$\langle E-1\rangle$} & 0.00559 & 0.00482 & 0.00445 & 0.00423 & 0.00397 & 0.00379 & 0.00370 & 0.00357 & 0.00342 & 0.00335 & 0.00802 & 0.00784 & 0.00783\tabularnewline
{$\sigma_{E-1}$} & 0.00047 & 0.00053 & 0.00053 & 0.00046 & 0.00047 & 0.00050 & 0.00054 & 0.00054 & 0.00059 & 0.00065 & 0.00057 & 0.00057 & 0.00057\tabularnewline
{max} & 0.00659 & 0.00598 & 0.00567 & 0.00548 & 0.00529 & 0.00515 & 0.00511 & 0.00493 & 0.00473 & 0.00518 & 0.00985 & 0.01013 & 0.01034\tabularnewline
{min} & 0.00401 & 0.00361 & 0.00301 & 0.00304 & 0.00245 & 0.00249 & 0.00202 & 0.00176 & 0.00176 & \textbf{0.00150} & 0.00572 & \textbf{0.00485} & 0.00523\tabularnewline
{$f_{E-1}^{*}$} & 0.08 & 0.53 & 0.77 & 0.93 & 0.99 & 0.98 & 0.98 & 0.98 & \textbf{1.00} & 0.99 & - - - & - - - & - - -\tabularnewline
\hline 
{$\langle R_{\epsilon}\rangle$ } & 67.41 & 72.06 & 74.25 & 75.58 & 77.10 & 78.23 & 78.75 & 79.63 & 80.44 & 80.93 & 50.75 & 51.98 & 52.03\tabularnewline
{$\sigma_{R_{\epsilon}}$} & 3.12 & 3.42 & 3.41 & 2.90 & 2.93 & 3.17 & 3.37 & 3.33 & 3.61 & 4.02 & 4.17 & 4.15 & 4.15\tabularnewline
{max} & 77.62 & 81.27 & 84.22 & 82.65 & 85.36 & 85.60 & 88.12 & 90.58 & 90.53 & \textbf{92.65} & 67.12 & \textbf{72.10} & 71.21\tabularnewline
{min} & 61.60 & 64.31 & 65.94 & 67.12 & 68.45 & 68.85 & 69.29 & 70.08 & 72.15 & 69.09 & 37.16 & 35.18 & 35.97\tabularnewline
{$1-f_{R\epsilon}^{*}$} & 0.07 & 0.50 & 0.73 & 0.90 & 0.97 & 0.97 & 0.98 & 0.99 & \textbf{1.00} & 0.99 & - - - & - - - & - - -\tabularnewline
\hline 
{$\langle H\rangle$ } & 90.99 & 92.04 & 92.61 & 92.90 & 93.32 & 93.57 & 93.77 & 93.95 & 94.19 & 94.33 & 88.38 & 88.50 & 88.50\tabularnewline
{$\sigma_{H}$} & 0.89 & 0.96 & 1.01 & 0.89 & 0.91 & 0.95 & 0.99 & 0.97 & 1.09 & 1.21 & 0.99 & 1.01 & 1.01\tabularnewline
{max} & 93.39 & 94.48 & 95.66 & 95.56 & 95.46 & 96.05 & 96.15 & 96.64 & 96.64 & \textbf{98.02} & 92.70 & \textbf{93.49} & 93.00\tabularnewline
{min} & 88.71 & 89.94 & 89.74 & 90.04 & 91.22 & 91.42 & 90.93 & 91.03 & 92.11 & 90.93 & 85.11 & 84.42 & 84.22\tabularnewline
{$1-f_{H}^{*}$} & 0.00 & 0.06 & 0.21 & 0.29 & 0.42 & 0.61 & 0.64 & 0.68 & 0.78 & \textbf{0.80} & - - - & - - - & - - -\tabularnewline
\hline 
\end{tabular}
\par\end{centering}
}
\caption{Statistics of Example 2 for $N=2029$ with a 100 runs of our hybrid method using $K=8000$.\label{tab: N=2029}}

\end{table}

\subsection{Square torus $\Omega$\label{subsec:5.3}}

We now work with $\Omega$ being the square torus and $N=\{n\times1000\}_{n=1}^{4}$ to investigate the behavior of our method when the shape effects are tightened.\\

\paragraph*{Example 3\label{par:Example-3}}

\begin{figure}[H]
\centering

\begin{tabular}{ccc}
\includegraphics[width=0.2\columnwidth]{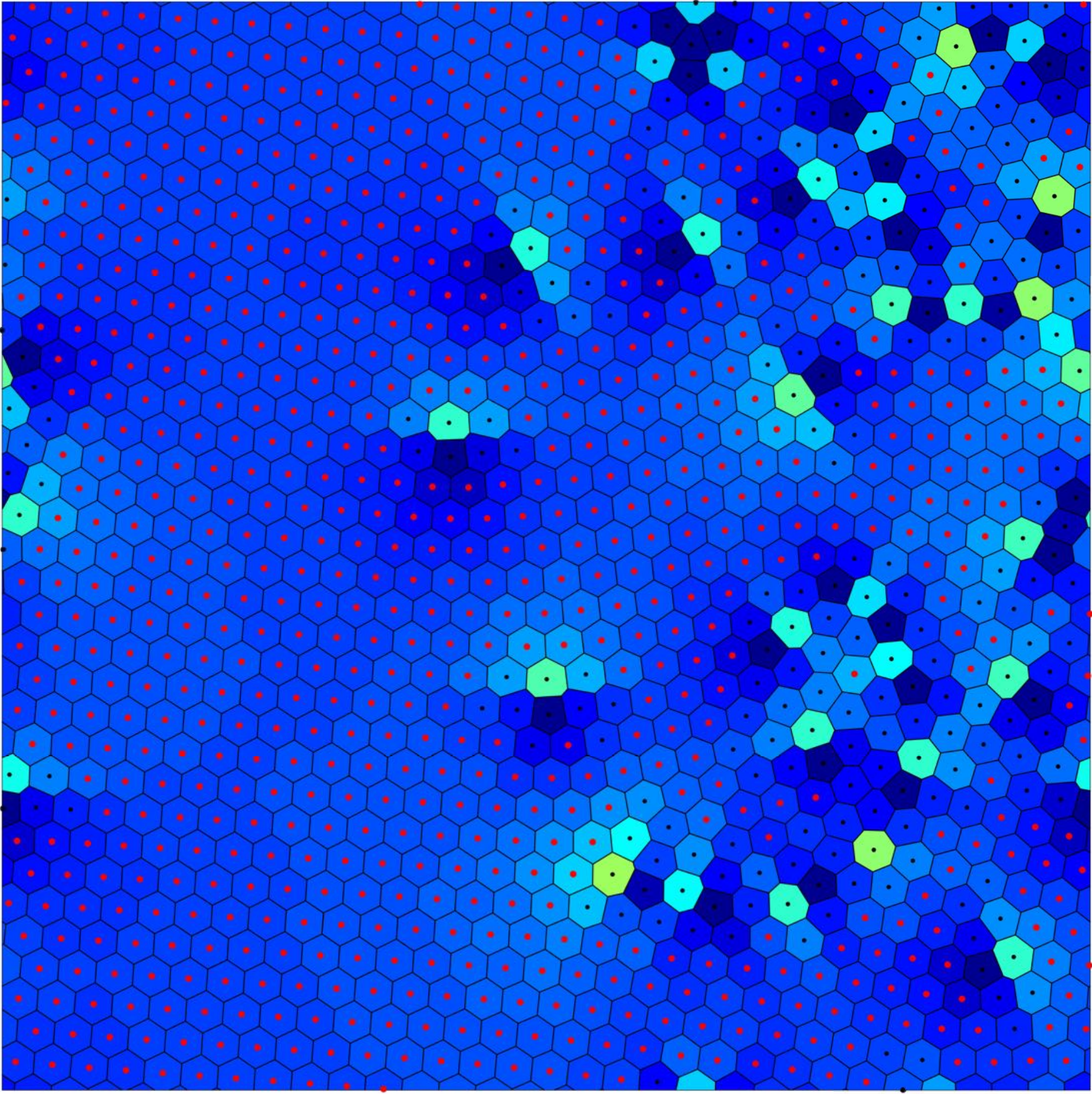} & \includegraphics[width=0.2\columnwidth]{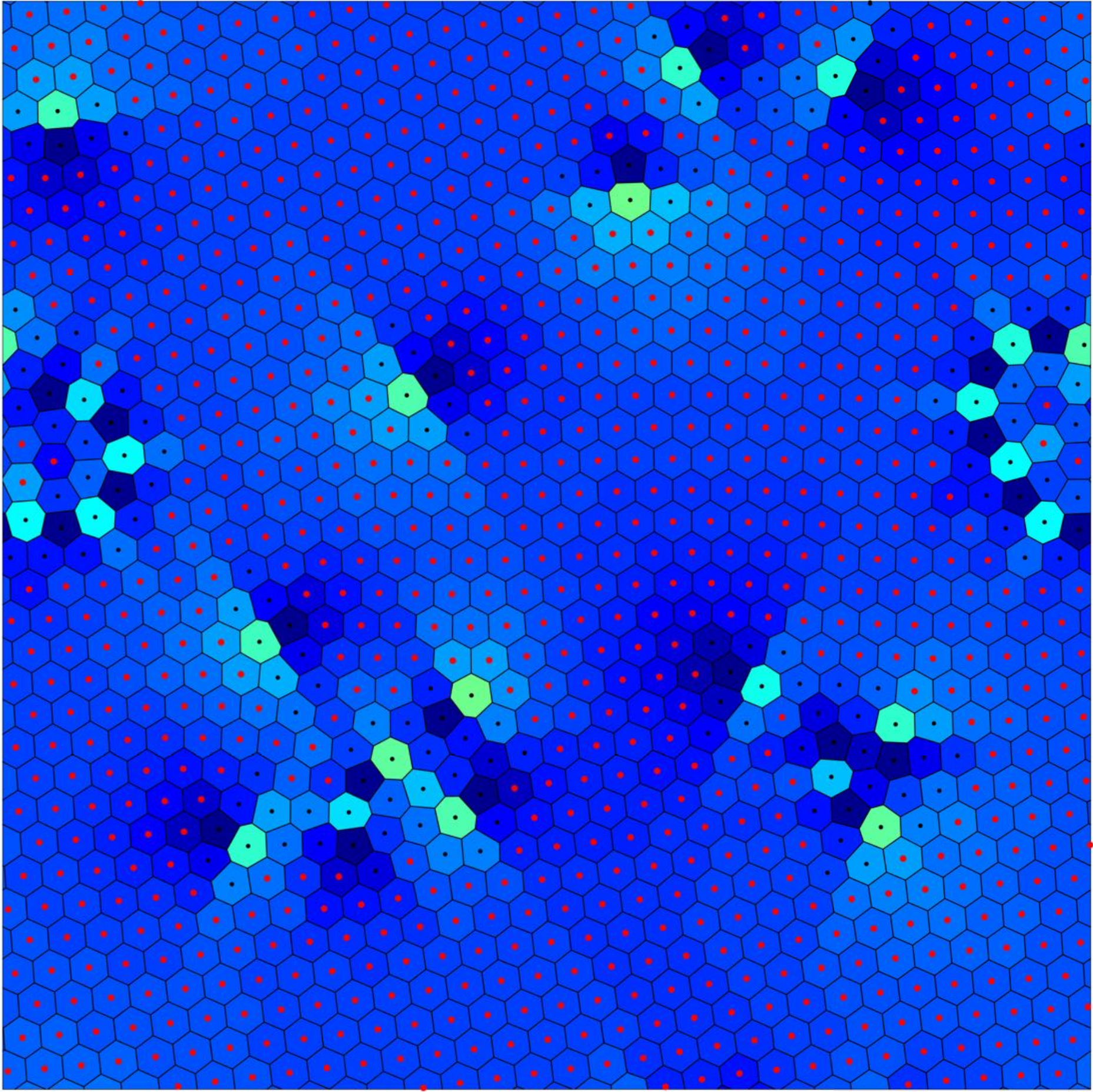} & \includegraphics[width=0.2\columnwidth]{\string"Graphics/N=1000_SquareTorus/MinimalEnergy_GLPLBFGS_N1000\string".pdf}\tabularnewline
{\footnotesize{}(a) Lloyd} & {\footnotesize{}(b) L-BFGS(7)} & {\footnotesize{}(c) P-L-BFGS(20,20)} \tabularnewline
\end{tabular}

\vspace{0.1cm}

\setlength\tabcolsep{0.005\columnwidth} 

\begin{tabular}{|ccc||ccc|}
\hline
{\footnotesize{} \textit{MACN-c}} & {\footnotesize{} Lloyd block} & {\footnotesize{} \textit{MACN-}$\delta$} & {\footnotesize{} \textit{MACN-c}} & {\footnotesize{} Lloyd block} & {\footnotesize{} \textit{MACN-}$\delta$} \tabularnewline
\includegraphics[width=0.15\columnwidth]{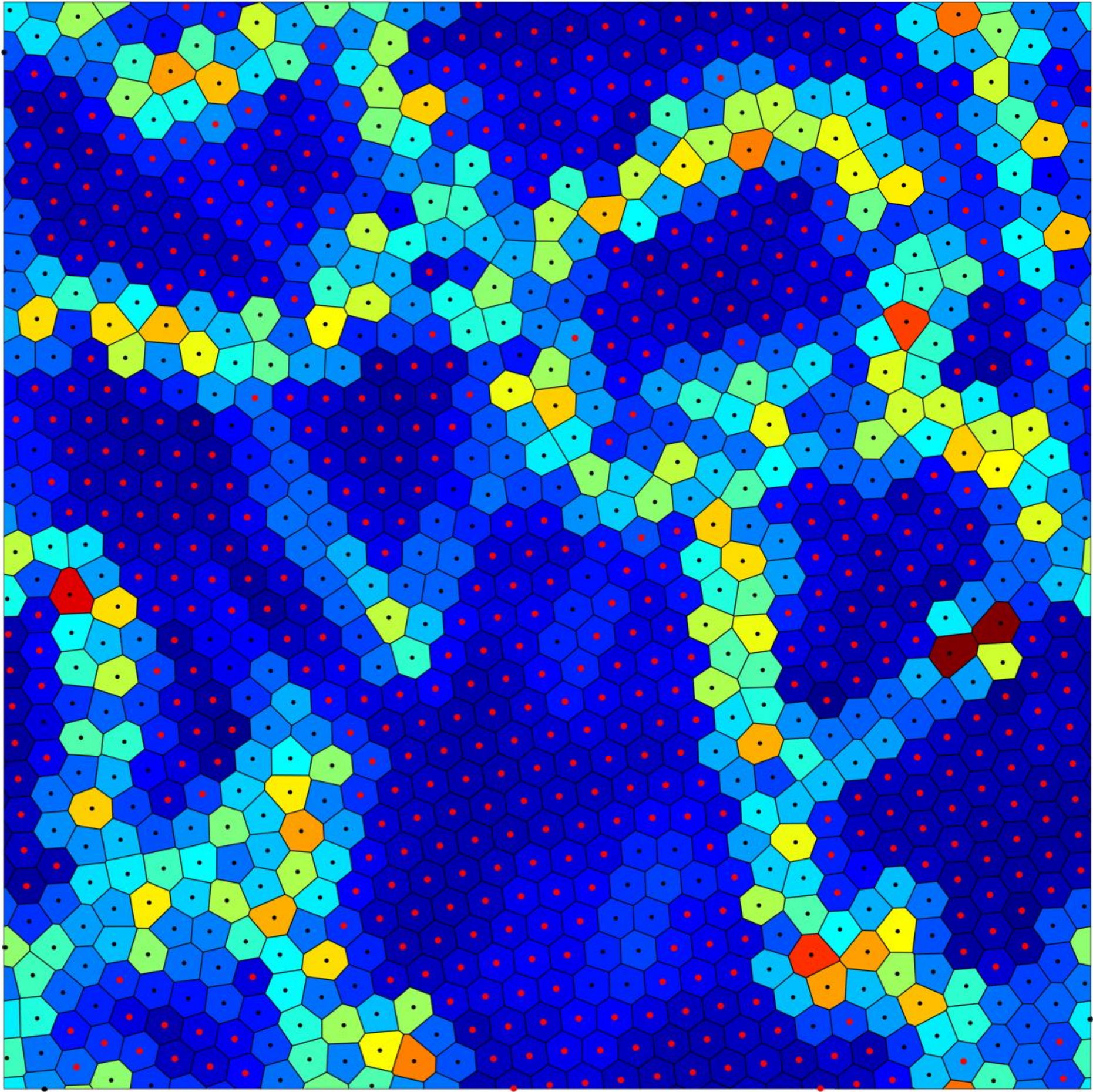} & \includegraphics[width=0.15\columnwidth]{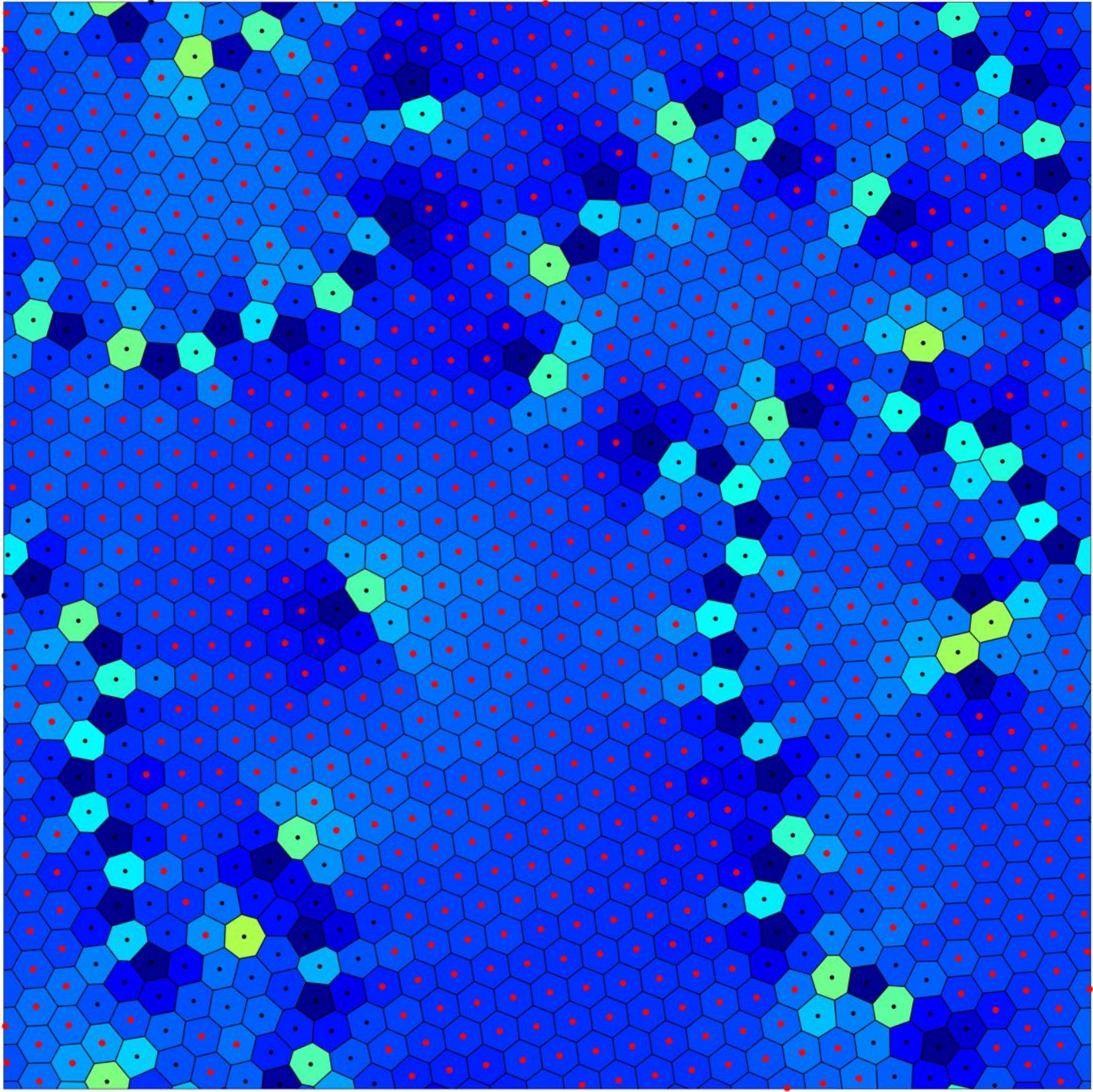} & \includegraphics[width=0.15\columnwidth]{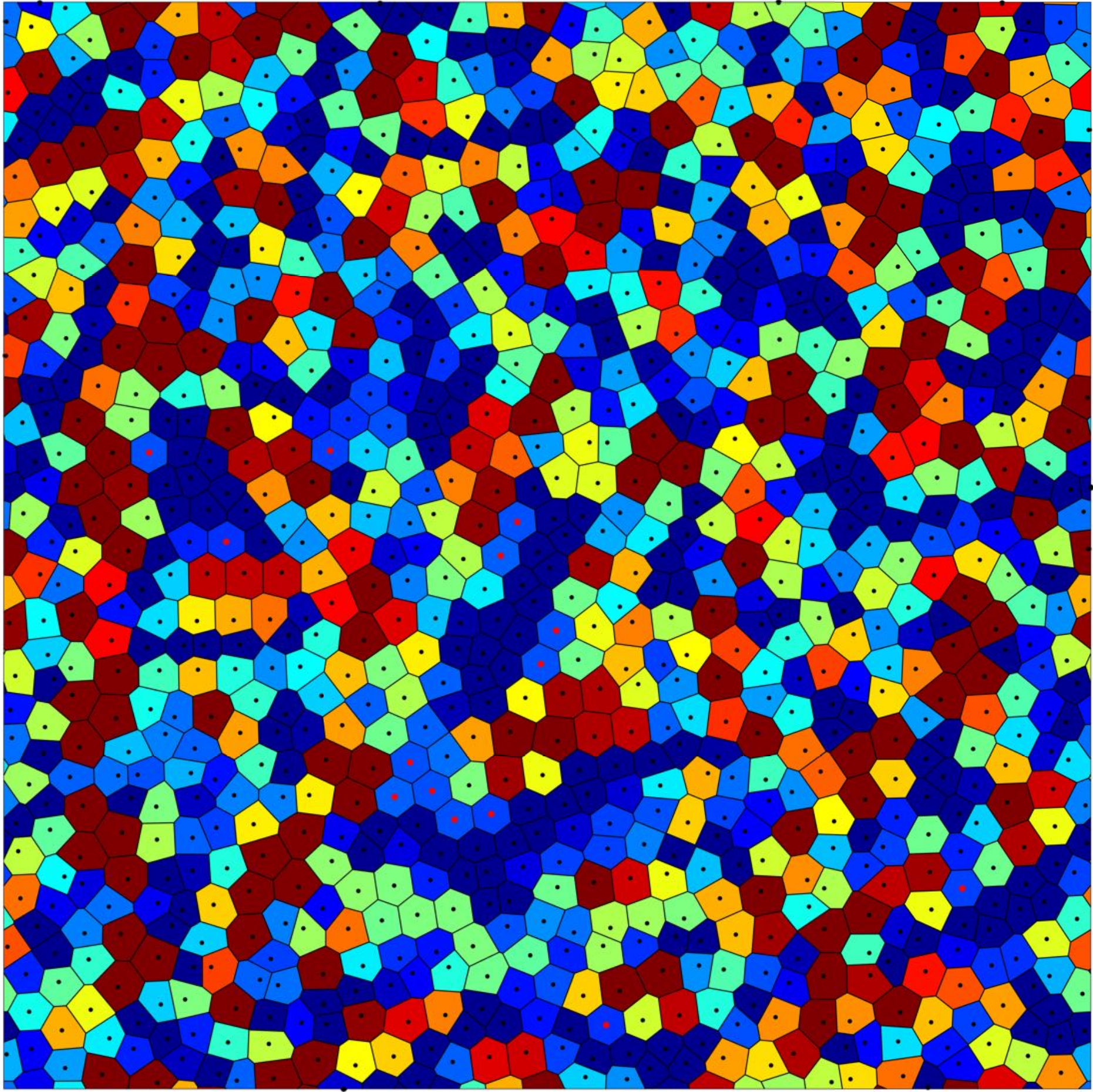} & \includegraphics[width=0.15\columnwidth]{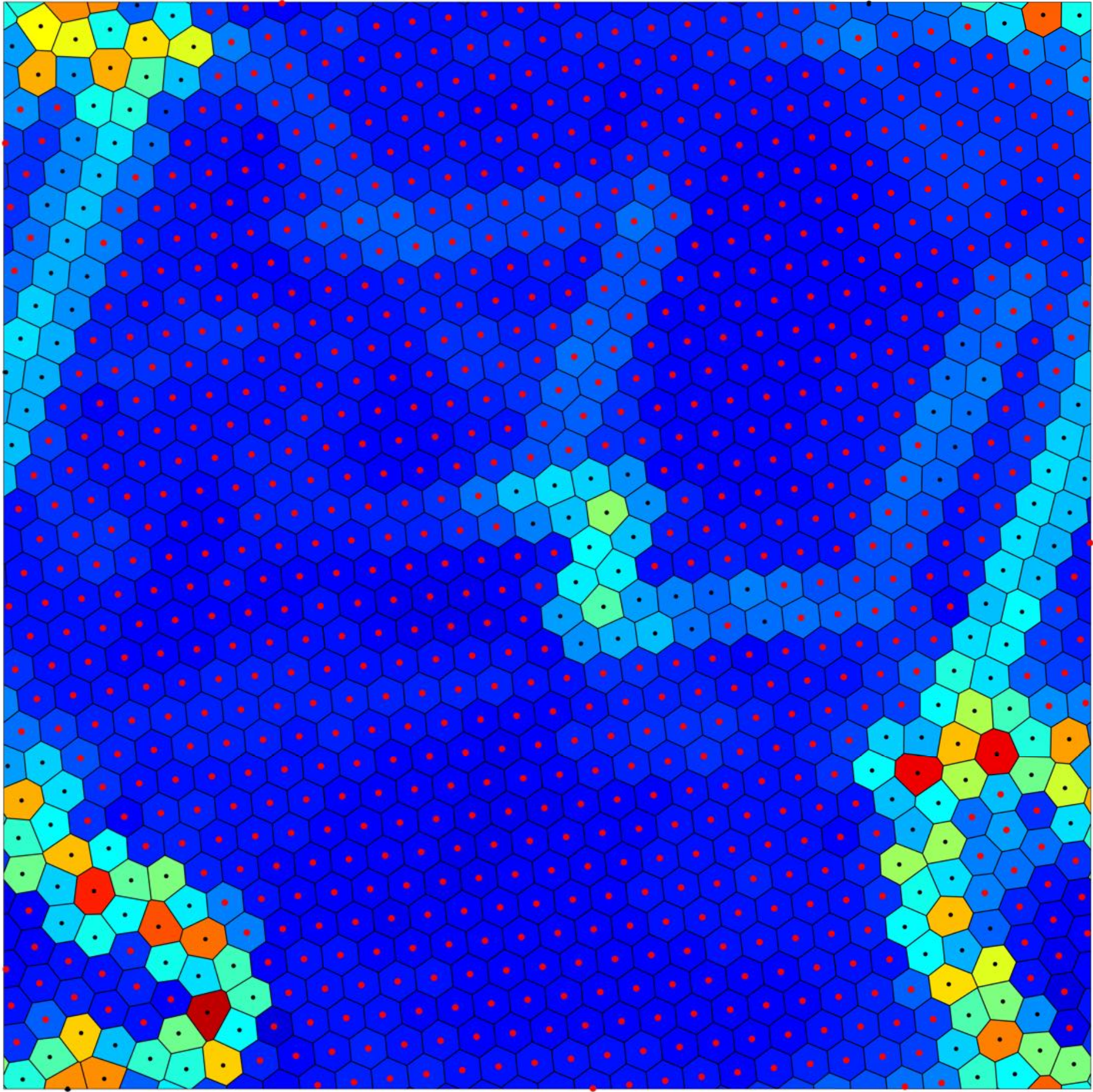} & \includegraphics[width=0.15\columnwidth]{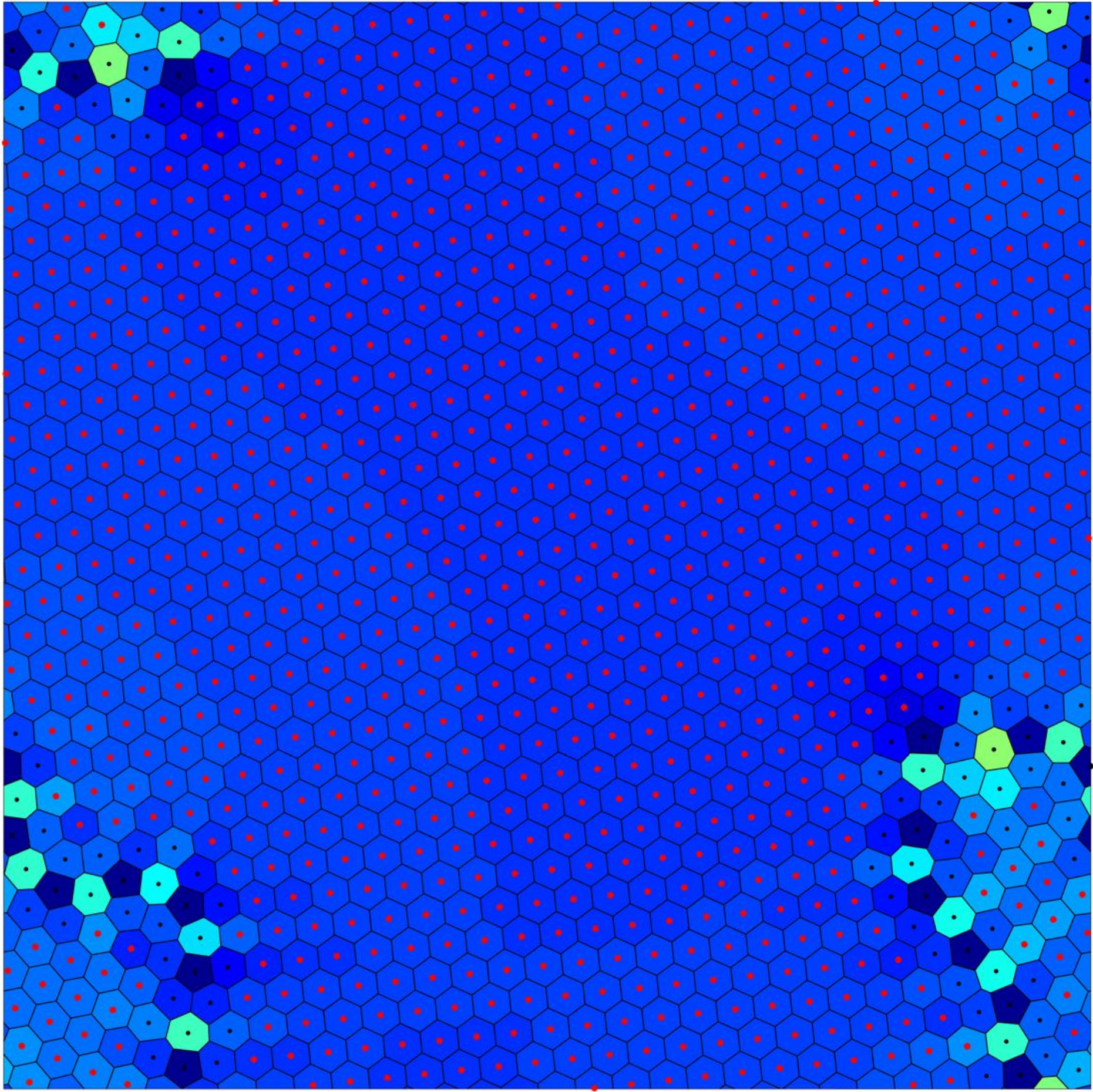} & \includegraphics[width=0.15\columnwidth]{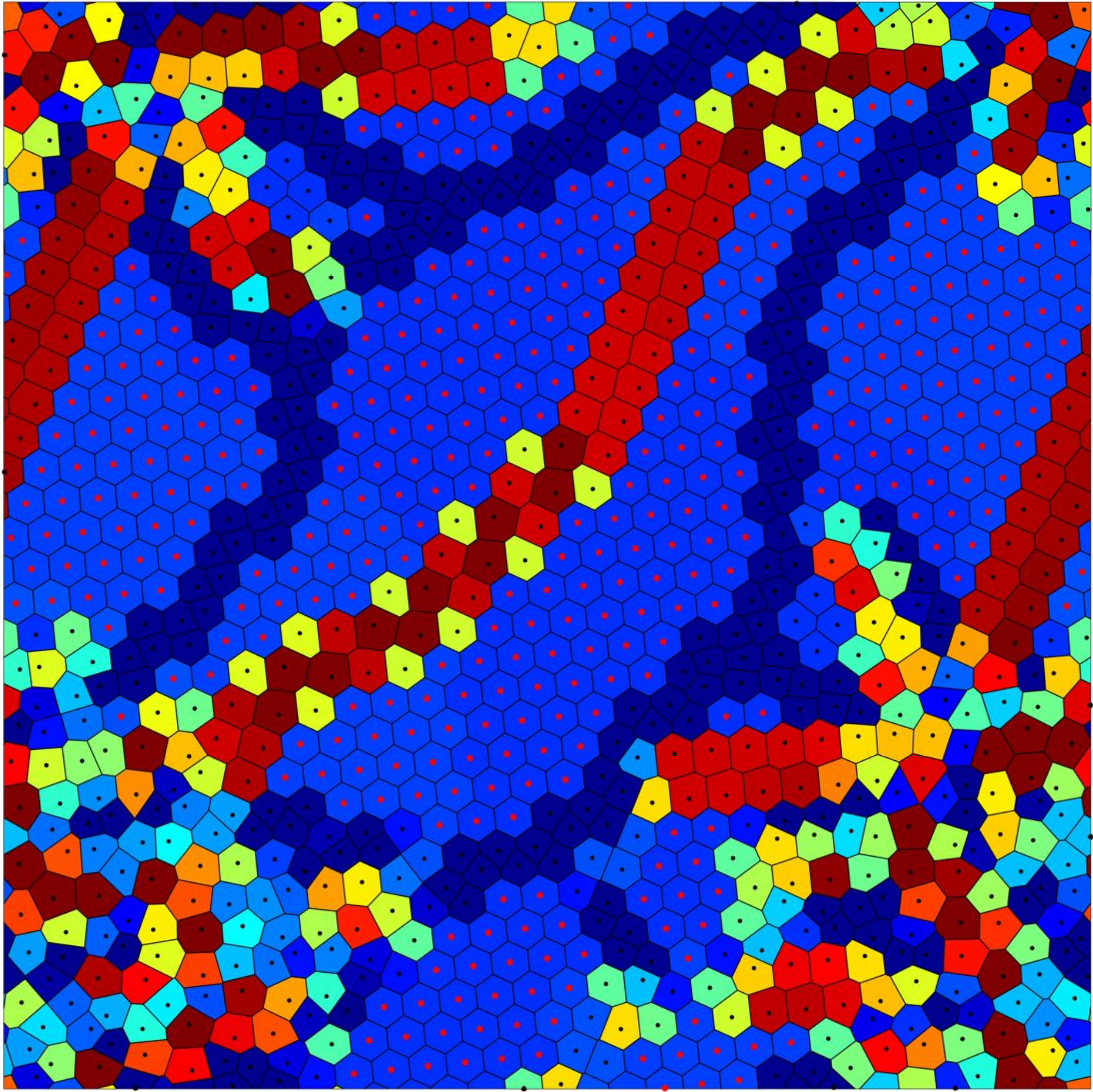}\tabularnewline
\includegraphics[width=0.15\columnwidth]{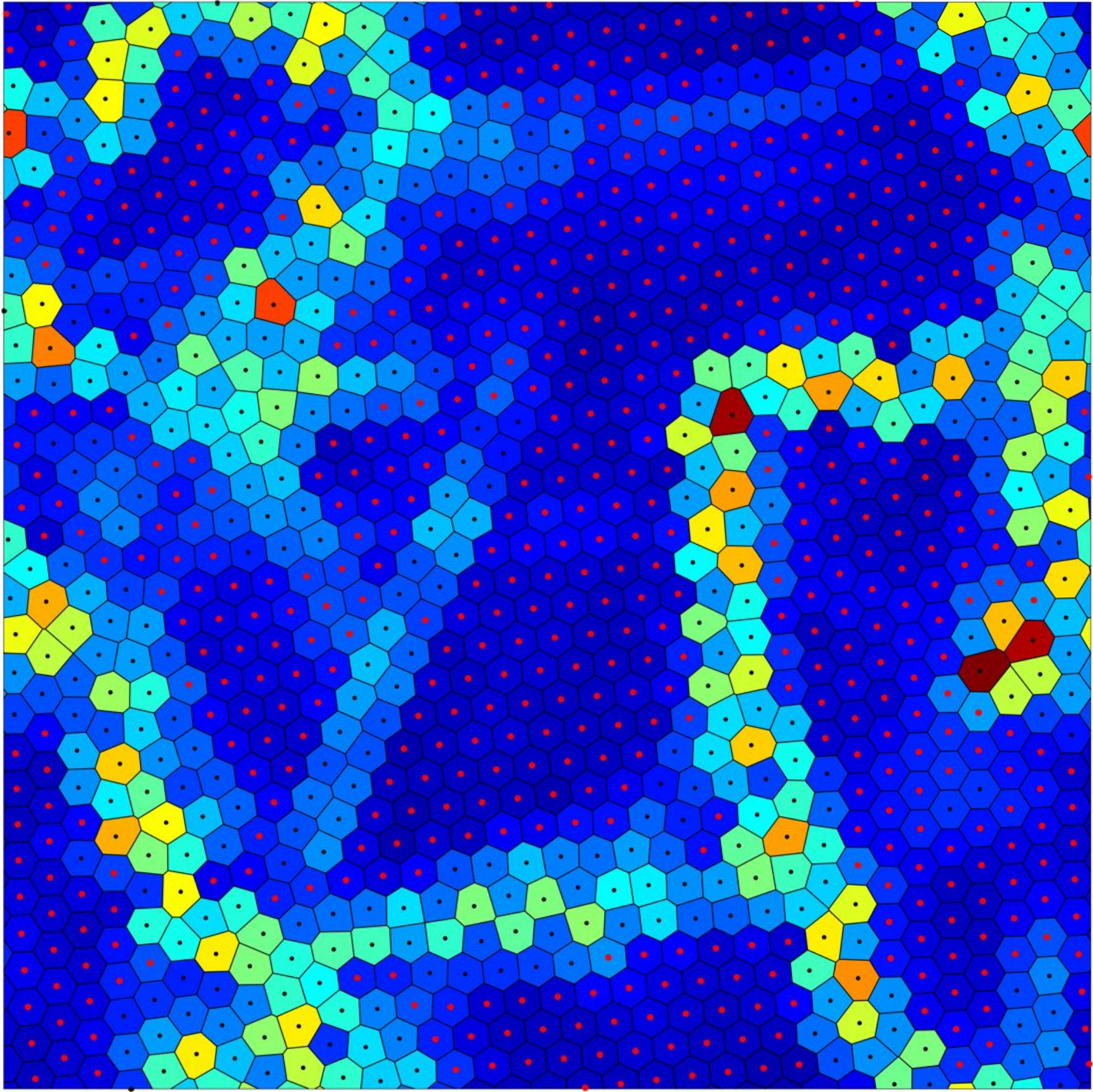} & \includegraphics[width=0.15\columnwidth]{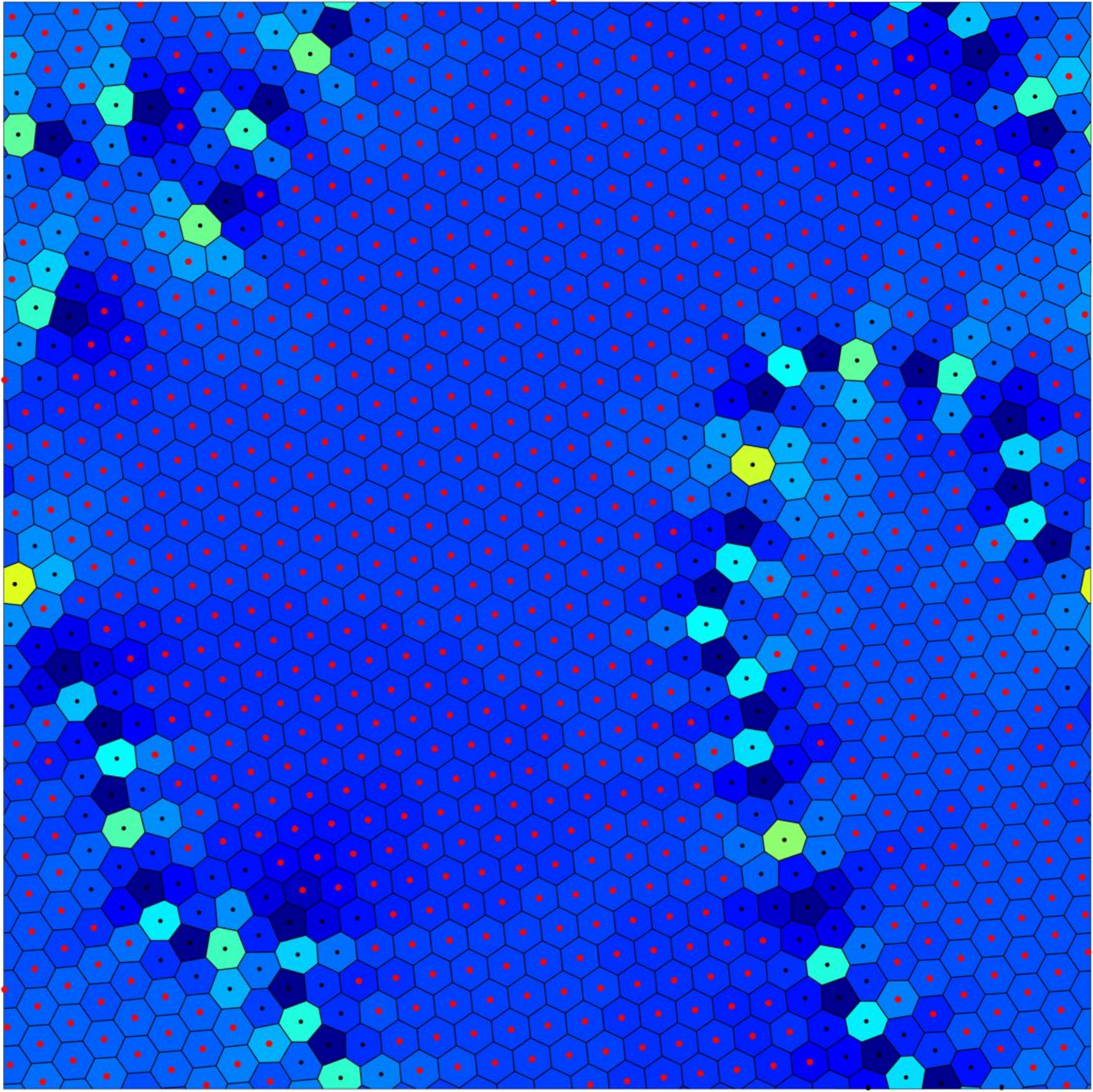} & \includegraphics[width=0.15\columnwidth]{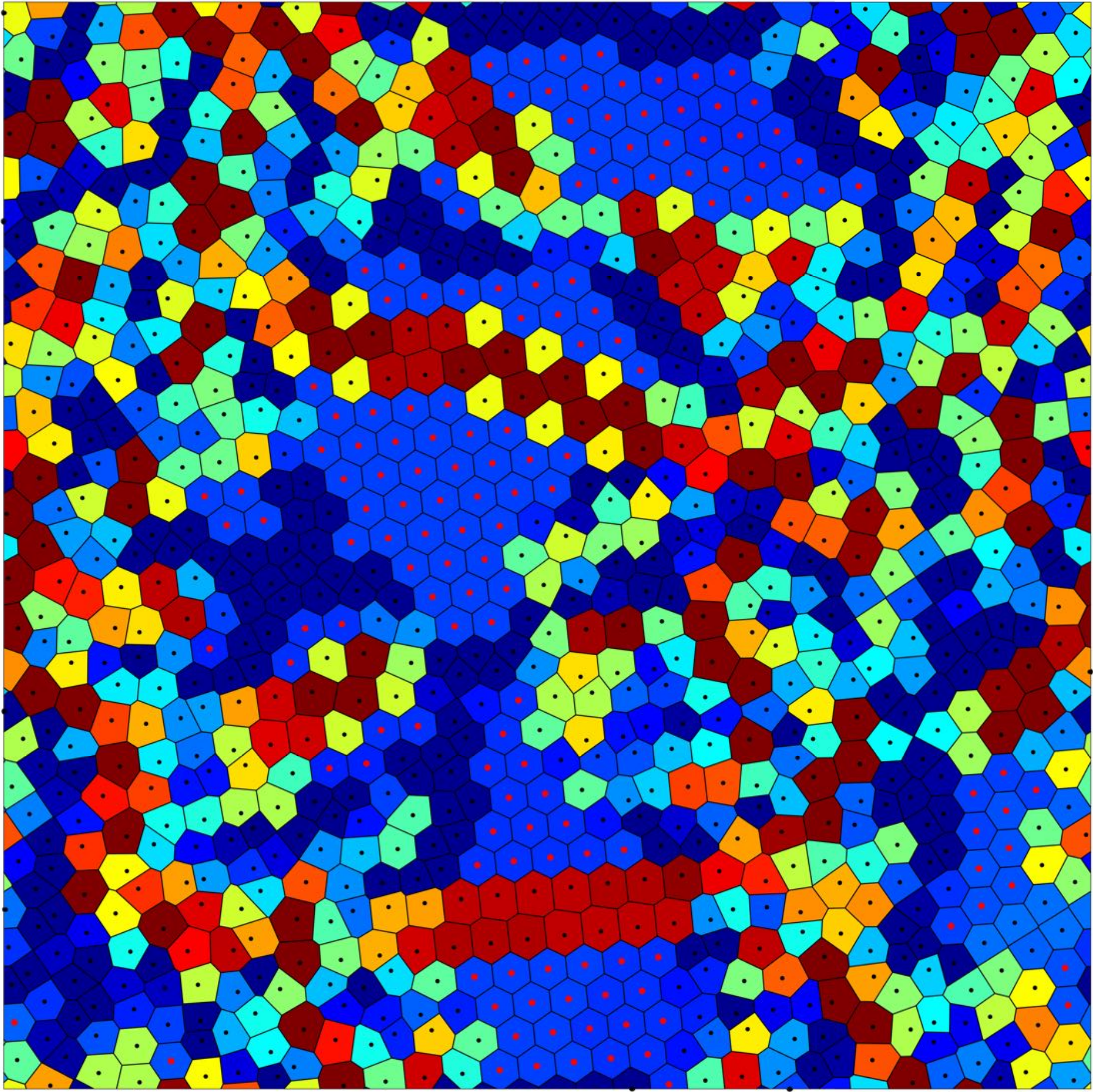} & \includegraphics[width=0.15\columnwidth]{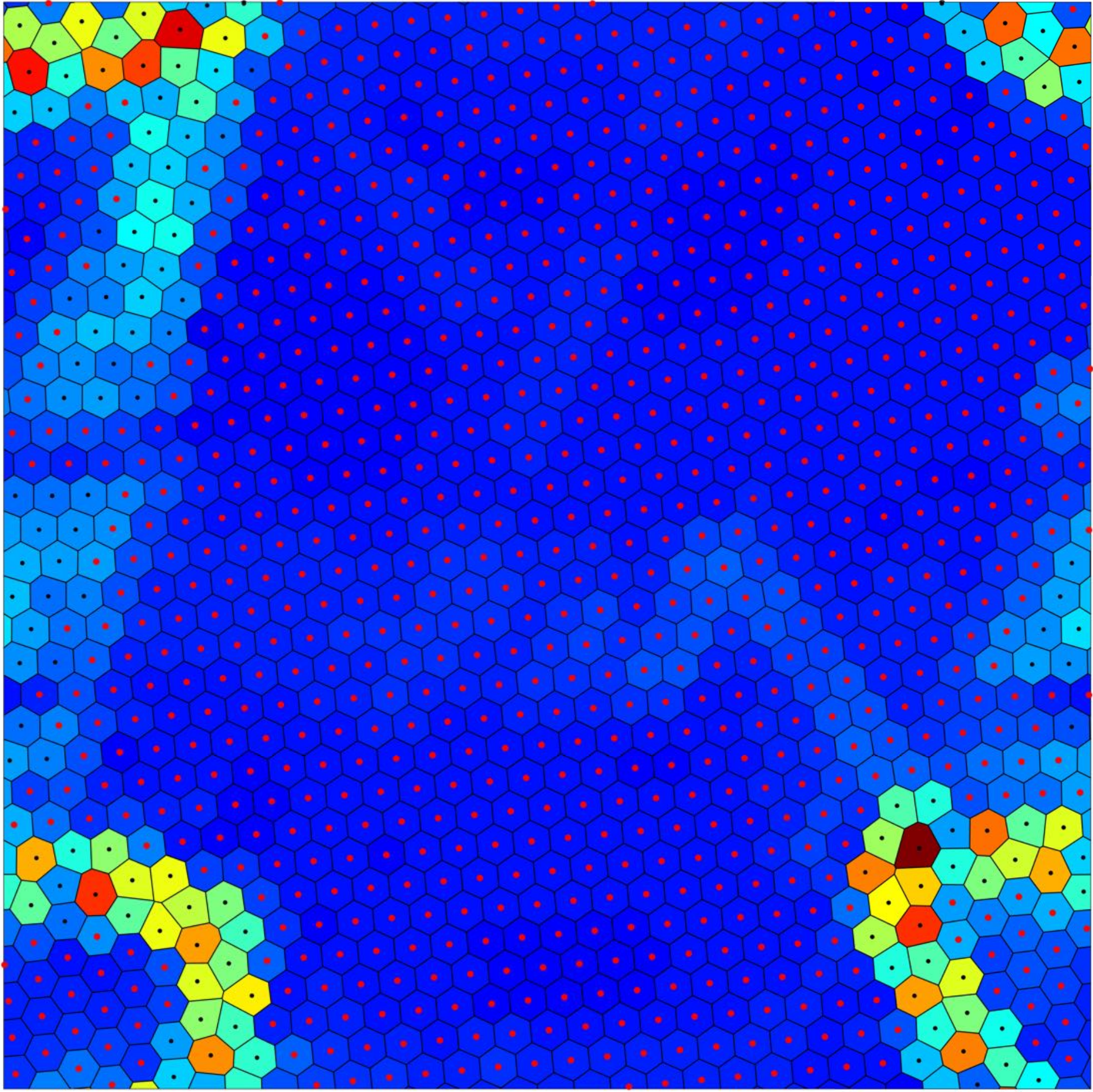} & \includegraphics[width=0.15\columnwidth]{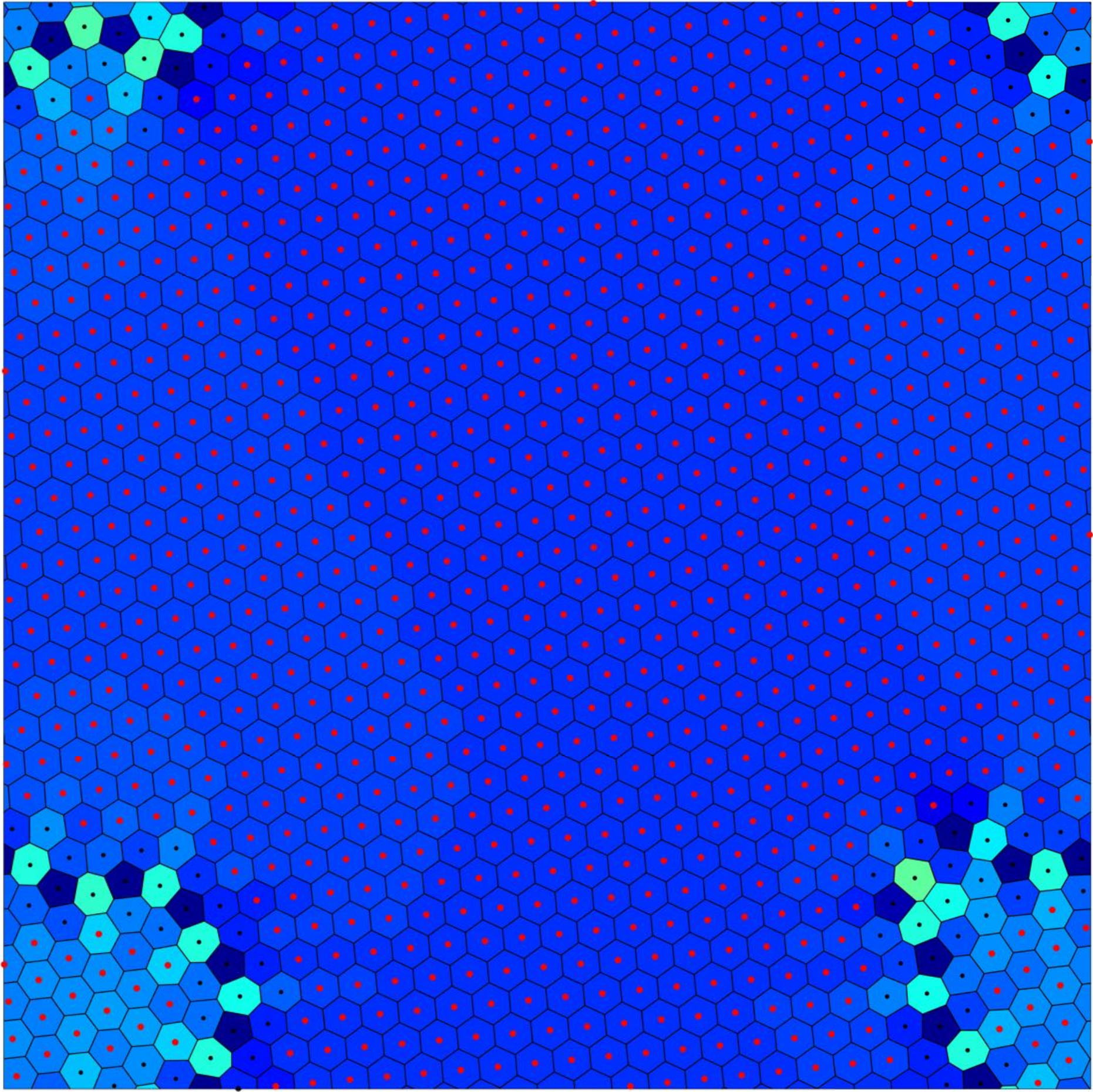} & \includegraphics[width=0.15\columnwidth]{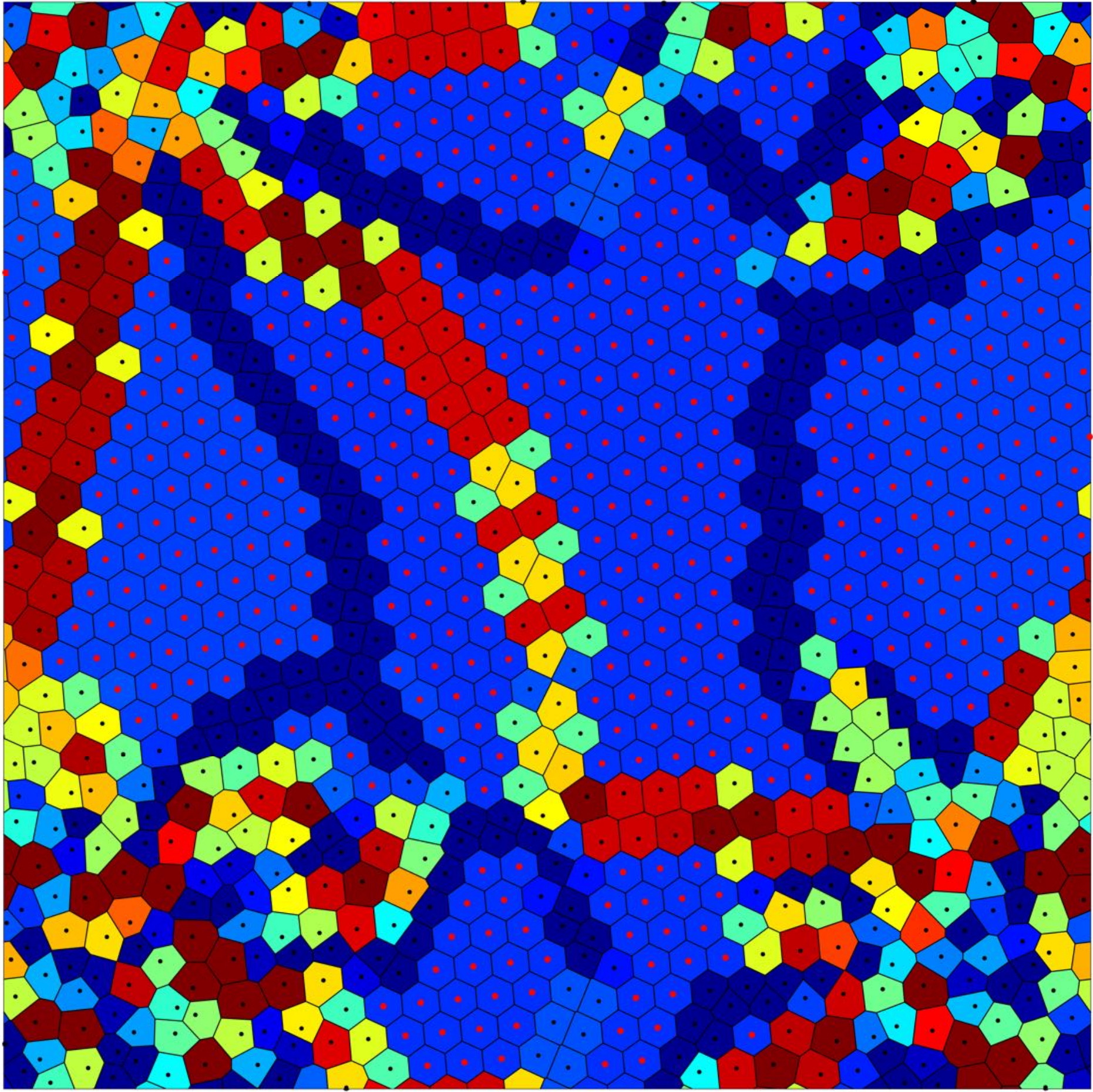}\tabularnewline
\includegraphics[width=0.15\columnwidth]{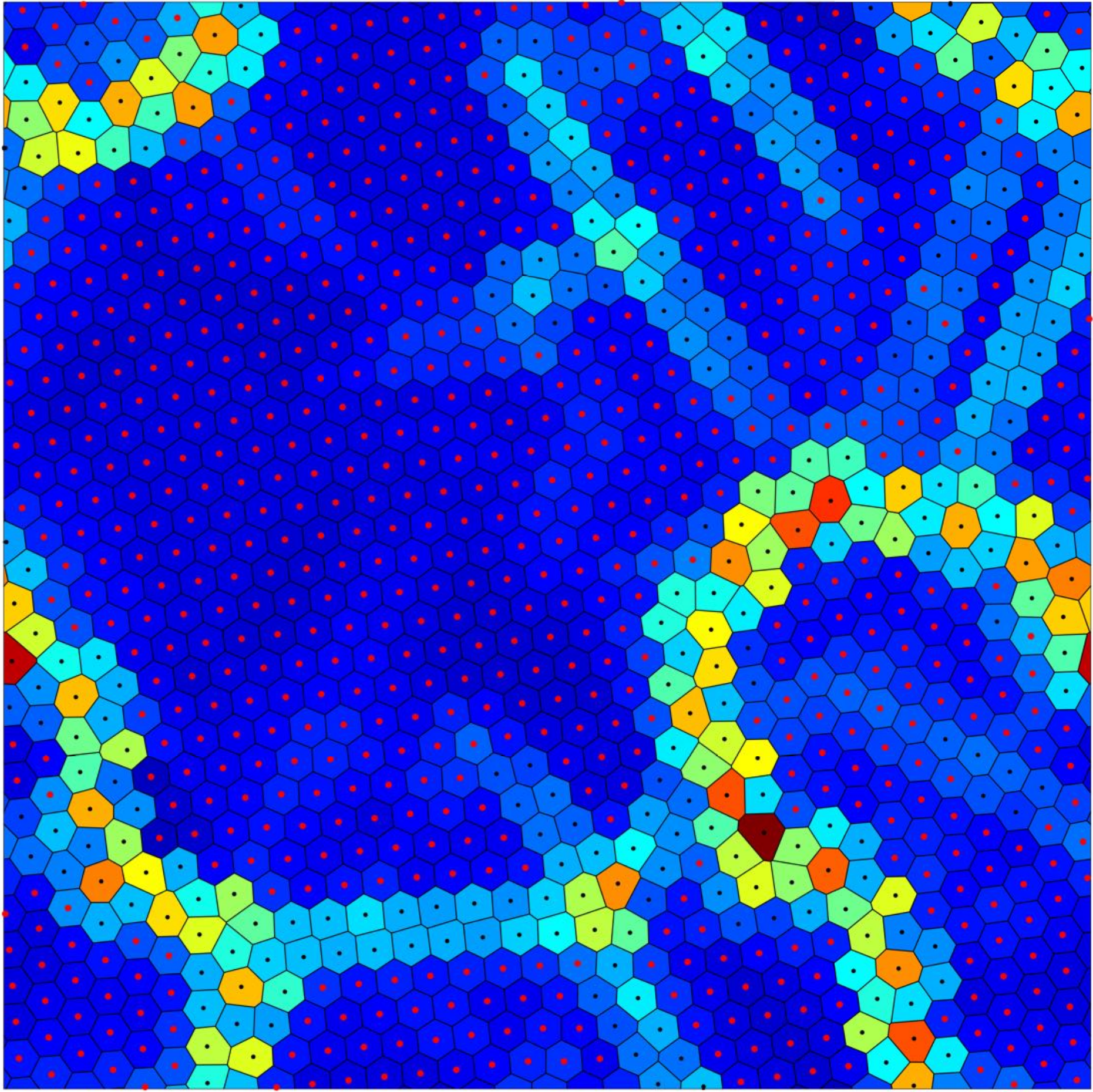} & \includegraphics[width=0.15\columnwidth]{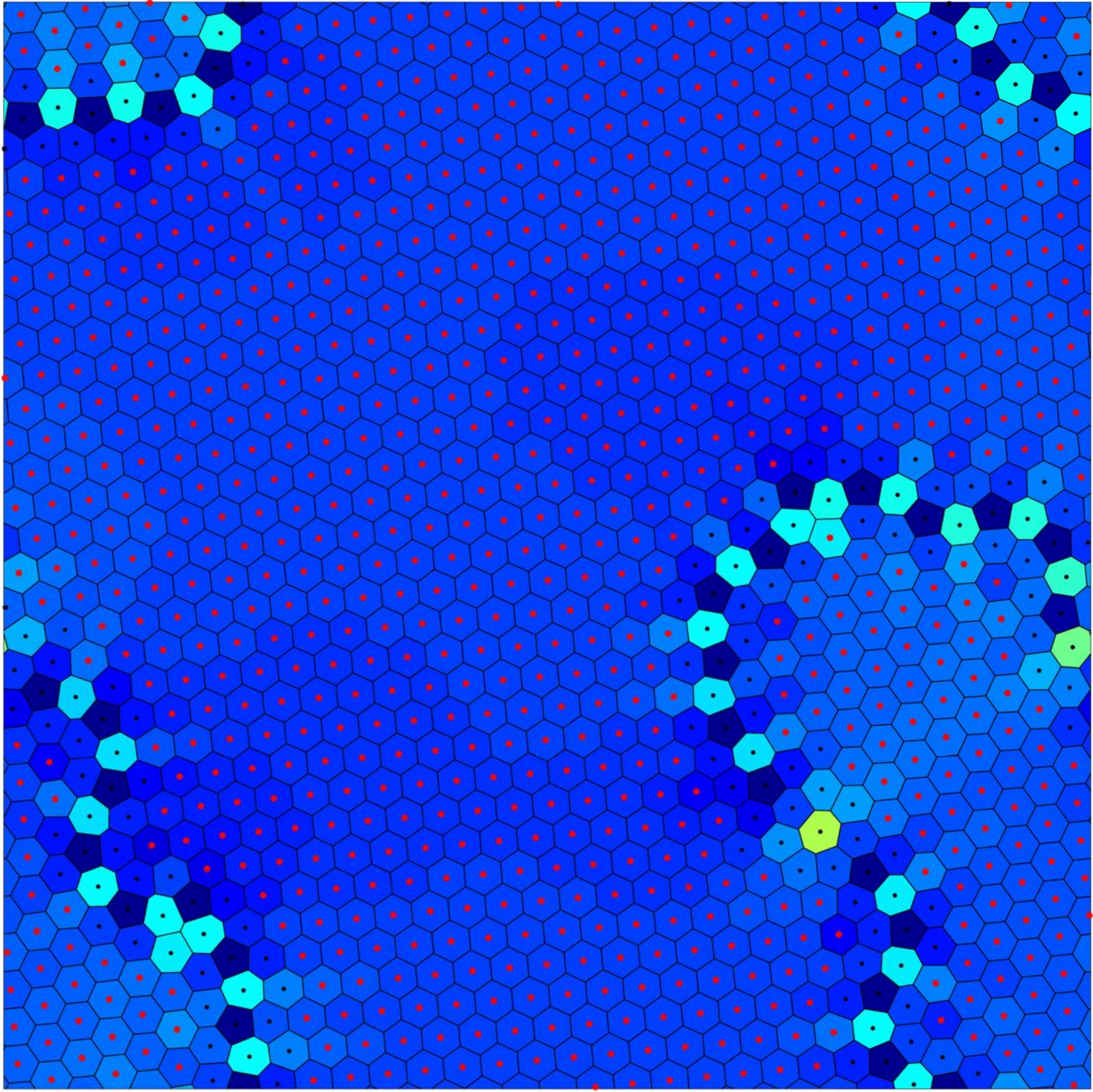} & \includegraphics[width=0.15\columnwidth]{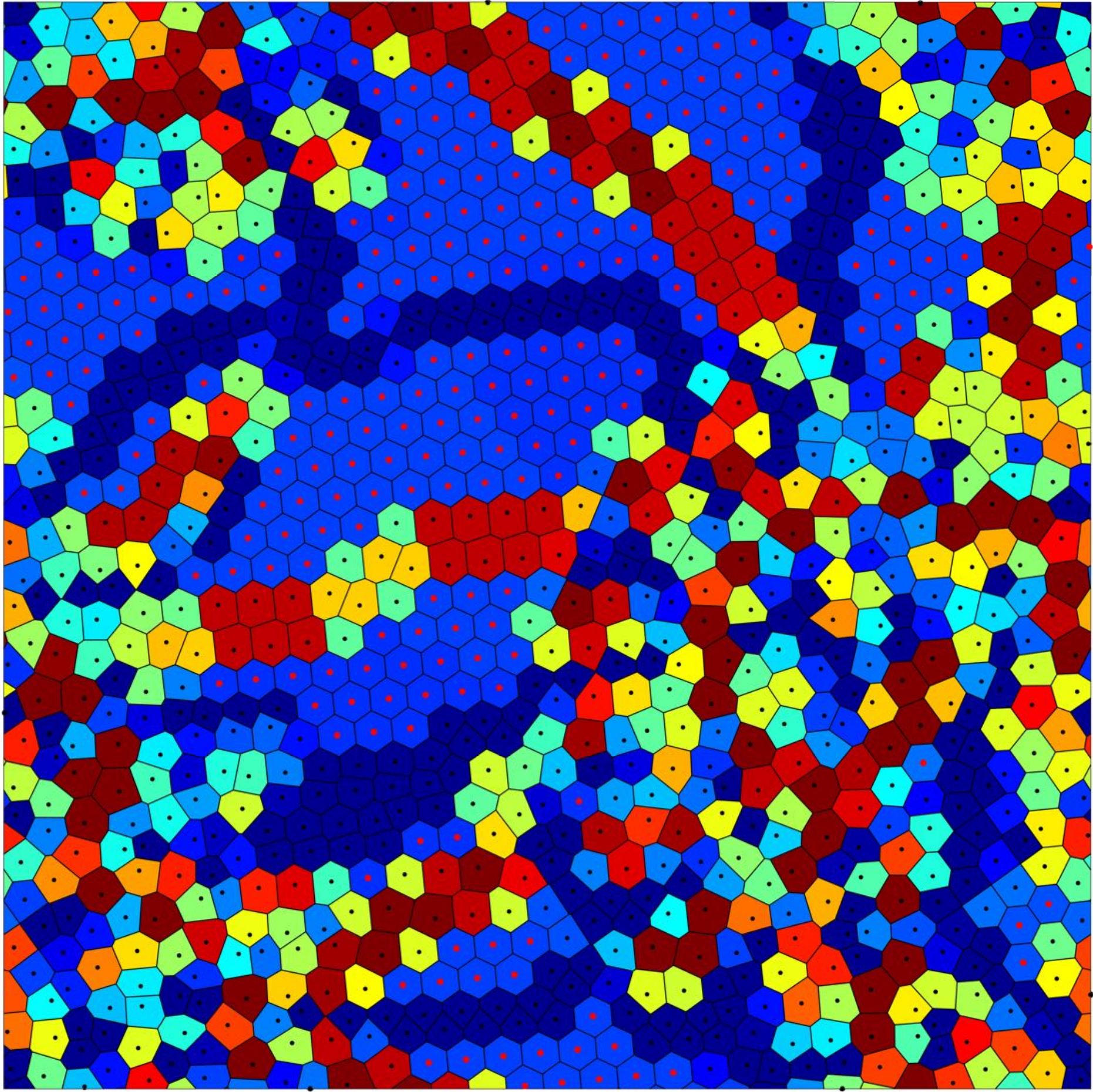} & \includegraphics[width=0.15\columnwidth]{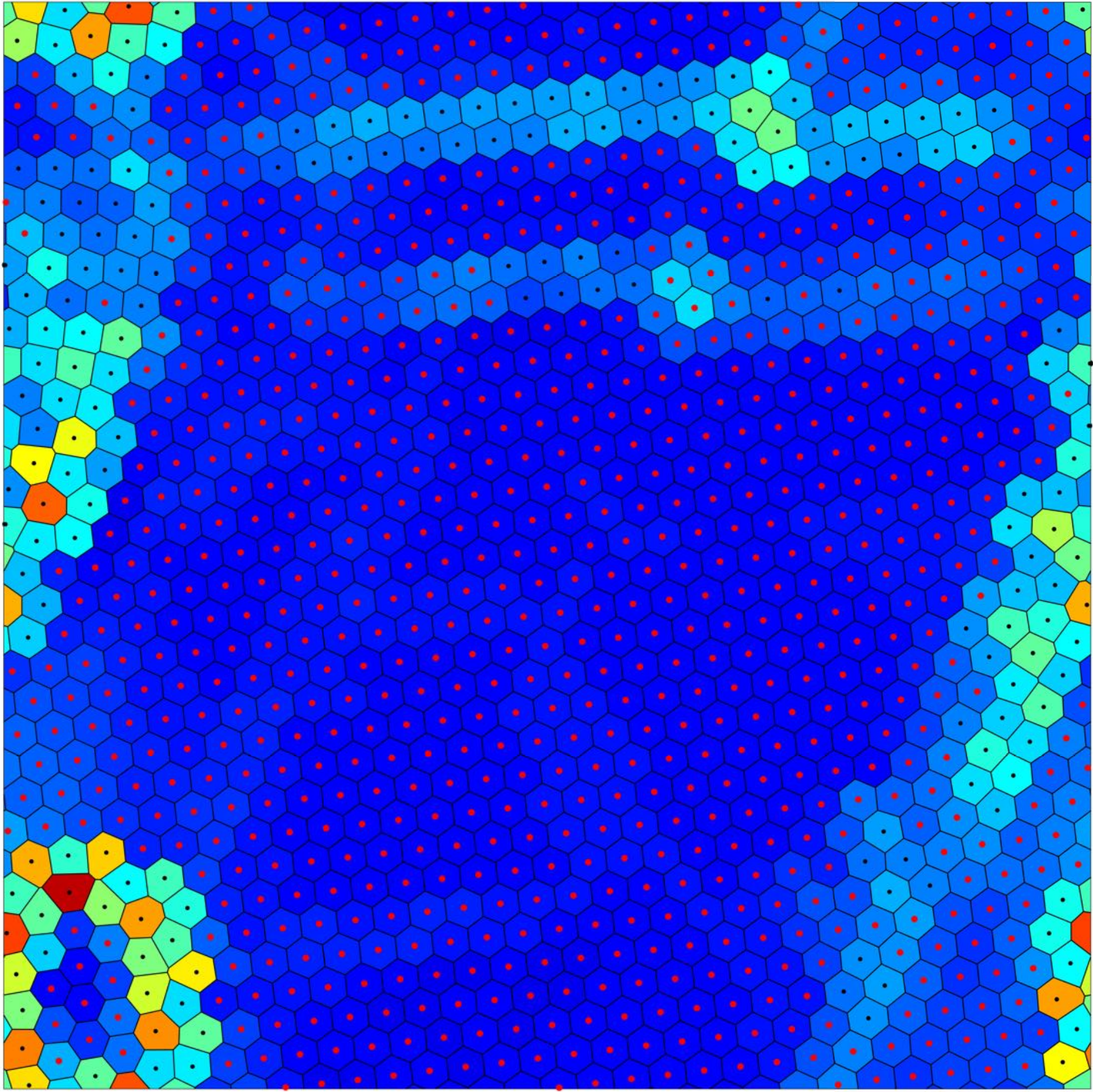} & \includegraphics[width=0.15\columnwidth]{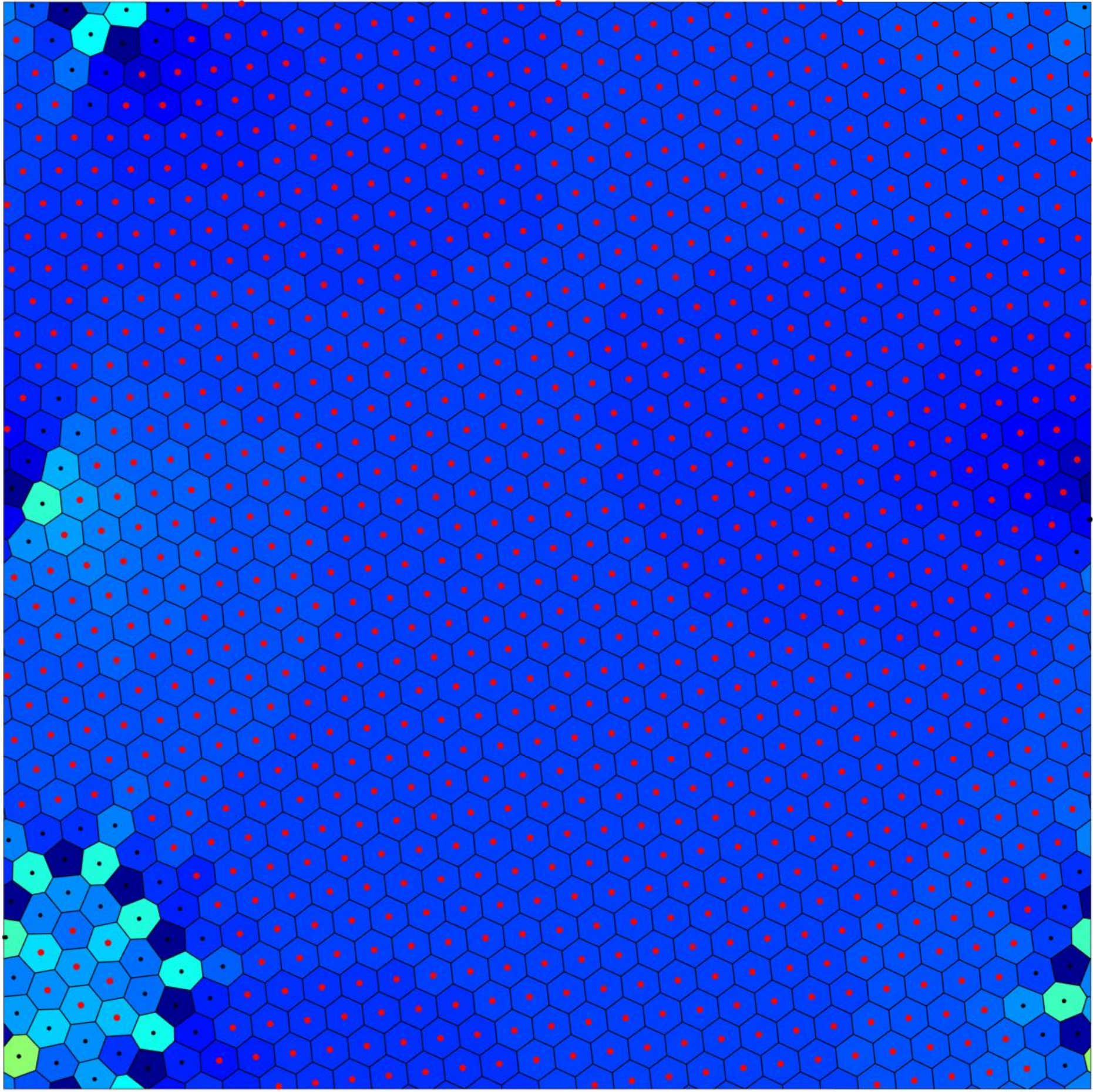} & \includegraphics[width=0.15\columnwidth]{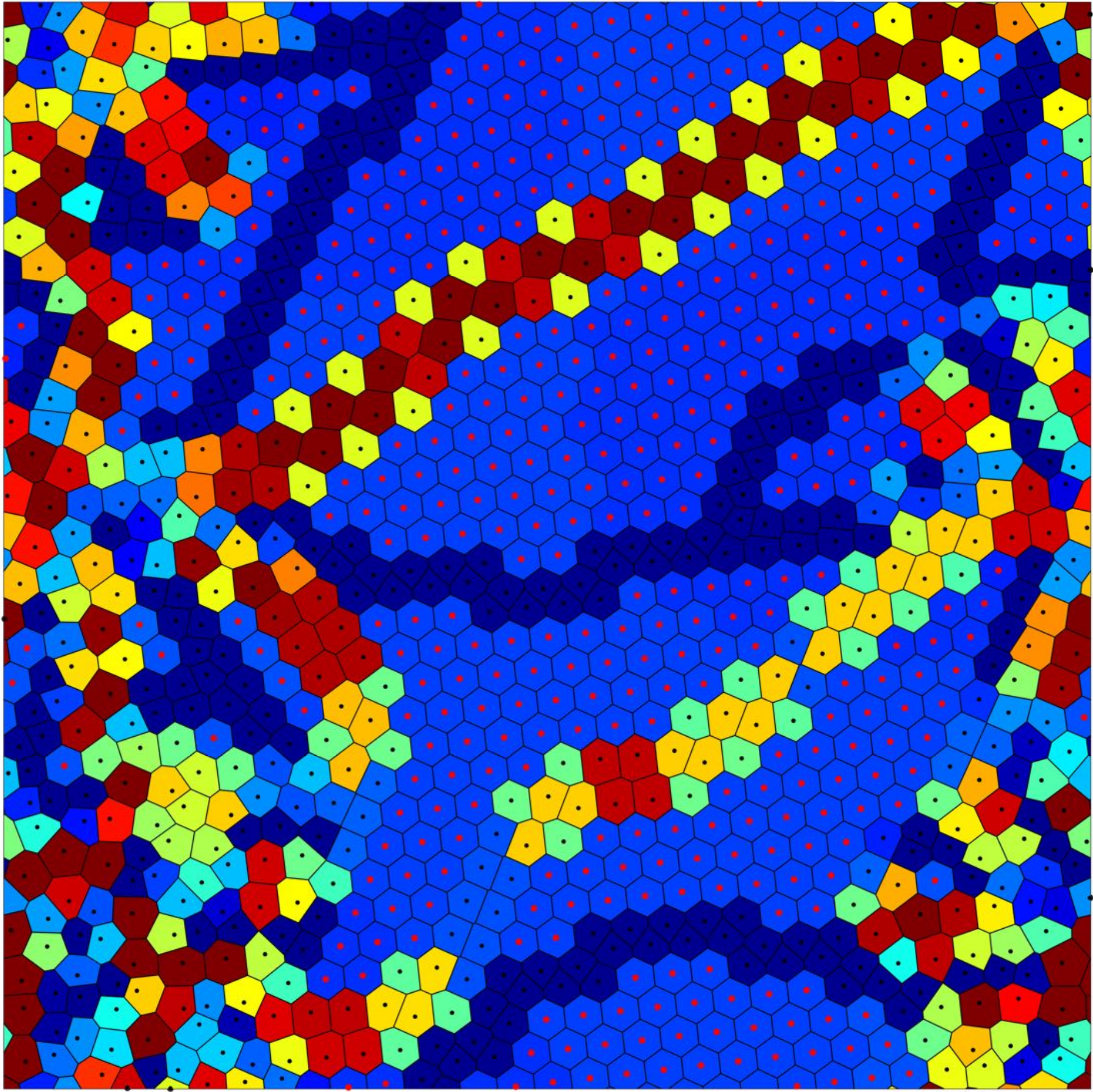}\tabularnewline
\includegraphics[width=0.15\columnwidth]{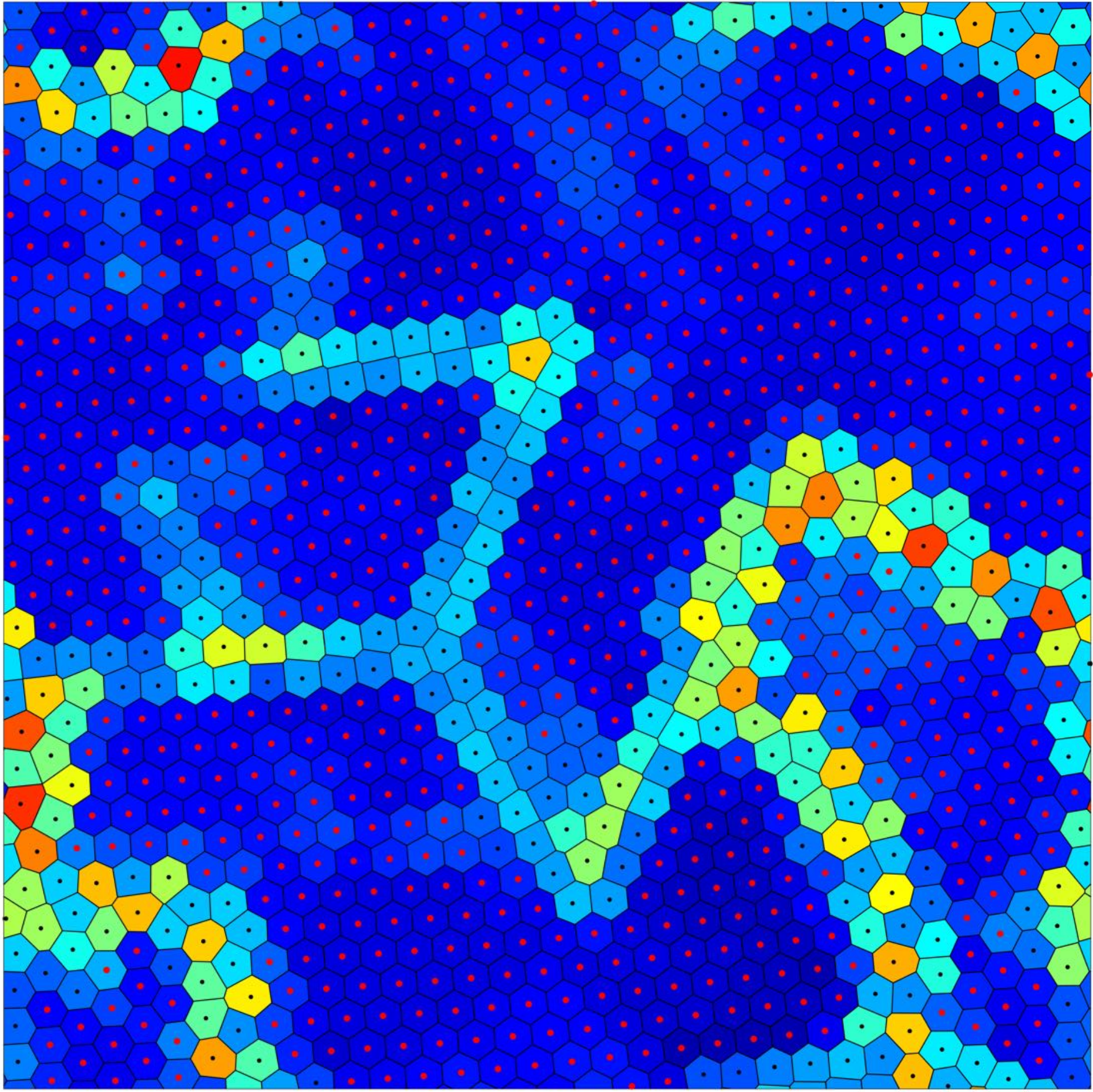} & \includegraphics[width=0.15\columnwidth]{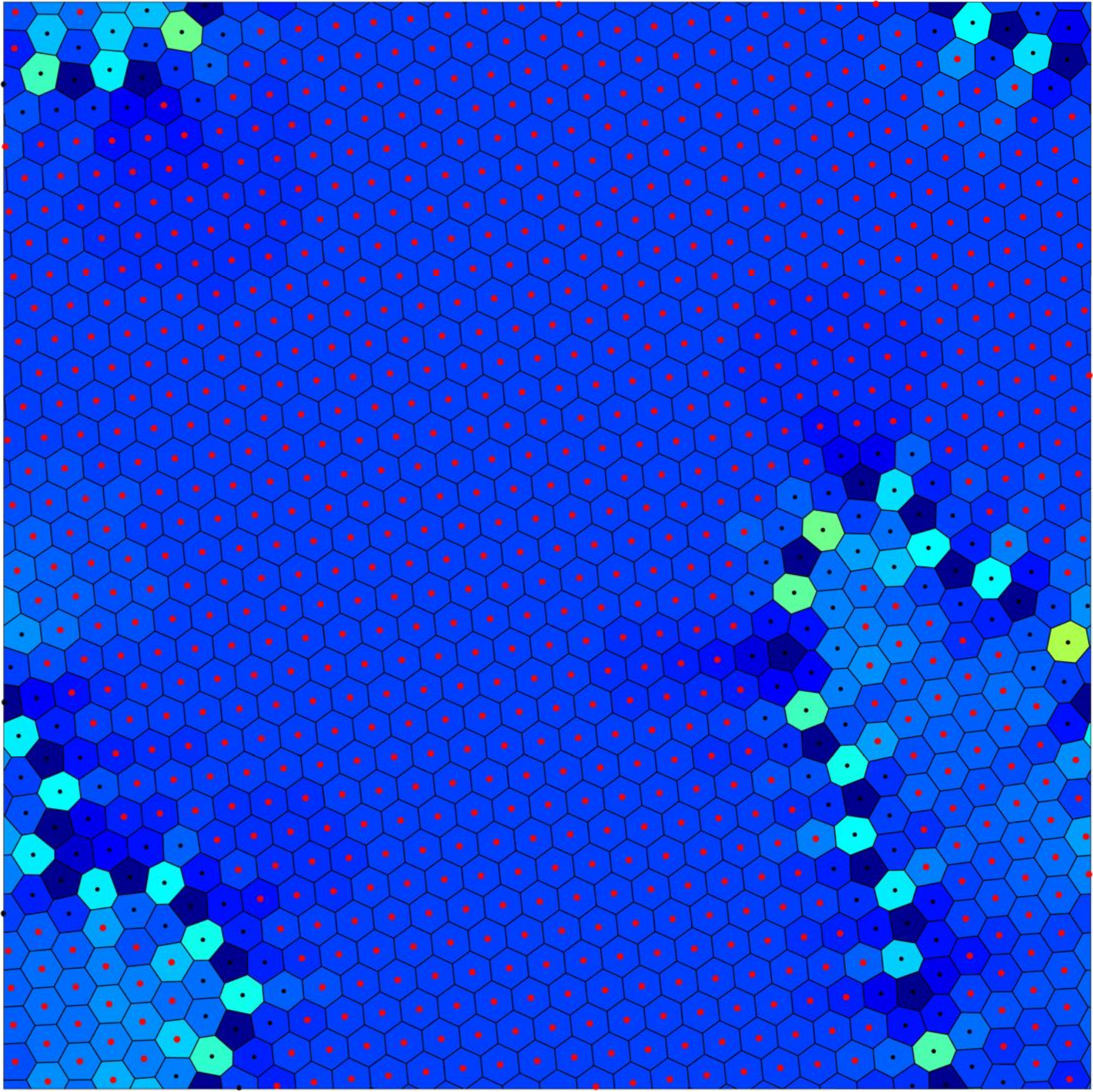} & \includegraphics[width=0.15\columnwidth]{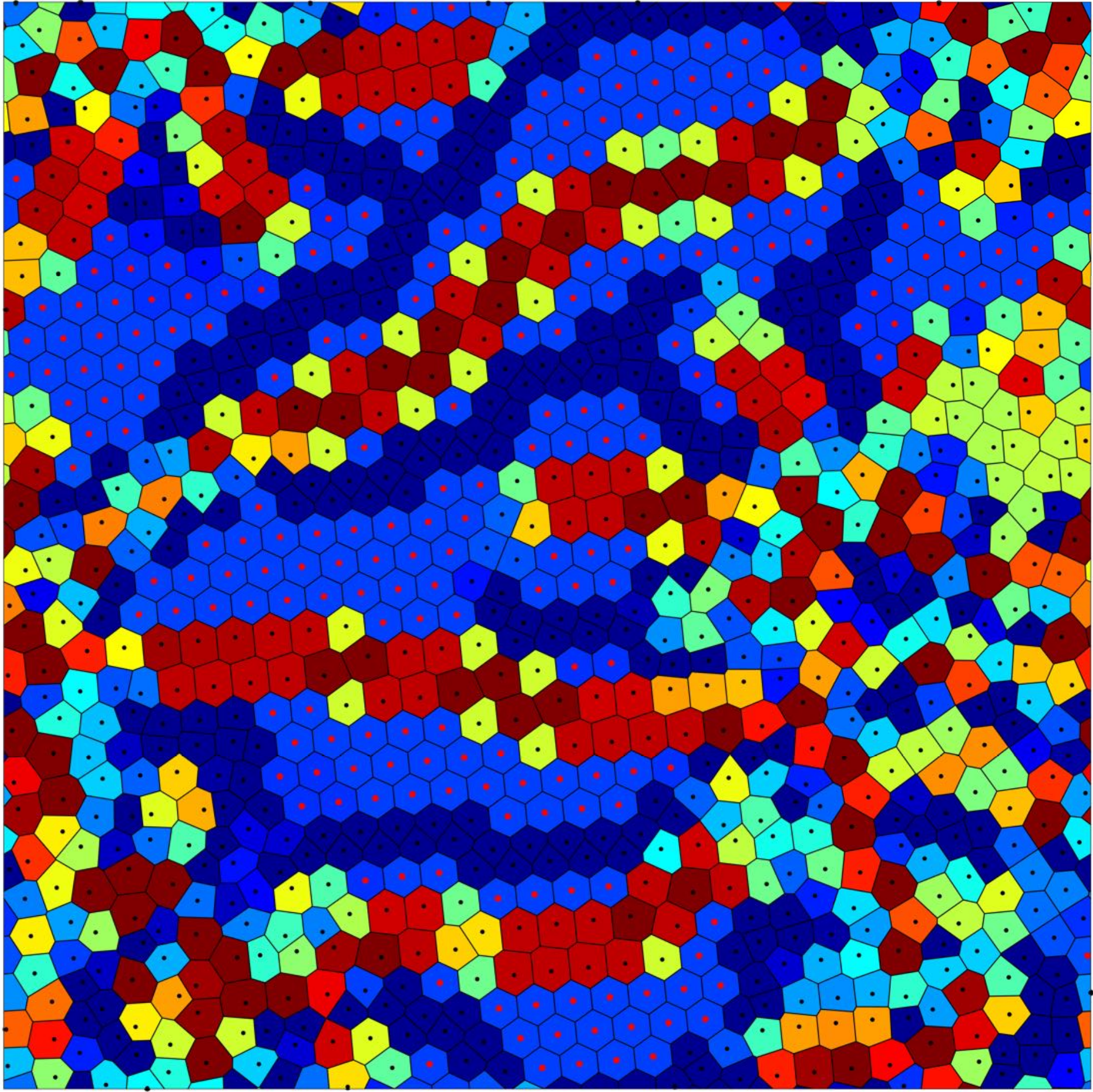} & \includegraphics[width=0.15\columnwidth]{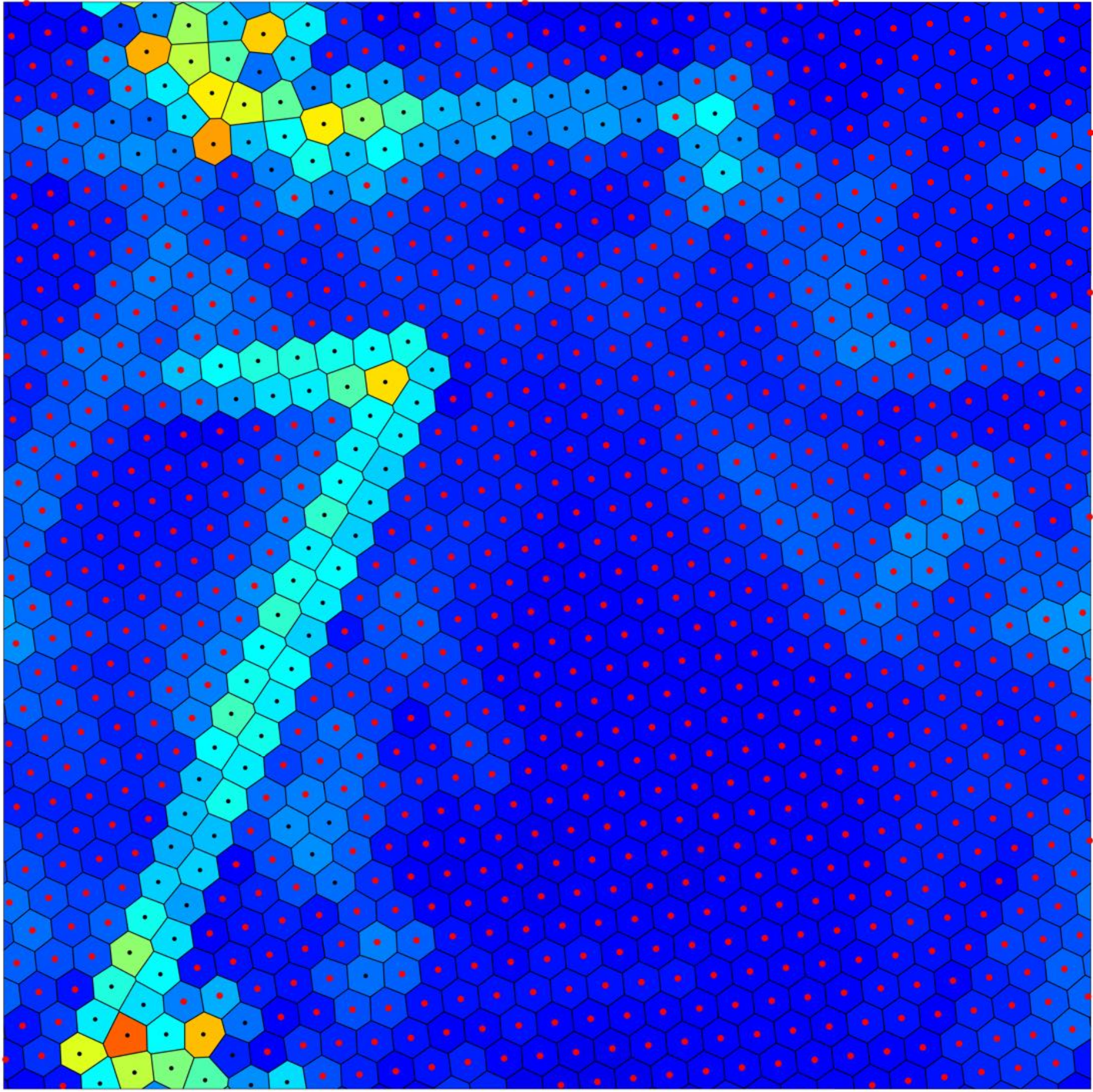} & %

\fcolorbox{white}{red}{\begin{minipage}[t]{0.15\columnwidth}%

\includegraphics[width=0.993\columnwidth]{\string"Graphics/N=1000_SquareTorus/L_9\string".pdf}%
\end{minipage}} & \includegraphics[width=0.15\columnwidth]{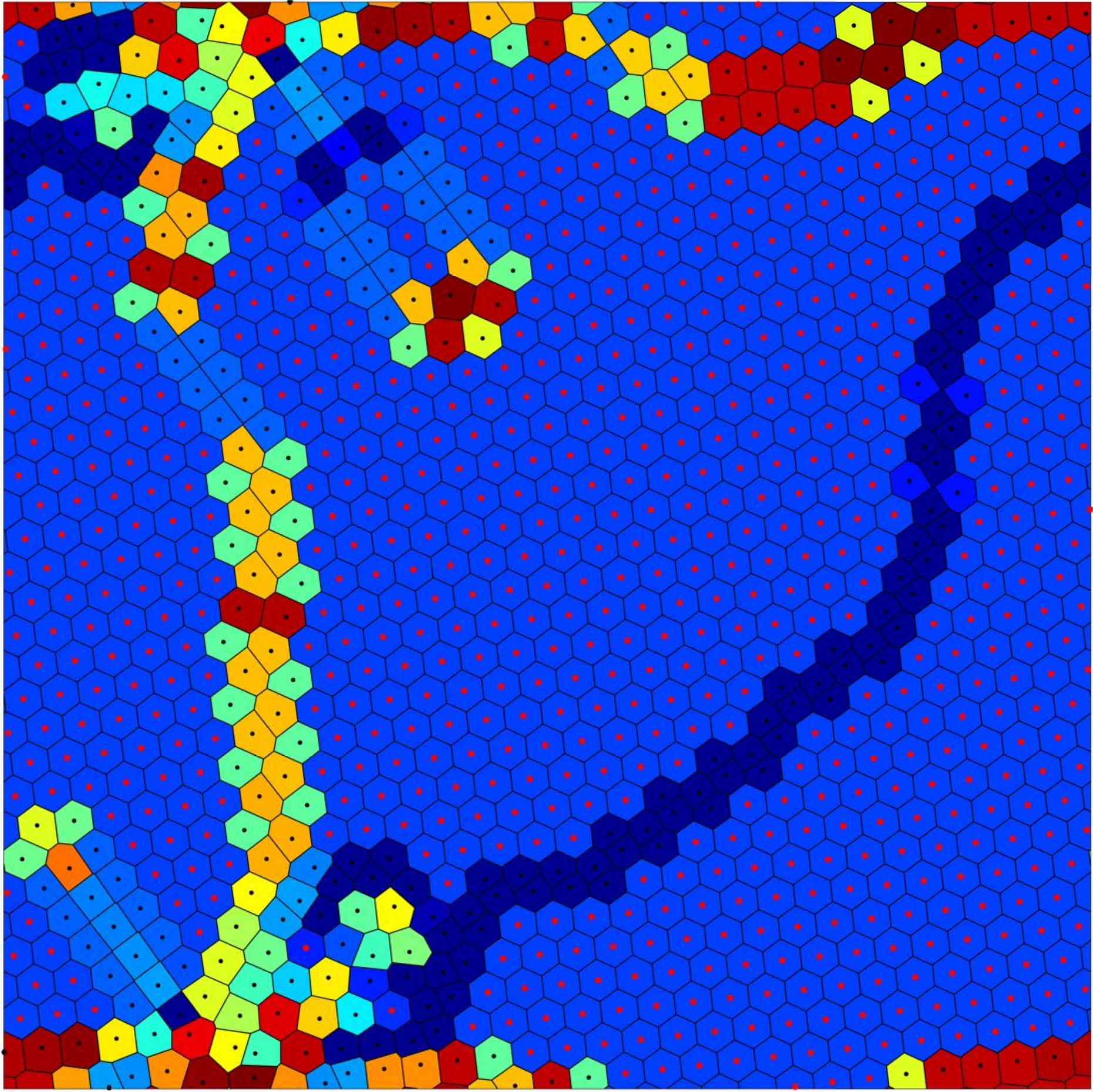}\tabularnewline
\cline{6-6} 
\includegraphics[width=0.15\columnwidth]{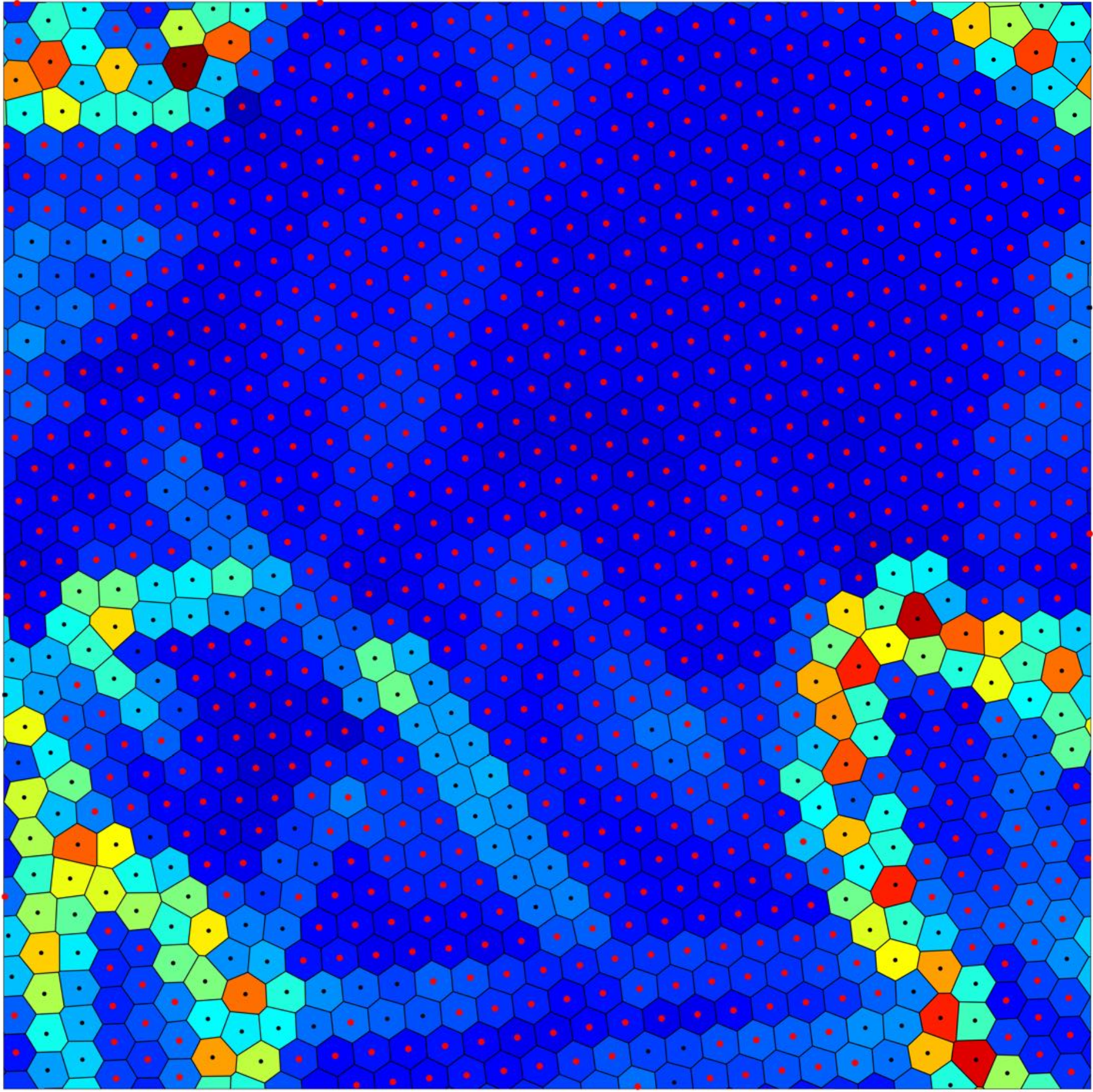} & \includegraphics[width=0.15\columnwidth]{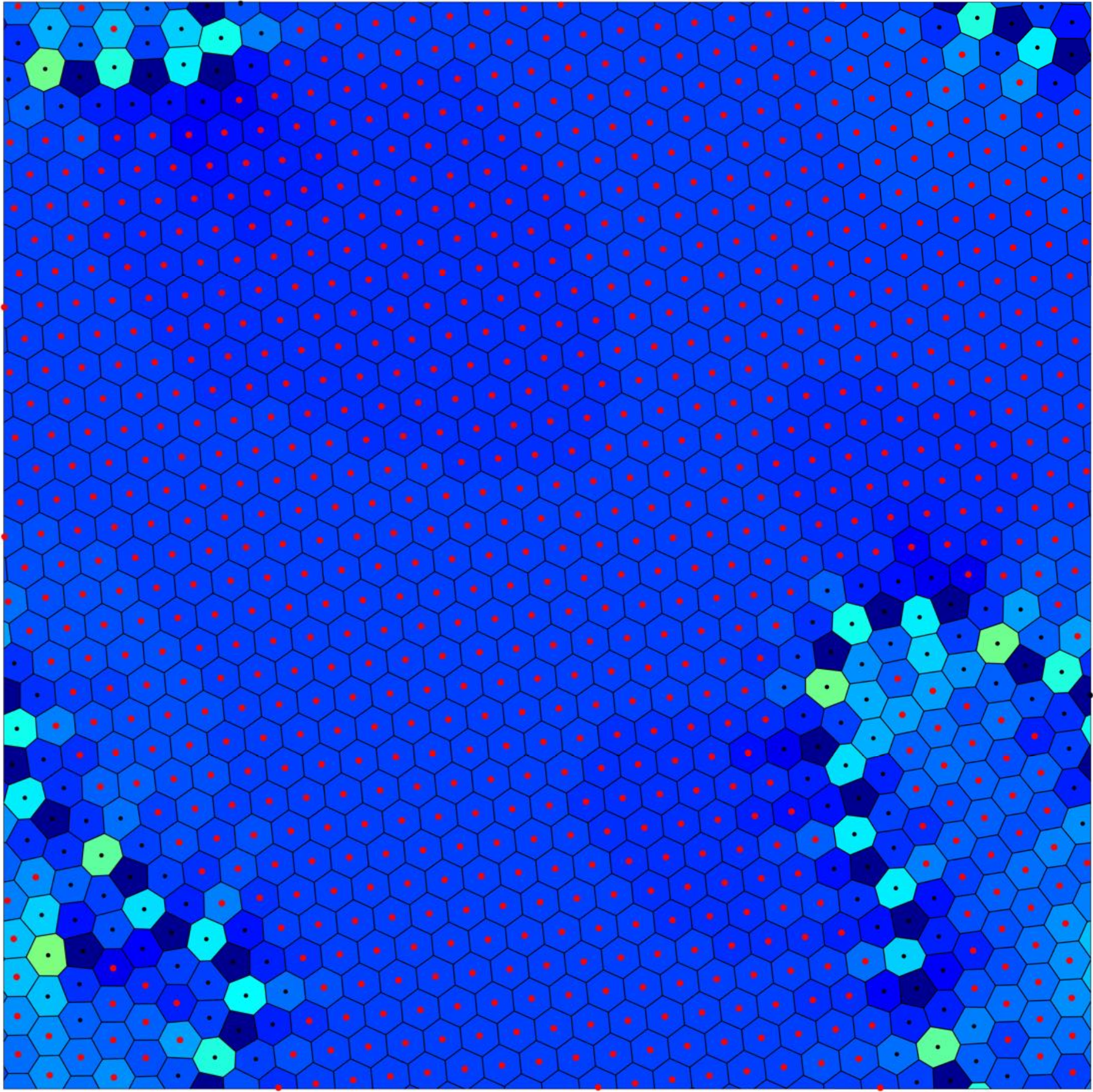} & \includegraphics[width=0.15\columnwidth]{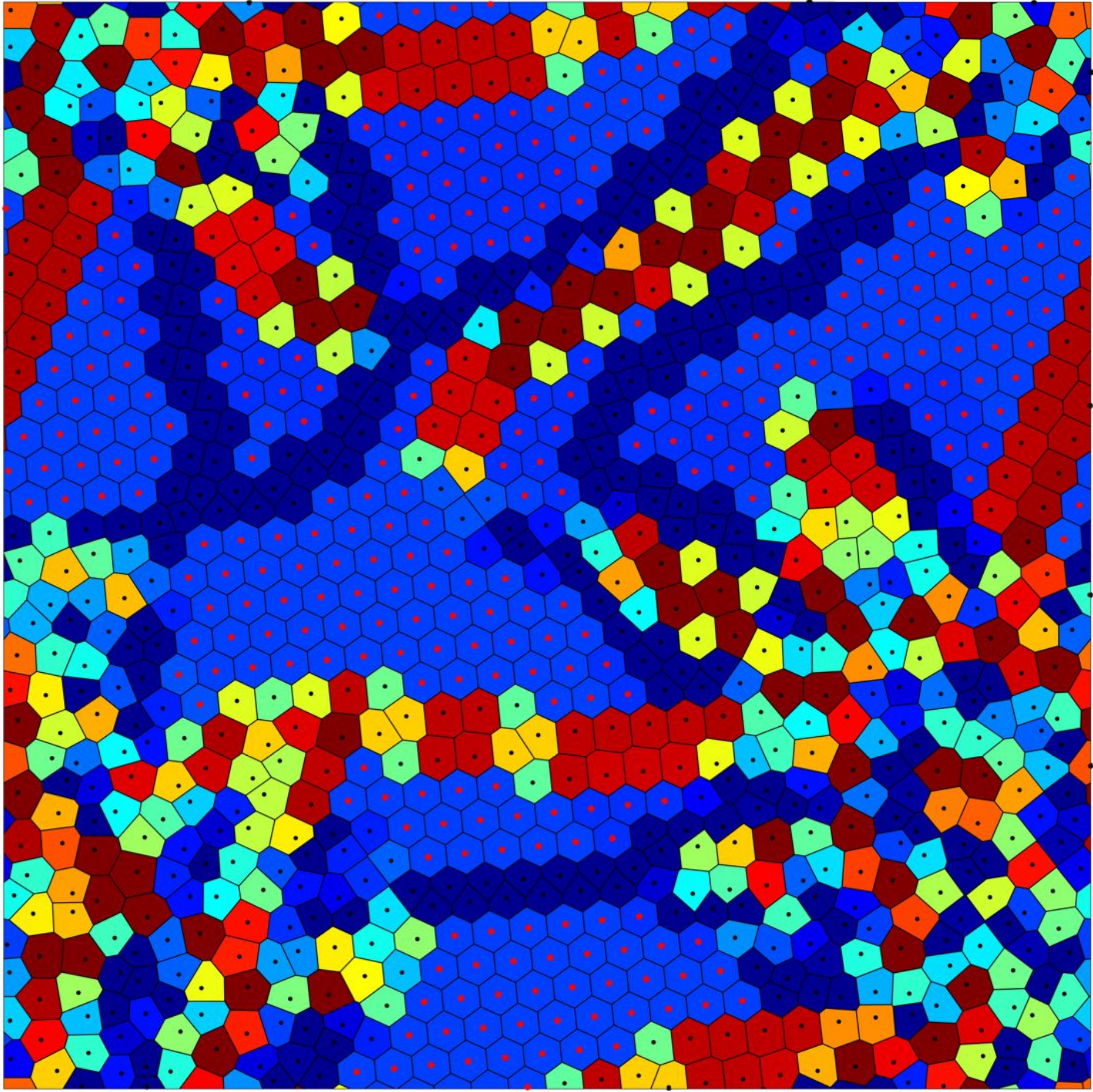} & 
 \includegraphics[width=0.15\columnwidth]{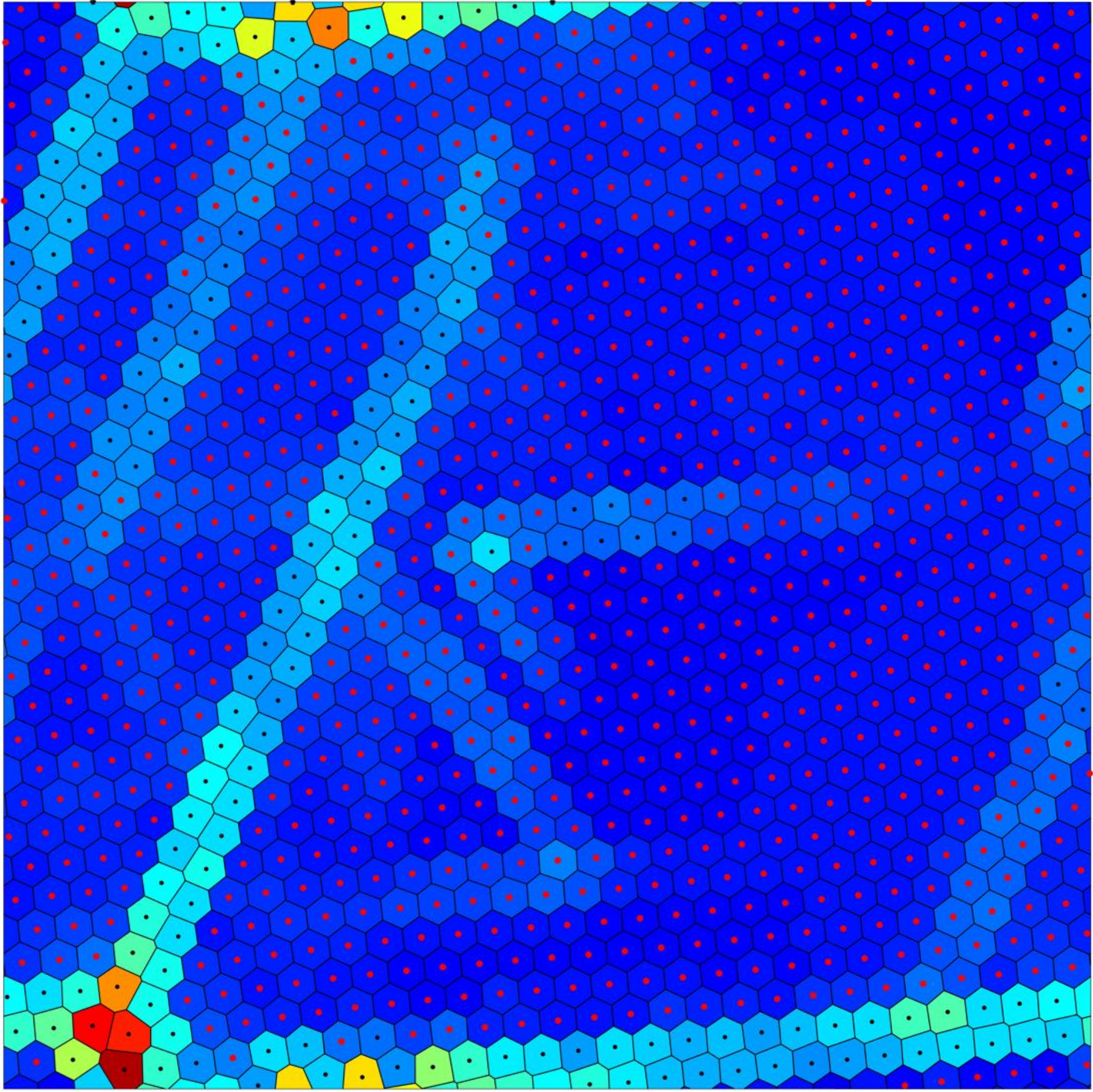}  & 
 \multicolumn{1} {c|}{\includegraphics[width=0.15\columnwidth]{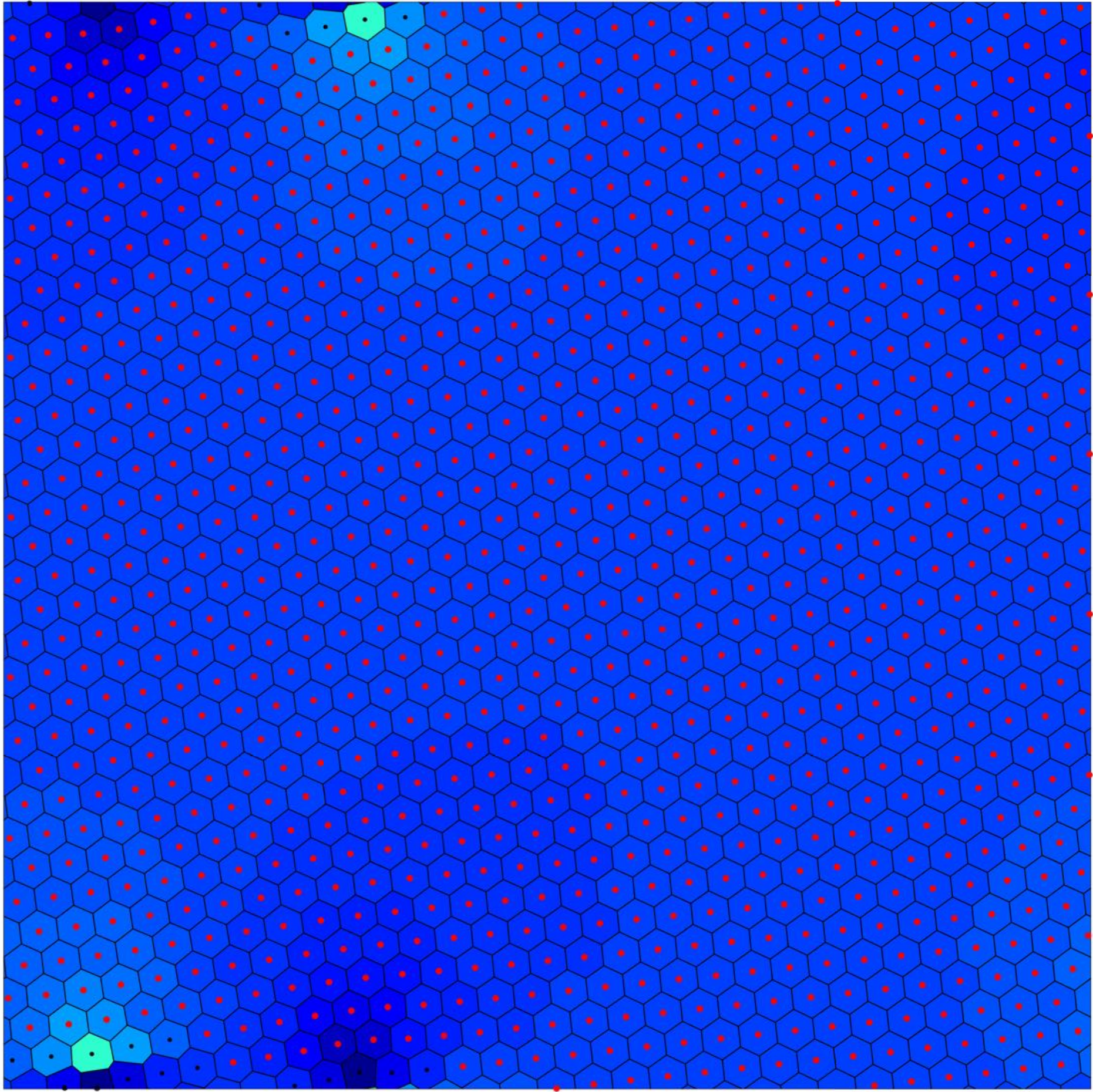}} &
 \multicolumn{1}{c}{}
\tabularnewline
\cline{1-5} 
\end{tabular}

\vspace{0.2cm}

\begin{tabular}{cc}
\includegraphics[width=0.45\columnwidth, height=0.25\columnwidth]{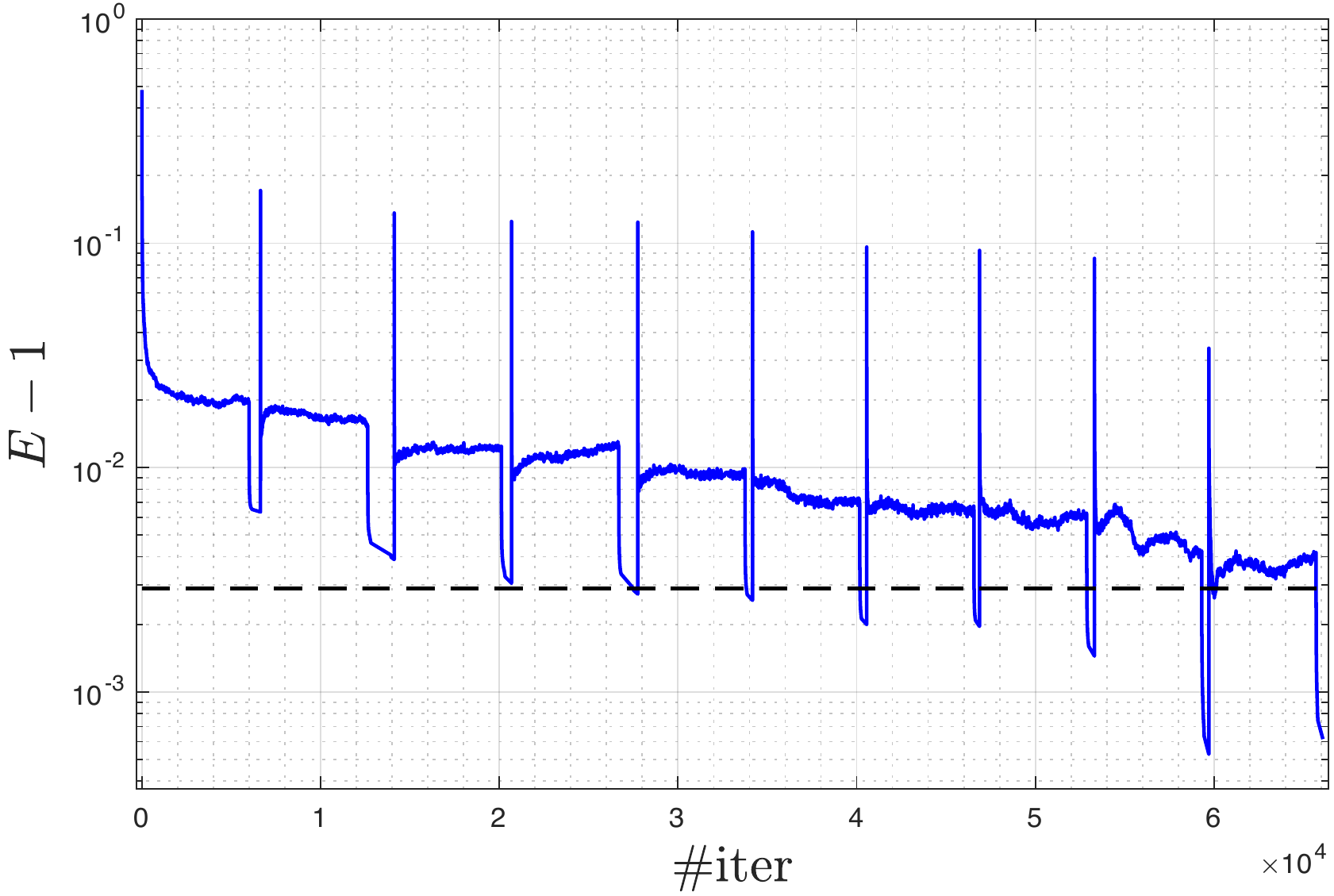} & \includegraphics[width=0.45\columnwidth, height=0.25\columnwidth]{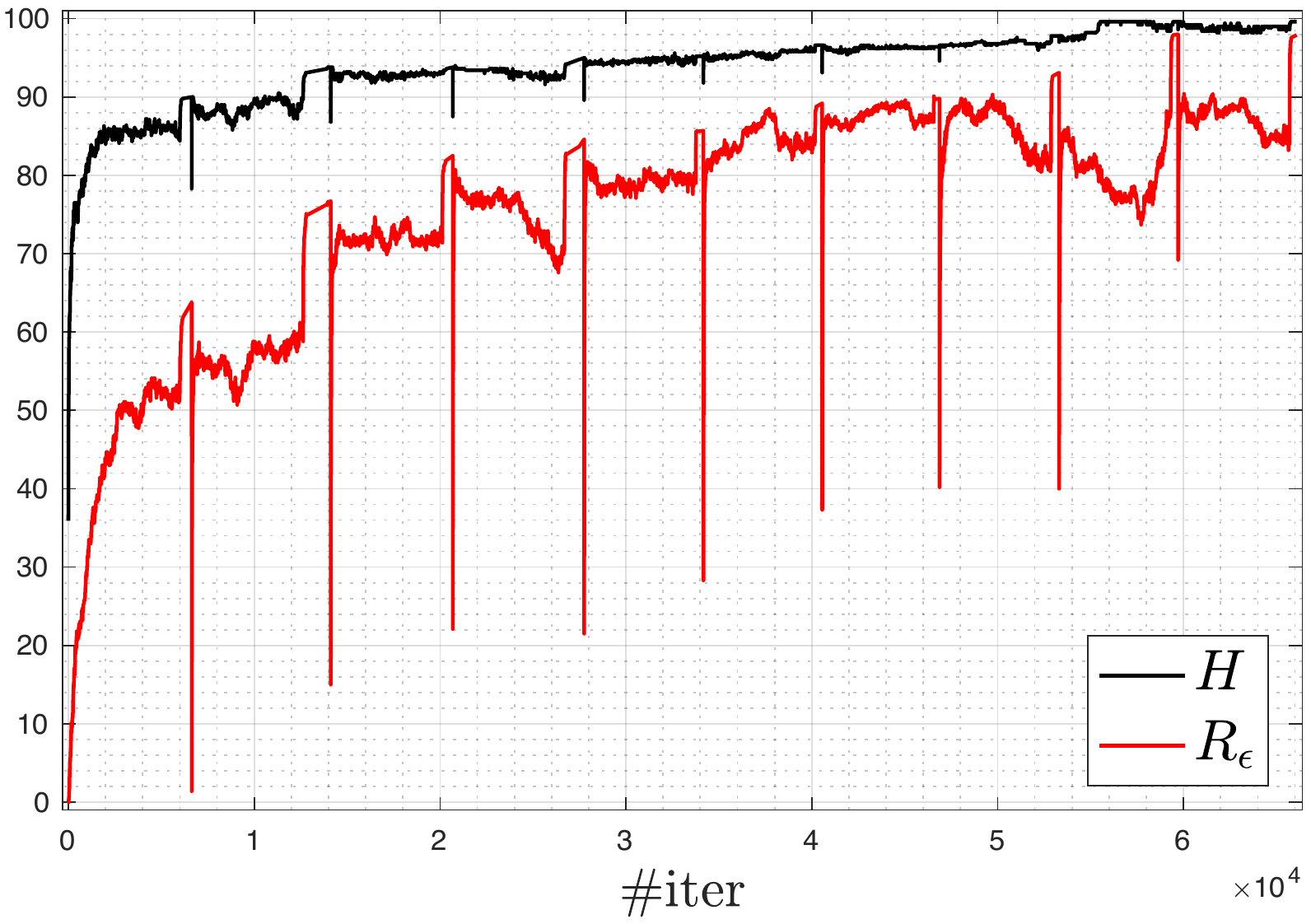}\tabularnewline
{\footnotesize{}(d) hybrid run achieving $E^{min}$ } & {\footnotesize{}(e) $R_{\epsilon},H$ profiles associated with (d)}\tabularnewline
\end{tabular}

\caption{Configurations with $N=1000$ from Example 3: (a) optimal PCVT obtained with Lloyd's
with $E-1=0.00466$, (b) optimal PCVT obtained with L-BFGS(7) having value $E-1=0.00349$, (c) minimal energy
configuration obtained amongst the three comparative methods --achieved by P-L-BFGS(20,20)--
with $E^{ref}-1=0.00289$.\protect \\
The mosaic shows the last iteration of each of the three blocks of
our hybrid method across the 10 stages of the run that achieved $E^{min}-1=5.27$e-$4$,
the PCVT carrying that minimal value is framed in a red box. These
images are to be read starting on the left side of the vertical double
black bar and then on the right, each row is for a stage $q=0,...,9$.
Both left columns are the last iteration of the \textit{MACN-c} blocks,
the middle columns contain the PCVTs from the Lloyd blocks and the
right columns are their respective \textit{MACN-}$\delta$ dislocations.
Finally (d) and (e) are the regularity measures' profiles of the run depicted in the mosaic, i.e. the one that achieved $E^{min}$ amongst the 1000 runs using $K=6000$.
\label{fig: mosaic N1000}}
\end{figure}

\begin{figure}[H]
\centering

\begin{tabular}{c}
\includegraphics[width=0.5\columnwidth,height=0.24\columnwidth]{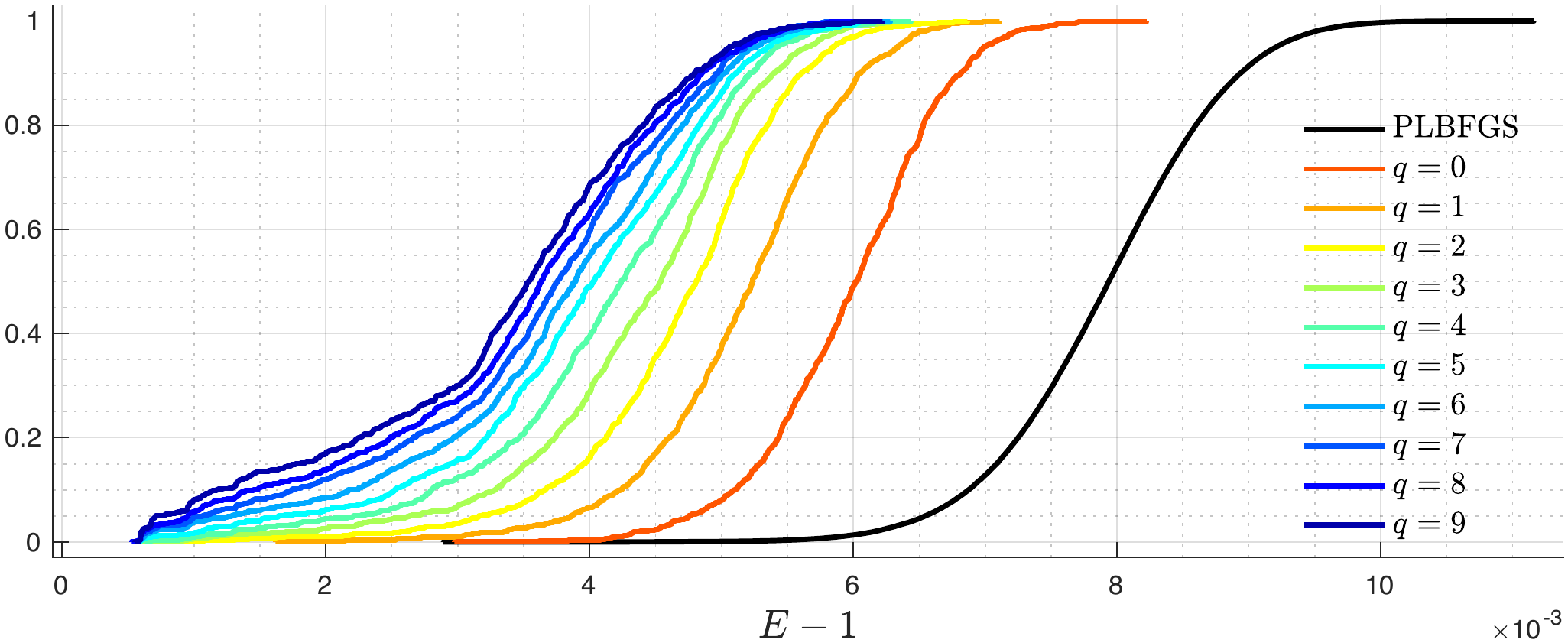}\tabularnewline
{\footnotesize{}(a) $f_{E-1}$}\tabularnewline
\includegraphics[width=0.53\columnwidth,height=0.24\columnwidth]{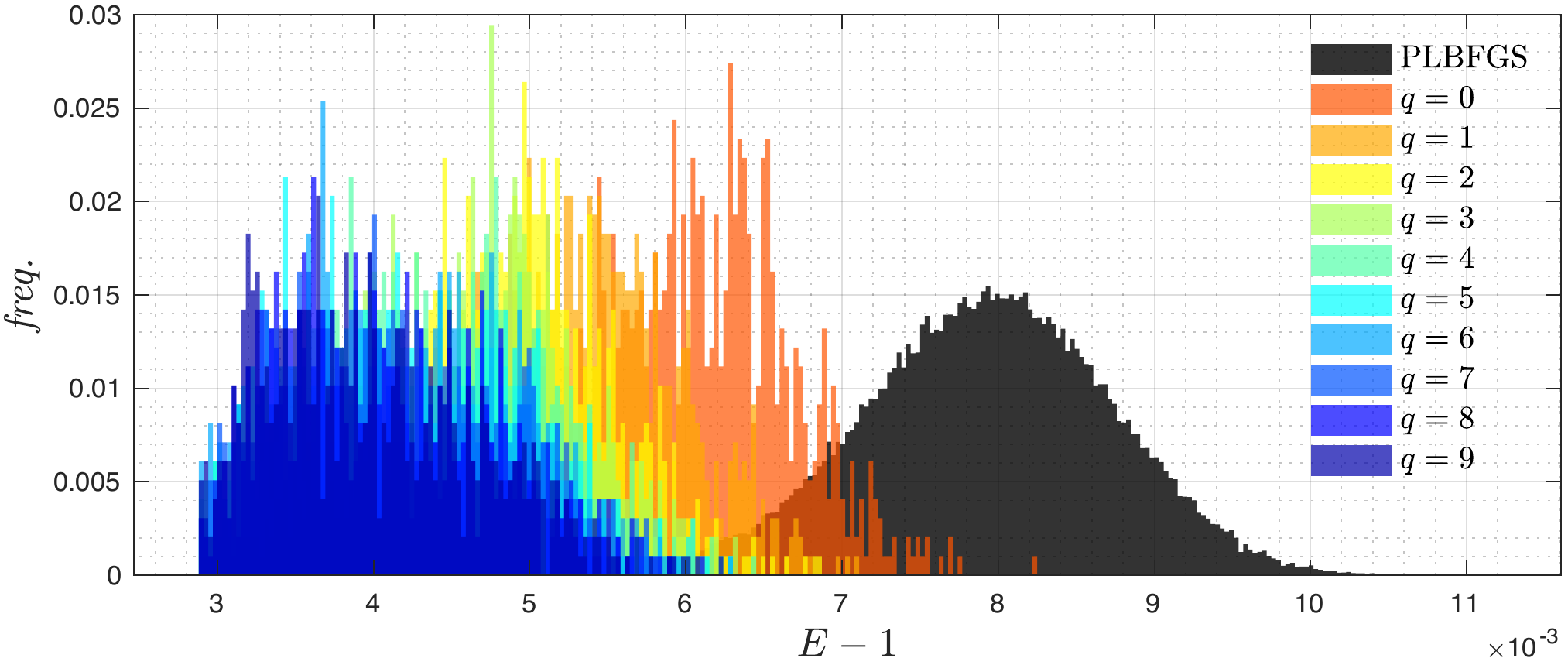}\tabularnewline
{\footnotesize{}(b) hist $E-1$}\tabularnewline

\end{tabular}

\begin{tabular}{cc}
\includegraphics[width=0.4\columnwidth, height=0.23\columnwidth]{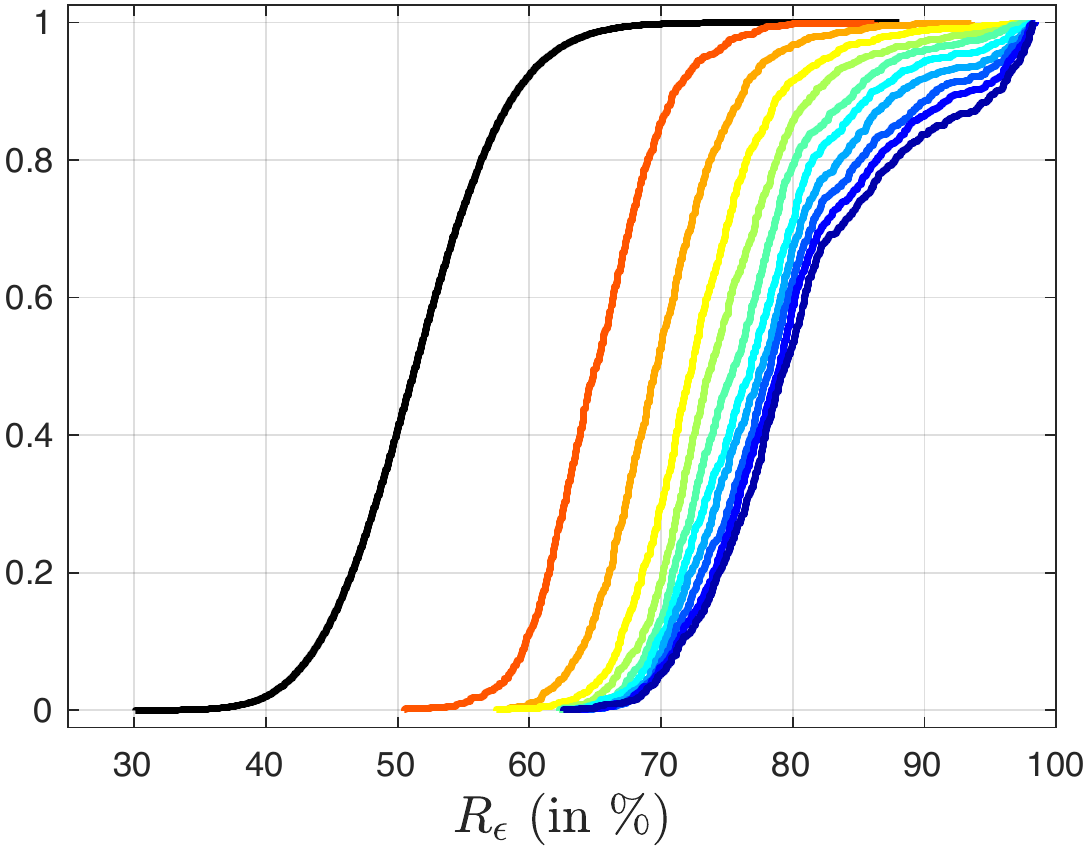} & \includegraphics[width=0.4\columnwidth, height=0.23\columnwidth]{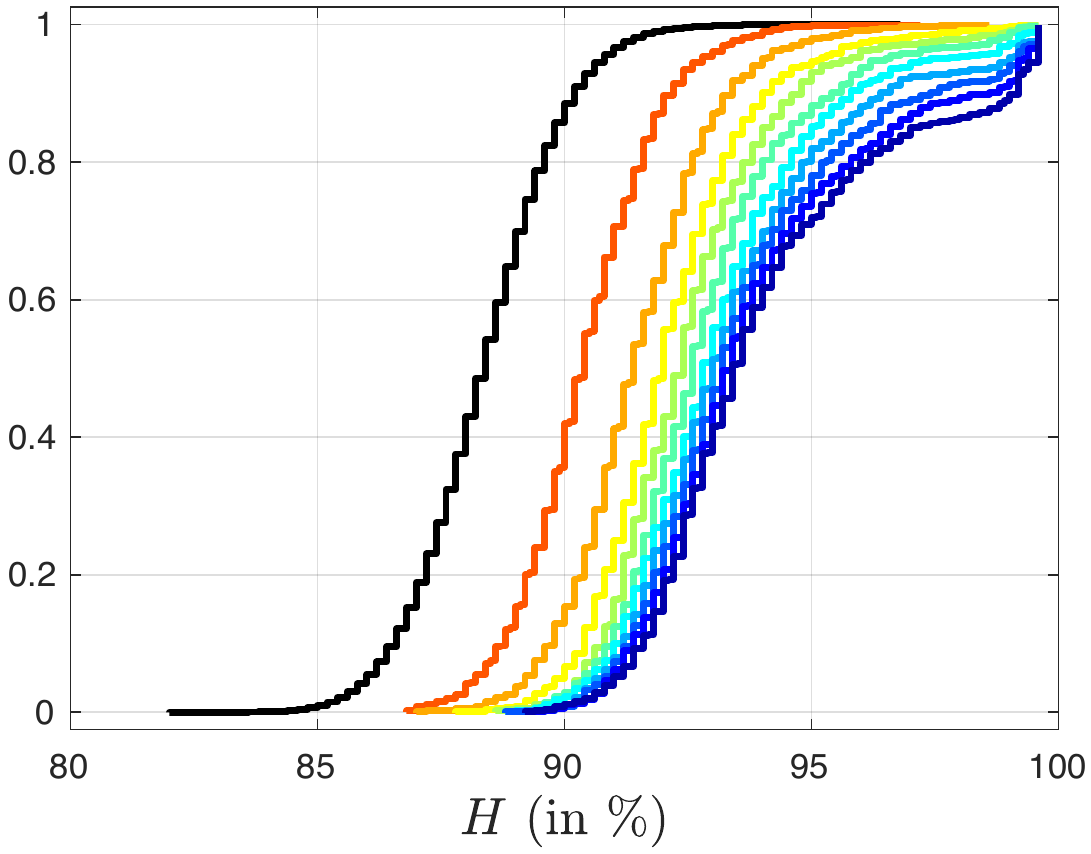}\tabularnewline
{\footnotesize{}(c) $f_{R_{\epsilon}}$} & {\footnotesize{}(d) $f_{H}$}\tabularnewline
\includegraphics[width=0.4\columnwidth, height=0.23\columnwidth]{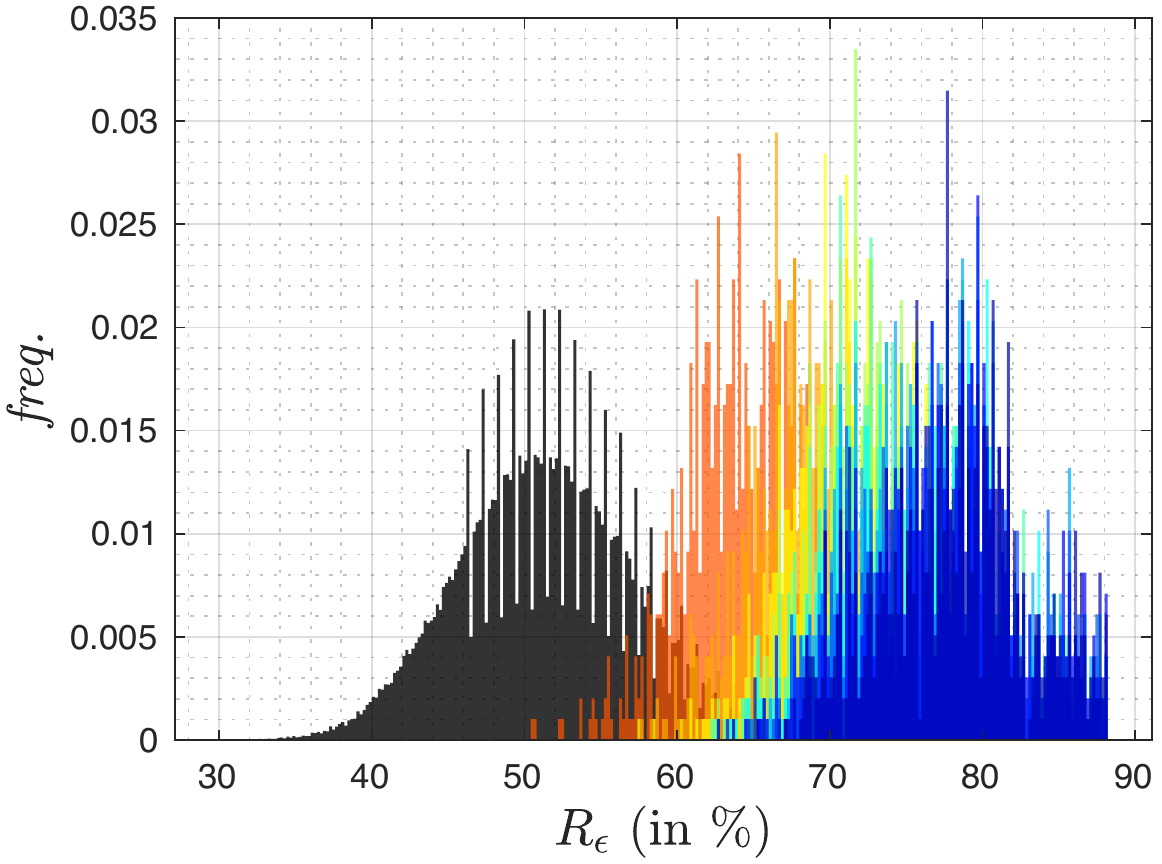} & \includegraphics[width=0.4\columnwidth, height=0.23\columnwidth]{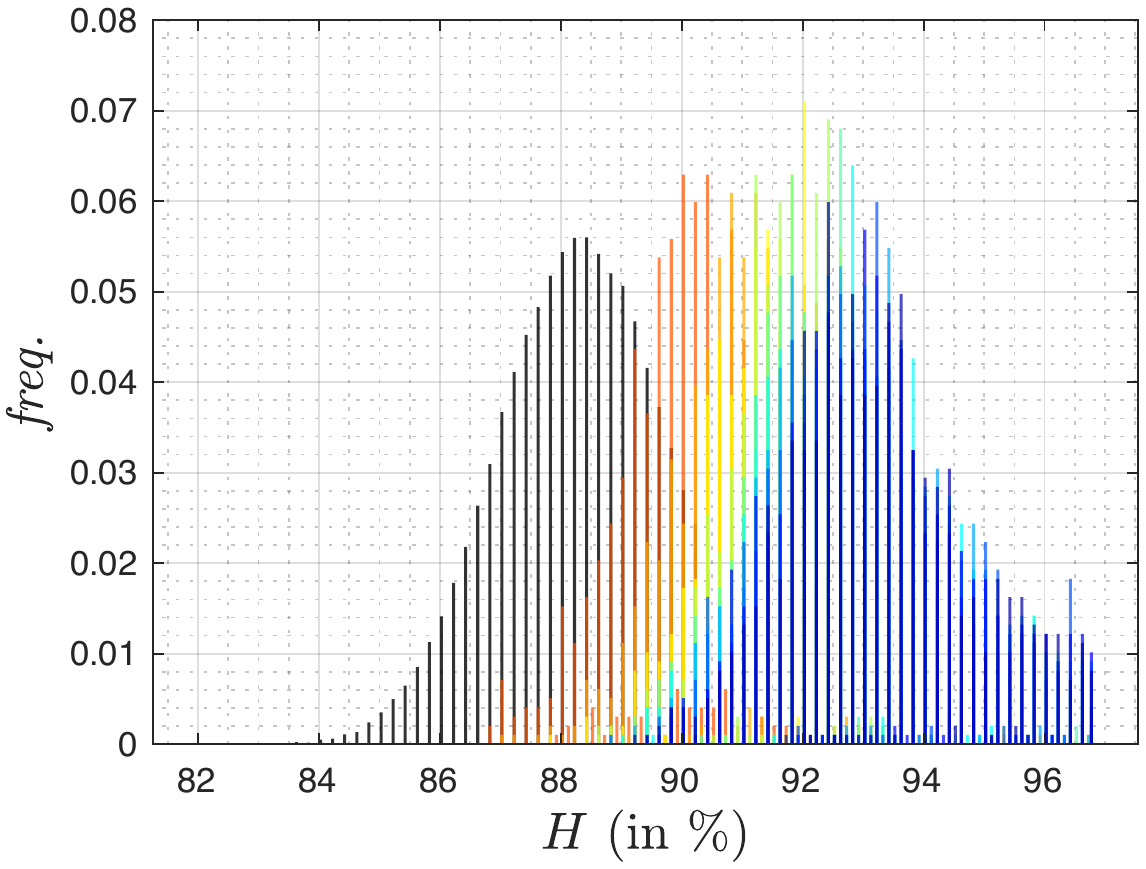}\tabularnewline
{\footnotesize{}(e) hist $R_{\epsilon}$} & {\footnotesize{}(f) hist $H$}\tabularnewline
\end{tabular}

\caption{Distributions from Example 3 for $N=1000$ with $K=6000$: ECDFs and histograms of regularity measures for the 1,000 hybrid runs as well as for the 100,000 runs of P-L-BFGS(20,20). Further detail on these measures are given in Table \ref{tab: N=1000}.
\label{fig: data N1000}}
\end{figure}

For $N=1000$ with $K=6000$ (as in Example 1) we implemented a large batch of 1,000 runs of our algorithm so as to have more robust statistics that are presented in Figure \ref{fig: data N1000} and in Table \ref{tab: N=1000}.
We further wish to emphasize, through the monotone skewness in $q$ of the ECDFs and histograms, the increasing performance tendency from one stage to another in this scenario where the contraction and relaxation phases are more constrained by the domain's shape.\\
The behavior is illustrated in the mosaic of PVTs from Figure \ref{fig: mosaic N1000} where snapshots of the run achieving $E^{min}-1=5.27$e-4 were taken. In particular one appreciates how stage after stage:

\begin{enumerate}
\item[i)]  the \textit{MACN-c} dynamics make defective (non-hexagonally
regular) regions of the PVTs ``communicate'' with each other
by creating a flow between regions with high average of individual
energy $E_i$ and ones with sensibly lower value, exhibiting then
a clear non-local behavior.
\item[ii)]  Lloyd's algorithm contracts the system while preserving the localization
of the defects.
\item[iii)]  \textit{MACN-$\delta$} preserves regions of hexagonal regularity
better and better as the energy of the perturbed PCVT diminishes whilst,
similarly to \textit{MACN-c,} creating a non-local ``communication''
between defects.
\end{enumerate}

\noindent At last, we point out that the progressive constriction of the defect
``interfaces'' we observe after the Lloyd block (the two middle columns
on the mosaic) is recurrent across all set ups that were tested, this is simply a consequence of the remarkable navigation of the energy that our hybrid method performs. The video animating
the iterations of this particular run can be found in the accompanying supplementary material, with it will appear in more
detail how the combination of these symbiotic blocks seem to have
a similar effect to a grain boundary evolution algorithm in polycrystals
when looking at the produced sequence of PCVTs.\\

\paragraph*{Example 4\label{par:Example-4} }

For $N=2000$ with $K=8000$ (as in Example 2) our results for 100 runs are summarized in Table \ref{tab: N=2000}.
Here P-L-BFGS(20,20) achieved $E^{ref}-1=0.00495$ which yields the ratio $\tau^{-1}\approx2.6$. We note that despite the non-negligible performance decrease of $\tau^{-1}$ when compared to Example
2; the values of $f_{E-1}^{*}$ above $90\%$ for $\mathbf{X}_{q\geq 4}^{*}$ indicate
that, statistically, the overall comparative efficiency is remarkably
similar to what we obtained in the hexagonal torus domain. This suggests
that the aforementioned change in $\tau$ could
be mainly attributed to the increased rigidity of the square torus.\\
On the other hand when restricting our attention to the distributions of $E$
gotten thus far from 100 runs or more, we make the point that there is no such
discrepancy since our energy scaling shows remarkable robustness for each method (e.g. $\langle E-1\rangle$ measurements seem to remain comparable regardless of $N$ and $\Omega$, this will be seen as well in the remaining cases).\\

\paragraph*{Example 5\label{par:Example-5:}}

When running our algorithm on a total of 100 initial configurations with $K=8000$ for $N=3000$
and comparing against P-L-BFGS(20,20), we get
$\tau^{-1}\approx3.1$ as well as values
of $f_{E-1}^{*}$ above $90\%$ for $\mathbf{X}_{q\geq 1}^{*}$,
see Table \ref{tab: N=3000}. The increasing relative performance with $N$ is once again apparent.\\

\paragraph*{Example 6\label{par:Example-6:} }

At last we consider $N=4000$ and 100 runs with $K=12000$, see Table \ref{tab: N=4000}. The reader
can appreciate how this set up comes to further corroborate the assertions
made thus far about the nature of our hybrid algorithm, namely:

\begin{enumerate}
\item[i)]  the monotone decrease of $\langle E-1\rangle$ in terms of $q$ for batches of 100 runs or more.
Additionally we've seen the robustness of $E$ in
terms of $N$ for the PCVTs produced by each method individually.
This suggests for example that, within
$Q=10$ stages, our hybrid method is able to get on average
twice as close to the non-achievable regular hexagonal configuration than the Quasi-Newton methods.
\item[ii)]  as the non-convexity of the problem increases with $N$ it is of
course harder to get closer to the ground state (the smallest energies
recorded increase regardless of the method), yet $f_{E-1}^{*}$ surpasses
90\% at earlier stages $q$ the larger the number of generators. This
shows that the manner our method probes the energy landscape manages
to overcome this stiffness remarkably better than gradient based methods
combined with random-like sampling.\\
\end{enumerate}

\noindent Graphics providing visual insight on the results of Examples 2,4,5 and 6 can be found in the accompanying supplementary material.

\begin{table}[H]

\resizebox{1\columnwidth}{0.13\columnwidth}{
\begin{centering}
\begin{tabular}{|>{\centering}p{13mm}||cccccccccc||ccc|}
\cline{2-14} 
\multicolumn{1}{>{\centering}p{9mm}|}{} &  &  &  &  & \textcolor{black}{Hybrid} &  &  &  &  &  & \multicolumn{1}{c|}{\multirow{2}{*}{\textcolor{black}{Lloyd}}} & \multicolumn{1}{c|}{\multirow{2}{*}{\textcolor{black}{LBFGS}}} & \multirow{2}{*}{\textcolor{black}{PLBFGS}}\tabularnewline
\cline{2-11} 
\multicolumn{1}{>{\centering}p{9mm}|}{} & \textcolor{black}{$\mathbf{X}_{0}^{*}$} & \textcolor{black}{$\mathbf{X}_{1}^{*}$} & \textcolor{black}{$\mathbf{X}_{2}^{*}$} & \textcolor{black}{$\mathbf{X}_{3}^{*}$} & \textcolor{black}{$\mathbf{X}_{4}^{*}$} & \textcolor{black}{$\mathbf{X}_{5}^{*}$} & \textcolor{black}{$\mathbf{X}_{6}^{*}$} & \textcolor{black}{$\mathbf{X}_{7}^{*}$} & \textcolor{black}{$\mathbf{X}_{8}^{*}$} & \textcolor{black}{$\mathbf{X}_{9}^{*}$} &  &  & \tabularnewline
\hline 
\textcolor{black}{$\langle E-1\rangle$} & \textcolor{black}{0.00597} & \textcolor{black}{0.00516} & \textcolor{black}{0.00469} & \textcolor{black}{0.00437} & \textcolor{black}{0.00412} & \textcolor{black}{0.00392} & \textcolor{black}{0.00374} & \textcolor{black}{0.00358} & \textcolor{black}{0.00346} & \textcolor{black}{0.00333} & \textcolor{black}{0.00848} & \textcolor{black}{0.00790} & \textcolor{black}{0.00791}\tabularnewline
\textcolor{black}{$\sigma_{E-1}$} & \textcolor{black}{0.00068} & \textcolor{black}{0.00075} & \textcolor{black}{0.00083} & \textcolor{black}{0.00092} & \textcolor{black}{0.00099} & \textcolor{black}{0.00106} & \textcolor{black}{0.00113} & \textcolor{black}{0.00117} & \textcolor{black}{0.00121} & \textcolor{black}{0.00126} & \textcolor{black}{0.00080} & \textcolor{black}{0.00081} & \textcolor{black}{0.00081}\tabularnewline
\textcolor{black}{max} & \textcolor{black}{0.00822} & \textcolor{black}{0.00710} & \textcolor{black}{0.00686} & \textcolor{black}{0.00631} & \textcolor{black}{0.00643} & \textcolor{black}{0.00628} & \textcolor{black}{0.00628} & \textcolor{black}{0.00592} & \textcolor{black}{0.00610} & \textcolor{black}{0.00622} & \textcolor{black}{0.01106} & \textcolor{black}{0.01142} & \textcolor{black}{0.01116}\tabularnewline
\textcolor{black}{min} & \textcolor{black}{0.00297} & \textcolor{black}{0.00161} & \textcolor{black}{7.40e-4} & \textcolor{black}{5.96e-4} & \textcolor{black}{6.26e-4} & \textcolor{black}{5.91e-4} & \textcolor{black}{5.82e-4} & \textcolor{black}{5.50e-4} & \textbf{\textcolor{black}{5.27e-4}} & \textcolor{black}{5.50e-4} & \textcolor{black}{0.00466} & \textcolor{black}{0.00349} & \textbf{\textcolor{black}{0.00289}}\tabularnewline
\textcolor{black}{$f_{E-1}^{*}$} & 0.000 & 0.009 & 0.030 & 0.058 & 0.114 & 0.143 & 0.188 & 0.227 & 0.262 & \textbf{0.283} & \textcolor{black}{- - -} & \textcolor{black}{- - -} & \textcolor{black}{- - -}\tabularnewline
\hline 
\textcolor{black}{$\langle R_{\epsilon}\rangle$ } & \textcolor{black}{65.32} & \textcolor{black}{70.05} & \textcolor{black}{72.93} & \textcolor{black}{74.76} & \textcolor{black}{76.26} & \textcolor{black}{77.43} & \textcolor{black}{78.47} & \textcolor{black}{78.48} & \textcolor{black}{80.24} & \textcolor{black}{81.03} & \textcolor{black}{47.77} & \textcolor{black}{51.57} & \textcolor{black}{51.51}\tabularnewline
\textcolor{black}{$\sigma_{R_{\epsilon}}$} & \textcolor{black}{4.67} & \textcolor{black}{4.96} & \textcolor{black}{5.42} & \textcolor{black}{5.93} & \textcolor{black}{6.33} & \textcolor{black}{6.78} & \textcolor{black}{7.19} & \textcolor{black}{7.50} & \textcolor{black}{7.75} & \textcolor{black}{8.06} & \textcolor{black}{5.79} & \textcolor{black}{5.90} & \textcolor{black}{5.87}\tabularnewline
\textcolor{black}{max} & \textcolor{black}{86.20} & \textcolor{black}{93.60} & \textcolor{black}{98.20} & \textcolor{black}{98.20} & \textcolor{black}{98.20} & \textcolor{black}{98.20} & \textcolor{black}{98.20} & \textcolor{black}{98.20} & \textbf{\textcolor{black}{98.40}} & \textcolor{black}{98.30} & \textcolor{black}{74.59} & \textcolor{black}{82.50} & \textbf{\textcolor{black}{88.09}}\tabularnewline
\textcolor{black}{min} & \textcolor{black}{50.50} & \textcolor{black}{57.60} & \textcolor{black}{57.50} & \textcolor{black}{62.30} & \textcolor{black}{62.30} & \textcolor{black}{63.50} & \textcolor{black}{62.70} & \textcolor{black}{64.40} & \textcolor{black}{63.30} & \textcolor{black}{62.60} & \textcolor{black}{27.89} & \textcolor{black}{28.30} & \textcolor{black}{30.10}\tabularnewline
\textcolor{black}{$1-f_{R\epsilon}^{*}$} & 0.000 & 0.004 & 0.017 & 0.034 & 0.048 & 0.074 & 0.103 & 0.145 & 0.166 & \textbf{0.192} & \textcolor{black}{- - -} & \textcolor{black}{- - -} & \textcolor{black}{- - }\tabularnewline
\hline 
\textcolor{black}{$\langle H\rangle$ } & \textcolor{black}{90.42} & \textcolor{black}{91.50} & \textcolor{black}{92.15} & \textcolor{black}{92.56} & \textcolor{black}{92.91} & \textcolor{black}{93.22} & \textcolor{black}{93.50} & \textcolor{black}{93.73} & \textcolor{black}{93.95} & \textcolor{black}{94.16} & \textcolor{black}{87.71} & \textcolor{black}{88.38} & \textcolor{black}{88.37}\tabularnewline
\textcolor{black}{$\sigma_{H}$} & \textcolor{black}{1.38} & \textcolor{black}{1.50} & \textcolor{black}{1.64} & \textcolor{black}{1.74} & \textcolor{black}{1.89} & \textcolor{black}{2.04} & \textcolor{black}{2.19} & \textcolor{black}{2.32} & \textcolor{black}{2.42} & \textcolor{black}{2.53} & \textcolor{black}{1.42} & \textcolor{black}{1.44} & \textcolor{black}{1.44}\tabularnewline
\textcolor{black}{max} & \textcolor{black}{97.20} & \textcolor{black}{98.60} & \textbf{\textcolor{black}{99.60}} & \textbf{\textcolor{black}{99.60}} & \textbf{\textcolor{black}{99.60}} & \textbf{\textcolor{black}{99.60}} & \textbf{\textcolor{black}{99.60}} & \textbf{\textcolor{black}{99.60}} & \textbf{\textcolor{black}{99.60}} & \textbf{\textcolor{black}{99.60}} & \textcolor{black}{94.00} & \textcolor{black}{95.79} & \textbf{\textcolor{black}{96.79}}\tabularnewline
\textcolor{black}{min} & \textcolor{black}{86.80} & \textcolor{black}{87.00} & \textcolor{black}{87.80} & \textcolor{black}{88.60} & \textcolor{black}{89.20} & \textcolor{black}{88.80} & \textcolor{black}{89.20} & \textcolor{black}{88.80} & \textcolor{black}{89.40} & \textcolor{black}{89.20} & \textcolor{black}{82.59} & \textcolor{black}{82.00} & \textcolor{black}{82.00}\tabularnewline
\textcolor{black}{$1-f_{H}^{*}$} & 0.001 & 0.003 & 0.019 & 0.029 & 0.042 & 0.061 & 0.079 & 0.112 & 0.138 & \textbf{0.161} & \textcolor{black}{- - -} & \textcolor{black}{- - -} & \textcolor{black}{- - }\tabularnewline
\hline 
\end{tabular}
\par\end{centering}
}
\caption{Statistics from Example 3 with $N=1000$ and $K=6000$: a 1000 runs of our hybrid method versus the comparative algorithms, the values $E^{min}-1,\,R_{\epsilon}^{min} ,\,H^{min}$ as well as $E^{ref}-1,\,R_{\epsilon}^{ref} ,\,H^{ref}$ are bold faced.\label{tab: N=1000}}

\resizebox{1\columnwidth}{0.13\columnwidth}{
\begin{centering}
\begin{tabular}{|>{\centering}p{13mm}||cccccccccc||ccc|}
\cline{2-14} 
\multicolumn{1}{>{\centering}p{9mm}|}{} &  &  &  &  & Hybrid &  &  &  &  &  & \multicolumn{1}{c|}{\multirow{2}{*}{Lloyd}} & \multicolumn{1}{c|}{\multirow{2}{*}{LBFGS}} & \multirow{2}{*}{PLBFGS}\tabularnewline
\cline{2-11} 
\multicolumn{1}{>{\centering}p{9mm}|}{} & $\mathbf{X}_{0}^{*}$ & $\mathbf{X}_{1}^{*}$ & $\mathbf{X}_{2}^{*}$ & $\mathbf{X}_{3}^{*}$ & $\mathbf{X}_{4}^{*}$ & $\mathbf{X}_{5}^{*}$ & $\mathbf{X}_{6}^{*}$ & $\mathbf{X}_{7}^{*}$ & $\mathbf{X}_{8}^{*}$ & $\mathbf{X}_{9}^{*}$ &  &  & \tabularnewline
\hline 
{$\langle E-1\rangle$} & 0.00572 & 0.00500 & 0.00462 & 0.00433 & 0.00418 & 0.00402 & 0.00395 & 0.00377 & 0.00368 & 0.00355 & 0.00842 & 0.00784 & 0.00786\tabularnewline
{$\sigma_{E-1}$} & 0.00043 & 0.00044 & 0.00053 & 0.00057 & 0.00056 & 0.00060 & 0.00059 & 0.00061 & 0.00060 & 0.00060 & 0.00057 & 0.00057 & 0.00057\tabularnewline
{max} & 0.00671 & 0.00582 & 0.00584 & 0.00586 & 0.00552 & 0.00581 & 0.00523 & 0.00533 & 0.00520 & 0.00551 & 0.01051 & 0.01013 & 0.01034\tabularnewline
{min} & 0.00471 & 0.00345 & 0.00256 & 0.00196 & 0.00226 & 0.00195 & 0.00197 & 0.00193 & \textbf{0.00191} & 0.00204 & 0.00585 & 0.00511 & \textbf{0.00495}\tabularnewline
{$f_{E-1}^{*}$} & 0.04 & 0.42 & 0.75 & 0.88 & 0.93 & 0.97 & \textbf{0.99} & \textbf{0.99} & \textbf{0.99} & \textbf{0.99} & - - - & - - - & - - -\tabularnewline
\hline 
{$\langle R_{\epsilon}\rangle$ } & 66.70 & 71.17 & 73.43 & 75.27 & 76.12 & 76.95 & 77.46 & 78.60 & 79.03 & 79.85 & 48.21 & 51.95 & 51.82\tabularnewline
{$\sigma_{R_{\epsilon}}$} & 2.97 & 2.80 & 3.47 & 3.49 & 3.51 & 3.67 & 3.65 & 3.73 & 3.72 & 3.69 & 4.12 & 4.15 & 4.18\tabularnewline
{max} & 73.70 & 81.65 & 86.65 & 89.85 & 87.90 & 89.90 & 90.00 & \textbf{90.10} & 89.85 & 89.45 & 65.80 & \textbf{72.54} & 71.95\tabularnewline
{min} & 60.00 & 66.05 & 64.70 & 65.85 & 66.80 & 65.90 & 69.25 & 69.65 & 70.45 & 67.75 & 33.55 & 35.39 & 35.75\tabularnewline
{\scriptsize{}$1-f_{R\epsilon}^{*}$} & 0.02 & 0.25 & 0.63 & 0.83 & 0.91 & 0.90 & 0.94 & 0.96 & 0.97 & \textbf{0.99} & - - - & - - - & - - -\tabularnewline
\hline 
{$\langle H\rangle$ } & 90.75 & 91.77 & 92.40 & 92.85 & 93.15 & 93.37 & 93.46 & 93.74 & 93.88 & 94.07 & 87.84 & 88.48 & 88.45\tabularnewline
{$\sigma_{H}$} & 0.84 & 0.91 & 1.01 & 1.03 & 1.03 & 1.09 & 1.05 & 1.10 & 1.08 & 1.06 & 1.01 & 1.01 & 1.02\tabularnewline
{max} & 92.40 & 94.70 & 95.95 & 96.40 & 96.50 & 96.40 & 96.40 & 96.50 & \textbf{96.65} & 97.00 & 91.70 & \textbf{93.50} & 93.30\tabularnewline
{min} & 88.90 & 89.80 & 90.00 & 90.10 & 90.00 & 90.40 & 90.50 & 91.10 & 91.90 & 91.10 & 83.90 & 84.00 & 83.90\tabularnewline
{$1-f_{H}^{*}$} & 0.00 & 0.04 & 0.10 & 0.23 & 0.33 & 0.38 & 0.44 & 0.58 & 0.56 & \textbf{0.72} & - - - & - - - & - - -\tabularnewline
\hline 
\end{tabular}
\par\end{centering}
}
\caption{Statistics out of 100 hybrid method runs from Example 4 with $N=2000$ and $K=8000$.\label{tab: N=2000}}

\resizebox{1\columnwidth}{0.13\columnwidth}{
\begin{centering}
\begin{tabular}{|>{\centering}p{13mm}||cccccccccc||ccc|}
\cline{2-14} 
\multicolumn{1}{>{\centering}p{9mm}|}{} &  &  &  &  & Hybrid &  &  &  &  &  & \multicolumn{1}{c|}{\multirow{2}{*}{Lloyd}} & \multicolumn{1}{c|}{\multirow{2}{*}{LBFGS}} & \multirow{2}{*}{PLBFGS}\tabularnewline
\cline{2-11} 
\multicolumn{1}{>{\centering}p{9mm}|}{} & $\mathbf{X}_{0}^{*}$ & $\mathbf{X}_{1}^{*}$ & $\mathbf{X}_{2}^{*}$ & $\mathbf{X}_{3}^{*}$ & $\mathbf{X}_{4}^{*}$ & $\mathbf{X}_{5}^{*}$ & $\mathbf{X}_{6}^{*}$ & $\mathbf{X}_{7}^{*}$ & $\mathbf{X}_{8}^{*}$ & $\mathbf{X}_{9}^{*}$ &  &  & \tabularnewline
\hline 
{$\langle E-1\rangle$} & 0.00577 & 0.00499 & 0.00462 & 0.00443 & 0.00427 & 0.00412 & 0.00396 & 0.00380 & 0.00368 & 0.00361 & 0.00842 & 0.00782 & 0.00786\tabularnewline
{$\sigma_{E-1}$} & 0.00044 & 0.00044 & 0.00047 & 0.00049 & 0.00056 & 0.00057 & 0.00057 & 0.00380 & 0.00053 & 0.00060 & 0.00057 & 0.00047 & 0.00047\tabularnewline
{max} & 0.00678 & 0.00587 & 0.00566 & 0.00547 & 0.00528 & 0.00530 & 0.00529 & 0.00519 & 0.00533 & 0.00540 & 0.01051 & 0.00985 & 0.00986\tabularnewline
{min} & 0.00425 & 0.00318 & 0.00314 & 0.00245 & 0.00199 & \textbf{0.00179} & 0.00208 & 0.00231 & 0.00218 & 0.00219 & 0.00585 & 0.00583 & \textbf{0.00568}\tabularnewline
{$f_{E-1}^{*}$} & 0.44 & 0.98 & \textbf{1.00} & \textbf{1.00} & \textbf{1.00} & \textbf{1.00} & \textbf{1.00} & \textbf{1.00} & \textbf{1.00} & \textbf{1.00} & - - - & - - - & - - -\tabularnewline
\hline 
{$\langle R_{\epsilon}\rangle$ } & 66.41 & 71.15 & 73.42 & 74.52 & 75.518 & 76.44 & 77.38 & 78.32 & 79.03 & 79.53 & 48.21 & 52.12 & 51.89\tabularnewline
{$\sigma_{R_{\epsilon}}$} & 2.97 & 2.83 & 2.97 & 3.20 & 3.60 & 3.57 & 3.54 & 3.31 & 3.34 & 3.72 & 4.12 & 3.42 & 3.42\tabularnewline
{max} & 76.93 & 82.63 & 82.06 & 85.96 & 89.30 & \textbf{90.36} & 89.23 & 87.66 & 88.53 & 87.53 & 65.80 & 66.13 & \textbf{67.20}\tabularnewline
{min} & 60.40 & 65.10 & 66.40 & 67.83 & 69.00 & 69.40 & 69.53 & 70.33 & 69.23 & 67.80 & 33.55 & 38.36 & 38.46\tabularnewline
{$1-f_{R\epsilon}^{*}$} & 0.40 & 0.95 & 0.99 & \textbf{1.00} & \textbf{1.00} & \textbf{1.00} & \textbf{1.00} & \textbf{1.00} & \textbf{1.00} & \textbf{1.00} & - - - & - - - & - - -\tabularnewline
\hline 
{$\langle H\rangle$ } & 90.67 & 91.79 & 92.35 & 92.70 & 92.95 & 93.18 & 93.44 & 93.72 & 93.97 & 94.11 & 87.84 & 88.53 & 88.46\tabularnewline
{$\sigma_{H}$} & 0.94 & 0.86 & 0.89 & 0.98 & 1.09 & 1.12 & 1.10 & 1.03 & 1.05 & 1.14 & 1.01 & 0.83 & 0.83\tabularnewline
{max} & 93.76 & 95.40 & 95.00 & 96.53 & 97.20 & \textbf{97.53} & 97.26 & 96.73 & 96.93 & 96.73 & 91.70 & 92.06 & \textbf{92.73}\tabularnewline
{min} & 88.66 & 90.00 & 90.66 & 90.66 & 91.00 & 91.00 & 90.86 & 91.46 & 91.00 & 90.86 & 83.90 & 85.20 & 85.13\tabularnewline
{$1-f_{H}^{*}$} & 0.03 & 0.17 & 0.34 & 0.45 & 0.60 & 0.64 & 0.76 & 0.81 & \textbf{0.89} & \textbf{0.89} & - - - & - - - & - - -\tabularnewline
\hline 
\end{tabular}
\par\end{centering}
}
\caption{Statistics out of 100 hybrid method runs from Example 5 with $N=3000$ and $K=8000$.\label{tab: N=3000}}

\resizebox{1\columnwidth}{0.13\columnwidth}{
\begin{centering}
\begin{tabular}{|>{\centering}p{13mm}||cccccccccc||ccc|}
\cline{2-14} 
\multicolumn{1}{>{\centering}p{9mm}|}{} &  &  &  &  & Hybrid &  &  &  &  &  & \multicolumn{1}{c|}{\multirow{2}{*}{Lloyd}} & \multicolumn{1}{c|}{\multirow{2}{*}{LBFGS}} & \multirow{2}{*}{PLBFGS}\tabularnewline
\cline{2-11} 
\multicolumn{1}{>{\centering}p{9mm}|}{} & $\mathbf{X}_{0}^{*}$ & $\mathbf{X}_{1}^{*}$ & $\mathbf{X}_{2}^{*}$ & $\mathbf{X}_{3}^{*}$ & $\mathbf{X}_{4}^{*}$ & $\mathbf{X}_{5}^{*}$ & $\mathbf{X}_{6}^{*}$ & $\mathbf{X}_{7}^{*}$ & $\mathbf{X}_{8}^{*}$ & $\mathbf{X}_{9}^{*}$ &  &  & \tabularnewline
\hline 
{$\langle E-1\rangle$} & 0.00535 & 0.00474 & 0.00441 & 0.00424 & 0.00404 & 0.00396 & 0.00380 & 0.00376 & 0.00375 & 0.00368 & 0.00833 & 0.00781 & 0.00786\tabularnewline
{$\sigma_{E-1}$} & 0.00033 & 0.00036 & 0.00039 & 0.00038 & 0.00040 & 0.00042 & 0.00042 & 0.00046 & 0.00047 & 0.00048 & 0.00041 & 0.00040 & 0.00040\tabularnewline
{max} & 0.00626 & 0.00545 & 0.00517 & 0.00507 & 0.00502 & 0.00486 & 0.00486 & 0.00472 & 0.00479 & 0.00499 & 0.00989 & 0.00956 & 0.00946\tabularnewline
{min} & 0.00448 & 0.00389 & 0.00338 & 0.00321 & 0.00314 & 0.00300 & 0.00268 & 0.00255 & 0.00240 & \textbf{0.00228} & 0.00645 & \textbf{0.00584} & 0.00601\tabularnewline
{$f_{E-1}^{*}$} & 0.96 & \textbf{1.00} & \textbf{1.00} & \textbf{1.00} & \textbf{1.00} & \textbf{1.00} & \textbf{1.00} & \textbf{1.00} & \textbf{1.00} & \textbf{1.00} & - - - & - - - & - - -\tabularnewline
\hline 
{$\langle R_{\epsilon}\rangle$ } & 69.11 & 72.61 & 74.64 & 75.75 & 76.93 & 77.44 & 78.43 & 78.72 & 78.66 & 79.15 & 48.84 & 52.21 & 51.86\tabularnewline
{$\sigma_{R_{\epsilon}}$} & 2.24 & 2.36 & 2.51 & 2.49 & 2.50 & 2.65 & 2.66 & 2.93 & 2.83 & 2.94 & 2.96 & 2.96 & 2.95\tabularnewline
{max} & 75.25 & 78.15 & 81.15 & 82.50 & 82.55 & 84.15 & 84.60 & 86.87 & 87.27 & \textbf{88.22} & 61.65 & \textbf{66.02} & 65.40\tabularnewline
{min} & 61.75 & 67.85 & 69.75 & 70.12 & 71.17 & 71.80 & 71.52 & 73.17 & 72.40 & 71.57 & 38.25 & 40.12 & 40.25\tabularnewline
{$1-f_{R\epsilon}^{*}$} & 0.94 & \textbf{1.00} & \textbf{1.00} & \textbf{1.00} & \textbf{1.00} & \textbf{1.00} & \textbf{1.00} & \textbf{1.00} & \textbf{1.00} & \textbf{1.00} & - - - & - - - & - - -\tabularnewline
\hline 
{$\langle H\rangle$ } & 91.28 & 92.19 & 92.74 & 92.98 & 93.34 & 93.49 & 93.79 & 93.85 & 93.87 & 93.96 & 88.00 & 88.57 & 88.45\tabularnewline
{$\sigma_{H}$} & 0.64 & 0.73 & 0.77 & 0.78 & 0.75 & 0.77 & 0.76 & 0.79 & 0.77 & 0.83 & 0.72 & 0.72 & 0.72\tabularnewline
{max} & 92.90 & 93.90 & 94.60 & 95.30 & 95.10 & 95.52 & 95.40 & 95.95 & 96.30 & \textbf{96.60} & 91.30 & \textbf{92.00} & 91.60\tabularnewline
{min} & 89.35 & 90.65 & 91.35 & 91.15 & 91.45 & 91.75 & 91.80 & 92.20 & 91.90 & 91.65 & 85.50 & 85.59 & 85.60\tabularnewline
{$1-f_{H}^{*}$} & 0.16 & 0.53 & 0.80 & 0.91 & 0.98 & 0.99 & 0.99 & \textbf{1.00} & 0.99 & 0.98 & - - - & - - - & - - -\tabularnewline
\hline 
\end{tabular}
\par\end{centering}
}
\caption{Statistics out of 100 hybrid method runs from Example 6 with $N=4000$ and $K=12000$.\label{tab: N=4000}}
\end{table}

\subsection{Scope on $R_{\epsilon}$\label{subsec:5.4}}

\begin{figure}[H]
\begin{center}
\begin{tabular}{cc}
\includegraphics[width=0.4\columnwidth,height=0.35\columnwidth]{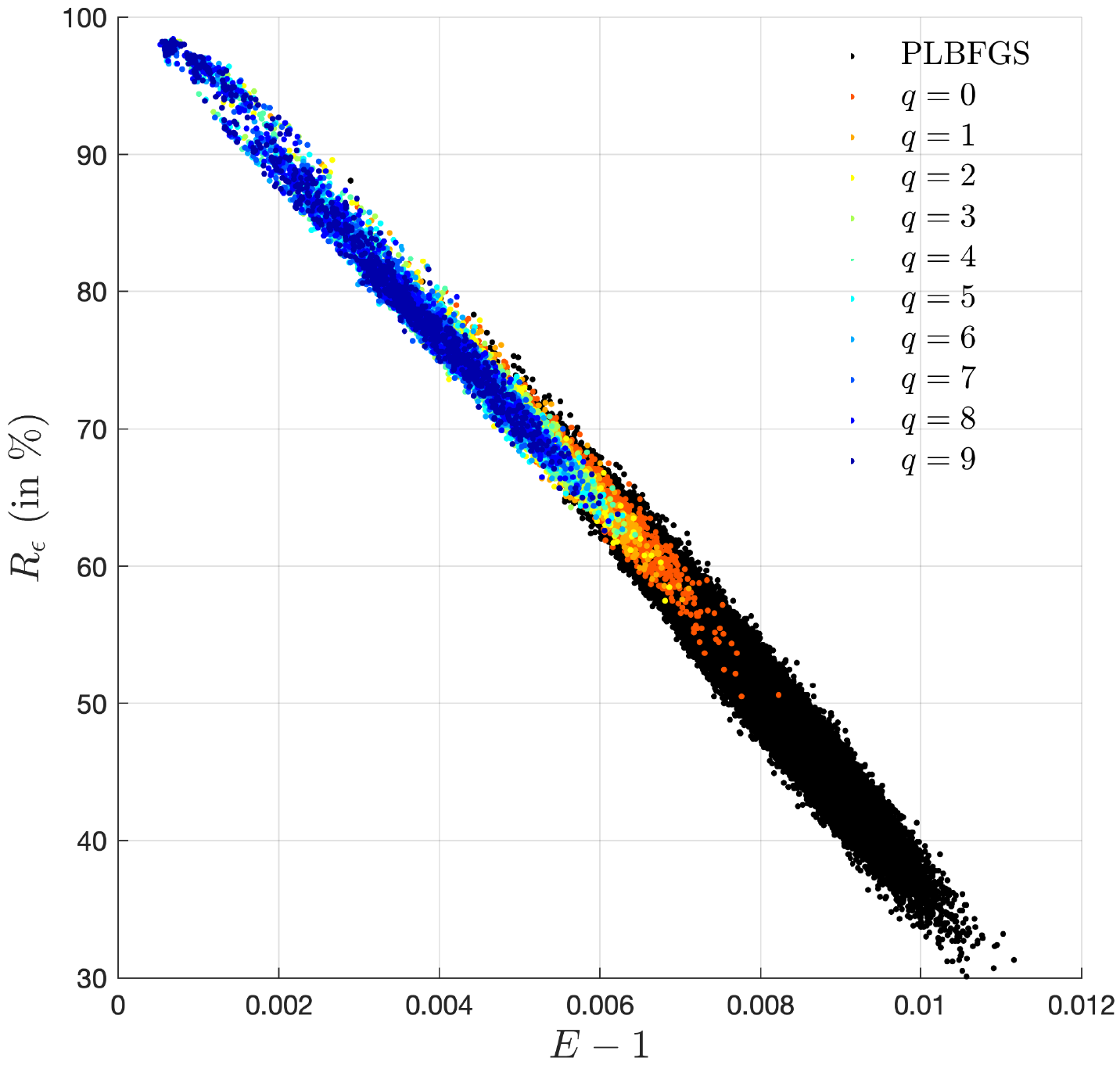} & \includegraphics[width=0.4\columnwidth,height=0.35\columnwidth]{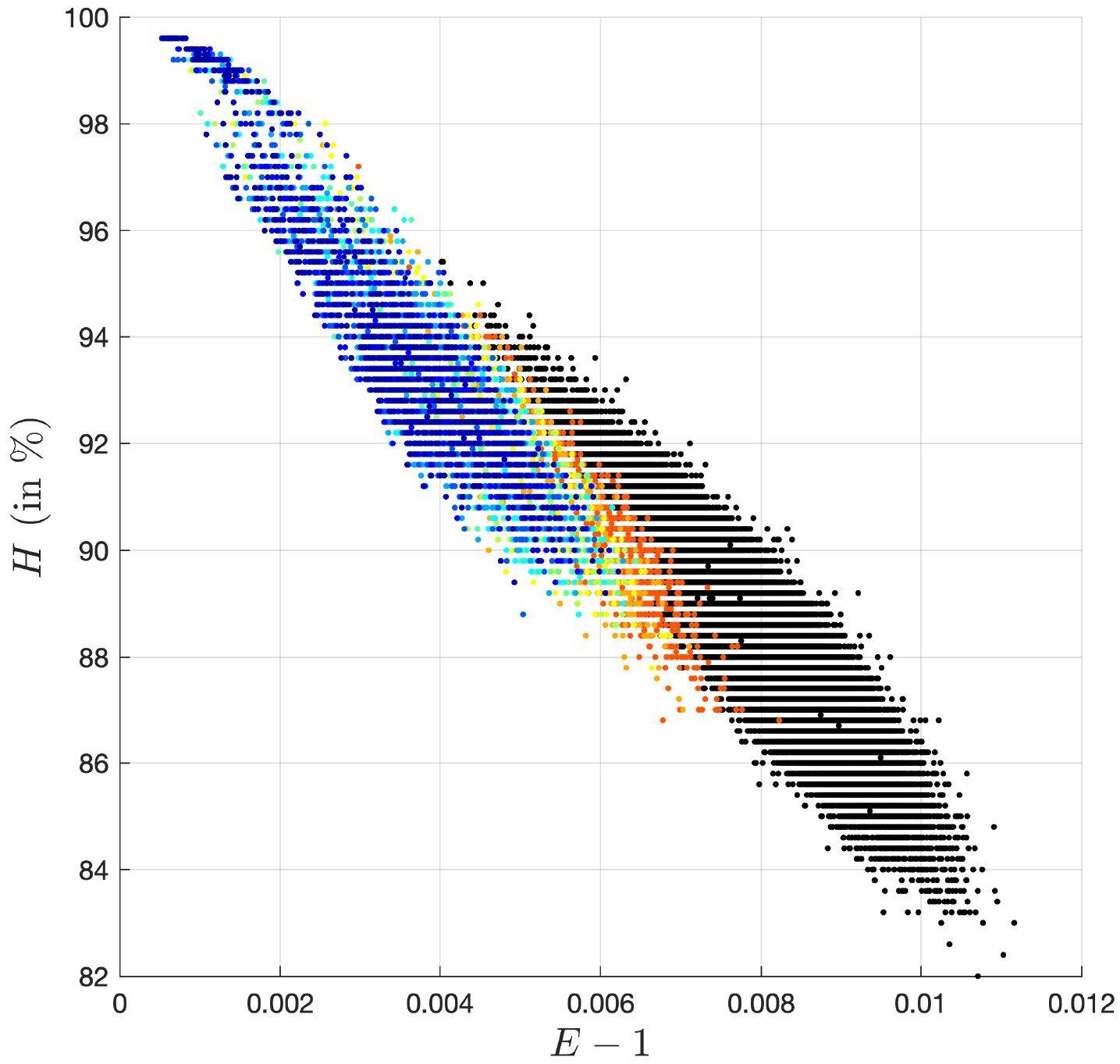}\tabularnewline
\end{tabular}
\end{center}
\caption{Scatter plots of $E$ vs. $H$ and $E$ vs. $R_{\epsilon}$ displaying
the same data as Figures \ref{fig: data N1000} (a) through (f)  from Example 3.}
\label{fig:R vs H_scatter plots}
\end{figure}

We have taken $\epsilon$ to be fixed at $0.5\%$ because this value
is bounded above by the deviation from $r_{hex}$ gotten from the
$\delta-$perturbation of a single generator in the honey comb PCVT,
yet it remains big enough so that the data clearly shows that:

\begin{itemize}
\item a higher variation $|\Delta E|$ between consecutive iterations in
our method results in a higher $|\Delta R_{\epsilon}|$ than $|\Delta H|$,
and this regardless of the block \textit{MACN-c}/Lloyd/\textit{MACN-}$\delta$
\item we have systematically that $|f_{E-1}^{*}-f_{R_{\epsilon}}^{*}|<|f_{E-1}^{*}-f_{H}^{*}|$ 
\item the ECDFs of $H$ present larger discontinuity jumps that the ones
of $R_{\epsilon}$; meaning that for given $\mathbf{X}^{*}$ the number
of computed PCVTS states sharing the value $H(\mathbf{X}^{*})$ is
higher than the one sharing $R_{\epsilon}(\mathbf{X}^{*})$.\\
\end{itemize}
These observations combined point out that $R_{\epsilon}$ is indeed a measure more faithful to $E$
and a better indicator of ``well distributed'' PCVTs than $H$ is.
We provide further insight on this matter in Table \ref{tab: Correlation ratios varrho} through the correlation ratio
\[
\varrho:=\frac{\sigma_{R_{\epsilon}}\text{cov}(E-1,H)}{\sigma_{H}\text{cov}(E-1,R_{\epsilon})}
\]
and in Figure \ref{fig:R vs H_scatter plots} through scatter plots of the data from Example 3.

\begin{table}[H]
\begin{center}
\resizebox{0.95\columnwidth}{0.085\columnwidth}{

\begin{tabular}{|>{\centering}p{9mm}||cccccccccc||ccc|}
\hline 
\multicolumn{1}{|>{\centering}p{9mm}|}{{\scriptsize{}Ex}} &  &  &  &  & {\scriptsize{}Hybrid} &  &  &  &  &  & \multicolumn{1}{c|}{\multirow{2}{*}{{\scriptsize{}Lloyd}}} & \multicolumn{1}{c|}{\multirow{2}{*}{{\scriptsize{}LBFGS}}} & \multirow{2}{*}{{\scriptsize{}PLBFGS}}\tabularnewline
\cline{2-11} 
\multicolumn{1}{|>{\centering}p{9mm}|}{} & $\mathbf{X}_{0}^{*}$ & $\mathbf{X}_{1}^{*}$ & $\mathbf{X}_{2}^{*}$ & $\mathbf{X}_{3}^{*}$ & $\mathbf{X}_{4}^{*}$ & $\mathbf{X}_{5}^{*}$ & $\mathbf{X}_{6}^{*}$ & $\mathbf{X}_{7}^{*}$ & $\mathbf{X}_{8}^{*}$ & $\mathbf{X}_{9}^{*}$ &  &  & \tabularnewline
\hline 
1 & 1.0621 & 1.0924 & 1.0838 & 1.0675 & 1.0550 & 1.0687 & 1.0674 & 1.0581 & 1.0523 & 1.0578 & 1.0750 & 1.0663 & 1.0666\tabularnewline
2 & 1.0691 & 1.0594 & 1.0732 & 1.1228 & 1.1250 & 1.1451 & 1.1594 & 1.1924 & 1.1515 & 1.1035 & 1.0817 & 1.0688 & 1.0683\tabularnewline
\hline 
\hline 
3 & 1.0636 & 1.0717 & 1.0746 & 1.0744 & 1.0667 & 1.0605 & 1.0561 & 1.0555 & 1.0476 & 1.0457 & 1.0822 & 1.0671 & 1.0678\tabularnewline
4 & 1.0995 & 1.0690 & 1.0646 & 1.0964 & 1.0759 & 1.0649 & 1.0743 & 1.0806 & 1.0960 & 1.0849 & 1.0869 & 1.0692 & 1.0673\tabularnewline
5 & 1.0644 & 1.0608 & 1.0752 & 1.0490 & 1.0382 & 1.0230 & 1.0389 & 1.0449 & 1.0455 & 1.0430 & 1.0820 & 1.0689 & 1.0679\tabularnewline
6 & 1.0566 & 1.0586 & 1.0742 & 1.0803 & 1.0724 & 1.0736 & 1.0602 & 1.0447 & 1.0508 & 1.0449 & 1.0842 & 1.0688 & 1.0685\tabularnewline
\hline 

\end{tabular}
} 
\end{center}
\vspace{3mm}
\caption{Correlation ratios $\varrho$ for each method from the PCVT data presented
in Examples 1 through 6 \label{tab: Correlation ratios varrho}}
\end{table}

\section{Deterministic vs. Stochastic}\label{sec:6}

In this section we present evidence of both positive and negative aspects of the fully deterministic nature of our method vis-\`a-vis of stochastic alternatives. For this we start by comparing our \textit{MACN} algorithm with the global Monte Carlo Method presented in \cite{Global_MCM_Lu} (hereafter referred to as \textit{MCM}) and we finish the section by considering variants of centroidal dislocations to ratify our use of \textit{MACN-$\delta$}.

\subsection{\textit{MACN} vs. \textit{MCM}}
\begin{figure}[H]
\begin{center}
\begin{tabular}{cc}
\includegraphics[width=0.45\columnwidth,height=0.31\columnwidth]{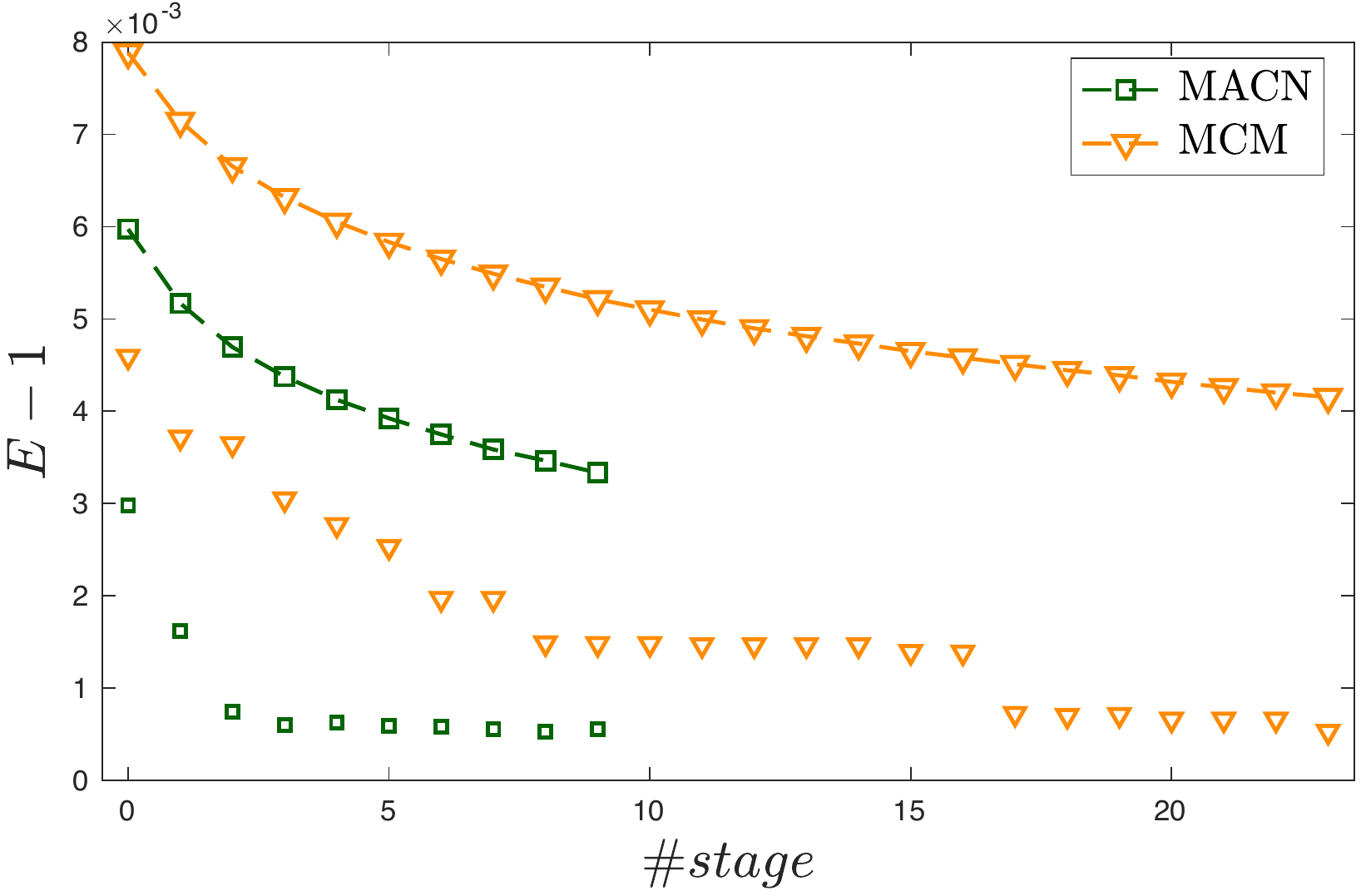} & \includegraphics[width=0.45\columnwidth,height=0.30\columnwidth]{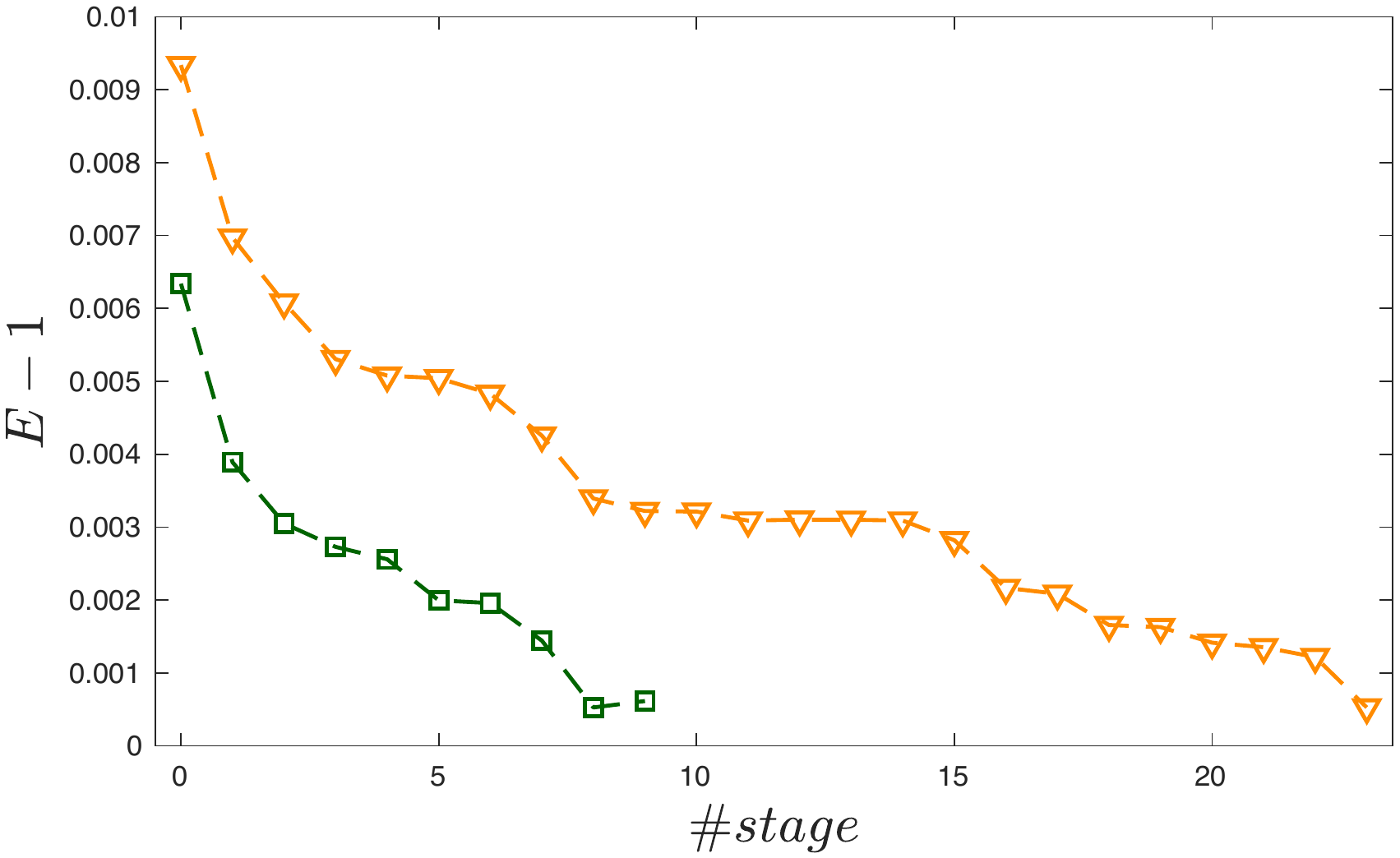}\tabularnewline
\end{tabular}
\end{center}
\caption{ Comparison of $Q=10$ stages of \textit{MACN} ($K=6000$) and $24$ stages of \textit{MCM} for $N=1000$: \textit{(left)} Averages $\langle E-1\rangle$ over 1000 runs are shown with joint markers while minimums achieved are isolated ones. \textit{(right)} The energy of the PCVTs obtained for the respective runs that first achieved the lowest recorded value $E^{min}-1=5.27$e-$4$}
\label{fig:MACN vs MCM}
\end{figure}

We implemented the \textit{MCM} method for $N=1000$ with the parameter values used for constant densities in \cite{Global_MCM_Lu}, the  results along with the \textit{MACN} data from Example 3 of \textsection \ref{sec:5} are summarized in Figure \ref{fig:MACN vs MCM}. We depict first the energy averages $\langle E-1\rangle$ and minimums obtained over 1000 runs, we also show the energy of each PCVT obtained during two runs respectively achieving $E^{min}-1=5.27$e-4 at the earliest iteration.\\
It becomes clear from these statistics that our method has a lower probing number than the comparative stochastic approach, e.g. Figure \ref{fig:MACN vs MCM} (\textit{right}) shows we need to compute 9 versus 24 PCVTs to achieve the same low energy basins of attraction as \textit{MCM}.\\
Moreover, as already noted in \textsection \ref{sec:3}, a crucial advantage of our method is the simple tuple $\{Q,K\}$ of parameters to adjust compared to the more complex set $\{K,T_0,T_k,h\}$ used by \textit{MCM} (note that $K$ serves a different purpose in each method).\\
On the other hand, \textit{MCM} benefits from having a substantially lower computation time: 1) since a Quasi-Newton method is used instead of Lloyd's and 2) due to the lack of preconditioning this alternative avoids a total of $K\times Q$ iterations of complexity $O(N\log(N))$ compared to \textit{MACN}. Thus making \textit{MCM} still desirable for applications where fast computation time is critical.\\

\subsection{Random direction versus \textit{MACN} motions}

As defined in (\ref{eq: original delta}), the value of $\delta$ is crucial since it has the peculiarity of making our perturbation stage preserve a certain regularity in the structure of the tessellation, however the direction of the perturbation
seems to be of primary importance compared to the step size when the
latter is fixed. To gain insight on this matter we define three variants
of our perturbation step:\\

\begin{enumerate}
\item Inspired by the relation between $\delta$ and the intrinsic length-scale
of the regular hexagonal lattice discussed in \textsection \ref{sec:4}, the first variant consists of moving
away from the closest neighbor by the length-scale proper to each
cell in the PCVT. Precisely, the perturbation follows
\begin{equation}
x_{i}\leftarrow x_{i}+{\displaystyle \frac{|V_{i}|}{|\partial V_{i}|}}\,\frac{x_{i}-x_{j_{i}^{*}}}{||x_{i}-x_{j_{i}^{*}}||}\qquad i=1,...,N\label{eq: variant of delta perturb (intrinsic)}
\end{equation}
\item The next variant contemplates $\delta$ as in (\ref{eq: original delta})
but choses a random neighbor $x_{j}\;,j\in\mathscr{N}_{i}$ to move
away from, thus not necessarily being the closest one.\\
\item Finally, we consider each generator moving by the distance $\delta$
and at a random angle $\theta_{i}\in[0,2\pi)$ taken from the uniform
distribution.\\
\end{enumerate}

\begin{figure}[H]
\begin{centering}
\includegraphics[width=0.5\columnwidth]{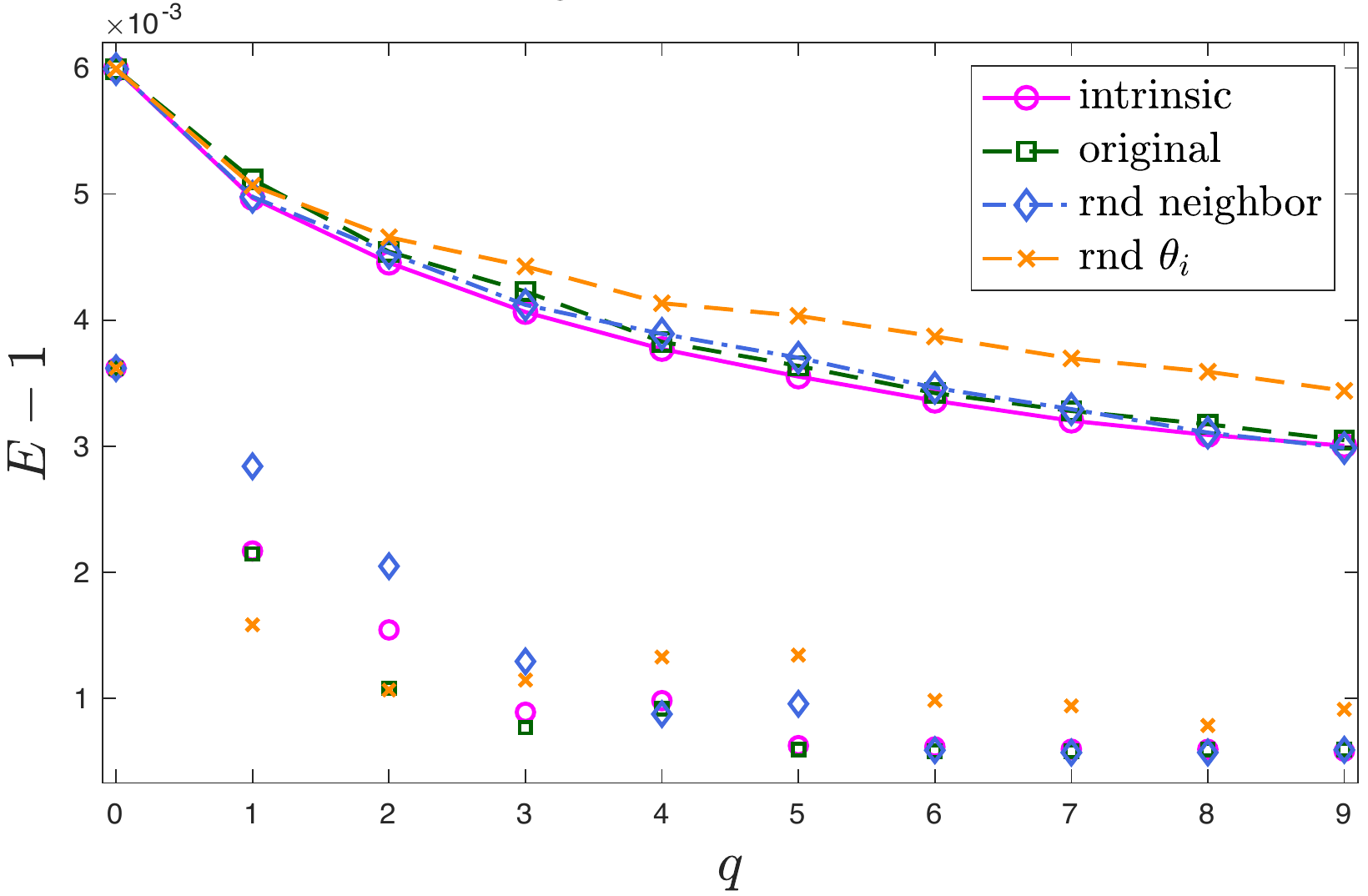}
\par\end{centering}
\caption{Comparison in the energy performance of the original \textit{MACN}-$\delta$
step and its three variants across stages of our hybrid method with
$N=1000$ and $K=6000$. The joint markers represent averages while
the isolated ones represent minima over 100 runs. \label{fig:Comparison of MACN-delta variants}}
\end{figure}

Figure \ref{fig:Comparison of MACN-delta variants} illustrates how
close the performance of our variants of the \textit{MACN}-$\delta$
step are from one another but with particular distinction of the random
angle $\theta_{i}$ perturbation, being then of some reassurance that
our original neighbor-guided dislocation of PCVTs is better suited than some
random search in a $\delta-$vicinity.

\section{Closing Remarks and future directions}\label{sec:7}

We have introduced and assessed a simple deterministic method for navigating the energy landscape of CVTs in two dimensions.
 This deterministic coupling algorithm: i) only has two degrees of freedom, ii) shows remarkable robustness with respect to the increasing non-convexity of the energy as $N$ grows larger, iii) statistically allows us to systematically obtain configurations closer to the ground state compared to the current state of the art deterministic methods and iv) finds energetically comparable results to the leading stochastic method while needing fewer CVT computations (lower probing number).
We also introduced the isoperimetric ratio via $R_{\epsilon}$ as an indicator of low energy CVTs. 

We point out that while we prioritized simplicity of the method in this paper, the algorithm's performance could be further improved, if needed be, by: i) using initial quasi-random distributions, ii) replacing Lloyd's with other gradient based descent methods that satisfy Wolfe conditions, and most importantly iii) by introducing a suitable decay in the sequence $\{K_q\}_{q=0}^{Q-1}$ (possibly adapted to $\{E(\textbf{X}_q^*)\}_{q=0}^{Q-1}$).
The point made is that even with the crude tunings made on $K_q\equiv K$ in \textsection \ref{sec:5}, the resulting regularity measures are remarkable.

It would be natural to explore our global method in two different settings. 
First,  explore our algorithm on the 3D cubic torus wherein we expect similar comparative results with the appearance of the BCC lattice and truncated octahedron Voronoi cells.
Second, explore our algorithm on the 2-sphere. As we previously mentioned in the Introduction, work in progress shows that that our algorithm works well on the 2-sphere, systematically obtaining configurations closer to the ground state.   Of  particular interesting on the sphere is the nature of the ground state. For certain values of $N$, for  example $N=32$, the folklore suggests that optimality is tied to {\it the soccer ball configuration}:  
regular hexagons except for exactly 12 regular pentagonal defects. Our results on the 2-sphere as well as structural results of the ground state are forthcoming. 
We note, however, that in the extension of our method to manifolds the increased complexity of the underlying computation of the Voronoi tessellation and centroids needs to be compensated by reducing the product $K\times Q$. 
In particular,   time tractability will be at the cost of energy efficiency.

A different related question is the inclusion of an underlying  inhomogeneous probability densities $\rho$ over $\Omega$ wherein the energy ({\ref{Vor-energy}) takes the form 
$
F(\mathbf{X})\, =\, \sum_{i=1}^N \int_{V_i} ||y-x_i||^2\, \rho (y) \, dy. 
$
However, here it is unclear how to choose the distance $\delta$ in the \textit{MACN} annealing step. In order to gain a better insight on how to tackle this generalized problem we strive to cast both our \textit{MACN} iterations as a low order approximation of some yet to be determined gradient flow involving the Voronoi energy.

\bibliographystyle{siam}
\bibliography{NavigatingCVTLandscape}

\newpage
\appendix

\section{Explicit pseudo code for P-L-BFGS method\label{App:C}}

\begin{algorithm}[H]
\caption{P-L-BFGS($M,T$)\label{PLFGS Algorithm-1}}

\begin{algorithmic}

\State \textbf{Prior definitions}: to ease notation we define at each iteration $k$

\[
\begin{array}{c}
 s_{k}:=\mathbf{X}^{(k+1)}-\mathbf{X}^{(k)} \, ;\, y_{k}:=DE(\mathbf{X}^{(k+1)})-DE(\mathbf{X}^{(k)}) \\
 \rho_k:=\dfrac{1}{y_k^{\top}s_k} \text{ and } H_{0}^{(k)}:=\dfrac{s_{k-1}^{\top}y_{k-1}}{y_{k-1}^{\top}y_{k-1}}I
\end{array}
\]
\State \textbf{Input:} i) initial iterate $\mathbf{X}^{(0)}$; ii) integer parameters $M$ and $T$; iii) tolerance $tol$ for convergence\\

\State set $k=0$
\State set \textit{diff}=Inf
 
	\While{\textit{diff > tol}}		
		\State set $q=DE(\mathbf{X}^{(k)})$\\
		\State $1^{\text{st}}$ \textit{L-BFGS update}
		\For {$i=k-1:-1:k-M$}
		\State $a_{i}=\rho_i s_i^{\top}q$
		\State $q\leftarrow q-a_iy_i$
		\EndFor\\
		
		\State \textit{Redirect search direction}
		\If {$k \mod T=0$}
		\State construct preconditioner matrix $\widetilde{A}_k$ and solve the system $\widetilde{A}_k \,r=q$
		
		\Else
		\State construct $H_{0}^{(k)}$ and set $r=H_{0}^{(k)}q$
		\EndIf\\
		
		\State $2^{\text{nd}}$ \textit{L-BFGS update}
		\For {$i=k-M:k-1$}
		
		\State $r\leftarrow r+s_i(a_i-\rho_i y_i^{\top}r)$
		\EndFor\\

		\State set descent direction $p^{(k)}=-r$
		\State update iterate $\mathbf{X}^{(k+1)}=\mathbf{X}^{(k)}+\alpha^{(k)}p^{(k)}$  where $				\alpha^{(k)}$ is a step length satisfying  \State the strong Wolfe conditions\\

		\If {$k>M$}
		\State erase the tuple $\{s_{k-M};y_{k-M}\}$
		\State compute and store $\{s_{k};y_{k}\}$
		\EndIf\\

		\State \textit{diff}$=||D\,E||/N$
		\State $k\leftarrow k+1$
	\EndWhile\\

\State \textbf{Output:} $\textbf{X}^{*}$, a stable local minimizer of \textit{E} and its corresponding PCVT and PDT

\end{algorithmic}
\end{algorithm}

\end{document}


\title{\textsc{A simple geometric method for navigating  the energy landscape of centroidal voronoi tessellations}}

\author{Ivan Gonzalez\thanks{Department of Mathematics and Statistics,
McGill University, Montreal, Quebec, H3A0B9
  (\email{ivan.gonzalez@mail.mcgill.ca}).}
\and Rustum Choksi\thanks{Department of Mathematics and Statistics,
McGill University, Montreal, Quebec, H3A0B9
 \newline  (\email{rustum.choksi@mcgill.ca}, \url{http://www.math.mcgill.ca/rchoksi/}).}
\and Jean-Christophe Navei\thanks{Department of Mathematics and Statistics,
McGill University, Montreal, Quebec, H3A0B9 
 \newline (\email{jean-christophe.nave@mcgill.ca}, \url{https://www.math.mcgill.ca/jcnave/}).}
}

\maketitle

\section{Additional visual material on \textsection 5}

We provide below further graphics illustrating our computations and the data obtained for Examples 2,4,5 and 6: namely the PCVTs of optimal configurations achieved as well as the ECDFs of the regularity measures. 

\vspace{0.5cm}

\begin{figure}[h!]
\centering

\begin{tabular}{cc}
\includegraphics[width=0.38\columnwidth]{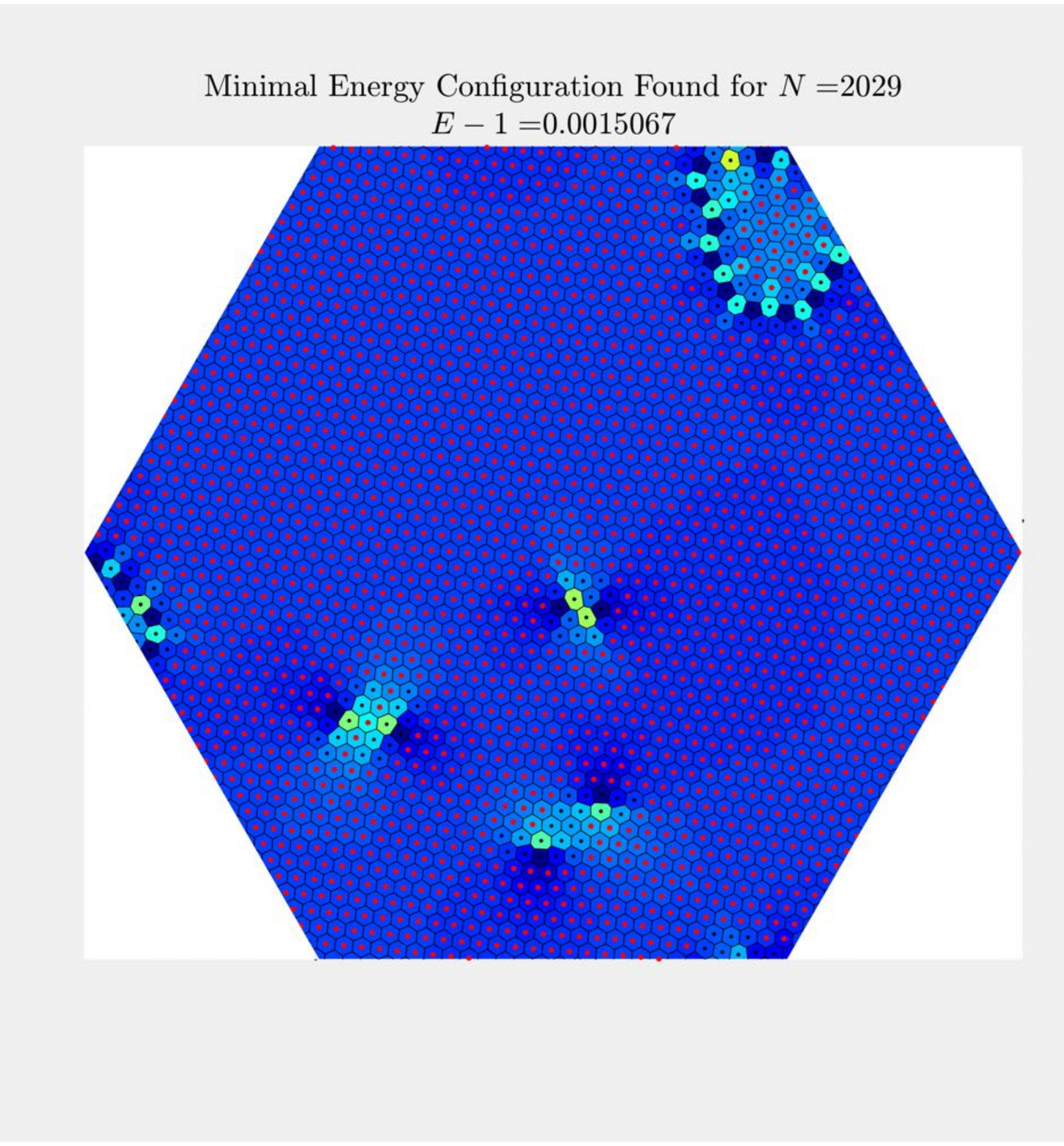} & \includegraphics[width=0.38\columnwidth]{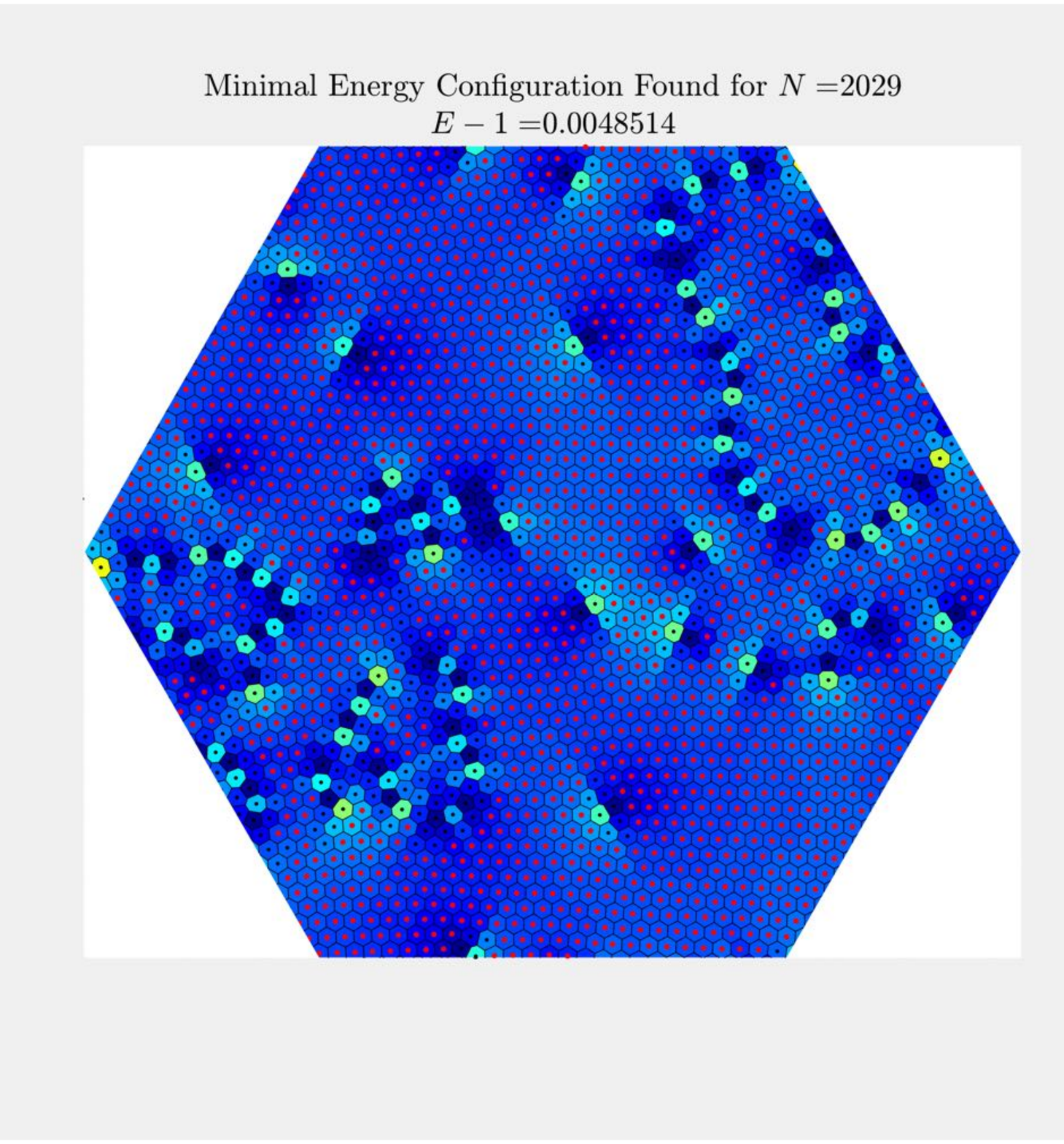}\tabularnewline
{\footnotesize{}(a) PCVT reaching $E^{min}$} & {\footnotesize{}(b) PCVT achieving $E^{ref}$}\tabularnewline
\end{tabular}

\begin{tabular}{cc}
\includegraphics[width=0.45\columnwidth, height=0.26\columnwidth]{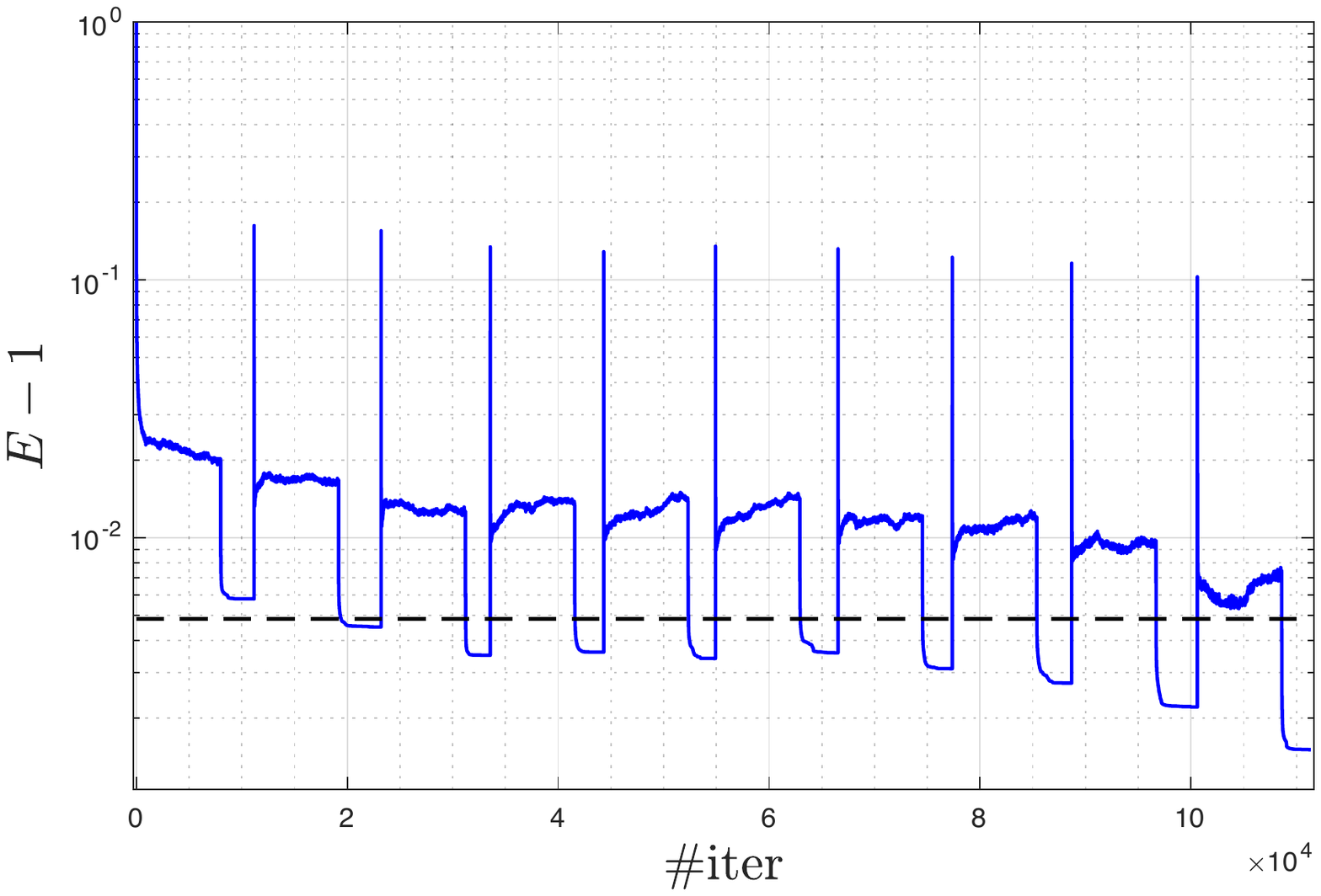} & \includegraphics[width=0.45\columnwidth, height=0.26\columnwidth]{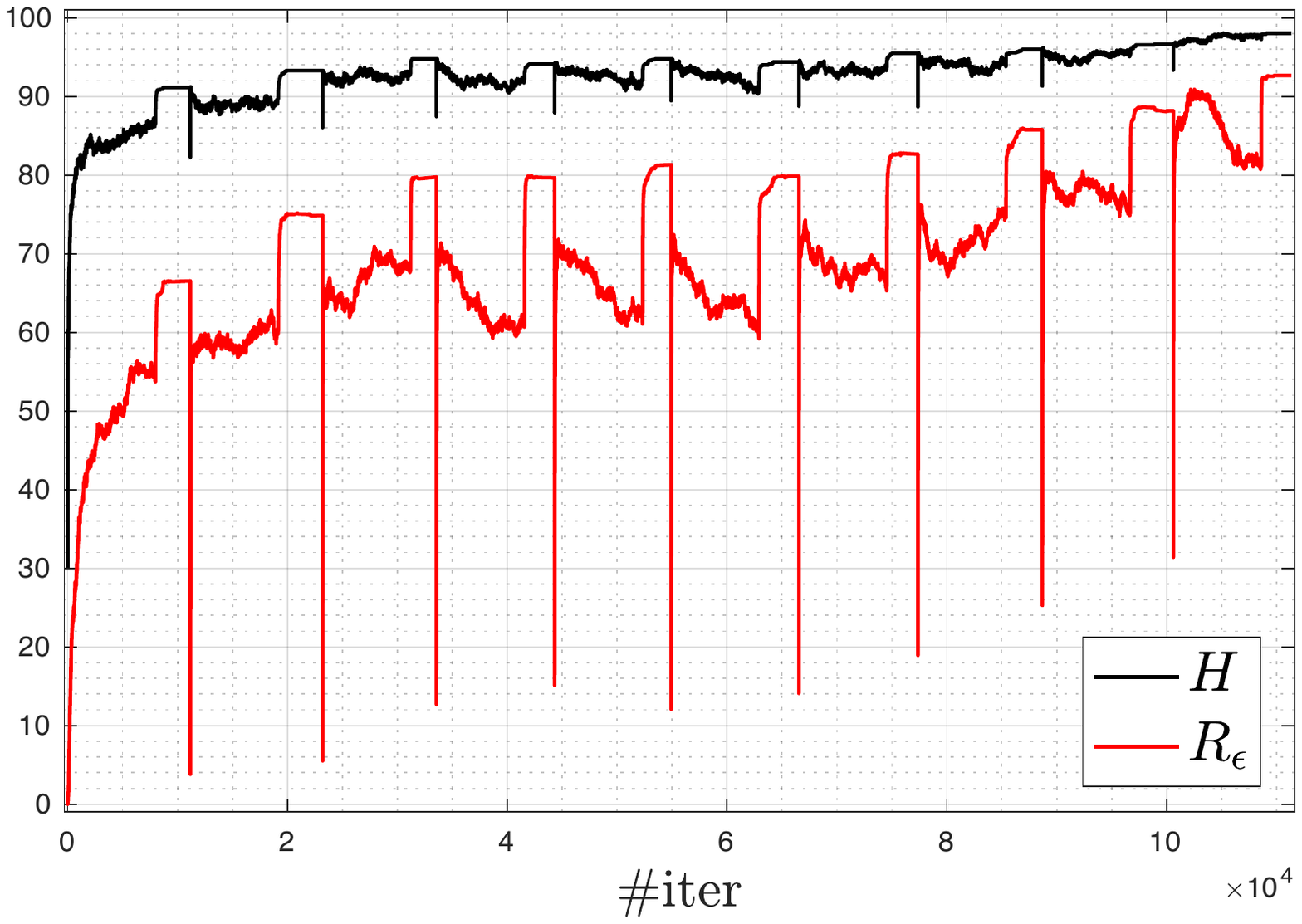}\tabularnewline
{\footnotesize{}(c) hybrid run reaching $E^{min}$ } & {\footnotesize{}(d) $R_{\epsilon},H$ profiles corresponding to (c)}\tabularnewline
\end{tabular}

\caption{Optimal PCVTs from Example 2 for $N=2029$ with $K=8000$: (a) lowest energy tessellation reached by the hybrid method out of 100 runs with $E^{min}-1=0.00150$. (b) tessellation carrying $E^{ref}-1=0.00485$ achieved with L-BFGS(7) out of 100,000 runs. (c) and (d) show the regularity measures as a function of the iterations of the particular hybrid run that achieved $E^{min}$.\label{fig: best configurations N2029}}

\end{figure}

\newpage

\begin{figure}[H]
\centering

\begin{tabular}{c}
\includegraphics[width=0.6\columnwidth, height=0.28\columnwidth]{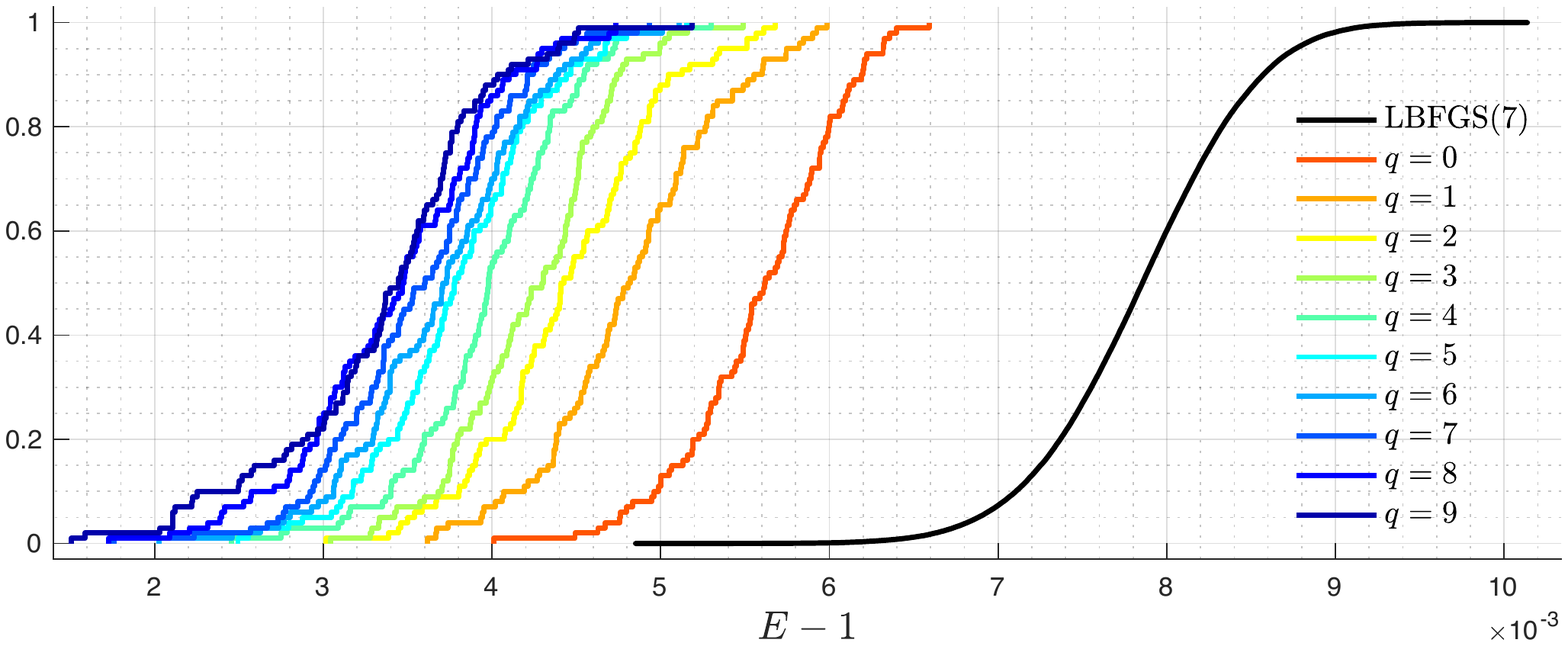}\tabularnewline
{\footnotesize{}(a) $f_{E-1}$}\tabularnewline
\end{tabular}

\begin{tabular}{cc}
\includegraphics[width=0.43\columnwidth, height=0.26\columnwidth]{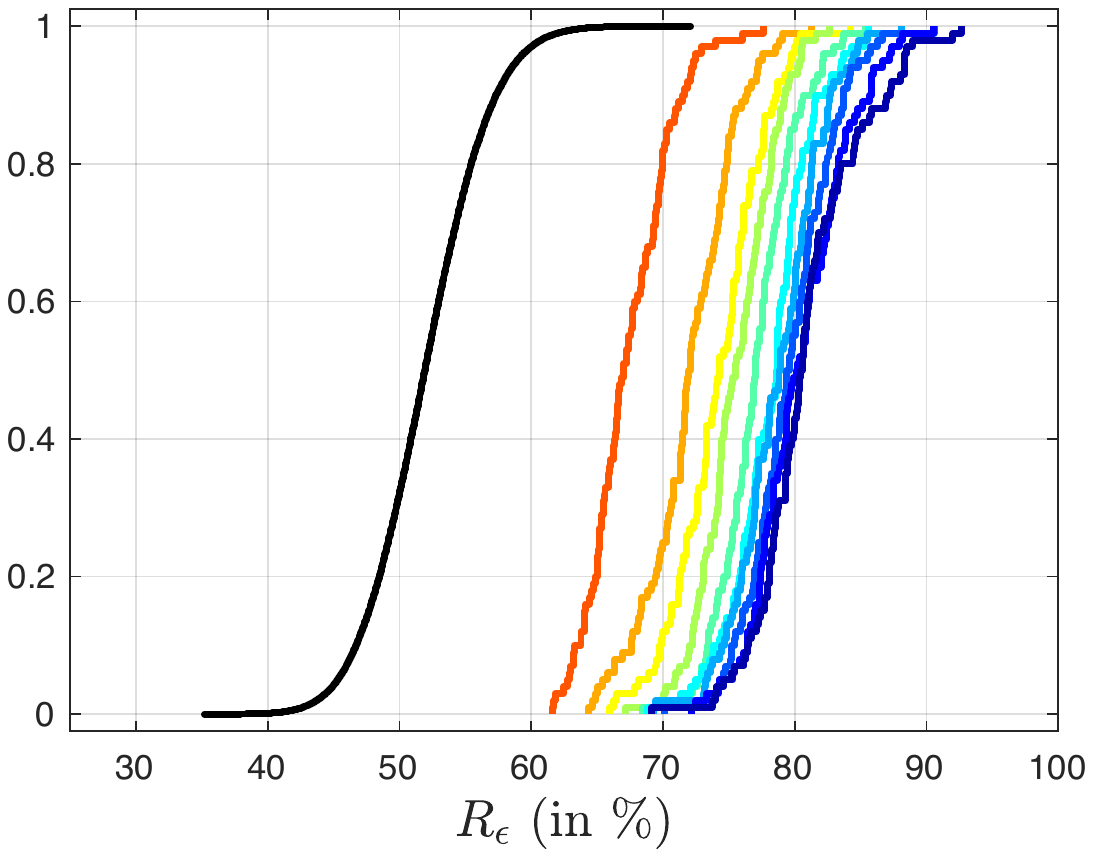} & \includegraphics[width=0.43\columnwidth, height=0.26\columnwidth]{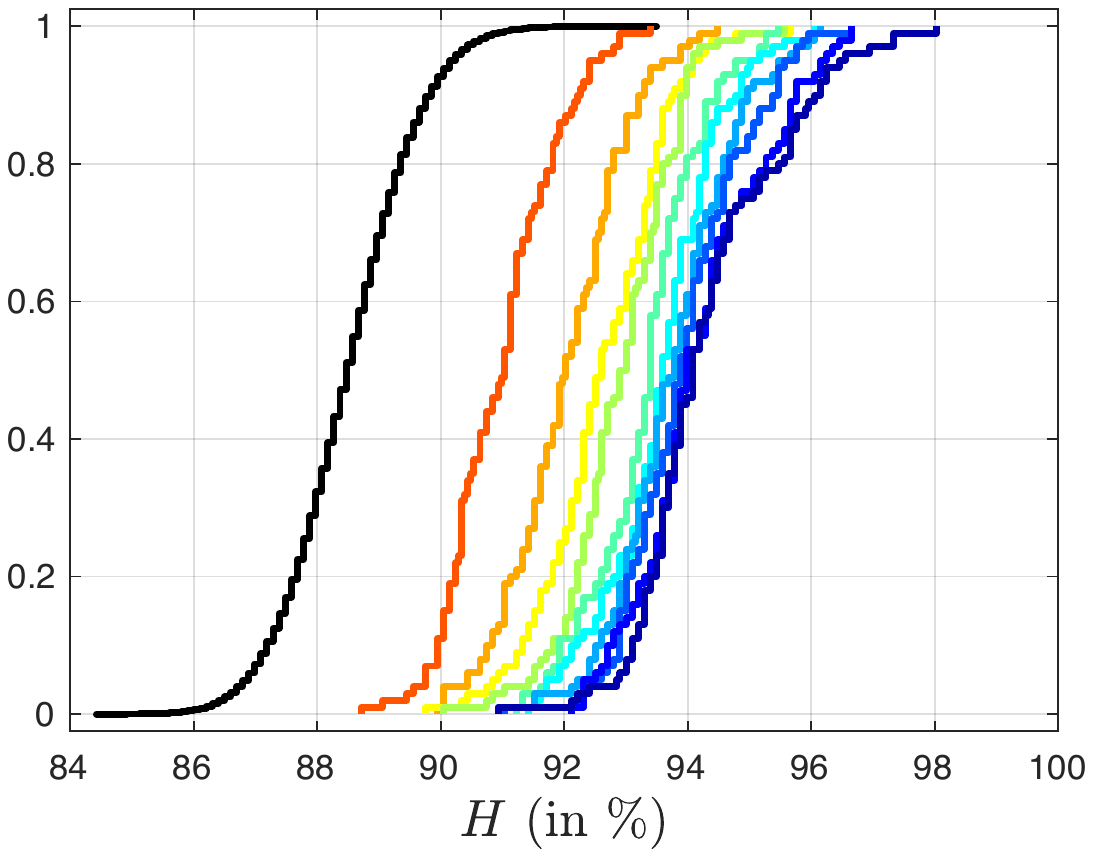}\tabularnewline
{\footnotesize{}(b) $f_{R_{\epsilon}}$} & {\footnotesize{}(c)$f_{H}$}\tabularnewline
\end{tabular}

\caption{Data from Example 2 with $N=2029$: ECDFs of the three regularity measures for 100 runs of the hybrid method with $K=8000$ along with 100,000 runs of L-BFGS(7)}

\end{figure}

\begin{figure}[H]
\centering

\begin{tabular}{cc}
\includegraphics[width=0.305\columnwidth]{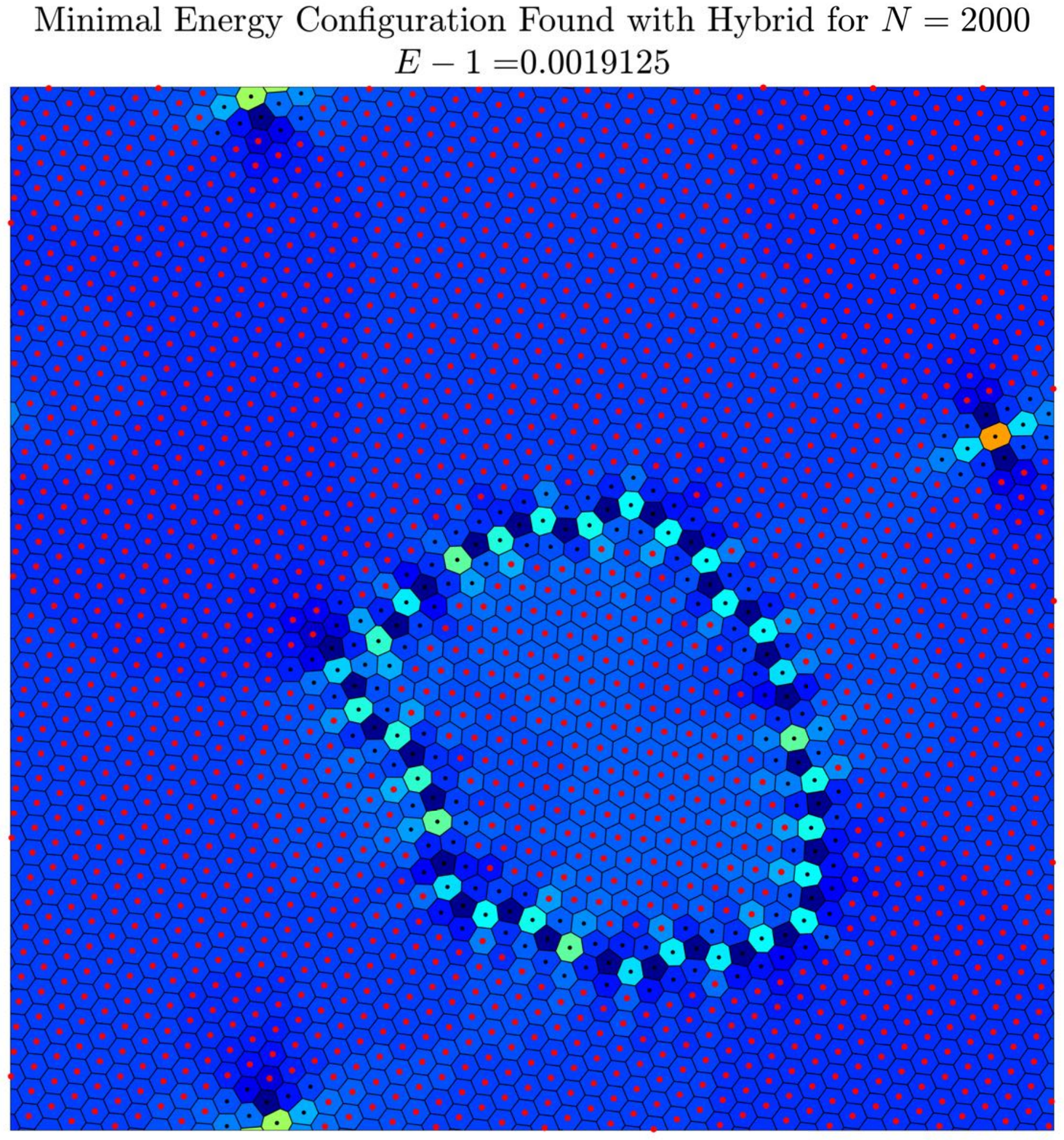} & \includegraphics[width=0.305\columnwidth]{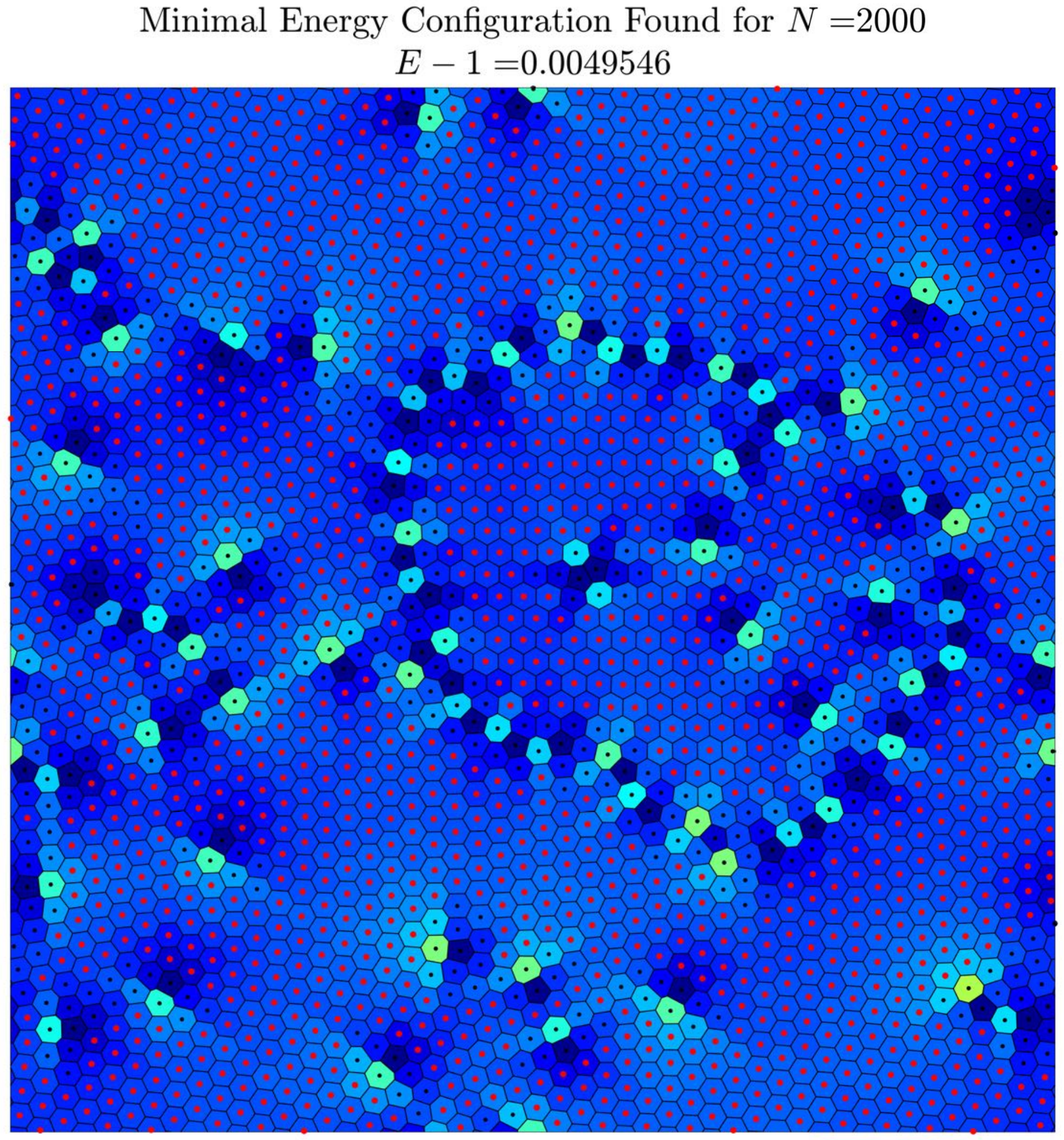}\tabularnewline
{\footnotesize{}(a) PCVT reaching $E^{min}$} & {\footnotesize{}(b) PCVT carrying $E^{ref}$}\tabularnewline
\end{tabular}

\begin{tabular}{cc}
\includegraphics[width=0.45\columnwidth, height=0.26\columnwidth]{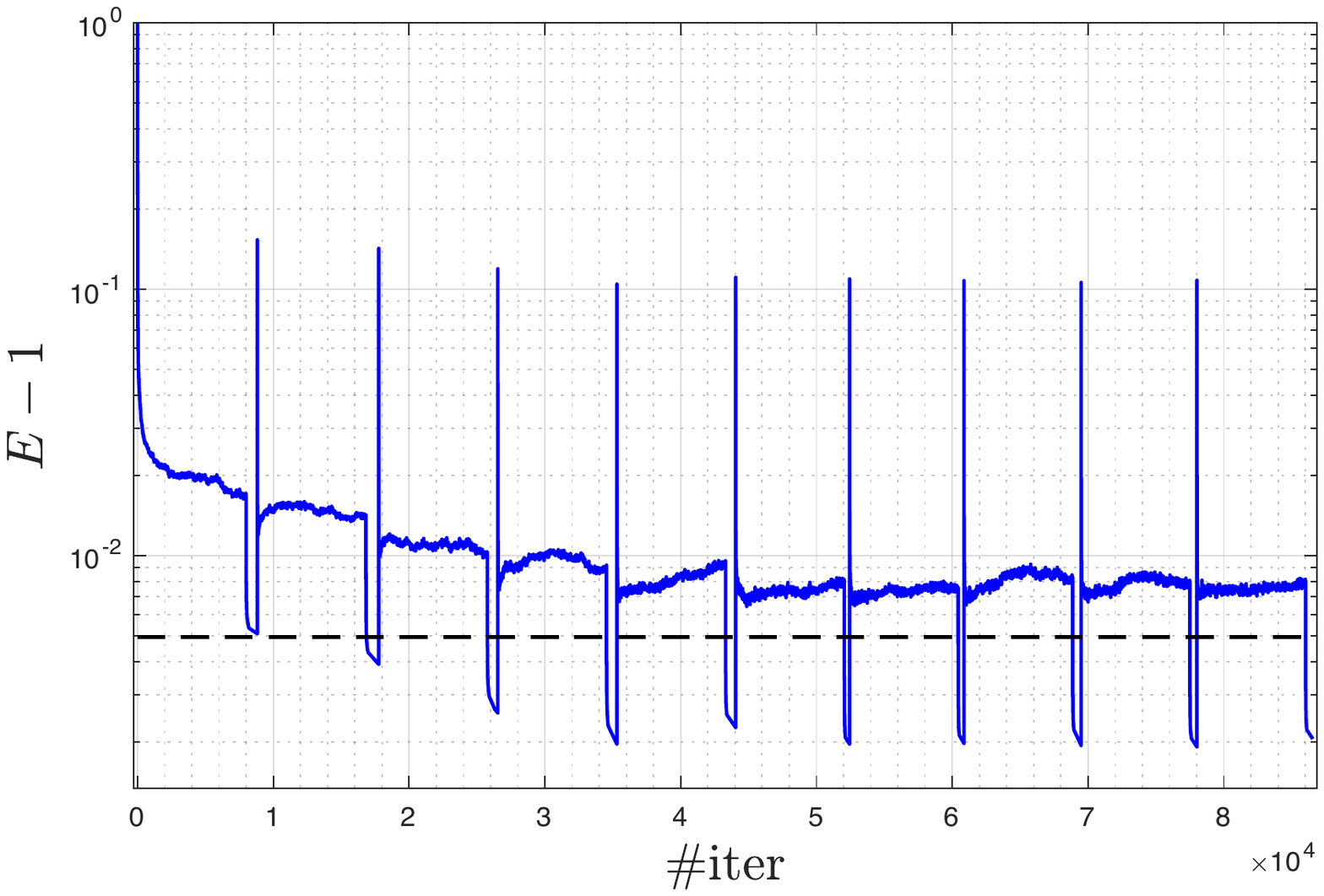} & \includegraphics[width=0.45\columnwidth, height=0.26\columnwidth]{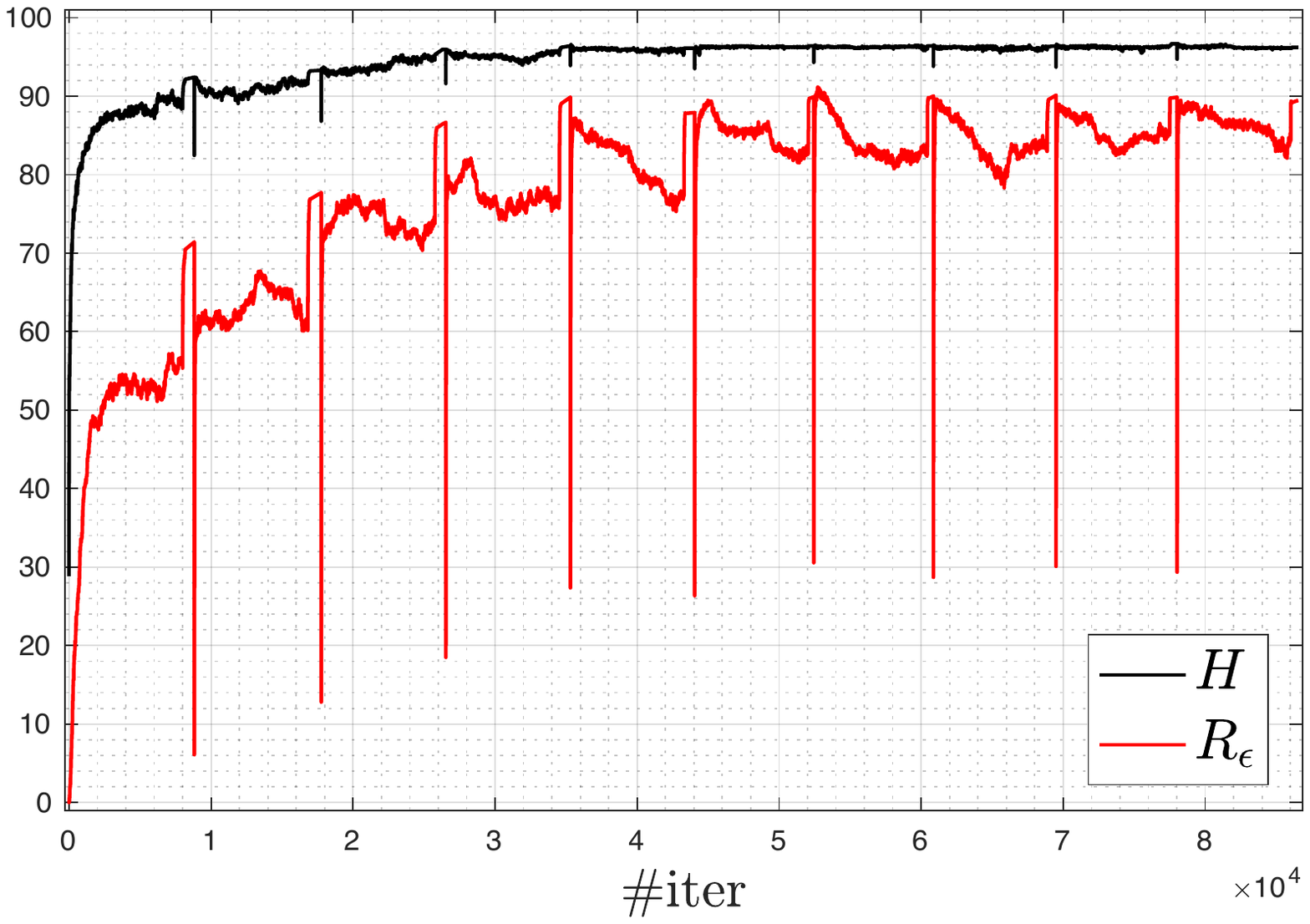}\tabularnewline
{\footnotesize{}(c) hybrid run achieving $E^{min}$ } & {\footnotesize{}(d) $R_{\epsilon},H$ profiles associated with (c)}\tabularnewline
\end{tabular}

\caption{Optimal PCVTs from Example 4 for $N=2000$ with $K=8000$: (a) lowest energy PCVT reached by the hybrid method out of 100 runs with $E^{min}-1=0.00191$. (b) tessellation carrying $E^{ref}-1=0.00495$ achieved with P-L-BFGS(20,20) out of 100,000 runs. (c) and (d) show the regularity measures as a function of the iterations of the particular hybrid run that achieved $E^{min}$.}

\end{figure}

\begin{figure}[H]
\centering

\begin{tabular}{c}
\includegraphics[width=0.6\columnwidth, height=0.28\columnwidth]{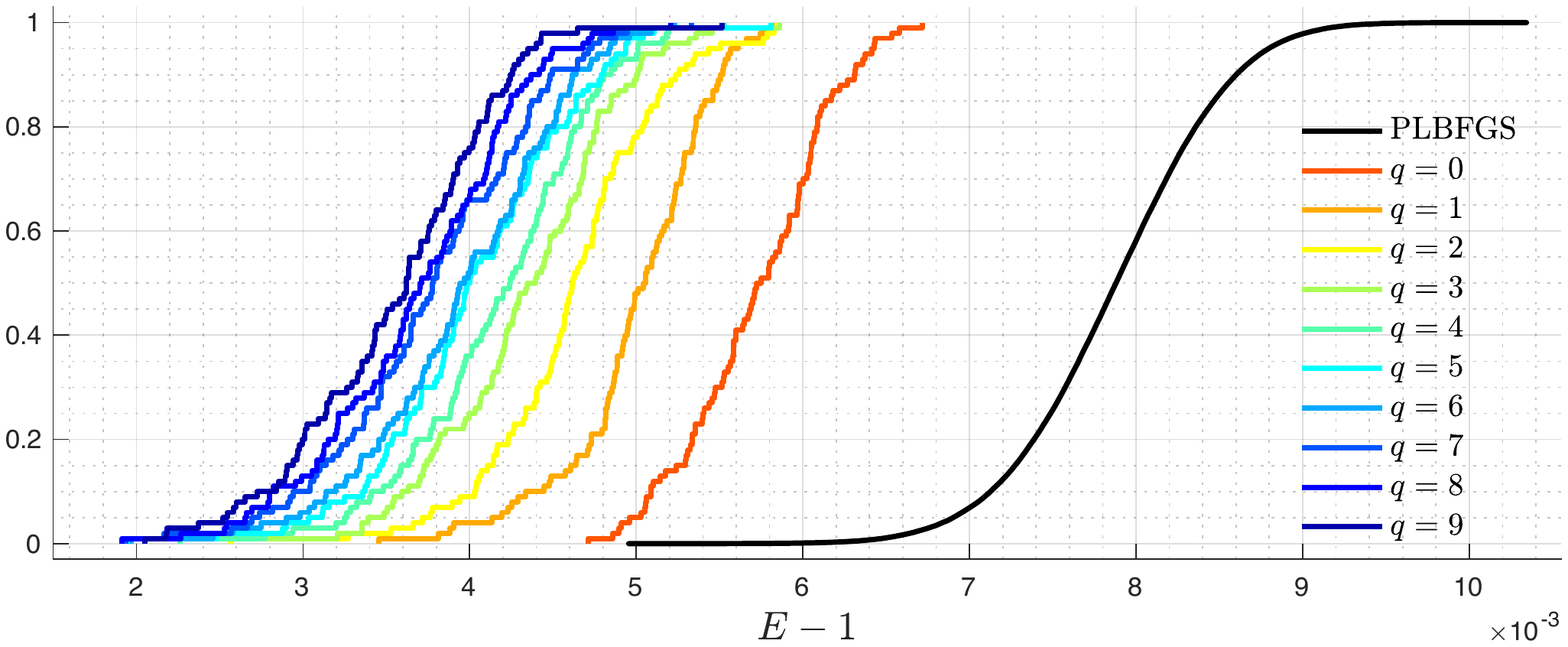}\tabularnewline
{\footnotesize{}(a) $f_{E-1}$}\tabularnewline
\end{tabular}

\begin{tabular}{cc}
\includegraphics[width=0.43\columnwidth, height=0.26\columnwidth]{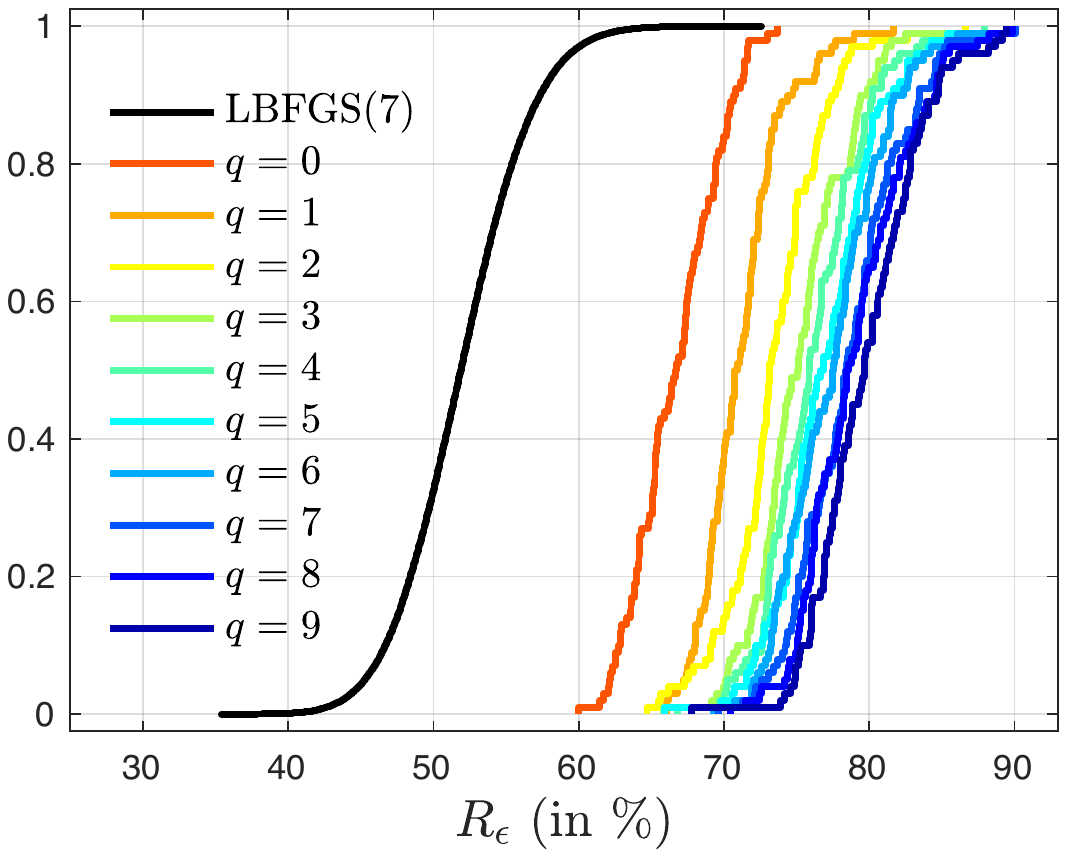} & \includegraphics[width=0.43\columnwidth, height=0.26\columnwidth]{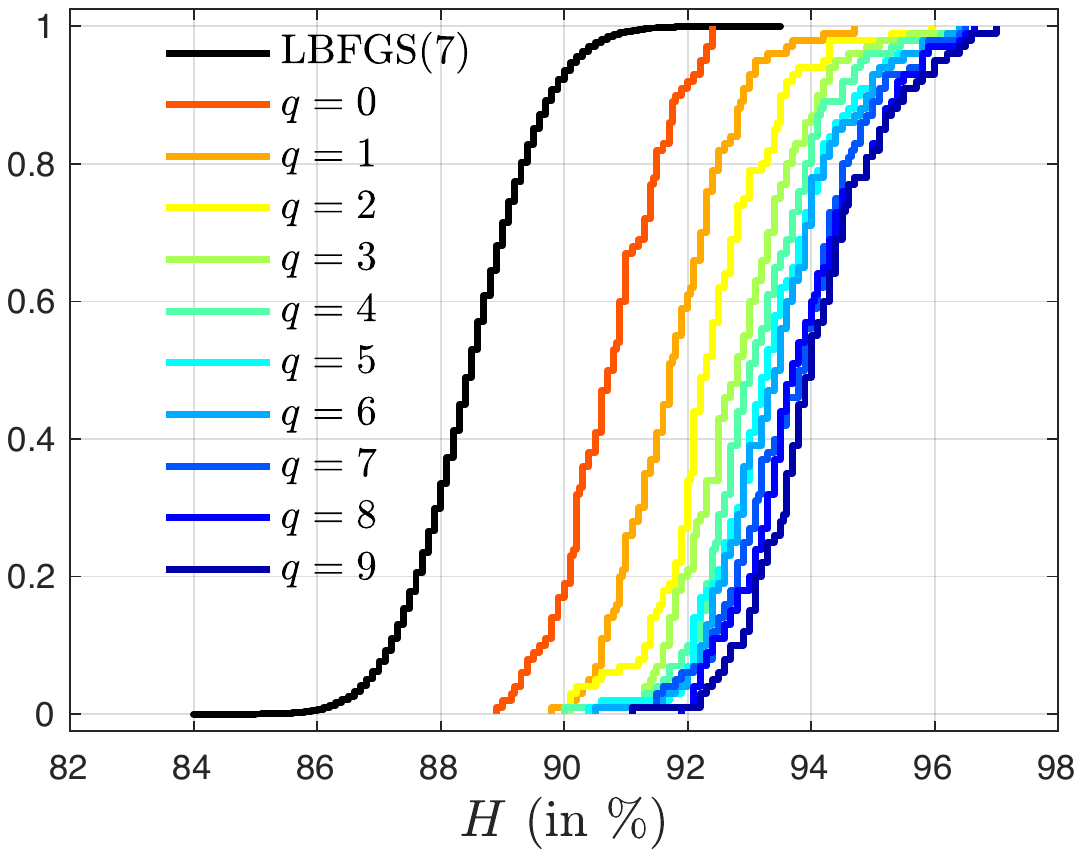}\tabularnewline
{\footnotesize{}(b) $f_{R_{\epsilon}}$} & {\footnotesize{}(c)$f_{H}$}\tabularnewline
\end{tabular}

\caption{Data from Example 4 with $N=2000$: ECDFs of our regularity
measures obtained with $K=8000$, for $H$ and $R_{\epsilon}$ we show the cumulative function
of L-BFGS(7) out of 100,000 runs rather than P-L-BFGS(20,20) since it is the former that
achieved the highest recorded values $H^{ref}=93.50\%$ and $R_{\epsilon}^{ref}=72.54\%$.}
\end{figure}

\begin{figure}[H]
\centering

\begin{tabular}{cc}
\includegraphics[width=0.305\columnwidth]{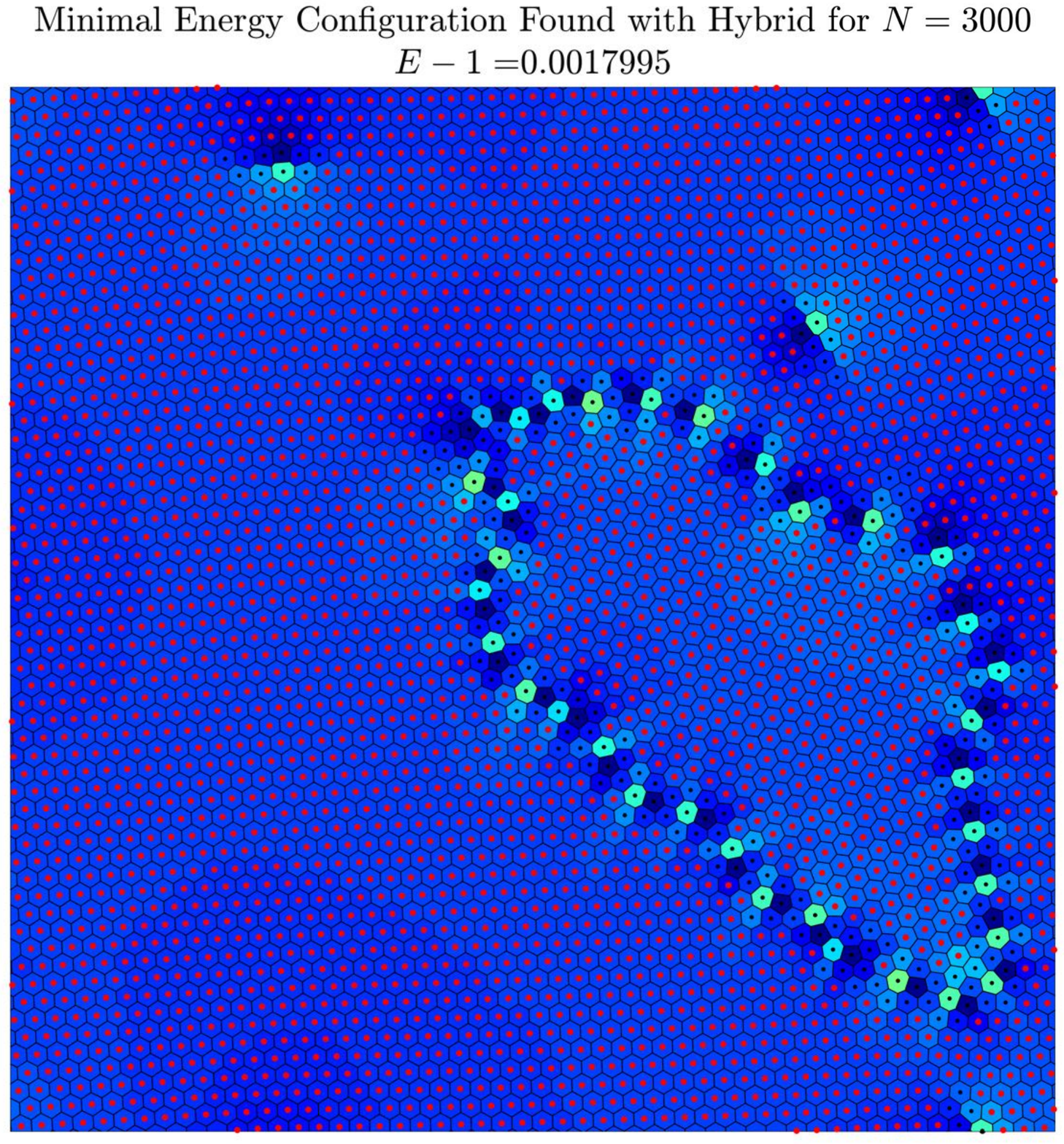} & \includegraphics[width=0.305\columnwidth]{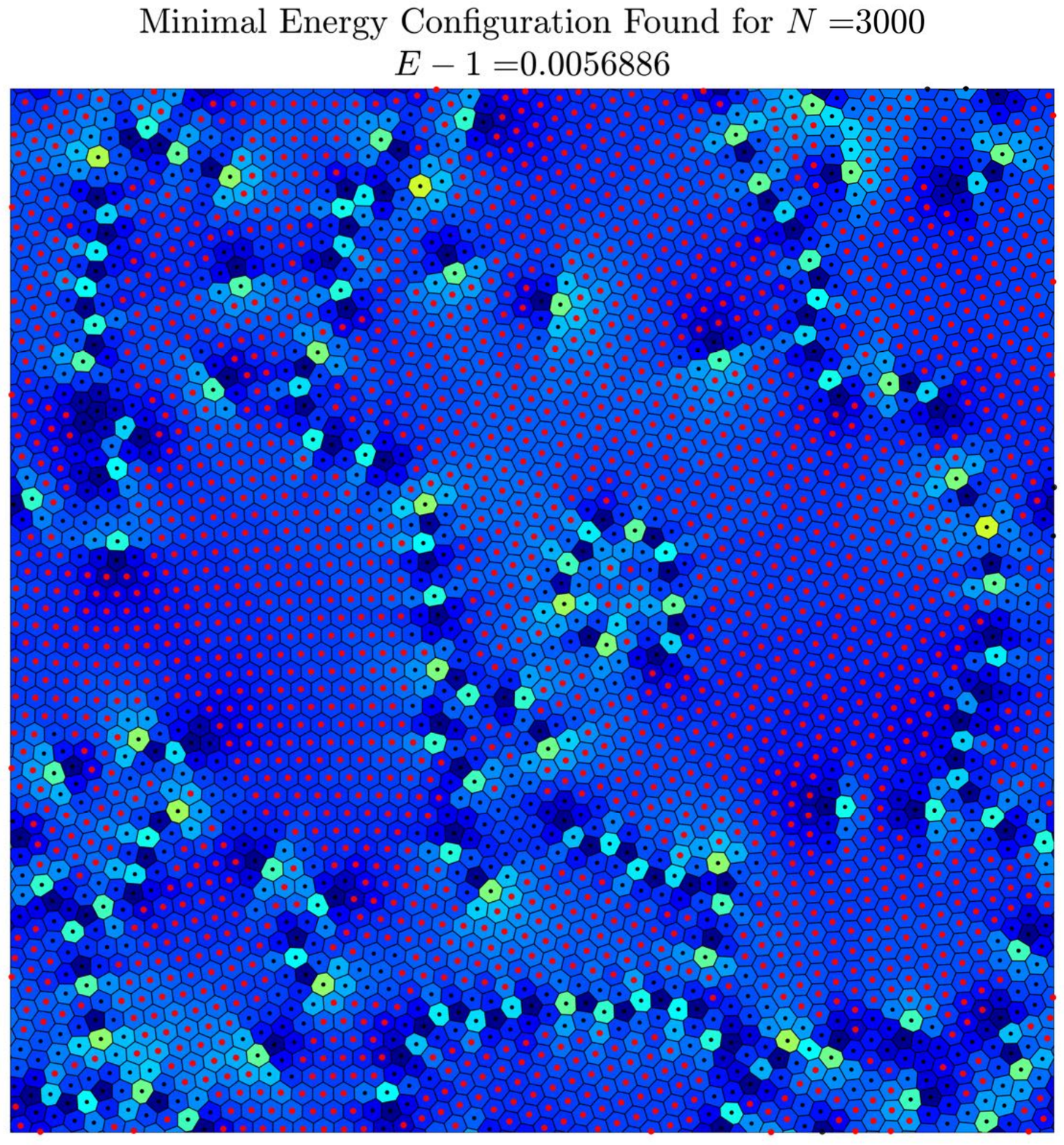}\tabularnewline
{\footnotesize{}(a) tessellation reaching $E^{min}$} & {\footnotesize{}(b) configuration with $E^{ref}$}\tabularnewline
\end{tabular}

\begin{tabular}{cc}
\includegraphics[width=0.45\columnwidth, height=0.26\columnwidth]{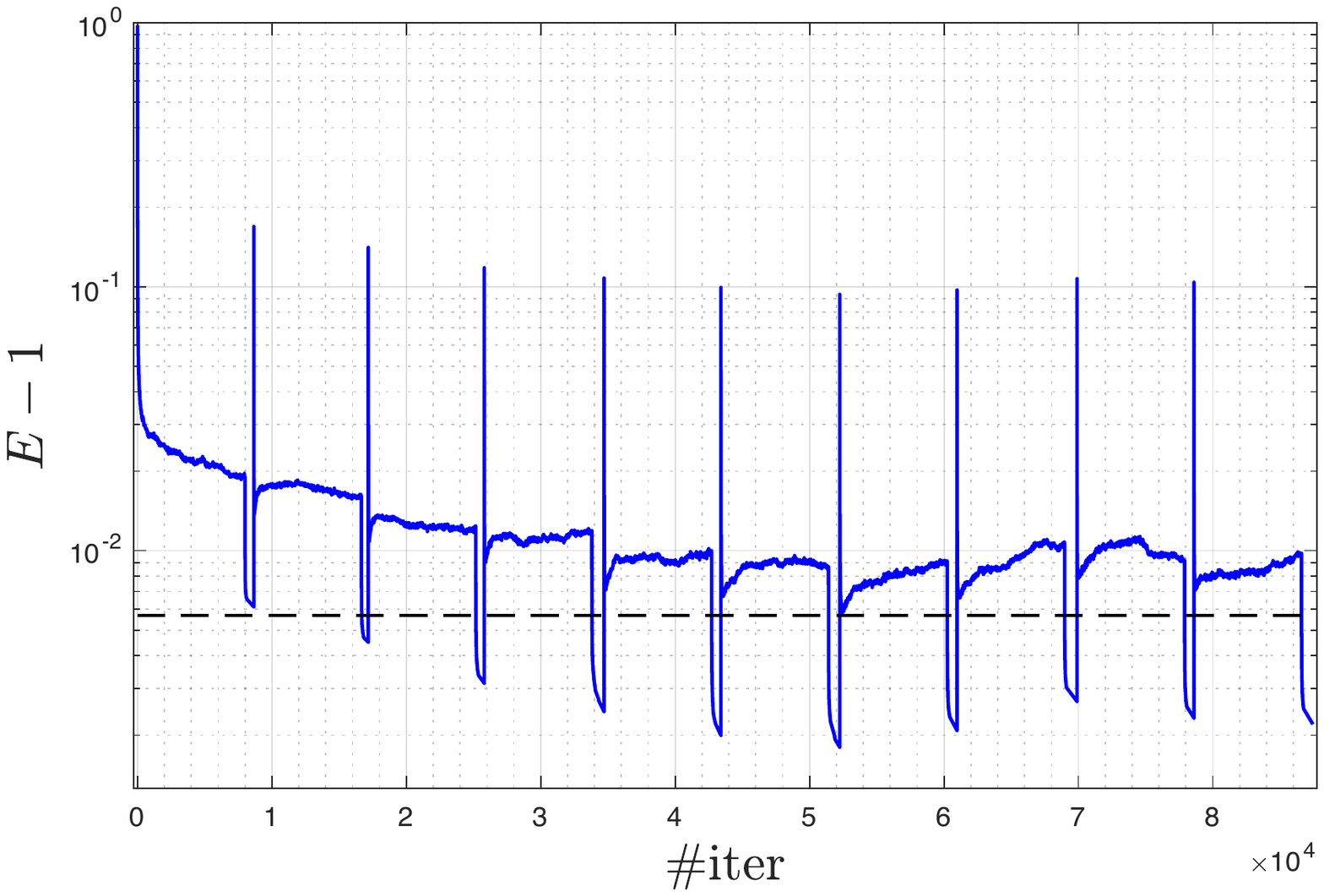} & \includegraphics[width=0.45\columnwidth, height=0.26\columnwidth]{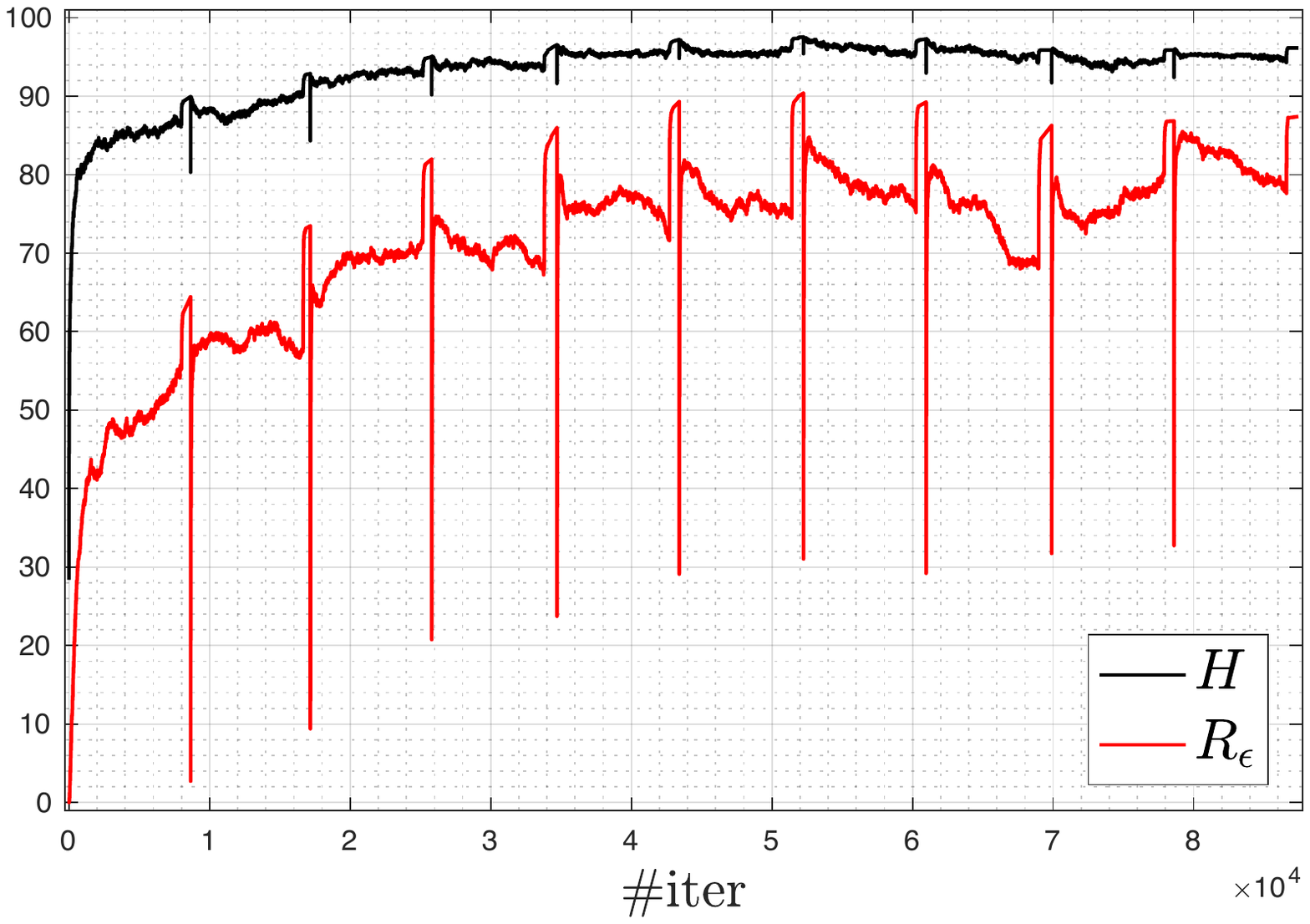}\tabularnewline
{\footnotesize{}(c) hybrid run achieving $E^{min}$ } & {\footnotesize{}(d) $R_{\epsilon},H$ profiles associated with (c)}\tabularnewline
\end{tabular}

\caption{Optimal PCVTs from Example 5 for $N=3000$ with $K=8000$: (a) lowest energy PCVT reached by the hybrid method out of 100 runs with $E^{min}-1=0.00179$. (b) tessellation carrying $E^{ref}-1=0.00568$ achieved with P-L-BFGS(20,20) out of 100,000 runs. (c) and (d) show the regularity measures as a function of the iterations of the particular hybrid run that achieved $E^{min}$.}

\end{figure}


\begin{figure}[H]
\centering

\begin{tabular}{c}
\includegraphics[width=0.6\columnwidth,height=0.28\columnwidth]{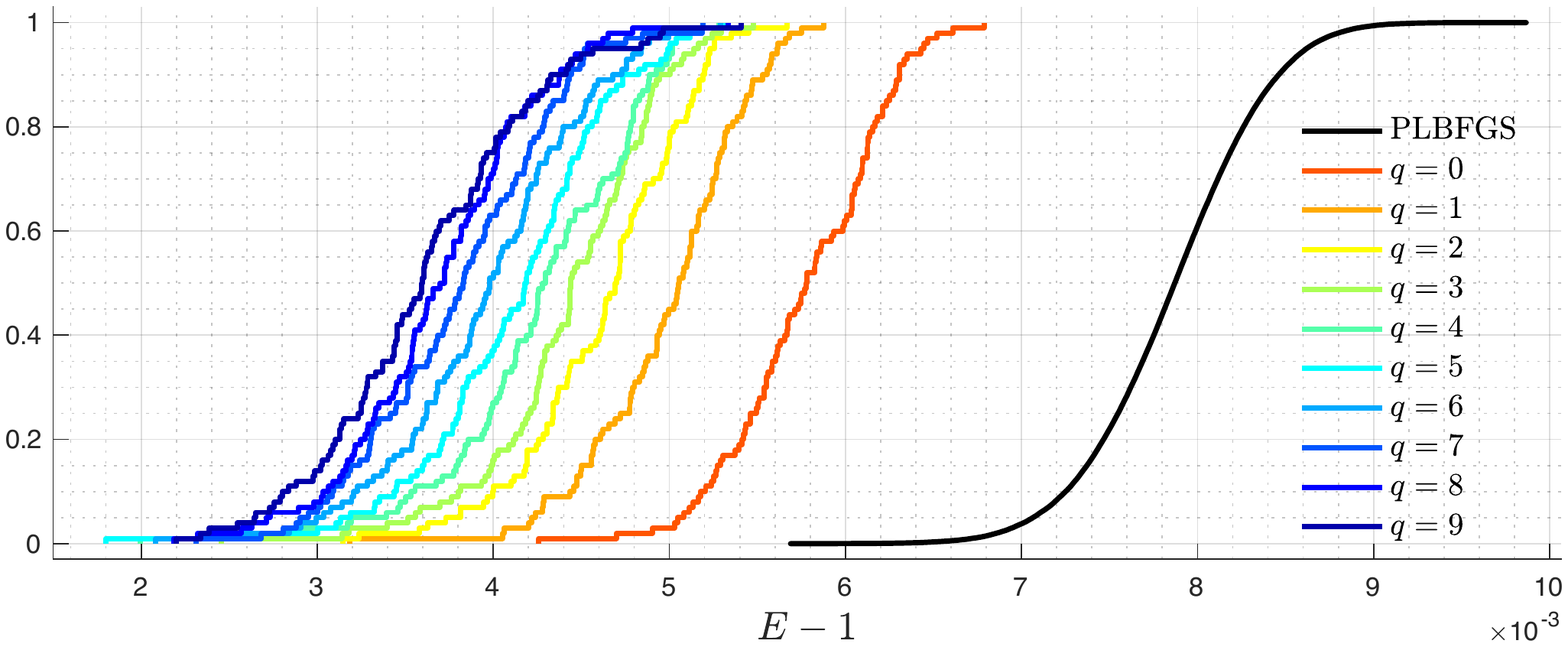}\tabularnewline
{\footnotesize{}(a) $f_{E-1}$}\tabularnewline
\end{tabular}

\begin{tabular}{cc}
\includegraphics[width=0.43\columnwidth, height=0.26\columnwidth]{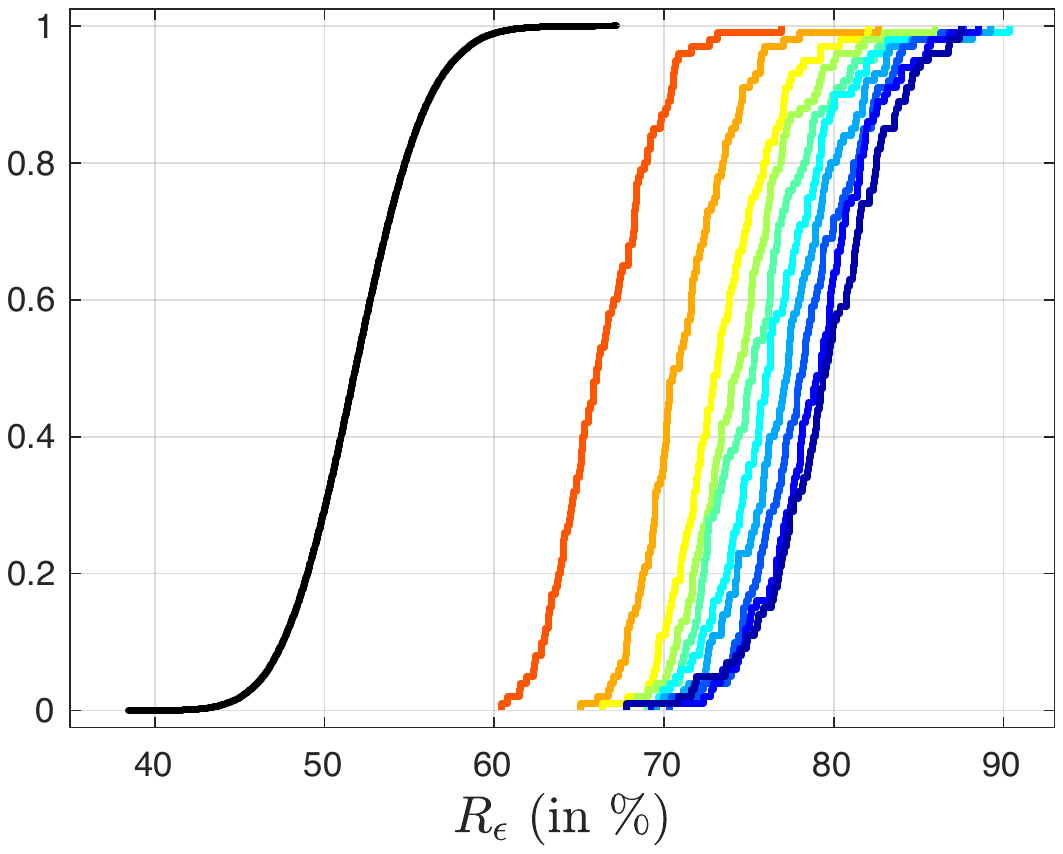} & \includegraphics[width=0.43\columnwidth, height=0.26\columnwidth]{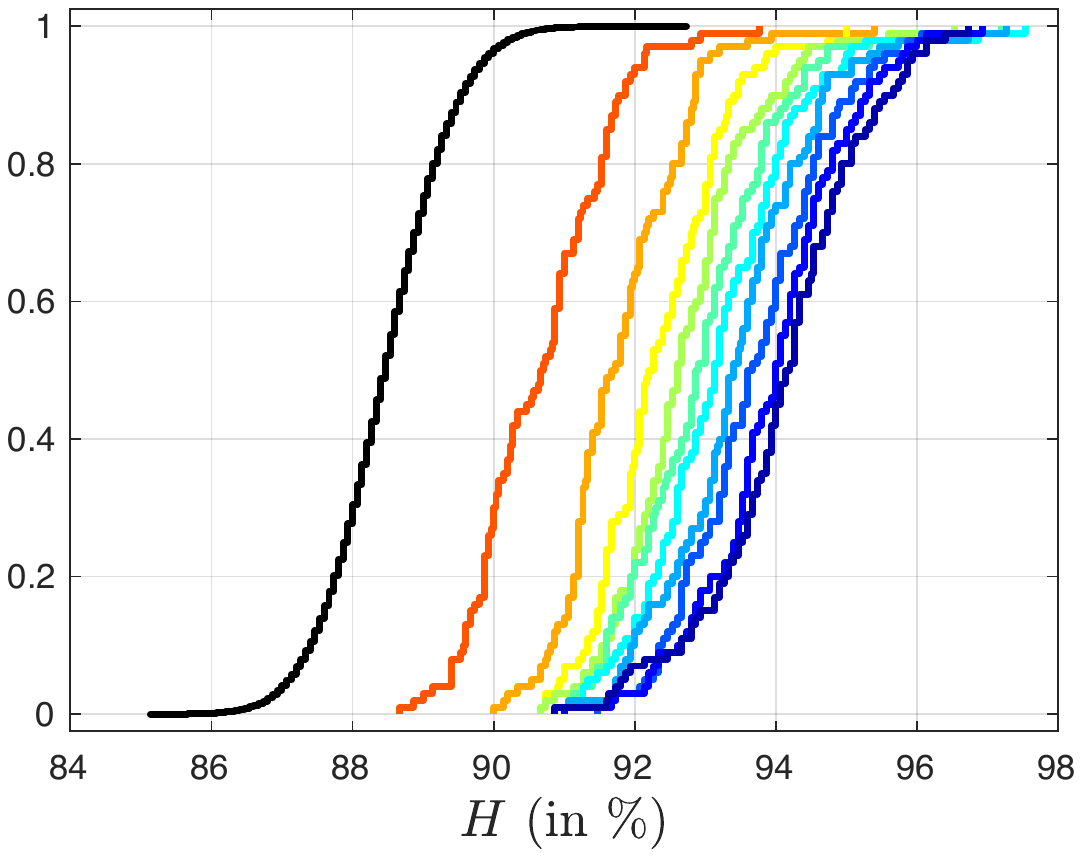}\tabularnewline
{\footnotesize{}(b) $f_{R_{\epsilon}}$} & {\footnotesize{}(c) $f_{H}$}\tabularnewline
\end{tabular}

\caption{Data from Example 5 with $N=3000$: ECDFs of the three regularity measures for 100 runs of the hybrid method with $K=8000$ along with 100,000 runs of P-L-BFGS(20,20)}

\end{figure}

\begin{figure}[H]
\centering

\begin{tabular}{cc}
\includegraphics[width=0.305\columnwidth]{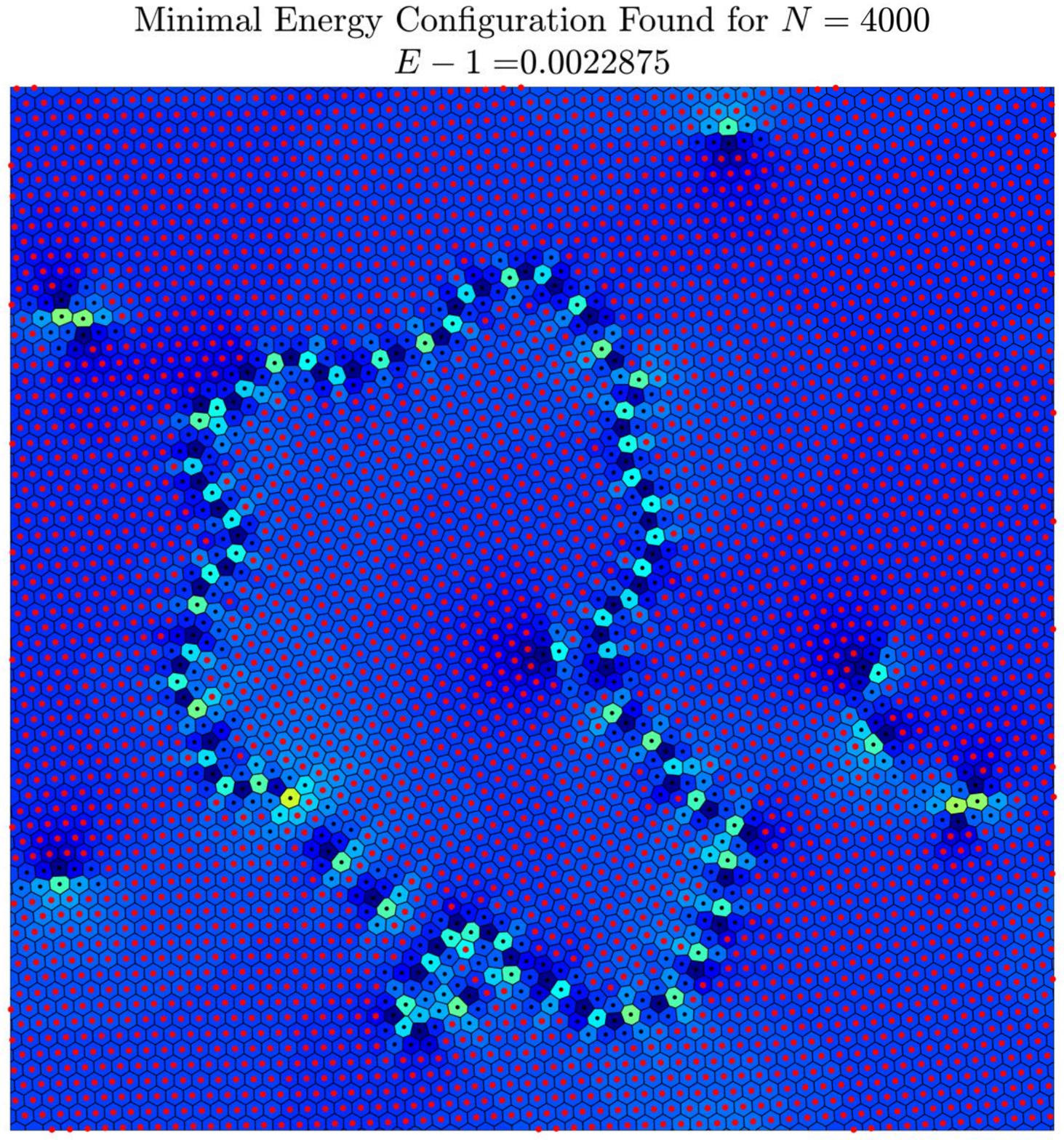} & \includegraphics[width=0.305\columnwidth]{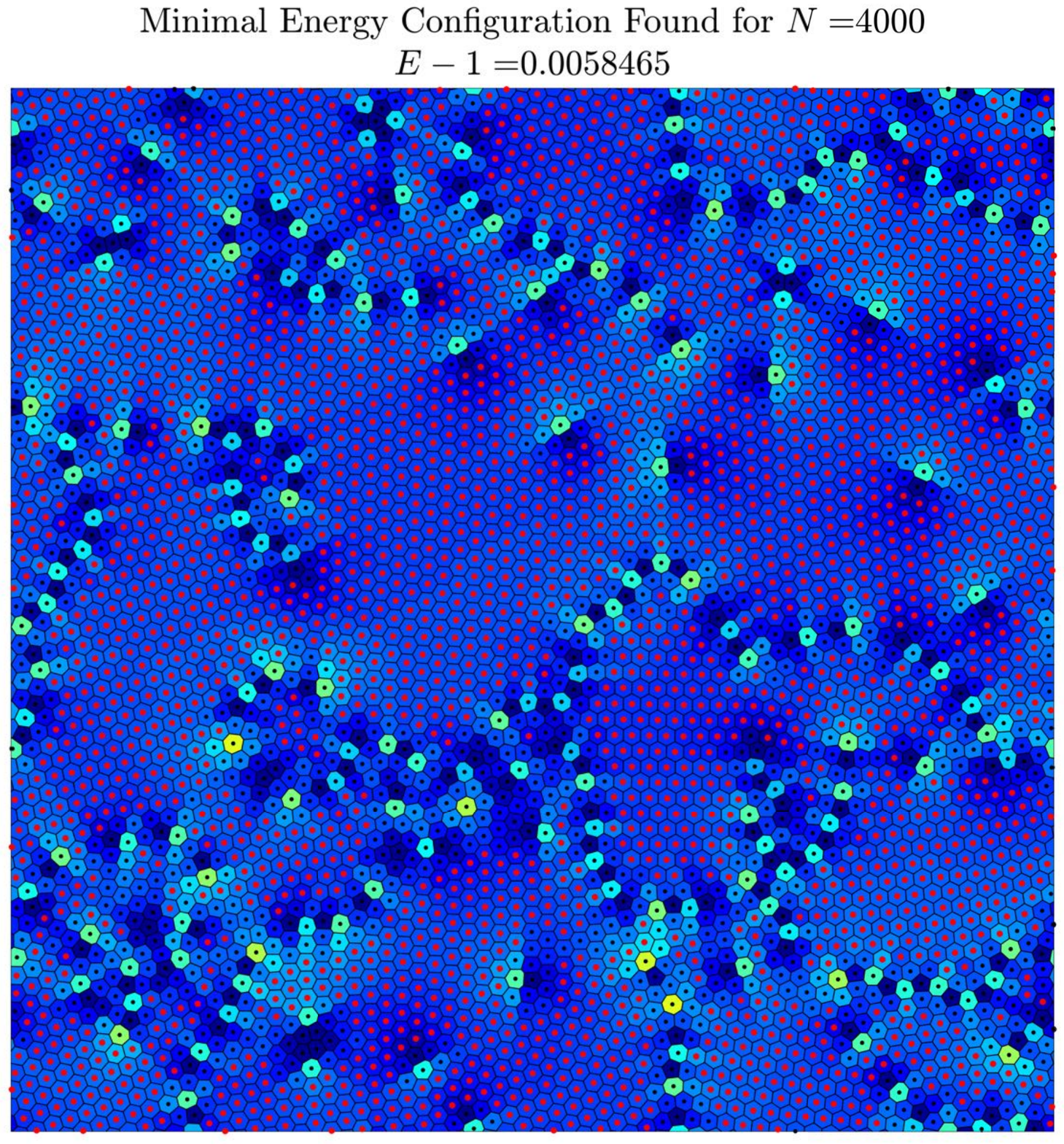}\tabularnewline
{\footnotesize{}(a) PCVT reaching $E^{min}$} & {\footnotesize{}(b) PCVT with $E^{ref}$}\tabularnewline
\end{tabular}

\begin{tabular}{cc}
\includegraphics[width=0.45\columnwidth,height=0.3\columnwidth]{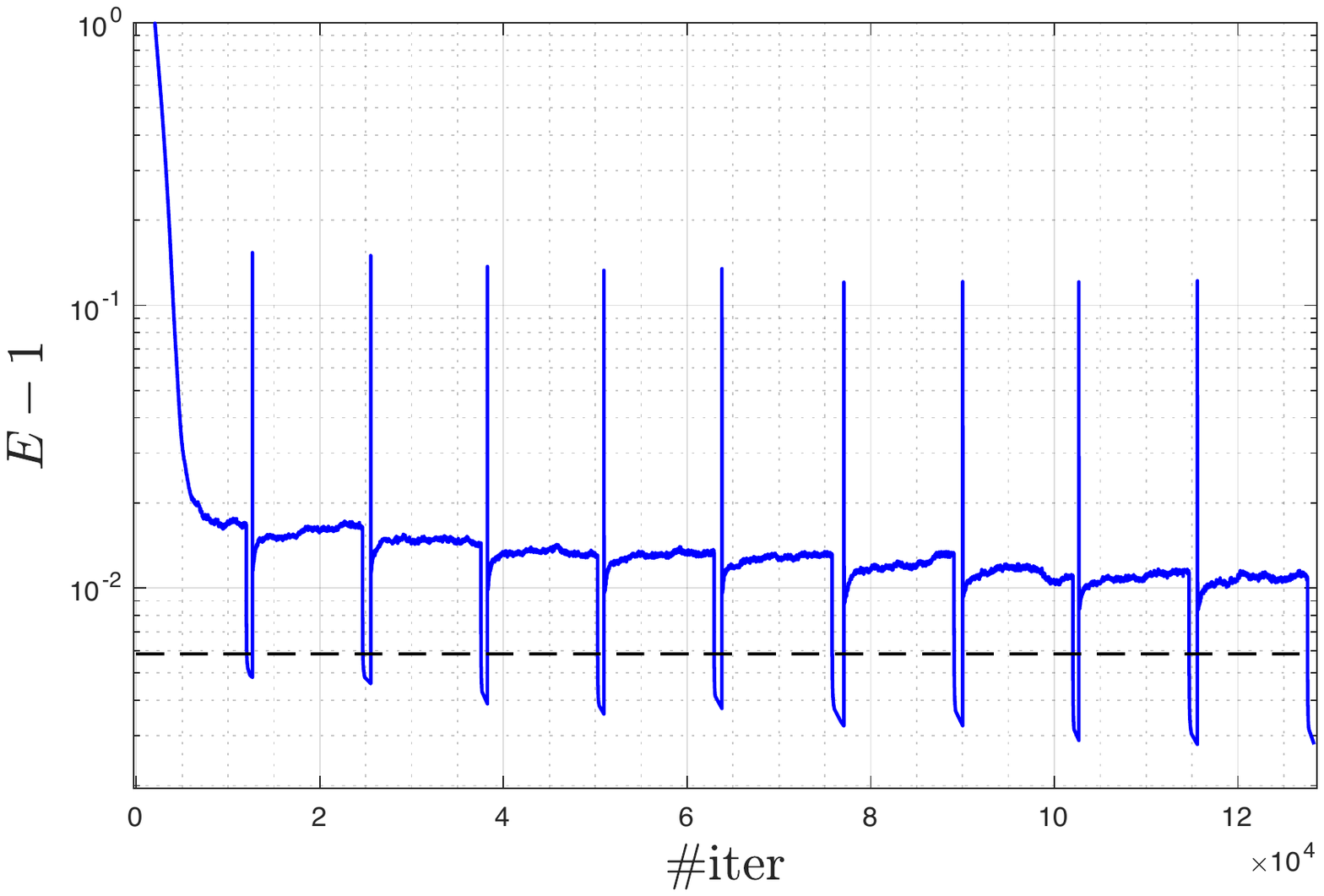} & \includegraphics[width=0.45\columnwidth,height=0.3\columnwidth]{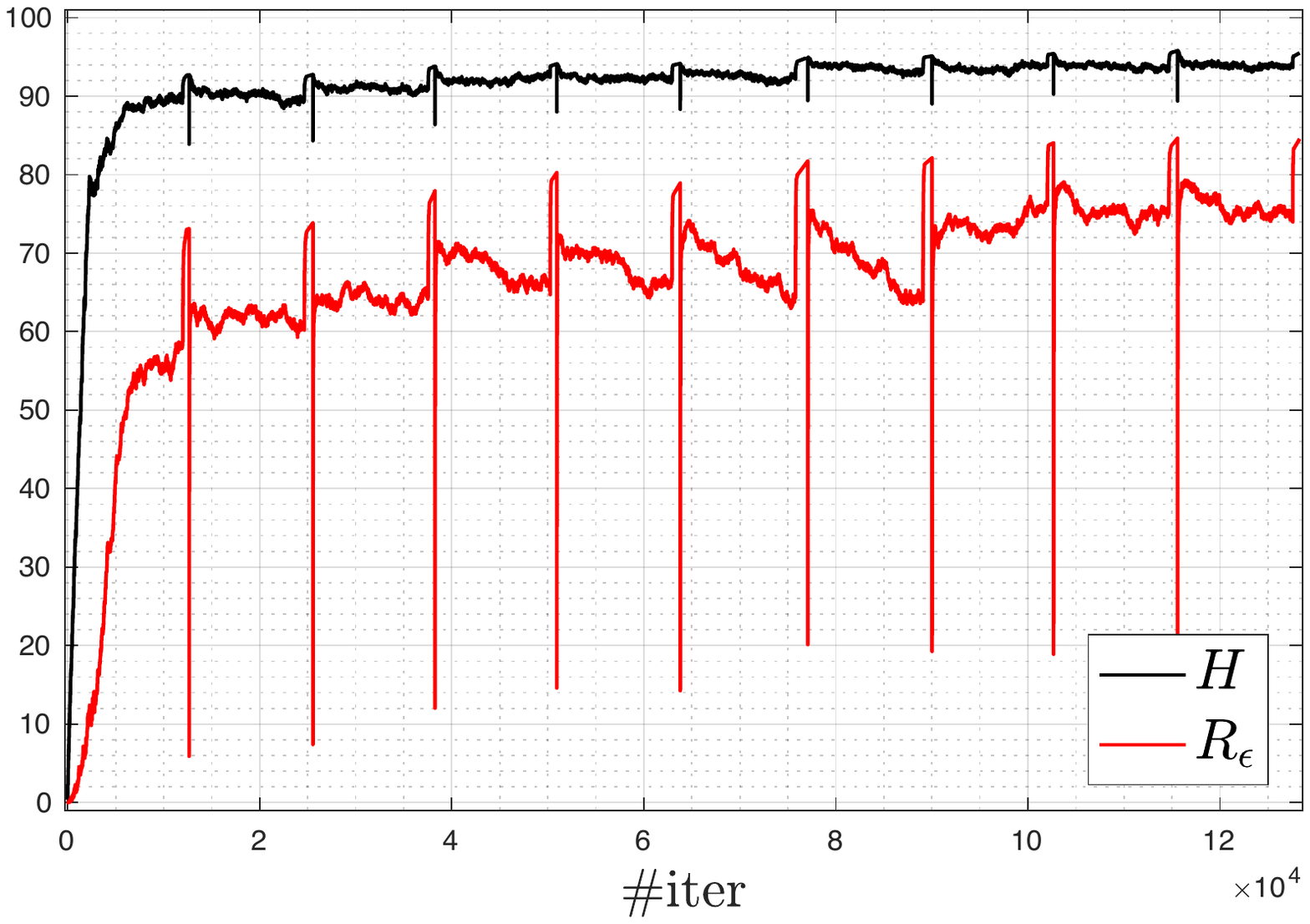}\tabularnewline
{\footnotesize{}(c) hybrid run achieving $E^{min}$ } & {\footnotesize{}(d) $R_{\epsilon},H$ profiles associated with (c)}\tabularnewline
\end{tabular}

\caption{Optimal PCVTs from Example 6 for $N=4000$ with $K=12000$: (a) lowest energy PCVT reached by the hybrid method out of 100 runs with $E^{min}-1=0.00228$. (b) tessellation carrying $E^{ref}-1=0.00584$ achieved with L-BFGS(7) out of 100,000 runs. (c) and (d) show the regularity measures as a function of the iterations of the particular hybrid run that achieved $E^{min}$.}

\end{figure}


\begin{figure}[H]
\centering
\begin{tabular}{c}
\includegraphics[width=0.6\columnwidth,height=0.28\columnwidth]{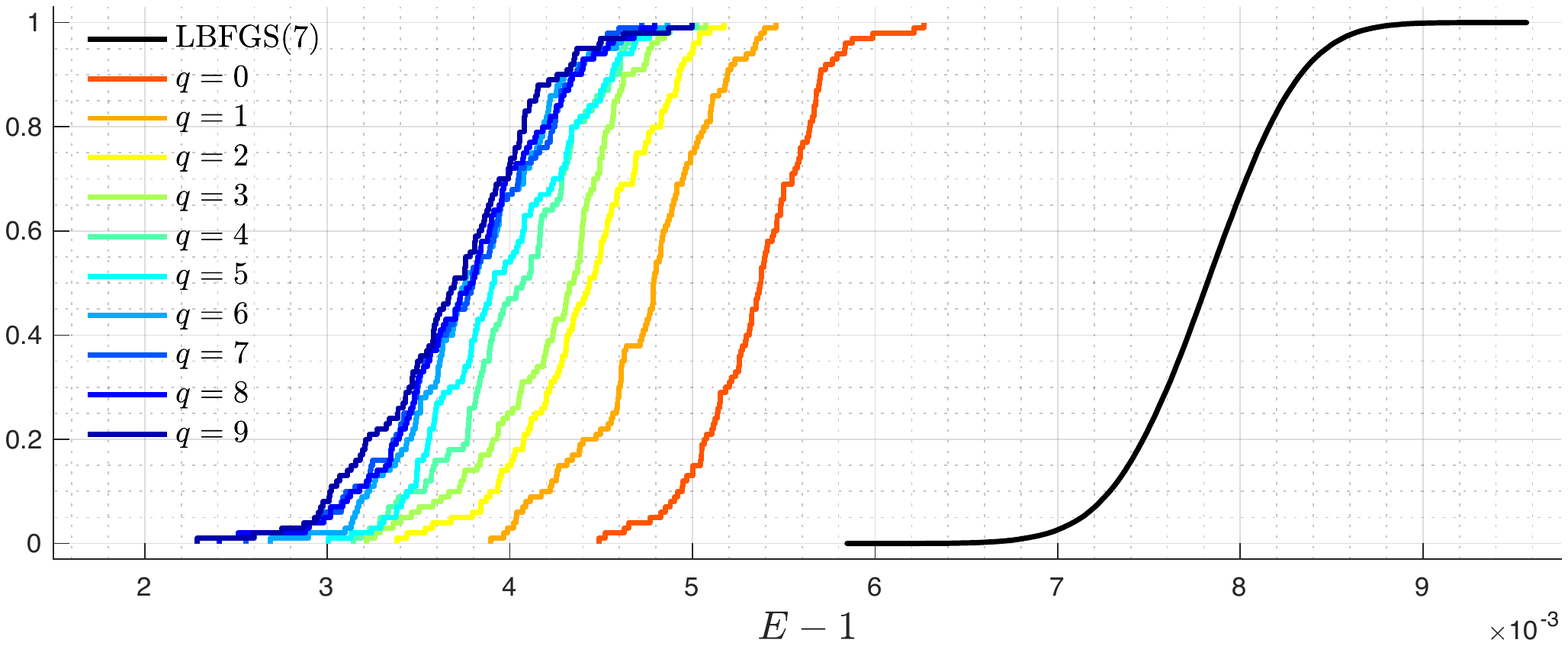}\tabularnewline
{\footnotesize{}(a) $f_{E-1}$}\tabularnewline
\end{tabular}

\begin{tabular}{cc}
\includegraphics[width=0.43\columnwidth, height=0.26\columnwidth]{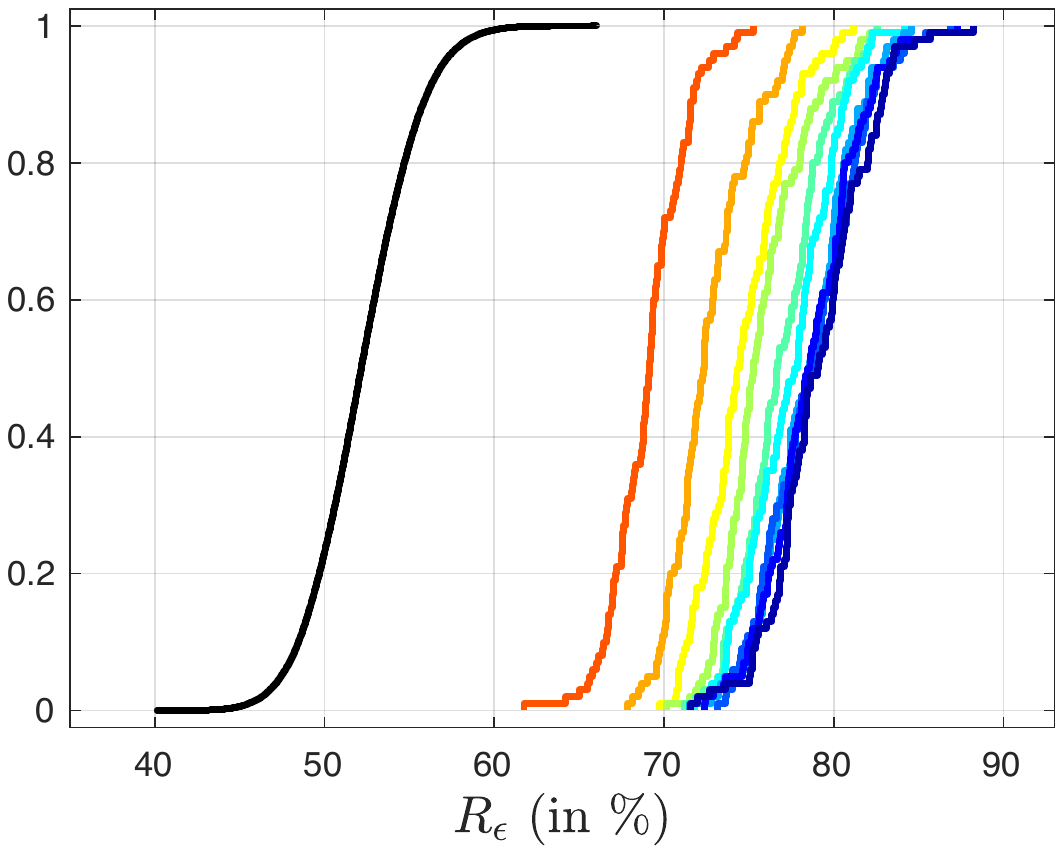} & \includegraphics[width=0.43\columnwidth, height=0.26\columnwidth]{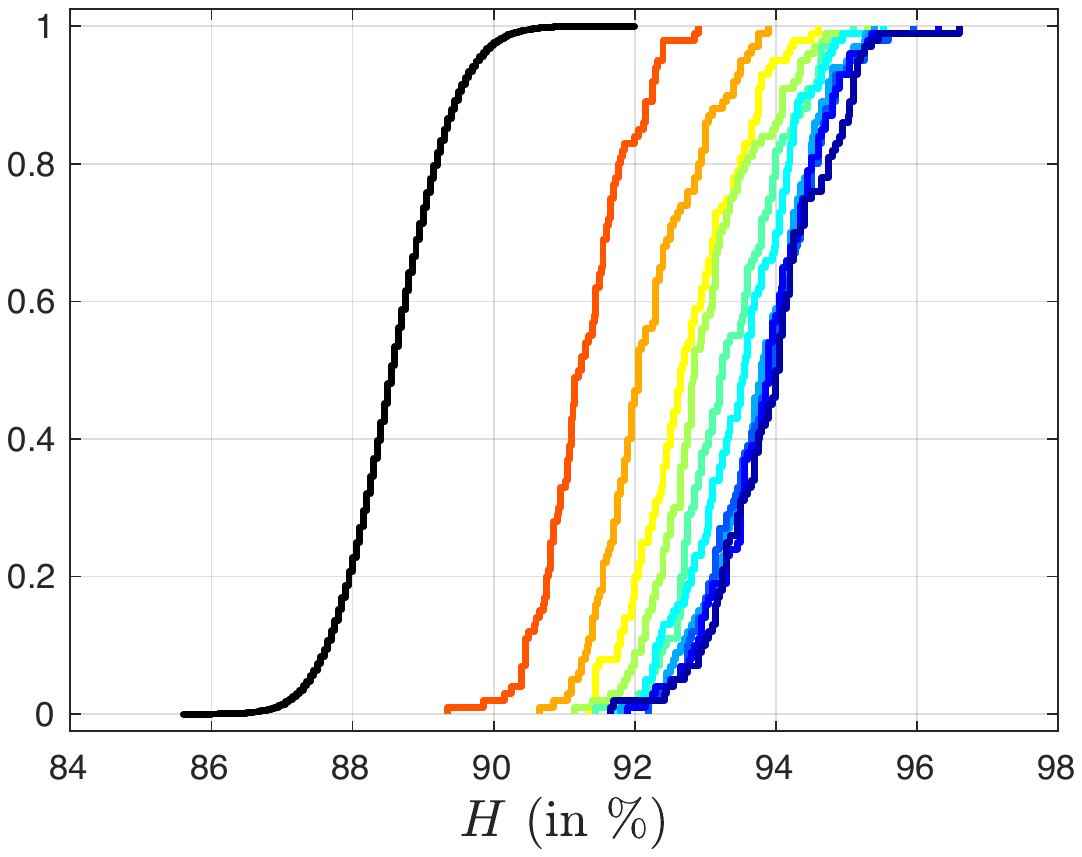}\tabularnewline
{\footnotesize{}(b) $f_{R_{\epsilon}}$} & {\footnotesize{}(c) $f_{H}$}\tabularnewline
\end{tabular}
\caption{Data from Example 6 with $N=4000$: ECDFs of the three regularity measures for 100 runs of the hybrid method with $K=12000$ along with 100,000 runs of L-BFGS(7).}
\end{figure}